\documentclass[12pt,a4paper]{article}

\usepackage[utf8]{inputenc}

\usepackage{mathrsfs,amsthm,graphicx,xcolor,verbatim,bbm,amsmath,amsfonts,amssymb,newclude,nicefrac,graphicx,hyperref,bm,geometry,mathtools,cite,mleftright}
\usepackage[shortlabels]{enumitem}
\usepackage[capitalise,sort]{cleveref} %

\crefname{enumi}{item}{items}
\crefname{equation}{}{}
\crefname{figure}{Figure}{Figures}
\crefname{subsection}{Subsection}{Subsections}

\usepackage{tikz}
\usetikzlibrary{matrix,chains,positioning,decorations.pathreplacing,arrows}
\usetikzlibrary{shapes,arrows}
\tikzset{
  font={\fontsize{9pt}{12}\selectfont}}
\usepackage{adjustbox}

\tikzset{
  font={\fontsize{9pt}{12}\selectfont}}
  
\usepackage{environ}
\makeatletter
\newsavebox{\measure@tikzpicture}
\NewEnviron{scaletikzpicturetowidth}[1]{%
  \def\tikz@width{#1}%
  \begin{lrbox}{\measure@tikzpicture}%
  \BODY
  \end{lrbox}%
  \pgfmathparse{#1/\wd\measure@tikzpicture}%
  \BODY
}
\makeatother

\usepackage{etoolbox} 
\makeatletter 
\apptocmd{\cref@getref}{\xdef\@lastusedlabel{#1}}{}{error}
\crefname{lemma}{Lemma\adddep{\loc}{\@lastusedlabel}}{Lemmas}
\crefname{prop}{Proposition\adddep{\loc}{\@lastusedlabel}}{Propositions}
\crefname{cor}{Corollary\adddep{\loc}{\@lastusedlabel}}{Corollaries}
\crefname{theorem}{Theorem\adddep{\loc}{\@lastusedlabel}}{Theorems}
\makeatother

\hypersetup{
    colorlinks,
    linkcolor={red!50!black},
    citecolor={green},
    urlcolor={blue!80!black}
}

\theoremstyle{plain}
\newtheorem{theorem}{Theorem}[section]

\theoremstyle{remark}

\theoremstyle{definition}
\newtheorem{definition}[theorem]{Definition}

\DeclareMathAlphabet{\mathpzc}{OT1}{pzc}{m}{it}

\DeclareFontEncoding{LS1}{}{}
\DeclareFontSubstitution{LS1}{stix}{m}{n}
\DeclareMathAlphabet{\mathscr}{LS1}{stixscr}{m}{n}

\newcommand{\R}{\mathbb{R}}
\newcommand{\N}{\mathbb{N}}

\newcommand{\Z}{\mathbb{Z}}

\newcommand{\smallsum}{\textstyle\sum}

\renewcommand{\emptyset}{\varnothing}

\newcommand{\fa}[1]{\forall\, #1 \in}

\DeclarePairedDelimiter{\norm}{\lVert}{\rVert\cfadd{def:Euclidean_norm}}
\DeclarePairedDelimiter{\vass}{\lvert}{\rvert}
\DeclarePairedDelimiter{\abs}{\lvert}{\rvert}
\DeclarePairedDelimiter{\pr}{(}{)}
\DeclarePairedDelimiter{\PR}{[}{]}
\DeclarePairedDelimiter{\pR}{\{}{\}}
\DeclarePairedDelimiter{\ceil}{\lceil}{\rceil\cfadd{def:ceiling}}
\DeclarePairedDelimiter{\mednorm}{\lVert}{\rVert\cfadd{def:Euclidean_norm}}
\newcommand{\bignorm}[1]{\mednorm[\big]{#1}}

\newcommand{\prb}[1]{\pr[\big]{ #1 }}
\newcommand{\PRb}[1]{\PR[\big]{ #1 }}
\newcommand{\pRb}[1]{\pR[\big]{ #1 }}
\newcommand{\prbbb}[1]{\pr[\bigg]{ #1 }}
\newcommand{\PRbbb}[1]{\PR[\bigg]{ #1 }}
\newcommand{\pRbbb}[1]{\pR[\bigg]{ #1 }}
\newcommand{\prbb}[1]{\pr[\Big]{ #1 }}
\newcommand{\PRbb}[1]{\PR[\Big]{ #1 }}

\newcommand{\PRbbbb}[1]{\PR[\Bigg]{ #1 }}

\newcommand{\andShort}{\text{ and }}

\newcommand{\radius}{R}

\newcommand{\CostLipwo}[2]{\cfadd{def:cost_Lipschitz_approx_set}\mathrm{Cost}_{#1}}
\newcommand{\CostLip}[5]{\cfadd{def:cost_Lipschitz_approx_set}\mathrm{Cost}_{#1}( #3,#4,#5)}
\newcommand{\CostLipA}[5]{\cfadd{def:cost_Lipschitz_approx_set}\mathrm{Cost}_{#1}\pr*{#3,#4,#5}}
\newcommand{\CostLipB}[5]{\cfadd{def:cost_Lipschitz_approx_set}\mathrm{Cost}_{#1}\prb{#3,#4,#5}}

\newcommand{\normmm}[1]{{\left\vert\kern-0.25ex\left\vert\kern-0.25ex\left\vert #1 
    \right\vert\kern-0.25ex\right\vert\kern-0.25ex\right\vert}} %

\newcommand{\qandq}{\qquad\text{and}\qquad}

\newcommand{\param}{\mathcal{P}}

\newcommand{\realisation}{\cfadd{def:ANNrealization}\mathcal{R}}

\newcommand{\indicator}[1]{\mathbbm{1}_{#1}}
\newcommand{\idMatrix}{\cfadd{def:identityMatrix}\operatorname{I}}

\newcommand{\ANNs}{\cfadd{def:ANN}\mathbf{N}}
\newcommand{\activation}{a}

\newcommand{\functionANN}{\cfadd{def:ANNrealization}\mathcal{R}_{\activation}}
\newcommand{\functionANNbis}[1]{\cfadd{def:ANNrealization}\pr*{\mathcal{R}_{\activation}\pr*{ #1 }} }
\newcommand{\paramANN}{\cfadd{def:ANN}\mathcal{P}}

\newcommand{\lengthANN}{\cfadd{def:ANN}\mathcal{L}}

\newcommand{\inDimANN}{\cfadd{def:ANN}\mathcal{I}}
\newcommand{\compANN}[2]{{#1 \bullet #2\cfadd{def:ANNcomposition}}}

\newcommand{\outDimANN}{\cfadd{def:ANN}\mathcal{O}}
\newcommand{\longerANN}[1]{\cfadd{def:ANNenlargement}\mathcal{E}_{#1}}

\newcommand{\idANNshort}[1]{\mathfrak{I}_{#1}}
\newcommand{\interpolatingNN}[2]{\cfadd{def:interpolatingDNN}\mathbf{I}^{#1}_{#2}} 
\newcommand{\ReLUidANN}[1]{\cfadd{def:ReLU_identity}\mathbb{I}_{#1}}
\newcommand{\clip}[2]{\cfadd{def:clipping_function}\mathfrak{c}_{ #1, #2 }}
\newcommand{\Clip}[3]{\cfadd{def:clip}\mathfrak{C}_{ #1, #2, #3 }}
\newcommand{\Mult}[1]{\cfadd{def:multidim_version}\mathfrak{M}_{#1}}
\newcommand{\imdim}[2]{\cfadd{def:smart_restriction}\mathfrak{d}_{#1,#2}}

\newcommand{\affineANN}{\cfadd{def:ANN:affine}\mathbf{A}}
\newcommand{\maxANN}{\cfadd{def:max_d}\mathcal{M}}

\newcommand{\dims}{\cfadd{def:ANN}\mathcal{D}}
\newcommand{\singledims}{\cfadd{def:ANN}\mathbb{D}}

\newcommand{\parallelization}{\cfadd{def:generalParallelization}\mathbf{P}}
\newcommand{\parallelizationSpecial}{\cfadd{def:simpleParallelization}\mathbf{P}}

\newcommand{\ReLU}{\mathfrak{r}\cfadd{def:ReLU}}
\newcommand{\functionReLUANN}{\cfadd{def:ANNrealization}\mathcal{R}_{\ReLU}}
\newcommand{\functionnbReLUANN}[1]{\cfadd{def:ANNrealization}\pr{\mathcal{R}_{\ReLU}\pr{{#1}}}}

\newcommand{\Dapprox}{\mathcal{V}\cfadd{def:polyD_mult}}
\newcommand{\Capprox}[2]{\cfadd{def:polyC}\mathcal{C}_{#1,#2}}

\newcommand{\fconstant}{\mathfrak{f}}
\newcommand{\gconstant}{\mathfrak{g}}

\newcommand{\gp}{\xi}

\renewenvironment{pmatrix}{\mleft(\begin{matrix}}{\end{matrix}\mright)}

\DeclareMathOperator*{\ssum}{\textstyle\sum}

\DeclareMathOperator{\ssssum}{\textstyle\sum}

\DeclareMathOperator*{\sprod}{\textstyle\prod}

\DeclareMathOperator{\sssprod}{\textstyle\prod}

\DeclareMathOperator{\sssbigtimes}{\textstyle\bigtimes}

\newcommand{\is}{\curvearrowleft}

\newcommand{\eps}{\varepsilon}

\ExplSyntaxOn

\seq_new:N \g_cflist_loaded
\seq_new:N \g_cflist_pending

\NewDocumentCommand{\cfadd}{ m }
{
  \seq_if_in:NnF \g_cflist_loaded { #1 } {
    \seq_if_in:NnF \g_cflist_pending { #1 } {
      \seq_gput_right:Nn \g_cflist_pending { #1 }
    }
  }
}

\NewDocumentCommand{\cfconsiderloaded}{ m }{
  \seq_gput_right:Nn \g_cflist_loaded {#1}
}

\NewDocumentCommand{\cfremove}{ m }
{
  \seq_gremove_all:Nn \g_cflist_pending { #1 }
}

\NewDocumentCommand{\cfload}{ o }
{
  \seq_if_empty:NTF \g_cflist_pending {\unskip} {
    (cf.\ \cref{\seq_use:Nn \g_cflist_pending {,}})\IfValueTF{#1}{#1~}{\unskip}
    \seq_gconcat:NNN \g_cflist_loaded \g_cflist_loaded \g_cflist_pending
    \seq_gclear:N \g_cflist_pending
  }
}

\NewDocumentCommand{\cfclear} {} {
  \seq_gclear:N \g_cflist_loaded
  \seq_gclear:N \g_cflist_pending
}

\NewDocumentCommand{\cfout}{ o }
{
  \seq_if_empty:NTF \g_cflist_pending {\unskip} {
    (cf.\ \cref{\seq_use:Nn \g_cflist_pending {,}})\IfValueTF{#1}{#1~}{\unskip}
    \seq_gclear:N \g_cflist_pending
  }
}

\NewDocumentCommand{\ifnocf} { m } {
  \seq_if_empty:NT \g_cflist_pending { #1 }
}

\ExplSyntaxOff

\ExplSyntaxOn

\bool_new:N \g_noteobserve

\NewDocumentCommand{\setnote}{}{
  \bool_gset_true:N \g_noteobserve
}

\NewDocumentCommand{\setobserve}{}{
  \bool_gset_false:N \g_noteobserve
}

\NewDocumentCommand{\nobs}{ o }{
  \IfValueT{#1}{
    \str_if_eq:noTF {note} {#1} {
      \bool_gset_true:N \g_noteobserve
    } {
      \str_if_eq:noTF {Note} {#1} {
        \bool_gset_true:N \g_noteobserve
      } {
        \bool_gset_false:N \g_noteobserve
      }
    }
  }
  \bool_if:nTF { \g_noteobserve } {
    \bool_gset_false:N \g_noteobserve 
    note
  } {
    \bool_gset_true:N \g_noteobserve 
    observe
  }
  \IfValueF{#1}{~}
}

\NewDocumentCommand{\Nobs}{ o }{
  \IfValueT{#1}{
    \str_if_eq:noTF {note} {#1} {
      \bool_gset_true:N \g_noteobserve
    } {
      \str_if_eq:noTF {Note} {#1} {
        \bool_gset_true:N \g_noteobserve
      } {
        \bool_gset_false:N \g_noteobserve
      }
    }
  }
  \bool_if:nTF { \g_noteobserve } {
    \bool_gset_false:N \g_noteobserve 
    Note
  } {
    \bool_gset_true:N \g_noteobserve 
    Observe
  }
  \IfValueF{#1}{~}
}

\bool_new:N \g_hencetherefore

\NewDocumentCommand{\hence}{ o }{
  \IfValueT{#1}{
    \str_if_eq:noTF {hence} {#1} {
      \bool_gset_true:N \g_hencetherefore
    } {
      \str_if_eq:noTF {Hence} {#1} {
        \bool_gset_true:N \g_hencetherefore
      } {
        \bool_gset_false:N \g_hencetherefore
      }
    }
  }
  \bool_if:nTF { \g_hencetherefore } {
    \bool_gset_false:N \g_hencetherefore 
    hence
  } {
    \bool_gset_true:N \g_hencetherefore 
    therefore
  }
  \IfValueF{#1}{~}
}

\NewDocumentCommand{\Hence}{ o }{
  \IfValueT{#1}{
    \str_if_eq:noTF {hence} {#1} {
      \bool_gset_true:N \g_hencetherefore
    } {
      \str_if_eq:noTF {Hence} {#1} {
        \bool_gset_true:N \g_hencetherefore
      } {
        \bool_gset_false:N \g_hencetherefore
      }
    }
  }
  \bool_if:nTF { \g_hencetherefore } {
    \bool_gset_false:N \g_hencetherefore 
    Hence
  } {
    \bool_gset_true:N \g_hencetherefore 
    Therefore
  }
  \IfValueF{#1}{~}
}

\int_new:N \g_furthermore

\NewDocumentCommand{\Moreover}{ o o }{
  \IfValueT{#1}{
    \str_case:nn {#1} {
      {Furthermore} {\int_set:Nn {\g_furthermore} {0}}
      {Moreover} {\int_set:Nn {\g_furthermore} {1}}
      {In~addition} {\int_set:Nn {\g_furthermore} {2}}
      {note} {\bool_gset_true:N \g_noteobserve}
      {observe} {\bool_gset_false:N \g_noteobserve}
    }
    \IfValueT{#2}{
      \str_case:nn {#2} {
        {Furthermore} {\int_set:Nn {\g_furthermore} {0}}
        {Moreover} {\int_set:Nn {\g_furthermore} {1}}
        {In~addition} {\int_set:Nn {\g_furthermore} {2}}
        {note} {\bool_gset_true:N \g_noteobserve}
        {observe} {\bool_gset_false:N \g_noteobserve}
      }
    }
  }
  \int_case:nn { \int_mod:nn {\g_furthermore} {3} } {
    { 0 } { Furthermore,~\nobs that}
    { 1 } { Moreover,~\nobs that}
    { 2 } { In~addition,~\nobs that}
  }
  \int_incr:N \g_furthermore
  \IfValueF{#1}{~}
}

\ExplSyntaxOff

\ExplSyntaxOn
\global\def\loc{dummy}

\NewDocumentEnvironment {athm} {m m} {%
\begin{#1}\label{#2}\global\def\loc{#2}%
}{%
\end{#1}%
}

\NewDocumentEnvironment{aproof} {} {%
\begin{proof}[Proof~of~\cref{\loc}]%
}{%
\global\def\loc{dummy}\end{proof}%
}

\ExplSyntaxOff

\newcommand{\eqqref}[1]{\eqref{eq:\loc.#1}}
\newcommand{\eqlabel}[1]{\label{eq:\loc.#1}}
\newcommand{\lref}[1]{\cref{\loc.#1}}
\newcommand{\Lref}[1]{\Cref{\loc.#1}}
\newcommand{\llabel}[1]{\label{\loc.#1}}

\newcommand{\finishproofthus}{The proof of \cref{\loc} is thus complete.}
\newcommand{\finishproofthis}{This completes the proof of \cref{\loc}.}

\ExplSyntaxOn

\seq_new:N \g_deplist
\seq_new:N \l_done
\tl_new:N \l_first_tl
\tl_new:N \l_second_tl
\int_new:N \l_count_int
\iow_new:N \l_logfile

\NewDocumentCommand{\adddep}{ m m }
{
	\seq_gput_right:Nx \g_deplist { #1 }
	\seq_gput_right:Nx \g_deplist { #2 }
}

\NewDocumentCommand{\listdeps}{}
{
	\int_set:Nn \l_count_int {\seq_count:N {\g_deplist}}
	\seq_clear:N \l_done
	\int_while_do:nNnn {\l_count_int} > {0} {
 	\iow_open:Nn \l_logfile {dependencies.dot}
	\seq_gpop_left:NN \g_deplist \l_first_tl
		\seq_gpop_left:NN \g_deplist \l_second_tl
		\cs_if_eq:NNTF { \l_first_tl } { \l_second_tl } {} {
			\tl_if_eq:NnTF { \l_first_tl } { dummy } {} {
				\seq_if_in:NxTF { \l_done } {\l_first_tl \l_second_tl} {} {
					\iow_now:Nn \l_logfile {\cref*{\tl_use:N \l_second_tl}}
					\iow_now:Nn \l_logfile {->}
					\iow_now:Nn \l_logfile {\cref*{\tl_use:N \l_first_tl}}
					{\cref{\tl_use:N \l_second_tl}}~
					~$\to$~
					{\cref{\tl_use:N \l_first_tl}}\\
					\seq_put_left:Nx {\l_done} {\l_first_tl \l_second_tl}
				}
			}
		}
	 	\int_decr:N {\l_count_int}
		\int_decr:N {\l_count_int}
	}
	\iow_close:N \l_logfile
}

\ExplSyntaxOff

\makeatletter
\newcommand{\vast}{\bBigg@{3.5}}
\newcommand{\Vast}{\bBigg@{5.5}}
\makeatother

\title{\vspace{-2cm}Deep neural network approximation\\theory for high-dimensional functions}

\author{
Pierfrancesco Beneventano$^{1,2,a}$,
Patrick Cheridito$^{1,b}$,\\
Robin Graeber$^{3,c}$,
Arnulf Jentzen$^{3,4,d}$,
and
Benno Kuckuck$^{3,e}$
\bigskip
\\
\small{$^1$\,Department of Mathematics, ETH Zurich, Switzerland}
\\
\small{$^2$\,Department of Operations Research and Financial Engineering,}
\\[-0.13cm]
\small{Princeton University, NJ, United States}
\\
\small{$^3$\,Applied Mathematics: Institute for Analysis and Numerics,}
\\[-0.13cm]
\small{Faculty of Mathematics and Computer Science,}
\\[-0.13cm]
\small{University of M{\"u}nster, Germany}
\\
\small{$^4$\,School of Data Science and Shenzhen Research Institute of Big Data,}
\\[-0.13cm]
\small{The Chinese University of Hong Kong, Shenzhen, China}
\\
\small{$^a$\,\textit{pierbene96@gmail.com}, \textit{pierb@princeton.edu}}
\\[-0.13cm]
\small{$^b$\,\textit{patrick.cheridito@math.ethz.ch}}
\\[-0.13cm]
\small{$^c$\,\textit{r\_grae02@uni-muenster.de}}
\\[-0.13cm]
\small{$^d$\,\textit{ajentzen@uni-muenster.de}, \textit{ajentzen@cuhk.edu.cn}}
\\[-0.13cm]
\small{$^e$\,\textit{bkuckuck@uni-muenster.de}}
}

\begin{document}

\maketitle

\begin{abstract}
The purpose of this article is to develop a machinery to study
the capacity of deep neural networks (DNNs) 
to approximate high-dimensional functions.
In particular, we show that DNNs have the expressive power to overcome
the curse of dimensionality in the approximation of a large class of 
functions. More precisely, we prove that these functions can be approximated by DNNs
on compact sets such that the number of parameters necessary to 
represent the approximating
DNNs grows at most polynomially in the reciprocal $\nicefrac1\eps$ of the
prescribed approximation error $\eps>0$ and in the 
input dimension $d\in\N$.
To this end, we introduce certain approximation spaces, consisting
of sequences of functions that can be efficiently approximated by
DNNs. We then establish closure properties which we combine with
known and new bounds on the number of parameters necessary to
approximate locally Lipschitz
continuous functions, maximum functions, and product functions by DNNs.
The main result of this article demonstrates
that DNNs have sufficient expressive power to approximate, without the curse of dimensionality, certain sequences
of functions which can be constructed by means of a finite number of compositions
using locally Lipschitz continuous functions, maxima, and products.
\end{abstract}

\begin{samepage}
\tableofcontents
\end{samepage}

\section{Introduction}
\label{sect:intro}

In recent years, deep learning has enjoyed tremendous success in many 
real-world application areas such as computer vision (e.g., image recognition,
image segmentation, medical image analysis), natural language processing 
(e.g., speech recognition,
machine translation, information retrieval), 
finance (e.g., fraud detection, risk management). 

In the basic setting of supervised learning, where the goal
is to find an approximation of a target function given a limited amount
of training data, a major appeal of deep learning methods is their
apparent ability to scale to very high-dimensional domains.
In these settings, more traditional approximation schemes frequently suffer from the
so-called \emph{curse of dimensionality}: The number of computational
steps necessary to achieve a given approximation accuracy grows
exponentially with the dimension of the domain of the
function that is to be approximated
(cf., e.g., Bellman \cite{Bellman1957}, Novak \& Wo\'zniakowski \cite{NovakWozniakowski2008,NovakWozniakowski2010},
and Novak \& Ritter \cite{novak1997curse}).
By contrast,
deep learning methods
appear to achieve good approximation accuracy in applications
even for functions
on very high dimensional domains in a manageable amount of time.
Empirical studies and simulations further corroborate
this impression that deep learning methods, in many settings,
are able to overcome the curse of dimensionality
in the sense that the number of computational steps
necessary to achieve a given approximation accuracy
grows at most polynomially with the dimension $d\in\N=\{1,2,3,\dots\}$ of the
domain of the target function and the reciprocal $\nicefrac1\eps$ of the prescribed approximation
error $\eps\in(0,\infty)$.
This latter property is sometimes called \emph{polynomial tractability}; 
cf., e.g., Novak \& Wo\'zniakowski \cite{NovakWozniakowski2008,NovakWozniakowski2010}.
Note, however, that lower bounds have also been established, proving that
general classes of algorithms (including deep learning methods)
cannot overcome the curse of dimensionality for all
reasonable classes of target functions; see, e.g.,
Heinrich \& Sindambiwe \cite{heinrich1999monte},
Heinrich \cite{heinrich2006randomized},
Grohs \& Voigtlaender \cite{grohs2021proof},
Petersen \& Voigtlaender \cite{petersen2018optimal},
and Yarotsky \cite{yarotsky2017error}. 

While there is not yet a comprehensive mathematical theory
explaining in a rigorous manner the encouraging empirical results
obtained so far, there
is by now a substantial body of literature -- of which we will provide a 
brief overview below --
shedding light on the \emph{expressive power}
of artificial neural networks (ANNs).
For a deep learning approximation scheme to be able to overcome the curse of dimensionality,
the class of ANNs used for approximation must be
sufficiently expressive in the sense that the number of parameters
needed to describe the approximating ANNs
grows at most polynomially
in the dimension $d$ of the domain and the reciprocal $\nicefrac1\eps$
of the prescribed approximation error $\eps$.

The present paper is a contribution to this line of research,
investigating the expressive power of deep ANNs and
in particular the question for which classes of target functions,
deep ANNs possess sufficient expressive power to
achieve approximations without the curse of dimensionality,
in the sense described above.

\paragraph{Previous research}
Initial research on the expressive power of ANNs
focused on universal approximation results, showing that 
even shallow ANNs, i.e., those with a single hidden layer,
can approximate very large classes of functions
to an arbitrary degree of accuracy as long as the
number of neurons is allowed to grow arbitrarily large. We refer
to \cite{gallant1988there,carroll1989construction,%
Cybenko1989,funahashi1989approximate,%
hornik1989multilayer} for universal approximation results 
using sigmoidal activation functions (see also
\cite{nguyen1999approximation,blum1991approximation}),
we refer to \cite{chen1995approximation,mhaskar1996neural,park1991universal}
for universal approximation results using radial basis functions
as activation functions, and
we refer to
\cite{hornik1990universal,hornik1991approximation,hornik1993some,leshno1993multilayer}
for universal approximation results using more general classes
of activation functions (see also \cite{irie1988capabilities,kidger2020universal}).
We also refer to \cite{Hanin2017,kidger2020universal}
for universal approximation results for deep ANNs
of bounded width.

First results which demonstrated that 
even shallow ANNs with sigmoidal activation functions
have the expressive power to break the curse of dimensionality
in the approximation of certain classes of functions
were obtained in \cite{barron1992neural,barron1993universal,%
barron1994approximation,jones1992simple} in the 1990s. This approach was subsequently
extended in several ways, see, e.g., \cite{donahue1997rates,KurKaiKre1997,%
makovoz1996random,makovoz1998uniform,kurkova2002comparison,%
BurgerNeubauer2001,KaiKurSang2012} for further results
in this mold using sigmoidal activation functions 
and \cite{breiman1993hinging,klusowski2018approximation,%
kainen2009complexity,MhaskarMicchelli1994,MhaskarMicchelli1995,%
caragea2020neural}
for similar results using other activation functions. The surveys
\cite{Ellacott1994,pinkus1999approximation} provide an overview
of the state-of-the-art of this research in the 1990s.

We also refer, e.g., to
\cite{lavretsky2002geometric,kurkova2008geometric,%
elbraechter2021deep,perekrestenko2018universal,wang2018exponential}
for upper bounds on the number of parameters necessary to achieve 
particular approximations with shallow ANNs that are polylogarithmic in the
reciprocal of the desired approximation error.
Furthermore,
\cite{maiorov2000near,devore1997approximation,guhring2019error,shen2019nonlinear}
prove upper bounds on the number
of parameters necessary to achieve particular
approximations with shallow ANNs
that suffer from the curse of dimensionality.

In the opposite direction, even in the 1990s, several authors
studied the limitations of shallow ANNs, see, e.g.,
\cite{ChuiMhaskar1994,candes1998ridgelets,maiorov2000near,schmitt2000lower}.
More recently, a large number of results have demonstrated that
deep ANNs overcome certain limitations
of shallow ANNs and in particular,
in various circumstances have the capacity to achieve
approximations
using significantly fewer parameters than shallow ANNs with the
same accuracy would need. We refer, e.g., to
\cite{daniely2017depth,almira2021negative,poggio2017why,%
chui2019deep,eldan2016power,safran17a} for results
comparing the expressive power of single-hidden layer
ANNs to that of ANNs with two or more
hidden layers and we refer, e.g., to \cite{elbraechter2021deep,%
GrohsIbrgimovJentzen2021,ChenWu2019,MhaskarPoggio2016,safran17a,%
cohen2016expressive}
for more general results demonstrating the superiority
of deeper ANNs in certain circumstances.

Perhaps surprisingly, it has been shown in Maiorov \& Pinkus \cite{MaioPinkus1999} 
that there
exists an analytic sigmoidal activation function such that any
continuous function on the unit cube in $\R^d$ can be approximated
with arbitrary precision by a two hidden layer network with $3d$ neurons
in the first hidden layer and $6d+3$ neurons in the second hidden layer
using this special activation function (for similar results,
see also, e.g., \cite{GuliIsm2018a,GULIYEV2018296}).
While this breaks the curse of dimensionality in a certain sense, the employed activation
function is pathological and not useful for practical purposes.

In recent years, much research has thus focused on finding
upper bounds on the number of parameters necessary to achieve a
particular approximation using deep ANNs
with various practically relevant activation functions.
We refer, e.g., to \cite{LuShenYangZhang2020,guhring2019error,%
shen2020deep,shen2019nonlinear,petersen2018optimal,voigtlaender2019approximation,%
mhaskar1993} for approximation results using deep ANNs
that suffer from the curse of dimensionality
(see also \cite{BeckJentzenKuckuck2019,%
BolcskeiGrohsKutyniokPetersen2019OptimalApproximation}).
For results showing that deep ANNs
have the expressive power to overcome the curse of dimensionality,
the approximation of solutions of various classes of PDEs
has been a particularly active area in recent years.
For results showing that deep ANNs
with rectified linear unit (ReLU) activations have the expressive power
to approximate solutions of certain PDEs without the
curse of dimensionality, we refer, e.g., to
\cite{beneventano2020highdimensional,BernerGrohsJentzen2018,%
ElbraechterSchwab2018,GononGrohsEtAl2019,gonon2020deep,%
GrohsHerrmann2020,GrohsWurstemberger2018,grohs2019space,%
GrohsJentzenSalimova2019,JentzenSalimovaWelti2018,%
HornungJentzenSalimova2020,HutzenthalerJentzenKruseNguyen2019,%
ReisingerZhang2019,SchwabZech2019}. We also refer to
\cite{EHanJentzen2017,HanJentzenE2018,Jiequn2020AlgorithmsPDEs,%
BeckBeckerCheriditoJentzenNeufeld2019} and the
surveys \cite{Kuckuck2020overview,Jiequn2020AlgorithmsPDEs} for more practically
oriented results giving empirical indications
that deep ANNs can approximate
solutions to certain PDEs without the curse of dimensionality.

There are also a number of recent results demonstrating that deep ANNs have the
expressive power to approximate more general classes of functions,
not related to PDEs, without the curse of dimensionality; 
we refer, e.g., to \cite{lee2017ability,bach2017breaking,%
cheridito2021efficient,ShahamCloningerCoifman2018,wang2018exponential,%
caragea2020neural} for such
results using deep ANNs with ReLU activations 
and we refer to \cite{li2020better,lee2017ability,JMLR:v25:23-0912} for results
using deep ANNs with other activation functions.

Lower bounds on the number of parameters necessary for a deep ANN
with ReLU activation functions to achieve a particular approximation
have been demonstrated, e.g., in \cite{almira2021negative,ChenWu2019,%
perekrestenko2018universal,yarotsky2017error,petersen2018optimal,GrohsIbrgimovJentzen2021}.
In particular, \cite{yarotsky2017error,petersen2018optimal} show that
there are natural classes of functions which deep ANNs with
ReLU activations cannot approximate without suffering from the curse of
dimensionality (see also \cite{grohs2021proof,heinrich1999monte,heinrich2006randomized,GrohsIbrgimovJentzen2021}).

Finally, we refer to the recent survey \cite{guehring2020expressivity}
for a much more complete overview of results on the expressivity
of deep ANNs than we can give here.

\paragraph{Our result}
It is the key purpose of this article to develop 
a machinery to study the high-dimensional approximation capacities of ANNs and, in particular, 
to show that deep ANNs have the expressive power to overcome
the curse of dimensionality in the approximation of a suitable large class of functions.
In that respect, our contribution fits in with some of the recent research
mentioned above, cf., e.g., \cite{lee2017ability,bach2017breaking,cheridito2021efficient,
ShahamCloningerCoifman2018,wang2018exponential}.
In particular, the character of our results is similar to those found in
Cheridito et al.\ \cite{cheridito2021efficient} in many respects.
However, our approach is quite different. 
While the setup in Cheridito et al.\ \cite{cheridito2021efficient} is more general 
and uses the concept of catalog networks, some of our arguments are more direct and 
therefore, allow us to derive approximation results for certain target functions 
that go beyond those shown in \cite{cheridito2021efficient}.
The present article also has certain ideas and methods in common with
Beneventano et al.~\cite{beneventano2020highdimensional}.

Before we present, in \cref{Theo:introduction} below,
a slightly simplified version of our main result,
let us briefly explain the statement and 
introduce some of the notions
used therein.
The class of approximating functions used throughout this
article consists of the realizations of (fully connected feed-forward) 
ANNs. These are alternating compositions of
affine linear functions and fixed, non-linear \emph{activation functions}.
In our case, the activation functions will always be chosen as the 
multi-dimensional
rectified linear unit (ReLU) function 
$A\colon\prb{\bigcup_{d\in\N}\R^d}\to\prb{\bigcup_{d\in\N}\R^d}$, which satisfies 
for all $d\in\N$ and all 
$x=(x_1,\dots,x_d)\in\R^d$ that 
$A(x)=(\max\{x_1,0\},\dots,\max\{x_d,0\})$.
Since the activation functions are fixed,
an ANN is determined by the matrices and vectors used to specify
all of the affine linear functions appearing in its realization.
More precisely, in our formalization (cf., e.g.,
Petersen \& Voigtlaender \cite[Definition~2.1]{petersen2018optimal}
and Beck et al.\ \cite[Definition~2.9]{BeckJentzenKuckuck2019}), 
the set of neural networks is given as
\begin{equation}
    \mathbf{N}
    =
    \bigcup_{L\in\N}\bigcup_{l_0,l_1,\dots,l_L\in\N}
    \prbbb{\bigtimes_{k=1}^L (\R^{l_k\times l_{k-1}}\times\R^{l_k})}
    .
\end{equation}
For every $L\in\N$, $l_0,l_1,\dots,l_L\in\N$, 
$\mathscr f\in\prb{\bigtimes_{k=1}^L(\R^{l_k\times l_{k-1}}\times\R^{l_k})}$,
we think of $L$ as the \emph{length}\footnote{What we call the length is 
also sometimes
called the depth of the neural network in the scientific literature.} of the 
neural network $\mathscr f$
and we think for every $k\in\{0,1,\dots,L\}$ of $l_k$ as the
\emph{dimension of} (or the \emph{number of neurons in}) \emph{the $k$-th layer}
of $\mathscr f$.
Furthermore, for every
$L\in\N$, $l_0,l_1,\dots,l_L\in\N$, 
$\mathscr f=((W_1,B_1),(W_2,B_2),\dots,(W_L,B_L))\in\prb{\bigtimes_{k=1}^L(\R^{l_k\times l_{k-1}}\times\R^{l_k})}$,
we define the \emph{number of parameters}
of $\mathscr f$ as $\mathcal P(\mathscr f)=\sum_{k=1}^Ll_k(l_{k-1}+1)$
(this is the total number of entries in the \emph{weight
matrices} $W_1,W_2,\dots,W_L$ and
\emph{bias vectors} $B_1,B_2,\dots,B_L$), 
we define for every
$k\in\{1,2,\dots,L\}$
the \emph{$k$-th layer affine transformation} as the function
$\mathfrak L_k^{\mathscr f}\colon \R^{l_{k-1}}\to\R^{l_k}$
which satisfies for all $x\in\R^{l_{k-1}}$ that $\mathfrak L_k^{\mathscr f}(x)=W_kx+B_k$,
and we define the \emph{realization} of
$\mathscr f$ as the composition
\begin{equation}
    \realisation(\mathscr f)
    =
    \mathfrak L_{L}^{\mathscr f}\circ A\circ \mathfrak L_{L-1}^{\mathscr f}
        \circ A\circ\dots\circ A\circ \mathfrak L_1^{\mathscr f}
    .
\end{equation}
For an illustration of such 
a neural network we refer to \cref{figure_1}.

\def\layersep{4cm}
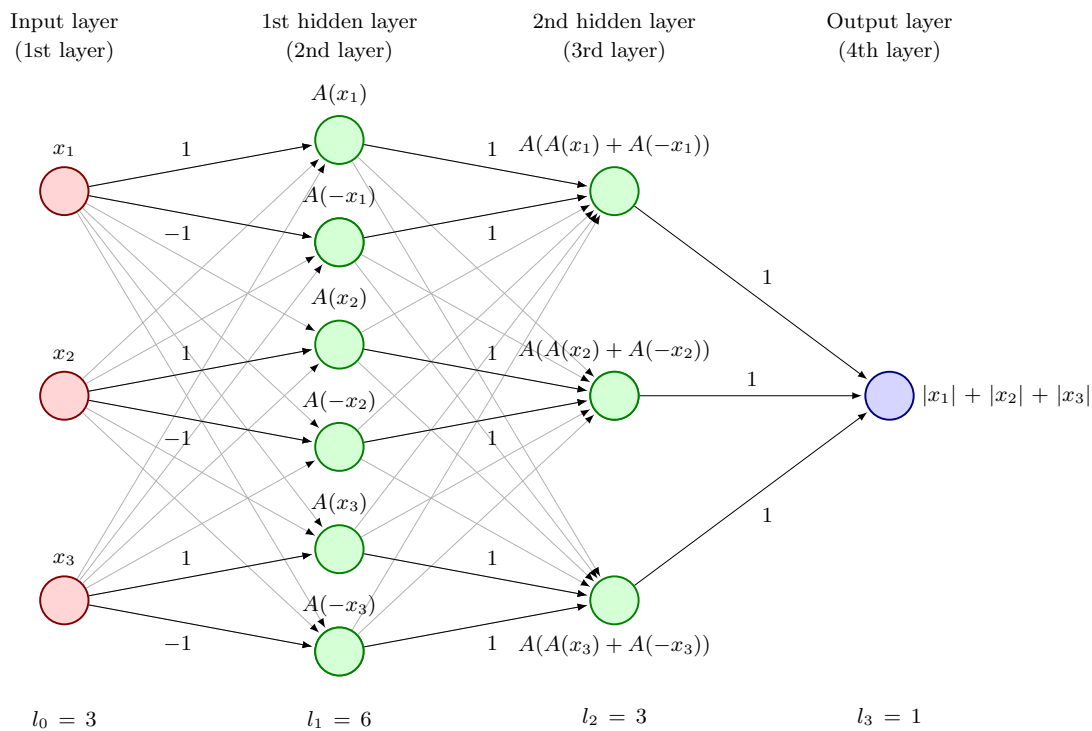
\begin{figure}
    \centering
    \begin{adjustbox}{width=\textwidth}
    \begin{tikzpicture}[shorten >=1pt,-latex,draw=black!100, node distance=\layersep,auto]
        \tikzstyle{every pin edge}=[<-,shorten <=1pt]
        \tikzstyle{neuron}=[circle,fill=black!25,draw=black!75,minimum size=20pt,inner sep=0pt,thick]
        \tikzstyle{input neuron}=[neuron,draw={rgb:red,1;black,1}, fill={rgb:red,1;white,5}];
        \tikzstyle{output neuron}=[neuron, draw={rgb:blue,1;black,1}, fill={rgb:blue,1;white,5}];
        \tikzstyle{hidden neuron}=[neuron,draw={rgb:green,1;black,1}, fill={rgb:green,1;white,5}];
        \tikzstyle{annot} = [text width=9em, text centered]
        \tikzstyle{annot2} = [text width=4em, text centered]
    
        \foreach \name / \y in {1,...,3}
            \node[input neuron,label=$x_{\y}$] (I-\name) at (0,-3*\y+2.75) {};
            \path[yshift=1.5cm]
								node[hidden neuron,label=${A}(x_1)$] (H-1) at (\layersep,-1 cm) {}
								node[hidden neuron] (Q-1) at (\layersep,-2.5 cm) {}
								node[hidden neuron,label=${A}(x_2)$] (H-2) at (\layersep,-4 cm) {}
								node[hidden neuron,label=${A}(-x_2)$] (Q-2) at (\layersep,-5.5 cm) {}
								node[hidden neuron,label=${A}(x_3)$] (H-3) at (\layersep,-7 cm) {}
								node[hidden neuron] (Q-3) at (\layersep,-8.5 cm) {};
                 \node[hidden neuron,label=${A}({A}(x_{1})+{A}(-x_{1}))$] (H2-1) at (2*\layersep,-3*1+2.75) {};
								 \node[hidden neuron] (H2-2) at (2*\layersep,-3*2+2.75) {};
				         \node[hidden neuron,label=below:${A}({A}(x_{3})+{A}(-x_{3}))$] (H2-3) at (2*\layersep,-3*3+2.75) {};
           \node[output neuron] (H3-1) at (3*\layersep,-3.25 cm) {};
					 \node[annot,right of=H3-1, node distance=1.7cm, align=center] () {$\vass{x_1}+\vass{x_2}+\vass{x_3}$};
                
				\foreach \y in {1,...,3}
							\foreach \target in {H,Q}
									\foreach \source in {1,...,3}
                \path (I-\source) edge [draw=black!30] (\target-\y);
				\foreach \y in {1,...,3}
							\foreach \source in {H,Q}
									\foreach \target in {1,...,3}
									\path (\source-\y) edge [draw=black!30] (H2-\target);
        \foreach \source in {1,...,3}
                \path (I-\source) edge node{$1$} (H-\source);
				\foreach \source in {1,...,3}
								\path (I-\source) edge [swap] node{$-1$} (Q-\source);
				\foreach \source in {1,...,3}
                \path (H-\source) edge node{$1$} (H2-\source);
				\foreach \source in {1,...,3}
								\path (Q-\source) edge [swap] node{$1$} (H2-\source);
								\path (H2-1) edge node{$1$} (H3-1);
								\path (H2-2) edge node{$1$} (H3-1);
								\path (H2-3) edge [swap] node{$1$} (H3-1);
    
				\path[yshift=1.5cm]
								node[hidden neuron,label=${A}(-x_1)$] (Q-1) at (\layersep,-2.5 cm) {}
								node[hidden neuron,label=${A}(-x_3)$] (Q-3) at (\layersep,-8.5 cm) {};
				\node[hidden neuron,label=${A}({A}(x_{2})+{A}(-x_{2}))$] (H2-2) at (2*\layersep,-3*2+2.75) {};
        \node[annot,above of=H-1, node distance=1.5cm, align=center] (hl) {1st hidden layer\\(2nd layer)};
        \node[annot,above of=H2-1, node distance=2.25cm, align=center] (hl2) {2nd hidden layer\\(3rd layer)};
        \node[annot,above of=H3-1, node distance=5.25cm, align=center] (hl3) {Output layer\\(4th layer)};
        \node[annot,left of=hl, align=center] {Input layer\\ (1st layer)};
				
				\node[annot2,below of=Q-3, node distance=1cm, align=center] (sl) {${l}_1=6$};
        \node[annot2,below of=H2-3, node distance=1.7cm, align=center] (sl2) {${l}_2=3$};
        \node[annot2,left of=sl, align=center] {${l}_0=3$};
        \node[annot2,right of=sl2, align=center] {${l}_3=1$};
    \end{tikzpicture}
    \end{adjustbox}
    \caption{\label{figure_1}Graphical illustration of an example neural network 
    which has as realizaton the $\ell^1$-norm on $\R^3$. The neural network
    has $2$ hidden layers and length $L=3$
    with $3$ neurons in the input layer (corresponding to ${l}_0 = 3$), 
    $6$ neurons in the first hidden layer (corresponding to ${l}_1 = 6$), 
    $3$ neurons in the second hidden layer (corresponding to ${l}_2 = 3$), 
    and one neuron in the output layer (corresponding to ${l}_3 = 1$). 
    In this situation we have an ANN with $39$ weights and $10$ biases 
    adding up to $\param(\mathscr{f})=49$ parameters overall of which $15$
    are nonzero (indicated by black arrows).
    The realization function 
    $\realisation(\mathscr f)\in C(\R^3,\R)$ 
    of the considered deep ANN maps each $3$-dimensional input vector 
    $x = ( x_1, x_2, x_3 ) \in \R^3$ to the 1-dimensional 
    output $(\realisation(\mathscr f))(x)=\vass{x_1}+\vass{x_2}+\vass{x_3}$.}
\end{figure}

Our goal is to show that ANNs have sufficient expressive power to
approximate certain sequences of functions
without the curse of dimensionality.
To formalize this, let
$\norm{\cdot}\colon\pr{\bigcup_{d\in\N}\R^d}\to\R$ be the function which
satisfies
for all $d\in\N$, $x\in\R^d$ that $\norm{x}=\PR{\sum_{i=1}^d (x_i)^2}^{\nicefrac12}$.
Given a sequence $(\mathfrak d_d)_{d\in\N}\subseteq\N$ of natural numbers,
a sequence $f=(f_d)_{d\in\N}\in\prb{\bigtimes_{d\in\N} C(\R^d,\R^{\mathfrak d_d})}$
of functions, and a radius $r\in[0,\infty)$,
we consider ANNs to have sufficient expressive power
to approximate the sequence $f$ on $\ell^\infty$-balls of
radius $r$ without the curse of dimensionality
if and only if 
there exists a constant $c\in\R$ such that for every 
dimension $d\in\N$ 
and every approximation error $\eps\in(0,1]$ there exists a neural network
$\mathscr f\in\ANNs$ such that
\begin{equation}
    \label{eq:approx_without_cod}
    \paramANN(\mathscr f)\leq cd^c\eps^{-c},
    \quad
    \realisation(\mathscr f)\in C(\R^d,\R^{\mathfrak d_{d}}),
    \quad\text{and}\quad
    \sup\nolimits_{[-r,r]^d}\norm{f_d(x)-(\realisation(\mathscr f))(x)}\leq \eps
    .
\end{equation}
The main result of this article,
\cref{Theo:example_multiple_composition_loclip}
in \cref{Section:7} below, demonstrates, roughly speaking, that
ANNs have sufficient expressive power to approximate, without the curse of dimensionality,
certain sequences of functions which can be constructed by
means of a finite number of compositions using
locally Lipschitz continuous functions, maxima, and products on arbitrarily large
$\ell^\infty$-balls uniformly with respect to the standard Euclidean norm.

In order to make this statement more precise, we present in
\cref{Theo:introduction} below a slightly simplified version
of our main result.

\begin{samepage}
\begin{athm}{theorem}{Theo:introduction}
	Let 
		$\norm{\cdot}\colon \prb{\bigcup_{d \in \N}\R^d} \to \R$ 
        and $A \colon \prb{\bigcup_{d \in \N}\R^d} \to \prb{\bigcup_{d \in \N}\R^d}$ 
        sa\-tis\-fy for all $d \in \N$, $x=(x_1,%
        \ldots, x_d) \in \R^d$ that
        \begin{equation}
            \norm{x}=  \PRb{\ssssum_{i=1}^d\vass{x_i}^2 }^{\nicefrac{1}{2}}
            \qquad\text{and}\qquad 
            A(x) = (\max\{x_1,0\},%
            \allowbreak \ldots, \max\{x_d,0\}),
        \end{equation}
	let
$
	\ANNs
	=
	\bigcup_{L \in \N}
	\bigcup_{l_0,%
    \ldots, l_L \in \N }
	\prb{
	\bigtimes_{k = 1}^L (\R^{l_k \times l_{k-1}} \times \R^{l_k})
	}
$,
let 
    $\realisation \colon \ANNs \to \prb{\bigcup_{k,l\in\N}\,C(\R^k,\R^l)}$ and
    $\paramANN \colon \ANNs \to \N$
satisfy for all 
	$ L\in\N$, 
	$l_0,l_1,
    \ldots, l_L \in \N$, 
	$
	\mathscr{f} 
	=
	((W_1, B_1),%
    \allowbreak \ldots, (W_L,\allowbreak B_L))
	\in  \allowbreak
	\prb{ \bigtimes_{k = 1}^L\allowbreak(\R^{l_k \times l_{k-1}} \times \R^{l_k})}
	$,
	$x_0 \in \R^{l_0},\, x_1 \in \R^{l_1}, 
    \ldots,\, \allowbreak x_{L} \in \R^{l_{L}}$ 
with 
	$\forall \, k \in \{1,%
    \ldots,L\} \colon\allowbreak x_k =A(W_k x_{k-1} + B_k)$  
that 
  \begin{align}
    \paramANN(\mathscr{f}) = \ssssum_{k = 1}^L l_k(l_{k-1} + 1),\;\;
	\realisation(\mathscr{f}) \in C(\R^{l_0},\R^{l_L}),\;\;\text{and}\;\;
	( \realisation(\mathscr{f}) ) (x_0) = W_L x_{L-1} + B_L,
	\end{align}
	let $r,R,n \in \N$, $a_1,a_2, 
    \ldots, a_n \in \N_0 \cup \{-1\}$, and let
	$f_{k,d} \colon \R \to \R$, $k, d \in \N$, and
	$F_{k} \colon \allowbreak \prb{\bigcup_{d \in \N} \R^d} \allowbreak\to \prb{\bigcup_{d \in \N} \R^d}$, $k \in \N_0 \cup \{-1\}$,
	satisfy for all 
        $k, d \in \N$, 
        $x,y \in \R$, 
        $v = (v_1,%
        \ldots, v_d) \in \R^d$ 
    that 
        $\vass{f_{k,d}(x)} \leq 1 \leq a_1$,
	    $\vass{f_{k,d}(x) - f_{k,d}(y)} \leq r(1 + \vass{x} + \vass{y})^{r} \vass{x-y}$,
	    $F_{-1}(v) = \prb{v_1, \allowbreak \max\{v_1,v_2\}, \allowbreak \ldots, \allowbreak \max\{v_1, %
        \ldots, v_d\} }$, 
	    $F_{0} (v) = \prb{v_1, \allowbreak v_1 v_2, \allowbreak \ldots, \allowbreak v_1 v_2\cdots v_d }$,
	and
	$F_{k}(v) = \prb{f_{k,1}(v_1), %
    \allowbreak \ldots, f_{k,d}(v_d)}$.
	Then there exist 
    $(\mathscr{F}_{d,\varepsilon})_{(d, \varepsilon) \in \N \times (0,1]} \subseteq \ANNs$ 
    and $c \in \R$ 
    such that for all 
    $d \in \N$, 
    $\varepsilon \in (0,1]$ 
    it holds that
	$\paramANN(\mathscr{F}_{d,\varepsilon}) \leq cd^c\varepsilon^{-c}$, 
    $\realisation(\mathscr{F}_{d,\varepsilon}) \in C(\R^d, \R^d)$, and
	\begin{equation}
	\label{intro:1}
    \textstyle
	    \sup_{x \in [-R,R]^d}\norm{ (F_{a_n} \circ \ldots %
        \circ F_{a_1})(x) - (\realisation(\mathscr{F}_{d,\varepsilon}))(x) } 
        \leq 
        \varepsilon
        .
	\end{equation}
\end{athm}
\end{samepage}

\cref{Theo:introduction} above is an immediate consequence of 
\cref{Coro:example_multiple_composition_loclip_2bis} 
in \cref{Subsection:7.2}. 
\cref{Coro:example_multiple_composition_loclip_2bis} follows from 
\cref{Theo:example_multiple_composition_loclip}
in \cref{Subsection:7.2}, the main result of this article.

As described above, \cref{Theo:introduction} concerns the approximation of
certain sequences of functions by (realizations of) ANNs.
For convenience of notation, these sequences of functions are
formalized as continuous functions from the union $\bigcup_{d\in\N}\R^d$
to itself, where $\bigcup_{d\in\N}\R^d$ is equipped with the topology of
the disjoint union. 
In this setting, the functions that are to be approximated
can be written as compositions of three types of functions.
The first type consists of the \emph{maximum function}
from $\bigcup_{d\in\N}\R^d$ to itself
which for every $d\in\N$ maps $v=(v_1,\dots,v_d)\in\R^d$ to
$(v_1,\max\{v_1,v_2\},\dots,\max\{v_1,v_2,\dots,v_d\})$;
the second type consists of the \emph{product function}
from $\bigcup_{d\in\N}\R^d$ to itself
which for every $d\in\N$ maps $v=(v_1,\dots,v_d)\in\R^d$ to
$(v_1,v_1v_2,\dots,v_1v_2\cdots v_d)$;
and the third type consists of functions
that are componentwise applications of
certain bounded, locally Lipschitz continuous
functions from $\R$ to $\R$.

In order to exhibit the scope of this result,
let us illustrate \cref{Theo:introduction} 
by means of several examples.
We note that \cref{Theo:introduction} implies that
deep ANNs have sufficient expressivity to
approximate all of the following sequences of functions\footnote{Note that
$\arcsin\in C([-1,1],\R)$, $\arctan\in C(\R,\R)$,
and
$\tanh\in C(\R,\R)$ are the unique continuous functions which 
satisfy for all $x\in(-\tfrac\pi2,\tfrac\pi2)$, $y\in\R$ that
$\arcsin(\sin(x))=x$,
$\arctan(\tan(x))=x$, and
$\tanh(y)=\frac{e^y-e^{-y}}{e^y+e^{-y}}$.}
without the curse of dimensionality on arbitrarily
large $\ell^\infty$-balls in the sense of \eqref{eq:approx_without_cod} above:
\begin{gather}
  \R^d \ni ( x_1, x_2, \dots, x_d ) \mapsto 
  \sin\prb{\textstyle{\prod\nolimits_{j=1}^{d}} \sin( x_j )} \in \R , \quad d \in \N , \label{eq:intro_ex1}
  \\
  \R^d \ni ( x_1, x_2, \dots, x_d ) \mapsto 
  \prb{\sssprod_{j=1}^{d} e^{-\nicefrac{\vass{x_j}^2}{2}}} \in \R , \quad d \in \N , \label{eq:intro_ex2}
  \\
  \R^d \ni ( x_1, x_2, \dots, x_d ) \mapsto 
  \tanh\prb{ \sssprod_{j=1}^{d} \arcsin\prb{ \tfrac{x_j}{1+\vass{x_j}^2} } } \in \R, \quad d \in \N, \label{eq:intro_ex3}
  \\
  \R^d \ni ( x_1, x_2, \dots, x_d ) \mapsto 
  \prb{ \sssprod_{j=1}^{\lfloor d/2 \rfloor} \PRb{ \tfrac{x_{2j} \arctan(x_{2j-1})}{1+\vass{x_{2j}}^2} } } \in \R, \quad d \in \N, \label{eq:intro_ex4}
	 \\
  \R^d \ni ( x_1, x_2, \dots, x_d ) \mapsto 
  \max\pR{\cos\pr{\vass{x_1}^2},\cos\pr{\vass{x_2}^2},\ldots,\cos\pr{\vass{x_d}^2}} \in \R, \quad d \in \N, \label{eq:intro_ex5}
	 \\
  \R^d \ni ( x_1, x_2, \dots, x_d ) \mapsto 
  \sssprod_{j=1}^{d}\max\{\sin(x_1),\sin(\tfrac{x_2}2),\ldots,\sin(\tfrac{x_j}j)\} \in \R, \quad d \in \N, \label{eq:intro_ex6}
	\\
  \R^d \ni ( x_1, x_2, \dots, x_d ) \mapsto 
  \prb{\sssprod_{j=1}^{d}(1+e^{(-1)^jx_j})}^{-1} \in \R, \quad d \in \N, \label{eq:intro_ex7}
  \\
  \R^d\ni(x_1,x_2,\dots,x_d)\mapsto
    \max\nolimits_{j\in\{1,2,\dots,d\}}\cos\prb{\sssprod_{i=1}^j \tanh(x_i)}
  \in \R,\quad d\in\N, \label{eq:intro_ex8}
  \\
  \R^d\ni(x_1,x_2,\dots,x_d)\mapsto
  \max\nolimits_{j\in\{1,2,\dots,d\}}\PRb{\sssprod_{i=1}^j\prb{\max\nolimits_{k\in\{1,2,\dots,i\}} \tfrac{1-\abs{x_k}^4}{1+\abs{x_k}^4}}}
  \in\R,\quad d\in\N, \label{eq:intro_ex9}
  \\
  \R^d\ni(x_1,x_2,\dots,x_d)\mapsto
  \sssprod_{j=1}^d\PRb{\max\nolimits_{i\in\{1,2,\dots,j\}}\prb{\sssprod_{k=1}^i \tfrac{(x_k)^3}{1+\abs{x_k}^3}}}
  \in\R,\quad d\in\N, \label{eq:intro_ex10}
  \\
  \R^d\ni(x_1,x_2,\dots,x_d)\mapsto
    \max\nolimits_{j\in\{1,2,\dots,d\}}\PRb{\pr{-1}^j\prb{\sssprod_{i=1}^{j} \cos(\vass{x_i}^4)}}
  \in \R,\quad d\in\N, \label{eq:intro_ex11}
  \\
  \R^d\ni(x_1,x_2,\dots,x_d)\mapsto
    \sssprod_{k=1}^{d}\PRb{\max\nolimits_{j\in\{1,2,\dots,k\}}\prb{\PRb{\sin\prb{\sssprod_{i=1}^{j} \tfrac{x_i}{1+\vass{x_i}}}}^3}}
  \in \R,\quad d\in\N \label{eq:intro_ex12}
  .
\end{gather}
For all of these examples, the fact that deep ANNs have sufficient expressivity to
approximate the sequence in question
without the curse of dimensionality is a direct consequence of
\cref{Coro:example_multiple_composition_loclip_proj} 
in \cref{Subsection:7.2} below,
which in
turn follows easily from \cref{Theo:introduction}.
We would like to emphasize that these sequences of functions, while illustrating by way of explicit instances how the abstract statement in \cref{Theo:introduction} may be applied concretely, are merely academic examples and do not necessarily reflect any real-world relevance. 

We leave it as a promising avenue for future research to explore the applications of the findings presented in this article to real-world problems. In this context, it is also important to point out that beyond \cref{Theo:introduction} above, the main contribution of our work in this article is to establish an extensible framework for proving similar results regarding the approximability of functions by ANNs without the curse of dimensionality (indeed, \cref{Theo:introduction} follows rather trivially from the more general construction principles we establish). In particular, we remark that extending \cref{Theo:introduction} to allow, e.g., for compositions with certain linear functions might be a useful addition to consider for future work.

The remainder of this article is organized as follows.
In \cref{Section:2} we introduce our formalization of ANNs and
operations on ANNs and prove or recall from the literature 
the relevant fundamental results needed later.
In \cref{Section:3} we define the cost of an ANN approximation
and ANN approximation spaces and we show how the cost
behaves under certain operations on functions, including composition.
In \cref{Section:4} we construct ANN approximations for locally Lipschitz
continuous function from $\R$ to $\R$.
In \cref{Section:5} we show how to represent the multi-dimensional maximum functions
efficiently with ReLU networks.
In \cref{Section:6} we develop efficient ANN approximations for the
multi-dimensional product functions.
In \cref{Section:7} we combine the results from 
\cref{Section:3,Section:4,Section:5,Section:6}
to obtain the central results of this article,
including \cref{Theo:introduction} above.

\section{Artificial neural network (ANN) calculus}
\label{Section:2}

In this section we present the formalism and theory of ANNs that we use 
throughout this article.
Though most of this material has appeared, sometimes in slightly
different form, in previous works (cf., e.g,
\cite{petersen2018optimal,grohs2019space,GrohsJentzenSalimova2019,%
BeckJentzenKuckuck2019,ElbraechterSchwab2018,cheridito2021efficient})
and the results are essentially elementary,
we recall them here, since these definitions and results are
used extensively throughout the
rest of this article.

In particular,
\cref{def:ANN} below is a slightly shortened version of, 
e.g., Grohs et al.\ \cite[Definition~2.1]{grohs2019space},
\cref{def:multidim_version} is, e.g., 
Grohs et al.\ \cite[Definition~2.2]{grohs2019space},
\cref{def:ANNrealization} is, e.g., 
Grohs et al.\ \cite[Definition~2.3]{grohs2019space},
\cref{def:ReLU} is, e.g., 
Beck et al.\ \cite[Definition~2.4]{BeckJentzenKuckuck2019},
\cref{subsec:parallelization} corresponds to
Grohs et al.~\cite[Subsection~2.3.1]{grohs2019space},
\cref{def:ANNcomposition} is, e.g., 
Petersen \& Voigtlaender \cite[Definition~2.2]{petersen2018optimal},
\cref{Lemma:CompositionAssociative} is, e.g., 
Grohs et al.\ \cite[Lemma~2.8]{grohs2019space},
\cref{Lemma:PropertiesOfCompositions_n2,lem:dimcomp} are an extension of, 
e.g., Grohs et al.\ \cite[Proposition~2.6]{grohs2019space},
\cref{Subsection:2.5} corresponds to 
Grohs et al.\ \cite[Subsection~2.2.4]{grohs2019space}, and
\cref{subsection:2.6} is based on 
Grohs et al.\ \cite[Subsection~2.3.2]{grohs2019space}.

\subsection{Definition of ANNs}

In this subsection we introduce the formalization of the set of ANNs and associated essential notations that will be used throughout this article. Note that, in \cref{def:ANN} below, we think of $\ANNs$ as the set of ANNs and for every $\mathscr{f}\in\ANNs$, $n\in\N_0$ we think of $\paramANN(\mathscr f)$, $\lengthANN(\mathscr f)$, $\inDimANN(\mathscr f)$, $\outDimANN(\mathscr f)$, $\dims(\mathscr f)$, and $\singledims_n(\mathscr f)$ as the number of parameters of $\mathscr f$, the length of $\mathscr f$, the input dimension of $\mathscr f$, the output dimension of $\mathscr f$, the dimension vector of $\mathscr f$, and the dimension of the $n$-th layer of $\mathscr{f}$, respectively.

\begin{definition}[Set of ANNs]
	\label{def:ANN}
	We denote by $\ANNs$ the set given by 
	\begin{equation}
	\begin{split}
	\ANNs
	&=
	\adjustlimits
	\bigcup_{L \in \N}
	\bigcup_{ (l_0,l_1,\ldots, l_L) \in \N^{L+1} }
	\prbbb{
	\bigtimes_{k = 1}^L \pr{\R^{l_k \times l_{k-1}} \times \R^{l_k}}
	},
	\end{split}
	\end{equation}
	we denote by
	$\paramANN \colon \ANNs \to \N$, 
	$\lengthANN \colon \ANNs \to \N$,
	$\inDimANN \colon \ANNs \to \N$,
	$\outDimANN \colon \ANNs \to \N$,
	$\dims\colon\ANNs\to  \prb{\bigcup_{L=1}^\infty\, \N^{L+1}}$, 
	and $\singledims_n \colon \ANNs \to \N_0$, $n \in \N_0$,
	the functions which satisfy
	for all $L\in\N$, $l_0,l_1,\ldots,\allowbreak l_L \in \N$, 
	$
	\mathscr{f} 
	\in  \allowbreak
	\prb{ \bigtimes_{k = 1}^L\allowbreak(\R^{l_k \times l_{k-1}} \times \R^{l_k})}$,
	$n \in \N_0$
	that
	$\paramANN(\mathscr{f}) =
	\sum_{k = 1}^L l_k(l_{k-1} + 1) $,
	$\lengthANN(\mathscr{f})=L$, 
	$\inDimANN(\mathscr{f})=l_0$,
	$\outDimANN(\mathscr{f})=l_L$,
	$\dims(\mathscr{f})= (l_0,l_1,\ldots, l_L)$, and
	\begin{equation}
	    \singledims_n(\mathscr{f})= \begin{cases}
	    l_n &\colon n\leq L\\
	    0   &\colon n>L.
	    \end{cases}
	\end{equation}
\end{definition}
\cfclear

\begin{definition}[ANNs]
	\label{def:neuralnetwork}
	We say that $\mathscr{f}$ is an artificial neural network (we say that $\mathscr{f}$
    is an ANN) if and only
	if it holds that $\mathscr{f}\in\ANNs$.
\end{definition}

\subsection{Realizations of ANNs}

In this subsection we recall the crucial notion of the realization of an ANN.

\begin{definition}[Multi-dimensional versions]\label{def:multidim_version}
	Let $a \colon \R \to \R$ be a function.
	Then we denote by 
        $\Mult{a}  \colon \prb{\bigcup_{d \in \N}\R^d} \to \prb{\bigcup_{d \in \N}\R^d}$ 
    the function which satisfies for all 
        $d \in \N$, 
        $ x = ( x_1, x_2, \dots, x_{d} ) \in \R^{d} $ 
    that
	\begin{equation}\label{multidim_version:Equation}
	\Mult{a}(x)
	=
	\pr*{
	a(x_1),
	a(x_2)
	,
	\ldots
	,
	a(x_d)
	}.
	\end{equation}
\end{definition}
\cfclear

\begin{definition}[Realizations associated to ANNs]
	\label{def:ANNrealization}
	\cfconsiderloaded{def:ANNrealization}
	Let $a\in C(\R,\R)$.
	Then we denote by 
	$
	\functionANN \colon \ANNs \to \prb{\bigcup_{k,l\in\N}\,C(\R^k,\R^l)}
	$
	the function which satisfies
	for all  $ L\in\N$, $l_0,l_1,\ldots, l_L \in \N$, 
	$
	\mathscr{f} 
	=
	((W_1, B_1),(W_2, B_2),\allowbreak \ldots, (W_L,\allowbreak B_L))
	\in  \allowbreak
	\prb{\bigtimes_{k = 1}^L\allowbreak(\R^{l_k \times l_{k-1}} \times \R^{l_k})}
	$,
	$x_0 \in \R^{l_0}, x_1 \in \R^{l_1}, \ldots, x_{L} \in \R^{l_{L}}$ 
	with $\forall \, k \in \{1,2,\ldots,L\} \colon x_k =\Mult{a}(W_k x_{k-1} + B_k)$  
	that
	\begin{equation}
	\label{ANNrealization:ass2}
	\functionANN(\mathscr{f}) \in C(\R^{l_0},\R^{l_L})\qandq
	( \functionANN(\mathscr{f}) ) (x_0) = W_L x_{L-1} + B_L
	\end{equation}
	\cfload.
\end{definition}
\cfclear

\begin{definition}[ReLU activation function]
\label{def:ReLU}
\cfconsiderloaded{def:ReLU}
	We denote by $\ReLU \colon \R \to \R$ the function which satisfies for all $x \in \R$ that $\ReLU(x) = \max\{x,0\}$.
\end{definition}
\cfclear

\subsection{Parallelizations of ANNs with the same length}
\label{subsec:parallelization}

In this subsection we recall the notion of parallelization of ANNs. Roughly speaking, the parallelization of ANNs corresponds, on the level of realizations, to the parallelization of functions, where we consider for every $n\in\N$, $k_1,k_2,\dots,k_n, l_1, l_2,\dots, l_n\in\N$ and $f_1\in C(\R^{k_1},\R^{l_1})$, $f_2\in C(\R^{k_2},\R^{l_2})$, $\dots$, $f_n\in C(\R^{k_n},\R^{l_n})$ the function $\R^{k_1+k_2+\dots+k_n}\cong \R^{k_1}\times \R^{k_2}\times \cdots \times \R^{k_n}\ni (x_1,x_2,\dots,x_n)\mapsto (f_1(x_1), f_2(x_2), \dots, f_n(x_n))\in \R^{l_1}\times \R^{l_2}\times \cdots \times \R^{l_n}\cong \R^{l_1+l_2+\dots + l_n}$ as the parallelization of $f_1, f_2, \dots, f_n$ (see \cref{Lemma:PropertiesOfParallelizationEqualLength} below for a precise statement). While we consider only parallelizations of ANNs with the same length here, a generalization without this restriction will be introduced in \cref{subsection:2.6} below.

\begin{definition}[Parallelization of ANNs with the same length]
	\label{def:simpleParallelization}
	\cfconsiderloaded{def:simpleParallelization}
	Let $n\in\N$. Then we denote by 
	\begin{equation}
		\parallelizationSpecial_{n}\colon \pRb{(\mathscr{f}_1,\mathscr{f}_2,\dots, \mathscr{f}_n)\in\ANNs^n\colon \lengthANN(\mathscr{f}_1)= \lengthANN(\mathscr{f}_2)=\ldots =\lengthANN(\mathscr{f}_n) }\to \ANNs
	\end{equation}
	the function which satisfies for all  $L\in\N$,
	$(l_{1,0},l_{1,1},\dots, l_{1,L}), (l_{2,0},l_{2,1},\dots, l_{2,L}),\dots,\allowbreak (l_{n,0},\allowbreak l_{n,1},\allowbreak\dots, l_{n,L})\in\N^{L+1}$, 
	$\mathscr{f}_1=((W_{1,1}, B_{1,1}),(W_{1,2}, B_{1,2}),\allowbreak \ldots, (W_{1,L},\allowbreak B_{1,L}))\in \prb{ \bigtimes_{k = 1}^L\allowbreak(\R^{l_{1,k} \times l_{1,k-1}} \times \R^{l_{1,k}})}$, 
		$\mathscr{f}_2=((W_{2,1}, B_{2,1}),(W_{2,2}, B_{2,2}),\allowbreak \ldots, (W_{2,L},\allowbreak B_{2,L}))\in \prb{ \bigtimes_{k = 1}^L\allowbreak(\R^{l_{2,k} \times l_{2,k-1}} \times \R^{l_{2,k}})}$,
\dots,  
		$\mathscr{f}_n=((W_{n,1}, B_{n,1}),(W_{n,2}, B_{n,2}),\allowbreak \ldots, (W_{n,L},\allowbreak B_{n,L}))\in \prb{ \bigtimes_{k = 1}^L\allowbreak(\R^{l_{n,k} \times l_{n,k-1}} \times \R^{l_{n,k}})}$
	that
	\begin{equation}\label{parallelisationSameLengthDef}
	\begin{alignedat}{2}
	\parallelizationSpecial_{n}(\mathscr{f}_1,\mathscr{f}_2,\dots,\mathscr{f}_n)
	&=
	\Vast(
	&&\pr*{\begin{pmatrix}
		W_{1,1}& 0& 0& \cdots& 0\\
		0& W_{2,1}& 0&\cdots& 0\\
		0& 0& W_{3,1}&\cdots& 0\\
		\vdots& \vdots&\vdots& \ddots& \vdots\\
		0& 0& 0&\cdots& W_{n,1}
		\end{pmatrix} ,\begin{pmatrix}B_{1,1}\\B_{2,1}\\B_{3,1}\\\vdots\\ B_{n,1}\end{pmatrix}},
	\\&&&
		\pr*{\begin{pmatrix}
		W_{1,2}& 0& 0& \cdots& 0\\
		0& W_{2,2}& 0&\cdots& 0\\
		0& 0& W_{3,2}&\cdots& 0\\
		\vdots& \vdots&\vdots& \ddots& \vdots\\
		0& 0& 0&\cdots& W_{n,2}
		\end{pmatrix} ,\begin{pmatrix}B_{1,2}\\B_{2,2}\\B_{3,2}\\\vdots\\ B_{n,2}\end{pmatrix}}
	,\dots,
	\\&&&
		\pr*{\begin{pmatrix}
	W_{1,L}& 0& 0& \cdots& 0\\
	0& W_{2,L}& 0&\cdots& 0\\
	0& 0& W_{3,L}&\cdots& 0\\
	\vdots& \vdots&\vdots& \ddots& \vdots\\
	0& 0& 0&\cdots& W_{n,L}
	\end{pmatrix} ,\begin{pmatrix}B_{1,L}\\B_{2,L}\\B_{3,L}\\\vdots\\ B_{n,L}\end{pmatrix}}
	\Vast)
	\end{alignedat}
	\end{equation}
	\cfout[.]
\end{definition}
\cfclear

\begin{athm}{lemma}{Lemma:ParallelizationElementary}
\cfconsiderloaded{Lemma:ParallelizationElementary}
	Let $n,L\in\N$, $\mathscr{f}_1,\mathscr{f}_2,\ldots,\mathscr{f}_n\in\ANNs$ satisfy for all $j\in\{1,2,\dots,n\}$ that $\lengthANN(\mathscr{f}_j)=L$
	\cfload[.]Then
	\begin{equation}\label{ParallelizationElementary:Display}
	\parallelizationSpecial_{n}(\mathscr{f}_1,\mathscr{f}_2,\dots,\mathscr{f}_n)\in 
	\prbbb{\bigtimes_{k = 1}^L\allowbreak\prb{\R^{\PR{\sum_{j=1}^n\singledims_k(\mathscr{f}_j)} \times \PR{\sum_{j=1}^n\singledims_{k-1}(\mathscr{f}_j)}} \times \R^{\PR{\sum_{j=1}^n\singledims_k(\mathscr{f}_j)}}}}
	\end{equation}
	\cfout[.]
\end{athm}

\begin{proof}[Proof of \cref{Lemma:ParallelizationElementary}]	
	Note that \eqref{parallelisationSameLengthDef} ensures that for all $k\in\{0,1,\dots, L\}$ it holds that  
	$\singledims_k(\parallelizationSpecial_{n}(\mathscr{f}_1,\mathscr{f}_2,\dots,\mathscr{f}_n))=\sum_{j=1}^n\singledims_k(\mathscr f_j)$ \cfload. Hence, we obtain \eqref{ParallelizationElementary:Display}.
	This completes the proof of \cref{Lemma:ParallelizationElementary}.
\end{proof}
\cfclear

\begin{athm}{prop}{Lemma:PropertiesOfParallelizationEqualLength}
	Let $a\in C(\R,\R)$, 
	$n,L\in\N$, 
	$\mathscr{f}=(\mathscr{f}_1,\mathscr{f}_2,\allowbreak\dots,\allowbreak \mathscr{f}_n)\in\ANNs^n$
	satisfy for all $j\in\{1,2,\dots,n\}$ that $\lengthANN(\mathscr{f}_j)=L$
  \cfload.
	Then
	\begin{enumerate}[(i)]
		\item\label{PropertiesOfParallelizationEqualLength:ItemOne} it holds that 
		\begin{equation}
\functionANN(\parallelizationSpecial_{n}(\mathscr{f}))\in C\prb{\R^{[\sum_{j=1}^n \inDimANN(\mathscr{f}_j)]},\R^{[\sum_{j=1}^n \outDimANN(\mathscr{f}_j)]}}
		\end{equation}
		and
		\item\label{PropertiesOfParallelizationEqualLength:ItemTwo} it holds for all    $x_1\in\R^{\inDimANN(\mathscr{f}_1)},x_2\in\R^{\inDimANN(\mathscr{f}_2)},\dots, x_n\in\R^{\inDimANN(\mathscr{f}_n)}$ that 
			\begin{equation}\label{PropertiesOfParallelizationEqualLengthFunction}
			\begin{split}
			&\prb{\functionANN\prb{\parallelizationSpecial_{n}(\mathscr{f}) }}(x_1,x_2,\dots, x_n) 
			\\&=\prb{(\functionANN(\mathscr{f}_1))(x_1), (\functionANN(\mathscr{f}_2))(x_2),\dots,
			(\functionANN(\mathscr{f}_n))(x_n) }
			\end{split}
			\end{equation}
	\end{enumerate}
	\cfout.
\end{athm}
\begin{proof}[Proof of \cref{Lemma:PropertiesOfParallelizationEqualLength}]
Note that \eqref{ANNrealization:ass2} and \cref{Lemma:ParallelizationElementary} imply 
\cref{PropertiesOfParallelizationEqualLength:ItemOne}. 
Observe that \eqref{ANNrealization:ass2} and \eqref{parallelisationSameLengthDef} establish 
\cref{PropertiesOfParallelizationEqualLength:ItemTwo}.
The proof of \cref{Lemma:PropertiesOfParallelizationEqualLength} is thus complete.
\end{proof}
\cfclear

\begin{athm}{prop}{Lemma:PropertiesOfParallelizationEqualLengthDims}
	Let 
	$n,L\in\N$, 
	$\mathscr{f}=(\mathscr{f}_1,\mathscr{f}_2,\allowbreak\dots,\allowbreak \mathscr{f}_n)\in\ANNs^n$,
	$(l_{1,0},l_{1,1},\dots, l_{1,L}),\allowbreak (l_{2,0},l_{2,1},\allowbreak\dots,\allowbreak l_{2,L}),\allowbreak \dots, (l_{n,0},l_{n,1},\dots, l_{n,L})
	\in\N^{L+1}$ satisfy for all $j\in\{1,2,\dots,n\}$ that 
	$\dims(\mathscr{f}_j)=(l_{j,0},l_{j,1},\dots,\allowbreak l_{j,L})$
	\cfload.
	Then
	\begin{enumerate}[(i)]
		\item \label{PropertiesOfParallelizationEqualLengthDims:Dims} it holds that 
		$\dims\prb{\parallelizationSpecial_{n}(\mathscr{f})}=\smallsum_{j=1}^n \dims(\mathscr{f}_j)=\prb{\smallsum_{j=1}^n l_{j,0}, \smallsum_{j=1}^n l_{j,1},\dots, \smallsum_{j=1}^n l_{j,L}}$ and
		\item \label{PropertiesOfParallelizationEqualLengthDims:Params}
		it holds that
		\begin{equation}
		\paramANN(\parallelizationSpecial_{n}(\mathscr{f}))\le \tfrac{1}{2} \PRb{\smallsum\nolimits_{j=1}^n \paramANN(\mathscr{f}_j)}^2
		\end{equation}
	\end{enumerate}
	\cfout.
\end{athm}

\begin{proof}[Proof of \cref{Lemma:PropertiesOfParallelizationEqualLengthDims}]	

Note that \cite[Proposition~2.20]{grohs2019space} establishes 
\cref{PropertiesOfParallelizationEqualLengthDims:Dims,PropertiesOfParallelizationEqualLengthDims:Params}. 
The proof of \cref{Lemma:PropertiesOfParallelizationEqualLengthDims} is thus complete.
\end{proof}
\cfclear

\begin{athm}{cor}{Lemma:ParallelizationImprovedBoundsOne}
	Let   $n\in\N$, $\mathscr{f}=( \mathscr{f}_1,\mathscr{f}_2,\dots,\allowbreak \mathscr{f}_n)\in\ANNs^n$ satisfy  
	$\dims(\mathscr{f}_1)=\dims(\mathscr{f}_2)=\ldots=\dims(\mathscr{f}_n)$
	\cfload.	
	Then $\paramANN(\parallelizationSpecial_{n}(\mathscr{f}))\le n^2 \paramANN(\mathscr{f}_1)$
\cfout.
\end{athm}

\begin{proof}[Proof of \cref{Lemma:ParallelizationImprovedBoundsOne}]
    
    Observe that
    \cite[Corollary~2.21]{grohs2019space} establishes
    $\paramANN(\parallelizationSpecial_{n}(\mathscr{f}))\le n^2 \paramANN(\mathscr{f}_1)$. The proof of
    \cref{Lemma:ParallelizationImprovedBoundsOne}
    is thus complete.
\end{proof}
\cfclear

\subsection{Compositions of ANNs}
\label{subsec:compositions}

In this subsection we recall the notion of the composition of ANNs, which, roughly speaking, corresponds to the composition of functions on the level of realizations (see \cref{Lemma:PropertiesOfCompositions_n2} below for a precise statement). After various technical lemmas on basic properties of the composition of ANNs, we provide an upper bound for the number of parameters of the composition of two ANNs in \cref{Lemma:PropertiesOfCompositions_n1,Lemma:PropertiesOfCompositions_n3} below, which will be important for our main task.

\begin{definition}[Standard composition of ANNs]
	\label{def:ANNcomposition}
	\cfconsiderloaded{def:ANNcomposition}
	We denote by $\compANN{(\cdot)}{(\cdot)}\colon\allowbreak \{(\mathscr{f}_1,\mathscr{f}_2)\allowbreak\in\ANNs\times \ANNs\colon \inDimANN(\mathscr{f}_1)=\outDimANN(\mathscr{f}_2)\}\allowbreak\to\ANNs$ the function which satisfies for all 
	$ L,\mathfrak{L}\in\N$, $l_0,l_1,\ldots, l_L, \allowbreak\mathfrak{l}_0,\mathfrak{l}_1,\allowbreak\ldots, \mathfrak{l}_\mathfrak{L} \in \N$, 
	$
	\mathscr{f}_1
	=
	((W_1, B_1),(W_2, B_2),\allowbreak \ldots, (W_L,\allowbreak B_L))
	\in  \allowbreak
	\prb{ \bigtimes_{k = 1}^L\allowbreak(\R^{l_k \times l_{k-1}} \allowbreak \times \R^{l_k})}
	$,
	$
	\mathscr{f}_2
	=
	((\mathfrak{W}_1, \mathfrak{B}_1),\allowbreak(\mathfrak{W}_2, \mathfrak{B}_2),\allowbreak \ldots, (\mathfrak{W}_\mathfrak{L},\allowbreak \mathfrak{B}_\mathfrak{L}))
	\in  \allowbreak
	\prb{ \bigtimes_{k = 1}^\mathfrak{L}\allowbreak(\R^{\mathfrak{l}_k \times \mathfrak{l}_{k-1}} \times \R^{\mathfrak{l}_k})}
	$ 
	with $l_0=\inDimANN(\mathscr{f}_1)=\outDimANN(\mathscr{f}_2)=\mathfrak{l}_{\mathfrak{L}}$
	that
	\begin{equation}\label{ANNoperations:Composition}
	\begin{split}
	&\compANN{\mathscr{f}_1}{\mathscr{f}_2}=\\&
	\begin{cases} 
			\begin{array}{r}
			\prb{(\mathfrak{W}_1, \mathfrak{B}_1),(\mathfrak{W}_2, \mathfrak{B}_2),\ldots, (\mathfrak{W}_{\mathfrak{L}-1},\allowbreak \mathfrak{B}_{\mathfrak{L}-1}),
			(W_1 \mathfrak{W}_{\mathfrak{L}}, W_1 \mathfrak{B}_{\mathfrak{L}}+B_{1}),\\ (W_2, B_2), (W_3, B_3),\ldots,(W_{L},\allowbreak B_{L})}
			\end{array}
	&: L>1<\mathfrak{L} \\[3ex]
	\prb{ (W_1 \mathfrak{W}_{1}, W_1 \mathfrak{B}_1+B_{1}), (W_2, B_2), (W_3, B_3),\ldots,(W_{L},\allowbreak B_{L}) }
	&: L>1=\mathfrak{L}\\[1ex]
	\prb{(\mathfrak{W}_1, \mathfrak{B}_1),(\mathfrak{W}_2, \mathfrak{B}_2),\allowbreak \ldots, (\mathfrak{W}_{\mathfrak{L}-1},\allowbreak \mathfrak{B}_{\mathfrak{L}-1}),(W_1 \mathfrak{W}_{\mathfrak{L}}, W_1 \mathfrak{B}_{\mathfrak{L}}+B_{1}) }
	&: L=1<\mathfrak{L}  \\[1ex]
	(W_1 \mathfrak{W}_{1}, W_1 \mathfrak{B}_1+B_{1}) 
	&: L=1=\mathfrak{L} 
	\end{cases}
	\end{split}
	\end{equation}
	\cfload.
\end{definition}
\cfclear

\begin{definition}
[ReLU identity networks]
\label{def:ReLU_identity}
\cfconsiderloaded{def:ReLU_identity}
We denote by 
    $(\ReLUidANN{d})_{d \in \N} \subseteq \ANNs$ 
the \cfadd{def:neuralnetwork}ANNs which satisfy 
    for all 
        $d \in \{2,3,\dots\}$ 
    that
\begin{equation}
\label{def:ReLU_identity_d_is_one}
\begin{split}
	\ReLUidANN{1} = \pr*{ \!\pr*{\!\begin{pmatrix}
	1\\
	-1
	\end{pmatrix},
	\begin{pmatrix}
	0\\
	0
	\end{pmatrix}\! },
	\prbb{	\begin{pmatrix}
	1& -1
	\end{pmatrix}, 
	0 } \!
	 }  \in \prb{(\R^{2 \times 1} \times \R^{2}) \times (\R^{1 \times 2} \times \R^1) }
\end{split}
\end{equation}
and $\ReLUidANN{d} = \parallelizationSpecial_{d} (\ReLUidANN{1},\ReLUidANN{1},\dots, \ReLUidANN{1})$ \cfload.
\end{definition}
\cfclear

\begin{athm}{prop}{Prop:identity_representation}
Let $d \in \N$. Then
\begin{enumerate}[(i)]
    \item \label{identity_representation:1} it holds that $\dims(\ReLUidANN{d})= (d, 2d, d)$,
    \item \label{identity_representation:2} it holds for all $x \in \R^d$ that $(\realisation_\ReLU(\ReLUidANN{d}))(x) = x$, and
    \item \label{identity_representation:3} it holds that $\paramANN(\ReLUidANN{d}) = 4d^2+3d$
\end{enumerate}
\cfout.
\end{athm}

\begin{proof}
[Proof of \cref{Prop:identity_representation}]
Note that \eqref{def:ReLU_identity_d_is_one} and \cref{Lemma:PropertiesOfParallelizationEqualLengthDims} prove \cref{identity_representation:1}.
Observe that \eqref{def:ReLU_identity_d_is_one} and \cref{Lemma:PropertiesOfParallelizationEqualLength} establish \cref{identity_representation:2}.
Note that \cref{identity_representation:1} implies \cref{identity_representation:3}. The proof of \cref{Prop:identity_representation} is thus complete.
\end{proof}
\cfclear

\cfclear
\begingroup
\newcommand{\f}{\mathscr{f}}
\newcommand{\g}{\mathscr{g}}
\begin{athm}{lemma}{lem:comp2}
	Let $\f,\g\in\ANNs$ satisfy
	$\inDimANN(\f)=\outDimANN(\g)$ \cfload.
	Then
	\begin{enumerate}[(i)]
		\item \llabel{it:1}
        it holds that $\inDimANN(\compANN{\f}{\g})=\inDimANN(\g)$ and
		\item \llabel{it:2}
        it holds that $\outDimANN(\compANN{\f}{\g})=\outDimANN(\f)$
	\end{enumerate}
	\cfout.
\end{athm}
\begin{aproof}
  \Nobs that \cref{ANNoperations:Composition} establishes
  \lref{it:1} and \lref{it:2}.
  \finishproofthus
\end{aproof}
\endgroup

\cfclear
\begin{athm}{lemma}{Lemma:CompositionAssociative}
	Let 
	$\mathscr{f}_1,\mathscr{f}_2,\mathscr{f}_3\in\ANNs$
	satisfy
	$\inDimANN(\mathscr{f}_1)=\outDimANN(\mathscr{f}_2)$ and
	$\inDimANN(\mathscr{f}_2)=\outDimANN(\mathscr{f}_3)$ 
	\cfload.
	Then $
	\compANN{(\compANN{\mathscr{f}_1}{\mathscr{f}_2})}{\mathscr{f}_3}=\compANN{\mathscr{f}_1}{(\compANN{\mathscr{f}_2}{\mathscr{f}_3})} $
	\cfadd{lem:comp2}\cfout.
\end{athm}
\begin{aproof}
    \Nobs that \cite[Lemma~2.8]{grohs2019space} establishes
    that $
	\compANN{(\compANN{\mathscr{f}_1}{\mathscr{f}_2})}{\mathscr{f}_3}=\compANN{\mathscr{f}_1}{\allowbreak(\compANN{\mathscr{f}_2}{\mathscr{f}_3})} $
    \finishproofthus
\end{aproof}

\cfclear
\begin{athm}{prop}{Lemma:PropertiesOfCompositions_n2}
   Let $n \in \{2,3,\dots\}$,
	$\mathscr{f}_1, \mathscr{f}_2, \ldots, \mathscr{f}_n \in \ANNs$ 
	satisfy for all $k\in \{2,3, \ldots, n\}$ that
	$\inDimANN(\mathscr{f}_{k-1})=\outDimANN(\mathscr{f}_{k})$
	\cfload.
	Then
\begin{enumerate}[(i)]
\item \label{PropertiesOfCompositions_n:Length}  it holds that
\begin{equation}
[\lengthANN(\compANN{\mathscr{f}_1}{{\mathscr{f}_2}}\bullet \ldots \bullet \mathscr{f}_n)-1]=[\lengthANN(\mathscr{f}_1)-1]+[\lengthANN(\mathscr{f}_2)-1]+ \ldots + [\lengthANN(\mathscr{f}_n)-1],
\end{equation}
\cfadd{Lemma:CompositionAssociative}

\item \label{PropertiesOfCompositions_n:Input} 
 it holds that
$\inDimANN(\mathscr{f}_1 \bullet \mathscr{f}_2\bullet \ldots \bullet \mathscr{f}_n)=\inDimANN(\mathscr{f}_n)$,

\item \label{PropertiesOfCompositions_n:Output} 
 it holds that
$\outDimANN(\mathscr{f}_1 \bullet \mathscr{f}_2\bullet \ldots \bullet \mathscr{f}_n)=\outDimANN(\mathscr{f}_1),
$
and
\item \label{PropertiesOfCompositions_n:Realization} 
 it holds for all  $\activation\in C(\R,\R)$ that
\begin{equation}\label{PropertiesOfCompositions:RealizationEquation_n}
\functionANN({\mathscr{f}_1}\bullet{\mathscr{f}_2}\bullet \ldots \bullet \mathscr{f}_n)=[\functionANN(\mathscr{f}_1)]\circ [\functionANN(\mathscr{f}_2)]\circ \ldots \circ [\functionANN(\mathscr{f}_n)]
\end{equation}
\end{enumerate}
\cfout.
\end{athm}

\begin{proof}
[Proof of \cref{Lemma:PropertiesOfCompositions_n2}]
Note that \cite[Proposition~2.6]{grohs2019space}
and induction prove \cref{PropertiesOfCompositions_n:Length,PropertiesOfCompositions_n:Input,PropertiesOfCompositions_n:Output,PropertiesOfCompositions_n:Realization}.
This completes the proof of \cref{Lemma:PropertiesOfCompositions_n2}.
\end{proof}

\cfclear
\begingroup
\newcommand{\f}{\mathscr f}
\newcommand{\g}{\mathscr g}
\begin{athm}{lemma}{lem:dimcomp}
	Let $\f,\g\in\ANNs$ satisfy $\inDimANN(\f)=\outDimANN(\g)$ \cfload. Then
	\begin{enumerate}[(i)]
		\item \label{it:dimcomp.1}
		it holds that
		\begin{equation}
			\dims(\compANN{\f}{\g})
			=
			(\singledims_0(\g),\singledims_1(\g),\dots,\singledims_{\lengthANN(\g)-1}(\g),
			\singledims_1(\f),\singledims_2(\f),\dots,\singledims_{\lengthANN(\f)}(\f))
		\end{equation}
		and
		\item \label{it:dimcomp.2}
		it holds that
		\begin{equation}
		\begin{split}
			&\dims(\compANN{\f}{\compANN{\ReLUidANN{\outDimANN(\g)}}\g})
			\\&=
			(\singledims_0(\g),\singledims_1(\g),\dots,\singledims_{\lengthANN(\g)-1}(\g),
			2\singledims_{\lengthANN(\g)}(\g),
			\singledims_1(\f),\singledims_2(\f),\dots,\singledims_{\lengthANN(\f)}(\f))
		\end{split}
		\end{equation}
	\end{enumerate}
	\cfout.
\end{athm}
\begin{aproof}
	Note that \cite[Proposition~2.6]{grohs2019space} proves \cref{it:dimcomp.1,it:dimcomp.2}.
	\finishproofthus
\end{aproof}
\endgroup

\cfclear
\begin{athm}{prop}{Lemma:PropertiesOfCompositions_n1}
	Let $n \in \{2,3,\dots\}$,
	$\mathscr{f}_1, \mathscr{f}_2, \ldots, \mathscr{f}_n \in \ANNs$ and let $l_{k,j}\in\N$, $j\in\{1,2,\ldots,\lengthANN({\mathscr{f}_k})\}$, $k\in\{1,2, \ldots, n\}$,  
	satisfy for all $k\in \{1,2, \ldots, n\}$ that
	$\inDimANN(\mathscr{f}_{\min\{k,n-1\}})=\outDimANN(\mathscr{f}_{\min\{k+1,n\}})$ and $\dims(\mathscr{f}_k)=(l_{k,0},l_{k,1},\dots, l_{k,\lengthANN(\mathscr{f}_k)})$
\cfload.
	Then 
\begin{enumerate}[(i)]

\item \label{PropertiesOfCompositions_n:Params} it holds that
\begin{equation}
\begin{split}
\paramANN(\compANN{\mathscr{f}_1}{\mathscr{f}_2}\bullet \ldots \bullet \mathscr{f}_n)&
= \PR*{\sum_{k=1}^n\paramANN(\mathscr{f}_k) }
+ \PR*{\sum_{k = 1}^{n-1}l_{k,1}( l_{k+1,\lengthANN(\mathscr{f}_{k+1})-1}+1) }
\\&\quad
-\PR*{\sum_{k = 1}^{n-1}l_{k,1}(l_{k,0} + 1) }
-\PR*{\sum_{k = 2}^{n}l_{k,\lengthANN(\mathscr{f}_{k})}(l_{k,\lengthANN(\mathscr{f}_{k})-1} + 1) }
\\&\leq
\PR*{ \sum_{k=1}^n\paramANN(\mathscr{f}_k) }
+ \PR*{\sum_{k = 1}^{n-1}l_{k,1}( l_{k+1,\lengthANN(\mathscr{f}_{k+1})-1}+1)},
\end{split}
\end{equation}
and
\item \label{PropertiesOfCompositions_id_n:Params} it holds that
\begin{equation}
\begin{split}
&\paramANN({\mathscr{f}_1} \bullet \ReLUidANN{\outDimANN(\mathscr{f}_2)} \bullet     {\mathscr{f}_2} \bullet \ReLUidANN{\outDimANN(\mathscr{f}_3)} \bullet
        \ldots \bullet \ReLUidANN{\outDimANN(\mathscr{f}_n)} \bullet \mathscr{f}_n)
\\&= \PR*{\sum_{k=1}^n\paramANN(\mathscr{f}_k) }
-\PR*{\sum_{k = 1}^{n-1}l_{k,1}(l_{k,0} + 1) }
-\PR*{\sum_{k = 2}^{n}l_{k,\lengthANN(\mathscr{f}_{k})}(l_{k,\lengthANN(\mathscr{f}_{k})-1} + 1) }
\\ &\quad
+ \PR*{\sum_{k = 1}^{n-1}l_{k,1}(2l_{k,0} + 1) }
+2\PR*{\sum_{k = 2}^{n}l_{k,\lengthANN(\mathscr{f}_{k})}(l_{k,\lengthANN(\mathscr{f}_{k})-1} + 1) }
\\&\leq
\PR*{ \sum_{k=1}^n\paramANN(\mathscr{f}_k) }
+\PR*{\sum_{k = 1}^{n-1}l_{k,1}(l_{k,0} + 1) }
+\PR*{\sum_{k = 2}^{n}l_{k,\lengthANN(\mathscr{f}_{k})}(l_{k,\lengthANN(\mathscr{f}_{k})-1} + 1) }
\end{split}
\end{equation}
\end{enumerate}
\cfout.
\end{athm}
\begin{proof}
[Proof of \cref{Lemma:PropertiesOfCompositions_n1}]
Observe that \cite[Proposition~2.6]{grohs2019space} and induction establish \cref{PropertiesOfCompositions_n:Params,PropertiesOfCompositions_id_n:Params}.
The proof of \cref{Lemma:PropertiesOfCompositions_n1} is thus complete.
\end{proof}
\cfclear

\begin{athm}{prop}{Lemma:PropertiesOfCompositions_n3}
    Let $n \in \{2,3,\dots\}$,
	$\mathscr{f}_1, \mathscr{f}_2, \ldots, \mathscr{f}_n \in \ANNs$
	satisfy for all $k\in \{2,3, \ldots, n\}$ that
	$\inDimANN(\mathscr{f}_{k-1})=\outDimANN(\mathscr{f}_{k})$
	\cfload.
	Then
\begin{enumerate}[(i)]
    \item \label{PropertiesOfCompositions_n:Params2}
	it holds that
\begin{equation}
\paramANN(\compANN{\mathscr{f}_1}{\mathscr{f}_2}\bullet \ldots \bullet \mathscr{f}_n)
\leq 2 \PR*{\sum_{k=1}^{n-1}\paramANN(\mathscr{f}_k)\paramANN(\mathscr{f}_{k+1})},
\end{equation} 
    \item \label{PropertiesOfCompositions_id_n:Params2} it holds that
\begin{equation}
\paramANN({\mathscr{f}_1} \bullet \ReLUidANN{\outDimANN(\mathscr{f}_2)} \bullet           {\mathscr{f}_2} \bullet \ReLUidANN{\outDimANN(\mathscr{f}_3)} \bullet
        \ldots \bullet \ReLUidANN{\outDimANN(\mathscr{f}_n)} \bullet \mathscr{f}_n)
\leq 3\PR*{ \sum_{k=1}^{n}\paramANN(\mathscr{f}_k) } - \paramANN(\mathscr f_1)-\paramANN(\mathscr f_n),
\end{equation}  
and
    \item \label{PropertiesOfCompositions_n:comparison} it holds that
\begin{equation}
 \sum_{k=1}^{n}\paramANN(\mathscr{f}_k) \leq  \sum_{k=1}^{n-1}\paramANN(\mathscr{f}_k)\paramANN(\mathscr{f}_{k+1}) 
\end{equation}
\end{enumerate}
\cfout.
\end{athm}

\begin{proof}[Proof of \cref{Lemma:PropertiesOfCompositions_n3}]
Throughout this proof let $l_{k,j}\in\N$, $j\in\{1,2,\ldots,\lengthANN({\mathscr{f}_k})\}$, $k\in\{1,2, \ldots, n\}$, 
satisfy for all $k\in\{1,2,\ldots,n\}$ that $\dims(\mathscr{f}_k)=(l_{k,0},l_{k,1},\dots, l_{k,\lengthANN(\mathscr{f}_k)})$ \cfload. 
Observe that the fact that for all $a,b \in [2, \infty)$ it holds that 
    $
        a+b \leq 2\max\{a,b\} \leq \min\{a,b\}\max\{a,b\}
        = ab
    $
    and the fact that 
	$\min_{i \in \{1,2, \ldots, n\}}\paramANN(\mathscr{f}_i) \geq 2$ ensure that 
	\begin{equation}
	\begin{split}
	     \sum_{k=1}^{n}\paramANN(\mathscr{f}_k)  
	    &\leq
	    \paramANN(\mathscr{f}_1)+\paramANN(\mathscr{f}_n)+2\PR*{ \sum_{k=2}^{n-1}\paramANN(\mathscr{f}_k) }
	    \\&= 
	   \PR*{\paramANN(\mathscr{f}_1)+\paramANN(\mathscr{f}_2)}
	    + \PR*{\paramANN(\mathscr{f}_2)+\paramANN(\mathscr{f}_3)}
	    + \ldots
	    + \PR*{\paramANN(\mathscr{f}_{n-1})+\paramANN(\mathscr{f}_n)}
	    \\&\leq
	     \sum_{k=1}^{n-1}\paramANN(\mathscr{f}_k)\paramANN(\mathscr{f}_{k+1}) \ifnocf.
	\end{split}
	\end{equation}
\cfload[. ]This establishes \cref{PropertiesOfCompositions_n:comparison}.
    Note that \cref{Lemma:PropertiesOfCompositions_n1} and \cref{PropertiesOfCompositions_n:comparison} imply that
\begin{equation}
    \begin{split}
	    \paramANN(\compANN{\mathscr{f}_1}{\mathscr{f}_2} \bullet \ldots \bullet \mathscr{f}_n) 
	    &\leq
        \PR*{ \sum_{k=1}^n\paramANN(\mathscr{f}_k) }
        + \PR*{ \sum_{k = 1}^{n-1}l_{k,1}( l_{k+1,\lengthANN(\mathscr{f}_{k+1})-1}+1)}
	    \\ & \leq 
	    \PR*{\sum_{k=1}^{n}\paramANN(\mathscr{f}_k)} + \PR*{
	    \sum_{k=1}^{n-1}\paramANN(\mathscr{f}_k)\paramANN(\mathscr{f}_{k+1}) }
	    \\ & \leq 2 \PR*{ 
	    \sum_{k=1}^{n-1}\paramANN(\mathscr{f}_k)\paramANN(\mathscr{f}_{k+1})}
	\end{split}
\end{equation}
\cfload. This establishes \cref{PropertiesOfCompositions_n:Params2}.
Observe that \cref{Lemma:PropertiesOfCompositions_n1} ensures that
\begin{equation}
    \begin{split}
	    &
	    \paramANN({\mathscr{f}_1} \bullet \ReLUidANN{\outDimANN(\mathscr{f}_2)} \bullet           {\mathscr{f}_2} \bullet \ReLUidANN{\outDimANN(\mathscr{f}_3)} \bullet
        \ldots \bullet \ReLUidANN{\outDimANN(\mathscr{f}_n)} \bullet \mathscr{f}_n)
				\\&\leq
        \PR*{ \sum_{k=1}^n\paramANN(\mathscr{f}_k) }
        +\PR*{\sum_{k = 1}^{n-1}l_{k,1}(l_{k,0} + 1)}
        +\PR*{\sum_{k = 2}^{n}l_{k,\lengthANN(\mathscr{f}_{k})}(l_{k,\lengthANN(\mathscr{f}_{k})-1} + 1) }
				\\&\leq
        \PR*{ \sum_{k=1}^n\paramANN(\mathscr{f}_k) }
        +\PR*{\sum_{k = 1}^{n-1}\paramANN(\mathscr{f}_k)}
        +\PR*{\sum_{k = 2}^{n}\paramANN(\mathscr{f}_k) }
        \\&=
				3\PR*{ \sum_{k=1}^{n}\paramANN(\mathscr{f}_k) }
				-\paramANN(\mathscr f_1)-\paramANN(\mathscr f_n)
								\ifnocf.
	\end{split}
\end{equation}
\cfload[.]%
    Hence, 
we obtain 
    \cref{PropertiesOfCompositions_id_n:Params2}.
The proof of \cref{Lemma:PropertiesOfCompositions_n3} is thus complete.
\end{proof}

\subsection{Powers and extensions of ANNs}
\label{Subsection:2.5}

In this section, we recall the notion of the extension of an ANN, which, under certain assumptions, provides a way of obtaining an ANN with the same realization but a larger length than a given ANN. For this purpose, we first recall the notion of powers of ANNs. Extensions of ANNs, in turn, will be used to define parallelizations of ANNs with different lengths in \cref{subsection:2.6}.

\begin{definition}
[Identity matrices]
\label{def:identityMatrix}
\cfconsiderloaded{def:identityMatrix}
	Let $d\in\N$. Then we denote by $\idMatrix_{d}\in \R^{d\times d}$ the identity matrix in $\R^{d\times d}$.
\end{definition}

\begin{definition}[Affine linear transformation ANNs]
\label{def:ANN:affine}
Let $m,n\in\N$, $W\in\R^{m\times n}$, $B\in\R^m$. Then we denote by 
$\affineANN_{W,B} \in (\R^{m\times n}\times \R^m) \subseteq \ANNs$ the \cfadd{def:neuralnetwork}ANN
given by $\affineANN_{W,B} = ((W,B))$ (cf.~\cref{def:ANN}).
\end{definition}

\cfclear
\begin{definition}
[Powers of ANNs]
\label{def:iteratedANNcomposition}
\cfconsiderloaded{def:iteratedANNcomposition}
	We denote by $(\cdot)^{\bullet n}\colon \{\mathscr{f}\in \ANNs\colon \inDimANN(\mathscr{f})=\outDimANN(\mathscr{f})\}\allowbreak\to\ANNs$, $n\in\N_0$, the functions
		 which satisfy for all $n\in\N_0$, $\mathscr{f}\in\ANNs$ with $\inDimANN(\mathscr{f})=\outDimANN(\mathscr{f})$ that 
	\begin{equation}\label{iteratedANNcomposition:equation}
		\begin{split}
		\mathscr{f}^{\bullet n}=
		\begin{cases} \affineANN_{\idMatrix_{\outDimANN(\mathscr{f})},0}
		&: n=0 \\
		\,\compANN{\mathscr{f}}{(\mathscr{f}^{\bullet (n-1)})} &: n\in\N
		\end{cases}
		\end{split}
	\end{equation}	
	\cfload.
\end{definition}

\cfclear
\begin{athm}{lemma}{lem:dimensionsiteratedcompositions}
Let $\mathscr{f}\in\ANNs$ satisfy $\inDimANN(\mathscr{f})=\outDimANN(\mathscr{f})$ \cfload. Then it holds for all $n\in\N_0$ that
\begin{equation}
\label{lem:dimensionsiteratedcompositions:equation}
\cfadd{def:iteratedANNcomposition}
\inDimANN(\mathscr{f}^{\bullet n})=\outDimANN(\mathscr{f}^{\bullet n})=\inDimANN(\mathscr{f})
\end{equation}
\cfout.
\end{athm}

\begin{proof}[Proof of \cref{lem:dimensionsiteratedcompositions}]
First, note that \eqref{iteratedANNcomposition:equation} assures that
\begin{equation}
\cfadd{def:iteratedANNcomposition}
\label{lem:dimensionsiteratedcompositions:eq1}
\inDimANN(\mathscr{f}^{\bullet 0})=\outDimANN(\mathscr{f}^{\bullet 0})=\outDimANN(\mathscr{f})=\inDimANN(\mathscr{f})
.
\end{equation}
\cfload[.]Furthermore, observe that \cref{Lemma:PropertiesOfCompositions_n2} and \eqref{iteratedANNcomposition:equation} ensure that for all $n\in\N$ it holds that $\outDimANN(\mathscr{f}^{\bullet n})=\outDimANN(\compANN{\mathscr{f}}{(\mathscr{f}^{\bullet (n-1)})})=\outDimANN(\mathscr{f})$ and 
\begin{equation}
\inDimANN(\mathscr{f}^{\bullet n})=\inDimANN(\compANN{\mathscr{f}}{(\mathscr{f}^{\bullet (n-1)})})=\inDimANN(\mathscr{f}^{\bullet (n-1)}).
\end{equation}
This, \eqref{lem:dimensionsiteratedcompositions:eq1}, and induction establish
\eqref{lem:dimensionsiteratedcompositions:equation}. The proof of \cref{lem:dimensionsiteratedcompositions} is thus complete.
\end{proof}

\cfclear
\begin{definition}[Extensions of ANNs]
\label{def:ANNenlargement}
\cfconsiderloaded{def:ANNenlargement}
	Let $L\in\N$, $\mathscr{g}\in \ANNs$ satisfy $\inDimANN(\mathscr{g})=\outDimANN(\mathscr{g})$ \cfload.
	Then
	we denote by $\longerANN{L,\mathscr{g}}\colon \{\mathscr{f}\in\ANNs\colon (\lengthANN(\mathscr{f})\le L \andShort \outDimANN(\mathscr{f})=\inDimANN(\mathscr{g})) \}\to \ANNs$ the function which satisfies for all $\mathscr{f}\in\ANNs$ with $\lengthANN(\mathscr{f})\le L$ and $\outDimANN(\mathscr{f})=\inDimANN(\mathscr{g})$ that
	\begin{equation}\label{ANNenlargement:Equation}
	\longerANN{L,\mathscr{g}}(\mathscr{f})=	 \compANN{\prb{\mathscr{g}^{\bullet (L-\lengthANN(\mathscr{f}))}}}{\mathscr{f}}
	\cfadd{def:iteratedANNcomposition}
	\cfadd{lem:dimensionsiteratedcompositions}
	\end{equation}
	\cfout.
\end{definition}
\cfclear

\begin{athm}{lemma}{Lemma:PropertiesOfANNenlargementGeometry}
	Let  $d,{i}\in\N$, $\mathscr{f},\mathscr{g}\in \ANNs$ satisfy $\outDimANN(\mathscr{f})=d$ and $\dims(\mathscr{g})=(d,{i},d)$ 
	\cfload.
	Then
	\begin{enumerate}[(i)]
	\item \label{PropertiesOfANNenlargementGeometry:BulletPower} it holds for all $n\in\N_0$ that
	$\lengthANN(\mathscr{g}^{\bullet n})=n+1$,  $\dims(\mathscr{g}^{\bullet n})\in \N^{n+2}$, and 
	\begin{equation}\label{BulletPower:Dimensions}
	\dims(\mathscr{g}^{\bullet n}) = 
	\begin{cases}
	(d,d) &: n=0\\
	(d,{i},{i},\dots,{i},d) &: n\in\N
	\cfadd{def:iteratedANNcomposition}
	\end{cases}
	\end{equation}
	 and
	\item \label{PropertiesOfANNenlargementGeometry:ItemLonger}
	it holds 
	for all
	$L\in\N\cap [\lengthANN(\mathscr{f}),\infty)$ 
	that $\lengthANN\prb{\longerANN{L,\mathscr{g}}(\mathscr{f})}=L$ and
	\begin{equation}
	\label{PropertiesOfANNenlargementGeometry:ParamsLonger}
	\begin{split}
	&\paramANN(\longerANN{L,\mathscr{g}}(\mathscr{f}))
\\&\le 
\begin{cases} \paramANN(\mathscr{f})
		&: \lengthANN(\mathscr{f})=L \\
		\PR*{\prb{\max\pRb{1,\tfrac{i}{d}}}\paramANN(\mathscr{f})+
		\prb{(L-\lengthANN(\mathscr{f})-1) \,{i}+d}({i}+1)
		} &: \lengthANN(\mathscr{f})<L
		\end{cases}
	\end{split}
	\end{equation}
\end{enumerate}
	\cfout.
\end{athm}

\begin{proof}[Proof of \cref{Lemma:PropertiesOfANNenlargementGeometry}]
    
    Note that \cite[Lemma~2.13]{grohs2019space} establishes
    \cref{PropertiesOfANNenlargementGeometry:BulletPower,PropertiesOfANNenlargementGeometry:ItemLonger}.
	The proof of \cref{Lemma:PropertiesOfANNenlargementGeometry} is thus complete.
\end{proof}

\cfclear
\begin{athm}{lemma}{Lemma:PropertiesOfANNenlargementRealization}
	Let $\activation\in C(\R,\R)$,  $ \mathscr f,\mathfrak{I}\in \ANNs$ satisfy 
	for all  $x\in\R^{\inDimANN(\mathfrak{I})}$ that 
	$\outDimANN(\mathscr f)=\inDimANN(\mathfrak{I})=\outDimANN(\mathfrak{I})$ 
	and $(\functionANN( \mathfrak{I}))(x)=x$
	\cfload.
	Then
	\begin{enumerate}[(i)]
		\item \label{PropertiesOfANNenlargementRealization:itemOne} it holds for all $n\in\N_0$, $x\in\R^{\inDimANN(\mathfrak{I})}$ that
		\begin{equation}\label{BulletPowerRealization:Dimensions}
			\functionANN(\mathfrak{I}^{\bullet n})\in C(\R^{\inDimANN(\mathfrak{I})},\R^{\inDimANN(\mathfrak{I})}) \qandq (\functionANN(\mathfrak{I}^{\bullet n}))(x)=x
			\cfadd{def:iteratedANNcomposition}
		\end{equation}
		 and	
		\item 	\label{PropertiesOfANNenlargementRealization:ItemIdentityLonger}
			it holds 
		for all
		$L\in\N\cap [\lengthANN(\mathscr{f}),\infty)$, 
		$x\in\R^{\inDimANN(\mathscr{f})}$ 
		that	
		\begin{equation}
		\functionANN(\longerANN{L,\mathfrak{I}}(\mathscr{f}))\in C(\R^{\inDimANN(\mathscr{f})},\R^{\outDimANN(\mathscr{f})})
		\qandq
		\prb{\functionANN(\longerANN{L,\mathfrak{I}}(\mathscr{f}))}(x)=\prb{\functionANN(\mathscr{f})}(x)
		\end{equation}
	\end{enumerate}
\cfout.
\end{athm}

\begin{proof}[Proof of \cref{Lemma:PropertiesOfANNenlargementRealization}]

Observe that \cite[Lemma~2.14]{grohs2019space} proves
\cref{PropertiesOfANNenlargementRealization:itemOne,PropertiesOfANNenlargementRealization:ItemIdentityLonger}.
	The proof of \cref{Lemma:PropertiesOfANNenlargementRealization} is thus complete.
\end{proof}

\subsection{Parallelizations of ANNs with different lengths}
\label{subsection:2.6}

In this section we recall how, under certain assumptions, the notion of parallelization of ANNs with the same length introduced in \cref{subsec:parallelization} above may be generalized to the case of ANNs with different lengths. In \cref{Lemma:PropertiesOfANNenlargementRealization} we establish an upper bound on the number of parameters of the parallelization of ANNs with different lengths, which will be crucial to our main results.

\cfclear
\begin{definition}[Parallelization of ANNs with different lengths]\label{def:generalParallelization}
\cfconsiderloaded{def:generalParallelization}
	Let $n\in\N$,  $\mathscr{g}=(\mathscr{g}_1,\mathscr{g}_2,\dots, \mathscr{g}_n)\in \ANNs^n$ satisfy for all $j\in\{1,2,\dots, n\}$ that $\lengthANN(\mathscr{g}_j)=2$ and $\inDimANN(\mathscr{g}_j)=\outDimANN(\mathscr{g}_j)$ \cfload.
	Then
	we denote by 
	\begin{equation}
		\parallelization_{n,\mathscr{g}}\colon \pRb{(\mathscr{f}_1,\mathscr{f}_2,\dots,\mathscr{f}_n)\in\ANNs^n\colon\smallsum_{j=1}^{n}\vass{\outDimANN(\mathscr{f}_j)-\inDimANN(\mathscr{g}_j)}=0}\to \ANNs
	\end{equation}
	the function which satisfies for all $\mathscr{f}=(\mathscr{f}_1,\mathscr{f}_2,\dots, \mathscr{f}_n)\in\ANNs^n$ with 
	$\sum_{j=1}^{n}\vass{\outDimANN(\mathscr{f}_j)- \inDimANN(\mathscr{g}_j)}=0$
	that
	\begin{multline}\label{generalParallelization:DefinitionFormula}
			\parallelization_{n,\mathscr{g}}(\mathscr{f})=\parallelizationSpecial_{n}\prb{\longerANN{\max_{k\in\{1,2,\dots,n\}}\lengthANN(\mathscr{f}_k),\mathscr{g}_1}({\mathscr{f}_1}),\longerANN{\max_{k\in\{1,2,\dots,n\}}\lengthANN(\mathscr{f}_k),\mathscr{g}_2}({\mathscr{f}_2}),\\\dots,\longerANN{\max_{k\in\{1,2,\dots,n\}}\lengthANN(\mathscr{f}_k),\mathscr{g}_n}({\mathscr{f}_n})}
	\end{multline}
	\cfadd{Lemma:PropertiesOfANNenlargementGeometry}
	\cfload.
\end{definition}

\cfclear
\begin{athm}{prop}{Lemma:PropertiesOfParallelizationRealization}
	Let $a\in C(\R,\R)$, $n\in\N$, $\mathfrak{I}=(\idANNshort{1},\idANNshort{2},\dots, \idANNshort{n})$, $\mathscr{f}=(\mathscr{f}_1,\mathscr{f}_2,\dots,\allowbreak \mathscr{f}_n)\in\ANNs^n$
	satisfy for all $j\in\{1,2,\dots, n\}$, $x\in\R^{\outDimANN(\mathscr{f}_j)}$ that $\lengthANN(\idANNshort{j}) =2$, $ \inDimANN(\mathfrak{I}_j)=\outDimANN(\mathfrak{I}_j)=\outDimANN(\mathscr{f}_j)$, and
	$(\functionANN(\idANNshort{j}))(x)=x$
	\cfload.
	Then
	\begin{enumerate}[(i)]
		\item \label{PropertiesOfParallelizationRealization:ItemOne} it holds that
		\begin{equation}
			\functionANN\prb{\parallelization_{n,\mathfrak{I}}(\mathscr{f})}\in C\prb{\R^{[\sum_{j=1}^n \inDimANN(\mathscr{f}_j)]},\R^{[\sum_{j=1}^n \outDimANN(\mathscr{f}_j)]}}
		\end{equation}
		and
		\item \label{PropertiesOfParallelizationRealization:ItemTwo}
		it holds for all $x_1\in\R^{\inDimANN(\mathscr{f}_1)},x_2\in\R^{\inDimANN(\mathscr{f}_2)},\dots, x_n\in\R^{\inDimANN(\mathscr{f}_n)}$ that 
		\begin{equation}\label{PropertiesOfParallelizationRealizationEqualLengthFunction}
		\begin{split}
		&\prb{ \functionANN\prb{\parallelization_{n,\mathfrak{I}}(\mathscr{f})} } (x_1,x_2,\dots, x_n) 
		\\&=\prb{(\functionANN(\mathscr{f}_1))(x_1), (\functionANN(\mathscr{f}_2))(x_2),\dots,
		(\functionANN(\mathscr{f}_n))(x_n) }
		\end{split}
		\end{equation}
	\end{enumerate}
\cfout.
\end{athm}

\begin{proof}[Proof of \cref{Lemma:PropertiesOfParallelizationRealization}]
    
    Note that \cite[Corollary~2.23]{grohs2019space} establishes
    \cref{PropertiesOfParallelizationRealization:ItemOne,PropertiesOfParallelizationRealization:ItemTwo}.
	The proof of \cref{Lemma:PropertiesOfParallelizationRealization} is thus complete.
\end{proof}

\cfclear
\begin{athm}{cor}{Lemma:PropertiesOfParallelization}
	Let $n\in\N$, ${l}_1,{l}_2,\dots, {l}_n\in\N$,  $\mathscr{g}=(\mathscr{g}_1,\mathscr{g}_2,\dots,\allowbreak \mathscr{g}_n)$, $\mathscr{f}=(\mathscr{f}_1,\mathscr{f}_2,\allowbreak\dots, \mathscr{f}_n)\in \ANNs^n$
	satisfy for all $j\in\{1,2,\dots, n\}$ that $\dims(\mathscr{g}_j) = (\outDimANN(\mathscr{f}_j),{l}_j,\allowbreak\outDimANN(\mathscr{f}_j))$ 
	\cfload.
	Then
	\begin{equation}
	\begin{split}
	\paramANN\prb{\parallelization_{n,\mathscr{g}}(\mathscr{f})}
	&\le 
    \frac{1}{2} \biggl[\ssum_{j=1}^n\Bigl(
	\PR[\big]{\max\bigl\{1,\tfrac{l_j}{\outDimANN(\mathscr{f}_j)}\bigr\}}\,\paramANN(\mathscr{f}_j)
    \\&\qquad
    +
	\prbb{\PRbb{\max_{k\in\{1,2,\dots,n\}} \lengthANN(\mathscr{f}_k)} {l}_j	+\outDimANN(\mathscr{f}_j)}({l}_j+1)\Bigl)\biggl]^2
	\end{split}
	\end{equation}
	\cfout.
\end{athm}
\begin{proof}[Proof of \cref{Lemma:PropertiesOfParallelization}]
    Throughout this proof 
        let 
            $L\in\N$ 
        satisfy
            $L=\max_{k\in\{1,2,\dots,n\}} \lengthANN(\mathscr{f}_k)$.
    Observe that 
        \eqref{generalParallelization:DefinitionFormula}, 
        \cref{PropertiesOfParallelizationEqualLengthDims:Params} in  \cref{Lemma:PropertiesOfParallelizationEqualLengthDims}, and
        \cref{PropertiesOfANNenlargementGeometry:ItemLonger} in \cref{Lemma:PropertiesOfANNenlargementGeometry} 
    assure that
    \begin{equation}
    \begin{split}
        \paramANN\prb{\parallelization_{n,\mathscr{g}}(\mathscr{f})}
        &=
        \paramANN\prb{\parallelizationSpecial_n\prb{\longerANN{L,\mathscr{g}_1}({\mathscr{f}_1}),\longerANN{L,\mathscr{g}_2}({\mathscr{f}_2}),\dots,\longerANN{L,\mathscr{g}_n}({\mathscr{f}_n})}}
        \\&\le 
        \tfrac{1}{2} \PR*{\smallsum\nolimits_{j=1}^n \paramANN(\longerANN{L,\mathscr{g}_j}({\mathscr{f}_j}))}^2
        \\&\le 
        \tfrac{1}{2} \prbb{\PR*{\smallsum\nolimits_{j=1}^n
            \PRb{\max\bigl\{1,\tfrac{{l}_j}{\outDimANN(\mathscr{f}_j)}\bigr\}}\,\paramANN(\mathscr{f}_j)
            \,\indicator{(\lengthANN(\mathscr{f}_j),\infty)}(L)}
            \\&\qquad+\PR*{\smallsum\nolimits_{j=1}^n\prb{
            (L-\lengthANN(\mathscr{f}_j)-1) \,{l}_j
            +\outDimANN(\mathscr{f}_j)}({l}_j+1)
            \,\indicator{(\lengthANN(\mathscr{f}_j),\infty)}(L)}
            \\&\qquad+\PR*{\smallsum\nolimits_{j=1}^n  \paramANN(\mathscr{f}_j)\,\indicator{\{\lengthANN(\mathscr{f}_j)\}}(L)}}^{\!2}
        \\&\le 
        \tfrac{1}{2} \prbb{\PR*{\smallsum\nolimits_{j=1}^n
            \PR[\big]{\max\bigl\{1,\tfrac{{l}_j}{\outDimANN(\mathscr{f}_j)}\bigr\}}\,\paramANN(\mathscr{f}_j)
            \,\indicator{(\lengthANN(\mathscr{f}_j),\infty)}(L)}
            \\&\qquad+\PR*{\smallsum\nolimits_{j=1}^n\prb{
            L {l}_j
            +\outDimANN(\mathscr{f}_j)}({l}_j+1)
            \,\indicator{(\lengthANN(\mathscr{f}_j),\infty)}(L)}
            \\&\qquad+\PR*{\smallsum\nolimits_{j=1}^n  \PR[\big]{\max\bigl\{1,\tfrac{{l}_j}{\outDimANN(\mathscr{f}_j)}\bigr\}}\paramANN(\mathscr{f}_j)\,\indicator{\{\lengthANN(\mathscr{f}_j)\}}(L)}}^{\!2}
        \\&\le 
        \tfrac{1}{2} \PR*{\smallsum\nolimits_{j=1}^n
            \PR[\big]{\max\bigl\{1,\tfrac{{l}_j}{\outDimANN(\mathscr{f}_j)}\bigr\}}\,\paramANN(\mathscr{f}_j)
            +\prb{L{l}_j+\outDimANN(\mathscr{f}_j)}({l}_j+1)}^2
    \end{split}
    \end{equation}	
    \cfload.
    This completes the proof of \cref{Lemma:PropertiesOfParallelization}.
\end{proof}

\subsection{Clipping functions as ANNs}
\label{subsection:2.7}

In this section we define clipping functions and show that they can be realized by ANNs with the ReLU activation function, a straightforward technical result that will be needed in subsequent sections.

\begin{definition}[Clipping function]
\label{def:clipping_function}
Let
$ u \in \R $,
$ v \in [ u, \infty ) $.
Then
we denote by
$ \clip{u}{v} \colon \R \to \R $
the function which satisfies for all
$ x \in \R $
that 
\begin{equation}
\clip{u}{v}( x )
=
\max\{ u, \min\{ x, v \} \}
.
\end{equation}
\end{definition}

\cfclear
\begin{definition}[Multi-dimensional clipping functions]
    \label{def:clip}
    \cfconsiderloaded{def:clip}
    Let 
        $n \in \N$,
        $u\in \R$,
        $v\in [u,\infty)$.
    Then we denote by 
        $ \Clip uvn \colon \R^{n} \to \R^{n} $ 
    the function which satisfies for all
        $ x \in \R^n $
    that 
    \begin{equation}
        \Clip uvn(x)
        =
        \Mult{\clip{u}{v}}(x)
    \end{equation}
    \cfload.
\end{definition}

\cfclear
\begin{athm}{lemma}{lem:clipping_function}
Let
    $n\in\N$,
    $ u \in \R $,
    $ v \in [ u, \infty ) $. 
Then there exists 
    $\mathscr{f}\in\ANNs$ 
such that
\begin{enumerate}[(i)]
    \item\label{lem:clipping_function:item_0}
        it holds that 
            $\realisation_\mathfrak{r}(\mathscr{f})\in C(\R^n,\R^n)$,
    \item\label{lem:clipping_function:item_1}
        it holds for all 
            $x\in\R^n$ 
        that 
            $(\realisation_\ReLU(\mathscr{f}))(x)=\Clip uvn( x )$,
    \item\label{lem:clipping_function:item_2}
        it holds that 
            $\dims(\mathscr{f})=(n,n,n,n)$, and
    \item\label{lem:clipping_function:item_3}
        it holds that 
            $\paramANN(\mathscr{f})=3n^2+3n$
\end{enumerate}
\cfout.
\end{athm}

\begin{proof}[Proof of \cref{lem:clipping_function}]
Throughout this proof let $\mathscr{g}_1,\mathscr{g}_2\in\ANNs$ satisfy
\begin{equation}
\begin{split}
\mathscr{g}_1=\compANN{\affineANN_{1,u}}{((\idMatrix_1,0),(\idMatrix_1,0))}\bullet\affineANN_{1,-u}\quad \text{and}\quad\mathscr{g}_2=\affineANN_{-1,v}\bullet((\idMatrix_1,0),(\idMatrix_1,0))\bullet\affineANN_{-1,v}
\end{split}
\end{equation}
\cfload. Observe that \cref{Lemma:PropertiesOfCompositions_n2} shows that for all $x\in\R$ it holds that
\begin{equation}
\begin{split}
\pr[\big]{\realisation_\ReLU(\mathscr{g}_1)}(x)=\max\{x-u,0\}+u=\max\{u,x\}
\end{split}
\end{equation}
and
\begin{equation}
\begin{split}
\pr[\big]{\realisation_\ReLU(\mathscr{g}_2)}(x)=-\mathord{\max}\{-x+v,0\}+v=\min\{v,x\}
\end{split}
\end{equation}
\cfload. This and \cref{Lemma:PropertiesOfCompositions_n2} imply that for all $x\in\R$ it holds that
\begin{equation}
\begin{split}
\pr[\big]{\realisation_\mathfrak{r}(\mathscr{g}_1\bullet\mathscr{g}_2)}(x)
&=
\pr[\big]{[\realisation_\mathfrak{r}(\mathscr{g}_1)]\circ[\realisation_\mathfrak{r}(\mathscr{g}_2)]}(x)=\max\{u,\min\{v,x\}\}=\clip{u}{v}( x )
\cfadd{def:clipping_function}
\end{split}
\end{equation}
\cfload. \cref{Lemma:PropertiesOfParallelizationEqualLength} hence ensures that for all $x=(x_1,x_2,\dots,x_n)\in\R^n$ it holds that
\begin{equation}
\begin{split}
\label{lem:clipping_function:eq1}
&\pr[\big]{\realisation_\mathfrak{r}(\parallelizationSpecial_{n}(\mathscr{g}_1\bullet\mathscr{g}_2,\mathscr{g}_1\bullet\mathscr{g}_2,\dots,\mathscr{g}_1\bullet\mathscr{g}_2))}(x)
\\&=\pr[\Big]{(\realisation_\mathfrak{r}(\mathscr{g}_1\bullet\mathscr{g}_2))(x_1),(\realisation_\mathfrak{r}(\mathscr{g}_1\bullet\mathscr{g}_2))(x_2),\dots,
			(\realisation_\mathfrak{r}(\mathscr{g}_1\bullet\mathscr{g}_2))(x_n)}=\mathfrak{C}_{u,v,n}(x)
			\cfadd{def:clip}
\end{split}
\end{equation}
\cfload.
Furthermore, note that \cref{Lemma:PropertiesOfCompositions_n2} 
and the fact that for all $i\in\{1,2\}$ it holds that $\dims(\mathscr{g}_i)=(1,1,1)$ 
imply that $\dims(\mathscr{g}_1\bullet\mathscr{g}_2)=(1,1,1,1)$. 
\cref{Lemma:PropertiesOfParallelizationEqualLengthDims} therefore shows 
that $\dims(\parallelizationSpecial_{n}(\mathscr{g}_1\bullet\mathscr{g}_2,\mathscr{g}_1\bullet\mathscr{g}_2,\dots,\mathscr{g}_1\bullet\mathscr{g}_2))=(n,n,n,n)$. 
Combining this with \eqref{lem:clipping_function:eq1} establishes 
\cref{lem:clipping_function:item_0,lem:clipping_function:item_1,lem:clipping_function:item_2,lem:clipping_function:item_3}. 
The proof of \cref{lem:clipping_function} is thus complete.
\end{proof}

\section{Approximability of functions}
\label{Section:3}
This section contains essential definitions and general results on
the approximation of high-dimensional functions by ANNs, which will
be crucial for the later results in this article and in particular
for proving \cref{Theo:introduction} in the introduction.

In \cref{Subsection:3.1} below we define the cost of an ANN approximation,
which is the number of parameters necessary to approximate
a function by a neural network up to a given accuracy and with the approximation
satisfying a specified Lip\-schitz condition.
This definition is inspired by Cheridito et al.~\cite[Definition~9]{cheridito2021efficient}.
In \cref{Subsection:3.3} we introduce ANN approximation spaces,
which serve as our central formalization of functions that can
be approximated by ANNs with the ReLU activation function without the 
curse of dimensionality.
\Cref{Subsection:3.4,Subsection:3.5} deal with ANN approximations of 
compositions of functions. The main goal in these two subsections is
to prove that the approximation
space for multi-dimensional functions defined in \cref{Subsection:3.3} is closed under composition.
Finally, \cref{Subsection:3.6} concerns parallelizations of
functions, i.e., functions which are componentwise applications
of univariate functions. Again we analyze how the cost of ANN approximations
behaves under this operation and then prove a certain result
about closedness of the approximation spaces defined in \cref{Subsection:3.3}.

We refer the reader to Beneventano et al.~\cite[Subsection~3.3]{beneventano2020highdimensional}
for related notions of approximation spaces of ANN approximable
functions. We also note that Gribonval et al.~\cite{gribonval2019approximation}
have introduced a concept of approximation spaces of functions in the context
of ANN approximations, which is, however, only distantly related
to the notion presented here.

\subsection{Costs for ANN approximations}
\label{Subsection:3.1}

In this section we define the cost of an ANN approximation, which is, roughly speaking the minimum number of parameters of an ANN whose realization approximates a given function up to a given approximation error and which is Lipschitz continuous with a given Lipschitz constant (see \cref{def:cost_Lipschitz_approx_set} below for the precise definition). The remaining results in this section are mostly of a technical nature. More specifically, \cref{lem_cost_of_Lip_approx_set_equivalence_0,lem_cost_of_Lip_approx_set_equivalence_1,lem_cost_of_Lip_approx_set_equivalence,cor_cost_of_Lip_approx_set_equivalence} serve as characterizations of the cost being finite and \cref{Lemma:Monotonicity_of_Cost} is a monotonicity result for the cost function. \cref{lemma_Lip_catalogs_are_Lip_set} shows, roughly speaking, that a function which has finite cost with respect to arbitrarily small prescribed approximation errors is necessarily Lipschitz continuous.

\begin{definition}
[Standard norm]
\label{def:Euclidean_norm}
	We denote by $\norm{\cdot} \colon \prb{\bigcup_{n \in \N} \R^n} \rightarrow [0,\infty)$ the function which satisfies for all $n \in \N$, $x = (x_1,x_2,\dots,x_n) \in \R^n$ that
$
	\norm{x} = \PRb{ \textstyle{\sum}_{j=1}^{n} \vass{x_j}^2 }^{\nicefrac{1}{2}}
$.
\end{definition}

\begin{definition}[Costs for ANN approximations] %
\label{def:cost_Lipschitz_approx_set}
	Let $a \in C(\R,\R)$.
	Then we denote by $\CostLipwo{a}{} \colon \prb{\bigcup_{m,n \in \N}\bigcup_{D \subseteq \R^m} C(D, \R^n)} \times [0, \infty]^2 \rightarrow [1,\infty]$ the function which satisfies for all $m,n \in \N$, $D \subseteq \R^m$, $f \in C(D,\R^n)$, $L,\varepsilon \in [0, \infty]$ that\footnote{Note that we consider a function from a set $A$ to a set $B$ to be a tuple $(A,B,G)$ where $G\subseteq A\times B$ is the graph of the function. In particular, we consider the codomain of a function to be a well-defined notion.}
\begin{equation}
\begin{split}
\label{Cost_Lipschitz_set}
	&\CostLip{a}{A}{f}{L}{\varepsilon} = \\
	&\min \pr*{ \pR*{ \mathfrak{p} \in \N \colon \PR*{\!\!
	\begin{array}{c}
	    \exists\, \mathfrak{m} \in \N, \mathscr{f} \in \ANNs \colon
	    \PRb{(\realisation_a(\mathscr{f}) \in C(\R^\mathfrak{m},\R^n))\land{}\\
			(\mathfrak{p}=\paramANN(\mathscr{f}))\land (D \subseteq \R^\mathfrak{m})\land (\fa{x,y} D \text{ with } x \ne y \colon\\
			\mednorm{ \functionANNbis{\mathscr{f}}(x)-\functionANNbis{\mathscr{f}}(y)} \leq L\norm{x-y})\land{}\\
        (\fa{x} D \colon 
        \mednorm{f(x) -  \functionANNbis{\mathscr{f}}(x)} \leq \varepsilon)}\\
    \end{array} \!\! }
} \cup \{\infty\} }
\end{split}
\end{equation}
(cf.\ \cref{def:ANN,def:ANNrealization,def:Euclidean_norm}).
\end{definition}

\cfclear
\begin{athm}{lemma}{lem_cost_of_Lip_approx_set_equivalence_0}
Let $a \in C(\R,\R)$, $m,n \in \N$, $D\subseteq\R^m$, $L,\varepsilon \in [0,\infty]$, $f \in C(D,\R^n)$. %
Then $\CostLip{a}{\emptyset}{f}{L}{\varepsilon}\geq 2n$ \cfload.
\end{athm}

\begin{proof}[Proof of \cref{lem_cost_of_Lip_approx_set_equivalence_0}.]
Observe that for all 
	$\mathfrak{m}\in\N$, 
	$\mathscr{f}\in\ANNs$ 
	with $\realisation_a(\mathscr{f}) \in C(\R^\mathfrak{m},\R^n)$ 
it holds that 
	$\singledims_{\lengthANN(\mathscr f)}(\mathscr{f})=n$
	\cfload. 
Hence, we obtain that for all 
	$\mathfrak{m}\in\N$, 
	$\mathscr{f}\in\ANNs$ 
	with $\realisation_a(\mathscr{f}) \in C(\R^\mathfrak{m},\R^n)$ 
it holds that 
	\begin{equation}
		\paramANN(\mathscr{f}) 
		\geq 
		[\singledims_{\lengthANN(\mathscr f)}(\mathscr{f})](\singledims_{\lengthANN(\mathscr f)-1}(\mathscr{f})+1)
		\geq
		n(1+1)= 2n.
	\end{equation}
The proof of \cref{lem_cost_of_Lip_approx_set_equivalence_0} is thus complete.
\end{proof}

\cfclear
\begin{athm}{lemma}{lem_cost_of_Lip_approx_set_equivalence_1}
Let $a \in C(\R,\R)$, $n \in \N$, $L,\varepsilon \in [0,\infty]$, $f \in C(\emptyset,\R^n)$. %
Then $\CostLip{a}{\emptyset}{f}{L}{\varepsilon}=2n$ \cfout.
\end{athm}

\begin{proof}[Proof of \cref{lem_cost_of_Lip_approx_set_equivalence_1}.]
Throughout this proof let $W\in\R^{n\times 1}$.
Note that \eqref{ANNrealization:ass2} shows that
\begin{enumerate}[(I)]
\item\label{lem_cost_Lip_app_set_equ_1_item_1_proof} it holds that $\realisation_a(\affineANN_{W,0}) \in C(\R,\R^n)$, 
\item\label{lem_cost_Lip_app_set_equ_1_item_2_proof} it holds that $\paramANN(\affineANN_{W,0}) = 2n$,
\item \label{lem_cost_Lip_app_set_equ_1_item_3_proof} it holds for all $x \in \emptyset$ that $ \mednorm{f(x) -  \functionANNbis{\affineANN_{W,0}}(x)} \leq \varepsilon$, and
\item \label{lem_cost_Lip_app_set_equ_1_item_4_proof} it holds for all $x,y \in \emptyset$ with $x \ne y$ that $\mednorm{ \functionANNbis{\affineANN_{W,0}}(x)-\functionANNbis{\affineANN_{W,0}}(y)} \leq L\norm{x-y}$
\end{enumerate}
\cfload[.]%
Hence, we obtain that $\CostLip aAfL\eps\leq 2n$.
Combining this with \cref{lem_cost_of_Lip_approx_set_equivalence_0} ensures that $\CostLip{a}{\emptyset}{f}{L}{\varepsilon}= 2n$ \cfload.
This completes the proof of \cref{lem_cost_of_Lip_approx_set_equivalence_1}.
\end{proof}

\cfclear
\begin{athm}{lemma}{lem_cost_of_Lip_approx_set_equivalence}
Let $a \in C(\R,\R)$, $m,n \in \N$, $D \subseteq \R^m $, $f \in C(D,\R^n)$, $L,\varepsilon \in [0,\infty]$\cfload. %
Then it holds that $\CostLip{a}{A}{f}{L}{\varepsilon} < \infty$ if and only if there exists $\mathscr{f} \in \ANNs$ such that
\begin{enumerate}[(i)]
\item\label{lem_cost_Lip_app_set_equ_item_1} there exists $\mathfrak{m}\in\N$ such that $D\subseteq\R^\mathfrak{m}$ and $\realisation_a(\mathscr{f}) \in C(\R^\mathfrak{m},\R^n)$,
\item\label{lem_cost_Lip_app_set_equ_item_2} it holds for all $x \in D$ that $\mednorm{f(x) - \functionANNbis{\mathscr{f}}(x)} \leq \varepsilon$,
\item\label{lem_cost_Lip_app_set_equ_item_3} it holds for all $x,y \in D$ with $x \ne y$ that 
$\mednorm{ \functionANNbis{\mathscr{f}}(x)-\functionANNbis{\mathscr{f}}(y)} \leq L \mednorm{x-y}$, and
\item\label{lem_cost_Lip_app_set_equ_item_4} it holds that $\paramANN(\mathscr{f}) = \CostLip{a}{A}{f}{L}{\varepsilon}$
\end{enumerate}
\cfout.
\end{athm}
\begin{proof}[Proof of \cref{lem_cost_of_Lip_approx_set_equivalence}.]
Observe that the fact that $\paramANN(\ANNs)\subseteq\N$ ensures that for all $\mathscr{f}\in\ANNs$ with $\paramANN(\mathscr{f})=\CostLip{a}{A}{f}{L}{\varepsilon}$ it holds that
$
\CostLip{a}{A}{f}{L}{\varepsilon}=\paramANN(\mathscr{f})<\infty
$
\cfload.
Combining this with \eqref{Cost_Lipschitz_set} completes the proof of \cref{lem_cost_of_Lip_approx_set_equivalence}.
\end{proof}

\cfclear
\begin{athm}{cor}{cor_cost_of_Lip_approx_set_equivalence}
Let $a \in C(\R,\R)$, $m,n \in \N$, $D \subseteq \R^m $, $f \in C(D,\R^n)$, $L,\varepsilon \in [0,\infty]$ satisfy $D \neq \emptyset$ \cfload.
Then it holds that $\CostLip{a}{A}{f}{L}{\varepsilon} < \infty$ if and only if there exists $\mathscr{f} \in \ANNs$ such that
\begin{enumerate}[(i)]
\item\label{cor_cost_Lip_app_set_equ_item_1} it holds that $\realisation_a(\mathscr{f}) \in C(\R^m,\R^n)$,
\item\label{cor_cost_Lip_app_set_equ_item_2} it holds for all $x \in D$ that $\mednorm{f(x) - \functionANNbis{\mathscr{f}}(x)} \leq \varepsilon$,
\item\label{cor_cost_Lip_app_set_equ_item_3} it holds for all $x,y \in D$ with $x \ne y$ that 
$\mednorm{ \functionANNbis{\mathscr{f}}(x)-\functionANNbis{\mathscr{f}}(y)} \leq L \mednorm{x-y}$, and
\item\label{cor_cost_Lip_app_set_equ_item_4} it holds that $\paramANN(\mathscr{f}) = \CostLip{a}{A}{f}{L}{\varepsilon}$
\end{enumerate}
\cfout.
\end{athm}
\begin{proof}[Proof of \cref{cor_cost_of_Lip_approx_set_equivalence}.]
Observe that for all $\mathscr{f}\in\ANNs$ with $\realisation_a(\mathscr{f})\in C(\R^m,\R^n)$ there exists $\mathfrak{m}\in\N$ such that 
\begin{equation}
\label{cor_cost_of_Lip_approx_set_equivalence:eq1}
D\subseteq\R^\mathfrak{m}\qquad\text{and}\qquad\realisation_a(\mathscr{f})\in C(\R^\mathfrak{m},\R^n)
\end{equation}
\cfload[.]Furthermore, note that the fact that $\emptyset\neq D\subseteq\R^m$ ensures that for all $\mathfrak{m}\in\N$ with $D\subseteq\R^\mathfrak{m}$ it holds that $\mathfrak{m}=m$. Combining this and \eqref{cor_cost_of_Lip_approx_set_equivalence:eq1} with \cref{lem_cost_of_Lip_approx_set_equivalence} completes the proof of \cref{cor_cost_of_Lip_approx_set_equivalence}.
\end{proof}

\cfclear
\begin{athm}{lemma}{lem_cost_of_Lip_approx_set_equivalence_incodomain}
Let 
    $m,n \in \N$, 
	$D\subseteq\R^m$,
    $R \in [0,\infty)$, 
    $f \in C(D,\R^n)$, 
    $L,\varepsilon \in [0,\infty]$ 
satisfy 
    $f(D)\subseteq [-R,R]^n$ and 
    $\CostLip{\ReLU}{A}{f}{L}{\varepsilon} < \infty$ 
\cfload.
Then there exists 
    $\mathscr{f} \in \ANNs$ 
such that
\begin{enumerate}[(i)]
    \item\label{lem_cost_Lip_app_set_equ_incodomain_item_1} 
    it holds that 
		$\realisation_{\ReLU}(\mathscr{f}) \in C(\R^m,\R^n)$,
    \item\label{lem_cost_Lip_app_set_equ_incodomain_item_1a} 
    it holds for all 
		$x \in D$ 
	that 
		$(\realisation_\ReLU(\mathscr{f}))(x) \in [-R,R]^n$,
    \item\label{lem_cost_Lip_app_set_equ_incodomain_item_2} 
    it holds for all 
		$x \in D$ 
	that 
		$\mednorm{f(x) -(\realisation_\ReLU(\mathscr{f}))(x)} \leq \varepsilon$,
    \item\label{lem_cost_Lip_app_set_equ_incodomain_item_3} 
    it holds for all 
		$x,y \in D$ 
		with $x \ne y$ 
	that 
		$\mednorm{ (\realisation_\ReLU(\mathscr{f}))(x)-(\realisation_\ReLU(\mathscr{f}))(y)} \leq L \mednorm{x-y}$, and
	\item\label{lem_cost_Lip_app_set_equ_incodomain_item_4} 
	it holds that 
		$\paramANN(\mathscr{f}) = \CostLip{\mathfrak{r}}{A}{f}{L}{\varepsilon}+2n(n+1)$
\end{enumerate}
\cfout.
\end{athm}
\begin{aproof}
\Nobs that 
	\cref{lem:clipping_function} 
ensures that there exists 
    $\mathscr{g}\in\ANNs$ 
which satisfies that
\begin{enumerate}[(I)]
    \item\label{proof:lem_cost_of_Lip_approx_set_equivalence_incodomain:item_0}
    it holds that 
        $\realisation_\mathfrak{r}(\mathscr{g})\in C(\R^n,\R^n)$,
    \item\label{proof:lem_cost_of_Lip_approx_set_equivalence_incodomain:item_1}
    it holds for all 
        $x\in\R^n$ 
    that 
        $(\realisation_\mathfrak{r}(\mathscr{g}))(x)=\Clip{-R}Rn( x )$, and
    \item\label{proof:lem_cost_of_Lip_approx_set_equivalence_incodomain:item_2}
    it holds that 
        $\dims(\mathscr{g})=(n,n,n,n)$
\end{enumerate}
\cfload. 
Furthermore, \nobs that 
    \cref{cor_cost_of_Lip_approx_set_equivalence} and 
    the assumption that 
        $\CostLip{\ReLU}{A}{f}{L}{\varepsilon} < \infty$ 
show that there exists 
    $\mathscr{h}\in\ANNs$ 
which satisfies that
\begin{enumerate}[(A)]
    \item\label{lem_cost_Lip_app_set_equ_incodomain_proof_item_1} 
    it holds that 
        $\realisation_\mathfrak{r}(\mathscr{h}) \in C(\R^m,\R^n)$,
    \item\label{lem_cost_Lip_app_set_equ_incodomain_proof_item_2} 
    it holds for all 
        $x \in D$ 
    that 
        $\mednorm{f(x) - (\realisation_\mathfrak{r}(\mathscr{h}))(x)} \leq \varepsilon$,
    \item\label{lem_cost_Lip_app_set_equ_incodomain_proof_item_3} 
    it holds for all 
        $x,y \in D$ 
        with $x \ne y$ 
    that 
        $\mednorm{ (\realisation_\mathfrak{r}(\mathscr{h}))(x)-(\realisation_\mathfrak{r}(\mathscr{h}))(y)} \leq L \mednorm{x-y}$, and
    \item\label{lem_cost_Lip_app_set_equ_incodomain_proof_item_4} 
    it holds that 
        $\paramANN(\mathscr{h}) = \CostLip{\mathfrak{r}}{A}{f}{L}{\varepsilon}$
\end{enumerate}
\cfload.
Next \nobs that 
    \cref{Lemma:PropertiesOfCompositions_n2}, 
    \cref{proof:lem_cost_of_Lip_approx_set_equivalence_incodomain:item_0},
    and \cref{lem_cost_Lip_app_set_equ_incodomain_proof_item_1} 
establish that
\begin{equation}
    \eqlabel{i1}
    \realisation_{\ReLU}(\compANN{\mathscr g}{\mathscr h})
    \in C(\R^m,\R^n)
    \ifnocf.
\end{equation}
\cfload[.]%
Moreover, \nobs that 
    \cref{Lemma:PropertiesOfCompositions_n2}, 
    \cref{proof:lem_cost_of_Lip_approx_set_equivalence_incodomain:item_1},
    and the fact that $\Clip{-R}Rn(\R^n)\subseteq[-R,R]^n$ 
imply that for all 
    $x\in D$ 
it holds that
\begin{equation}
    \eqlabel{i2}
    (\realisation_{\mathfrak{r}}(\compANN{\mathscr g}{\mathscr h}))(x)
    =
    \prb{[\realisation_{\mathfrak{r}}(\mathscr{g})]\circ[\realisation_{\mathfrak{r}}(\mathscr{h})]}(x)
    =
    \Clip{-R}Rn\prb{(\realisation_{\mathfrak{r}}(\mathscr{h}))(x)}
    \in[-R,R]^n
    .
\end{equation}
    \Cref{lem_cost_Lip_app_set_equ_incodomain_proof_item_2},
    the assumption that $f(D)\subseteq[-R,R]^n$, 
    and the fact that 
        for all 
            $x,y\in\R^n$ 
        it holds that 
            $\mednorm{\Clip{-R}Rn(x)-\Clip{-R}Rn(y)}\leq \mednorm{x-y}$ 
    hence
ensure that for all 
    $x\in D$ 
it holds that
\begin{equation}
\eqlabel{i3}
\begin{split}
    \mednorm{f(x)-(\realisation_{\mathfrak{r}}(\compANN{\mathscr g}{\mathscr h}))(x)}
    &=
    \norm[\big]{\Clip{-R}Rn(f(x))-\Clip{-R}Rn\prb{(\realisation_{\mathfrak{r}}(\mathscr{h}))(x)}}
    \\&\leq
    \mednorm{f(x)-(\realisation_{\mathfrak{r}}(\mathscr{h}))(x)}\leq \varepsilon
    .
\end{split}
\end{equation}
In the next step \nobs that 
    \cref{lem_cost_Lip_app_set_equ_incodomain_proof_item_3} 
    and the fact that 
        for all 
            $x,y\in\R^n$ 
        it holds that 
            $\norm{\Clip{-R}Rn(x)-\Clip{-R}Rn(y)}\leq\norm{x-y}$ 
imply that for all 
    $x,y\in D$ with $x\neq y$
it holds that
\begin{equation}
\eqlabel{i4}
\begin{split}
    &\mednorm{(\realisation_{\mathfrak{r}}(\compANN{\mathscr g}{\mathscr h}))(x)-(\realisation_{\mathfrak{r}}(\compANN{\mathscr g}{\mathscr h}))(y)}
    \\&=
    \norm[\big]{\Clip{-R}Rn\prb{(\realisation_{\mathfrak{r}}(\mathscr{h}))(x)}-\Clip{-R}Rn\prb{(\realisation_{\mathfrak{r}}(\mathscr{h}))(y)}}
    \\&\leq
    \mednorm{(\realisation_{\mathfrak{r}}(\mathscr{h}))(x)-(\realisation_{\mathfrak{r}}(\mathscr{h}))(y)}\leq L\mednorm{x-y}.
\end{split}
\end{equation}
Furthermore, \nobs that 
    \cref{Lemma:PropertiesOfCompositions_n1}, 
    \cref{proof:lem_cost_of_Lip_approx_set_equivalence_incodomain:item_2}, 
    and \cref{lem_cost_Lip_app_set_equ_incodomain_proof_item_4} 
prove that
\begin{equation}
\eqlabel{i5}
\begin{split}
&\paramANN(\mathscr{g}\bullet\mathscr{h})
\\&=
\paramANN(\mathscr{g})+\paramANN(\mathscr{h})+\PR{\singledims_1(\mathscr{g})}(\singledims_{\lengthANN(\mathscr{h})-1}(\mathscr{h})+1)
\\&\quad-\PR{\singledims_1(\mathscr{g})}(\singledims_{0}(\mathscr{g})+1)-\PR{\singledims_{\lengthANN(\mathscr{h})}(\mathscr{h})}(\singledims_{\lengthANN(\mathscr{h})-1}(\mathscr{h})+1)
\\&=3n(n+1)+\paramANN(\mathscr{h})+n(\singledims_{\lengthANN(\mathscr{h})-1}(\mathscr{h})+1)-n(n+1)-n(\singledims_{\lengthANN(\mathscr{h})-1}(\mathscr{h})+1)
\\&=\CostLip{\ReLU}{A}{f}{L}{\varepsilon}+2n(n+1).
\end{split}
\end{equation}
Combining this with \eqqref{i1}, \eqqref{i2}, \eqqref{i3}, and \eqqref{i4}
establishes \cref{lem_cost_Lip_app_set_equ_incodomain_item_1,lem_cost_Lip_app_set_equ_incodomain_item_1a,lem_cost_Lip_app_set_equ_incodomain_item_2,lem_cost_Lip_app_set_equ_incodomain_item_3,lem_cost_Lip_app_set_equ_incodomain_item_4}.
The proof of \cref{lem_cost_of_Lip_approx_set_equivalence_incodomain} is thus complete.
\end{aproof}

\cfclear
\begin{athm}{lemma}{Lemma:Monotonicity_of_Cost}
Let $a \in C(\R,\R)$, $m,n \in \N$, $D,E \subseteq \R^m$, $f \in C(D,\R^n)$, $\varepsilon, L \in [0,\infty]$, $\eta \in [\varepsilon, \infty]$, $\mathfrak{L} \in [L, \infty]$ satisfy $E\subseteq D$.
Then
$
    \CostLip{a}{B}{f|_{E}}{\mathfrak{L}}{\eta}
    \leq
    \CostLip{a}{B}{f}{L}{\varepsilon}
$
\cfout.
\end{athm}

\begin{proof}
[Proof of \cref{Lemma:Monotonicity_of_Cost}]
Throughout this proof assume w.l.o.g.\ that $\CostLip{a}{B}{f}{L}{\varepsilon}<\infty$ and $D \neq \emptyset$ \cfadd{lem_cost_of_Lip_approx_set_equivalence_1}\cfload.
Note that \cref{cor_cost_of_Lip_approx_set_equivalence} ensures that there exists $\mathscr{f} \in \ANNs$ which satisfies that
\begin{enumerate}[(i)]
\item it holds that $\realisation_a(\mathscr{f}) \in C(\R^m,\R^n)$,
\item it holds for all $x \in D$ that $\mednorm{f(x) - \functionANNbis{\mathscr{f}}(x)} \leq \varepsilon$,
\item it holds for all $x,y \in D$ with $x \ne y$ that 
$\mednorm{ \functionANNbis{\mathscr{f}}(x)-\functionANNbis{\mathscr{f}}(y)} \leq L \mednorm{x-y}$, and
\item it holds that $\paramANN(\mathscr{f}) = \CostLip{a}{D}{f}{L}{\varepsilon}$
\end{enumerate}
\cfload.
\Nobs that the assumption that $\varepsilon \leq \eta$, $L\leq\mathfrak{L}$, and $E\subseteq D$ demonstrates that
\begin{equation}\begin{split}
    &\CostLip{a}{B}{f|_E}{\mathfrak{L}}{\eta}
    \\&=
    \min \pr*{ \pR*{ \mathfrak{p} \in \N \colon \PR*{
	\begin{array}{cc}
	    \exists\, \mathscr{g} \in \ANNs, \mathfrak m\in\N \colon
	    \PRb{(\realisation_a(\mathscr{g}) \in C(\R^{\mathfrak m},\R^n))\land{}\\
			(\mathfrak{p}=\paramANN(\mathscr{g}))\land(\fa{x,y} E \text{ with } x \ne y \colon\\
			\mednorm{ \functionANNbis{\mathscr{g}}(x)-\functionANNbis{\mathscr{g}}(y)} \leq \mathfrak{L}\norm{x-y})\land{}\\
         (\fa{x} E \colon 
        \mednorm{f(x) -  \functionANNbis{\mathscr{g}}(x)} \leq \eta)}\\
    \end{array} }
} \cup \{\infty\} }
    \\&\leq 
    \paramANN(\mathscr{f})
    =
    \CostLip{a}{B}{f}{L}{\varepsilon}.
    \end{split}\end{equation}
This completes the proof of \cref{Lemma:Monotonicity_of_Cost}.
\end{proof}

\cfclear

\cfclear
\begin{athm}{lemma}{lemma_Lip_catalogs_are_Lip_set}
	Let $a \in C(\R,\R)$, $m,n \in \N$, $D \subseteq \R^m$, $f \in C(D,\R^n)$, $L \in [0,\infty]$, $c \in (0, \infty)$ satisfy for all $\varepsilon \in (0,c)$ that $\CostLip{a}{D}{f}{L}{\varepsilon} < \infty$ \cfload.
    Then it holds for all $x,y \in D$ with $x \ne y$ that $\mednorm{f(x)-f(y)} \leq L \mednorm{x-y}$ \cfout.
\end{athm}

\begin{proof}[Proof of \cref{lemma_Lip_catalogs_are_Lip_set}.]
Throughout this proof assume w.l.o.g.\ that $D \neq \emptyset$. Observe that the assumption that for all $\varepsilon \in (0,c)$ it holds that $\CostLip{a}{D}{f}{L}{\varepsilon} < \infty$ and \cref{cor_cost_of_Lip_approx_set_equivalence} imply that there exist $ (\mathscr{f}_{\varepsilon})_{\varepsilon \in (0,c)} \subseteq \ANNs$ %
which satisfy that
\begin{enumerate}[(i)]
\item \label{lemma_Lip_catalogs_are_Lip:1} 
    it holds for all 
        $\varepsilon \in (0,c)$ 
    that 
        $\realisation_{a}(\mathscr{f}_{\varepsilon}) \in C(\R^{m},\R^{n})$,

\item \label{lemma_Lip_catalogs_are_Lip:3} 
    it holds for all 
        $\varepsilon \in (0,c)$, 
        $x \in D$ 
    that 
        $\mednorm{f(x) - \functionANNbis{\mathscr{f}_{\varepsilon}}(x)} \leq \varepsilon$, and

\item \label{lemma_Lip_catalogs_are_Lip:2} 
    it holds for all 
        $\varepsilon \in (0,c)$, 
        $x,y \in D$ 
        with $x \ne y$ 
    that 
        $\mednorm{\functionANNbis{\mathscr{f}_{\varepsilon}}(x) - \functionANNbis{\mathscr{f}_{\varepsilon}}(y)} \leq L \mednorm{x-y}$
\end{enumerate}
\cfload.
Note that \cref{lemma_Lip_catalogs_are_Lip:3} and \cref{lemma_Lip_catalogs_are_Lip:2} ensure that for all $\varepsilon \in (0,c)$, $x,y \in D$ with $x \ne y$ it holds that
\begin{equation}\begin{split}
	&\mednorm{f(x) - f(y)} 
	\\&\leq 
    \norm{f(x) - \functionANNbis{\mathscr{f}_{\varepsilon}}(x)}
	+ \norm{\functionANNbis{\mathscr{f}_{\varepsilon}}(x) - \functionANNbis{\mathscr{f}_{\varepsilon}}(y)} + \norm{\functionANNbis{\mathscr{f}_{\varepsilon}}(y) - f(y)} 
    \\&\leq 
    \varepsilon + L \mednorm{x-y} + \varepsilon
    .
\end{split}\end{equation}
Hence, we obtain that for all $x,y \in D$ with $x \ne y$ it holds that
\begin{equation}
\begin{split}
	\mednorm{f(x) - f(y)} &= \limsup_{\varepsilon \searrow 0}  \mednorm{f(x) - f(y)} 
	\leq \limsup_{\varepsilon \searrow 0} \PR{ 2\varepsilon + L \mednorm{x-y} } = L \mednorm{x-y}.
\end{split}
\end{equation}
This completes the proof of \cref{lemma_Lip_catalogs_are_Lip_set}.
\end{proof}

\subsection{Spaces of ANN approximable functions}
\label{Subsection:3.3}

In this section we first introduce ANN approximation spaces for one-dimensional functions, which are, roughly speaking, sets of functions which can be approximated by ANNs with the ReLU activation function such that the number of parameters grows at most polynomially in the diameter of the compact set on which the function is approximated and in the reciprocal of the prescribed approximation error (see \cref{def:polyC} below for a precise definition). 
We then turn to multi-dimensional functions, where it will often be useful to represent sequences of functions $f_1\in C(\R,\R^{d_1})$, $f_2\in C(\R^2,\R^{d_2})$, $f_3\in C(\R^3, \R^{d_3})$, $\ldots$ with $d_1,d_2,d_3,\ldots\in\N$ as a single function $F\colon \bigcup_{n\in\N}\R^n\to\bigcup_{n\in\N}\R^{n}$ with $F(x)=f_n(x)$ for all $n\in\N$ and $x\in\R^n$, where $\bigcup_{n\in\N}\R^n$ is equipped with the topology of the disjoint union (considering the utility of this representation, note, e.g., that such functions can be composed using standard function composition). \cref{Lemma:codomain} is a technical result justifying this equivalence. Finally, we define the ANN approximation space for multi-dimensional functions, which serves as our central formalization of the set of sequences of functions that can be approximated by ANNs with the ReLU activation function without the curse of dimensionality. In more detail, the ANN approximation space for multi-dimensional functions can be considered a set of sequences of functions which can be approximated by ANNs with the ReLU activation function such that the number of parameters grows at most polynomially in the dimension of the domain, in the diameter of the compact set on which the function is approximated, and in the reciprocal of the prescribed approximation error (see \cref{def:polyD_mult} below for a precise definition).
Apart from giving these definitions, which are central to the rest of this article, we also supply some useful, but straightforward technical results.

\begin{definition}
[ANN approximation spaces for one-dimensional functions]
\label{def:polyC}
Let $c,r \in [0, \infty)$. Then we denote by $\Capprox{c}{r}$ the set given by
\begin{multline}
\label{polyC:1}
    \Capprox{c}{r}= 
    \pRb{
    f \in C(\R,\R) \colon \PRb{
    \forall \, \radius \in [r, \infty),\varepsilon \in (0,1] \colon
    \CostLip{\ReLU}{[-\radius, \radius]}{f|_{\PR{-R,R}}} {c\radius^{c}}{\varepsilon} \leq c\radius^c \varepsilon^{-c} }
    }
\end{multline}
(cf.\ \cref{def:ReLU,def:cost_Lipschitz_approx_set}).
\end{definition}

\cfclear
\begin{athm}{cor}{cor:polyC}
Let $c,r \in [0, \infty)$\cfload. Then
\begin{enumerate}[(i)]
\item\label{cor:polyC_item_1b} it holds for all $R\in[r,\infty)$ that $\Capprox{c}{r}\subseteq \Capprox{c}{R}$, 
\item\label{cor:polyC_item_1a} it holds for all $C\in[c,\infty)$ that $\Capprox{c}{r}\subseteq \Capprox{C}{\max\{1,r\}}$, and
\item\label{cor:polyC_item_2} it holds that $\Capprox{c}{r}\subseteq \Capprox{\max\{c,cr^c\}}{1}$
\end{enumerate}
\cfout.
\end{athm}

\begin{proof}
[Proof of \cref{cor:polyC}]
First, \nobs that \eqref{polyC:1} establishes \cref{cor:polyC_item_1b}.
Next \nobs that \cref{Lemma:Monotonicity_of_Cost} shows that for all $\rho\in[1,\infty)$, $C\in[c,\infty)$, $f\in \Capprox{c}{\rho}$, $\radius\in[\rho,\infty)$, $\varepsilon\in(0,1]$ it holds that
\begin{equation}
\begin{split}
\CostLip{\ReLU}{[-\radius,\radius] }{f|_{[-R,R]}}{C\radius^C}{\varepsilon}\leq\CostLip{\ReLU}{[-\radius,\radius] }{f|_{[-R,R]}}{c\radius^c}{\varepsilon}\leq c\radius^c\varepsilon^{-c}\leq C\radius^C\varepsilon^{-C}
\end{split}
\end{equation}
\cfload. Hence, we obtain that for all $\rho\in[1,\infty)$, $C\in[c,\infty)$ it holds that $\Capprox{c}{\rho}\subseteq\Capprox{C}{\rho}$. Combining this with \cref{cor:polyC_item_1b} ensures that for all $C\in[c,\infty)$ it holds that
\begin{equation}
\begin{split}
\Capprox{c}{r}\subseteq\Capprox{c}{\max\{1,r\}}\subseteq\Capprox{C}{\max\{1,r\}}.
\end{split}
\end{equation}
This establishes \cref{cor:polyC_item_1a}. 
\Nobs that 
    the assumption that $c\geq 0$
implies that
for all $\radius \in [1, \infty)$ it holds that $\radius\leq \max\{1,r\}\radius\in [r, \infty)$ and
\begin{equation}
c(\max\{1,r\}\radius)^c= \max\{c,cr^c\}\radius^{c}\leq \max\{c,cr^c\} \radius^{\max\{c,cr^c\}}.
\end{equation}
This, \eqref{polyC:1}, and \cref{Lemma:Monotonicity_of_Cost} (applied with 
$a \curvearrowleft \ReLU$,
$m \curvearrowleft 1$,
$n \curvearrowleft 1$,
$D \curvearrowleft [-\mathord{\max}\{1,r\}\radius,\allowbreak\max\{1,r\}\radius]$,
$E \curvearrowleft [-\radius,\radius]$,
$f\is f|_{[-\mathord{\max}\{1,r\}\radius,\max\{1,r\}\radius]}$,
$\varepsilon \curvearrowleft \varepsilon$,
$L \curvearrowleft c\pr{\max\{1,r\}\radius}^c$,
$\eta \curvearrowleft \varepsilon$,
$\mathfrak{L} \curvearrowleft \max\{c,cr^c\}\radius^{\max\{c,cr^c\}}$
for $f\in\Capprox{c}{r}$, $\radius\in[1,\infty)$, $\varepsilon\in(0,1]$
in the notation of \cref{Lemma:Monotonicity_of_Cost}) demonstrate that for all $f\in \Capprox{c}{r}$, $\radius \in [1, \infty)$, $\varepsilon \in (0,1]$ it holds that
\begin{equation}
\begin{split}
    &\CostLipB{\ReLU}{[-\radius,\radius] }{f|_{[-R,R]}}{\max\{c,cr^c\} \radius^{\max\{c,cr^c\}}}{\varepsilon}
    \\&\leq
    \CostLipB{\ReLU}{[-\mathord{\max}\{1,r\}\radius,\max\{1,r\}\radius] }{ f|_{[-\mathord{\max}\{1,r\}\radius,\max\{1,r\}\radius]}}{c(\max\{1,r\}\radius)^c}{\varepsilon}
    \\&\leq c(\max\{1,r\}\radius)^c \varepsilon^{-c}\leq \max\{c,cr^c\}\radius^{\max\{c,cr^c\}}\varepsilon^{-\mathord{\max}\{c,cr^c\}}.
\end{split}
\end{equation}
Therefore, we obtain \cref{cor:polyC_item_2}.
The proof of \cref{cor:polyC} is thus complete.
\end{proof}
\cfclear

\begin{athm}{lemma}{Lemma:codomain}
Let\footnote{Here and in the rest of the article, $\bigcup_{n\in\N}\R^n$ is equipped with the
topology of the disjoint union.} $f \in C\prb{\bigcup_{n \in \N}\R^n, \bigcup_{n \in \N}\R^n}$, $d \in \N$. Then there exists a unique $\delta \in \N$ such that $f(\R^d) \subseteq \R^{\delta}$.
\end{athm}
\begin{proof}
[Proof of \cref{Lemma:codomain}]
Throughout this proof let $\delta \in \N$ satisfy $f(1,2, \ldots, d) \allowbreak  \in \R^{\delta}$. 
Observe that 
    the assumption that 
        $f \in C\prb{\bigcup_{n \in \N}\R^n, \bigcup_{n \in \N}\R^n}$ 
    and the fact that 
        $\R^d$ is connected 
imply that 
    $f(\R^d)$ is connected. 
    The fact that
    \begin{equation}
    \label{Lemma:codomain:disjoint_decomposition}
    \begin{split}
        f\pr{\R^d}
        =
        \pr{f\pr{\R^d}\cap\R^{\delta}}\cup\prb{f\pr{\R^d}\cap\prb{\textstyle\bigcup_{n \in \N\setminus\{\delta\}}\R^n}},
    \end{split}
    \end{equation}
    the fact that 
        $\pr{f\pr{\R^d}\cap\R^{\delta}}$, $\prb{f\pr{\R^d}\cap\prb{\bigcup_{n \in \N\setminus\{\delta\}}\R^n}}\subseteq f\pr{\R^d}$ are open sets with respect to the subspace topology on $f(\R^d)$, 
    the fact that 
        $\pr{f\pr{\R^d}\cap\R^{\delta}}\cap\prb{f\pr{\R^d}\cap\prb{\bigcup_{n \in \N\setminus\{\delta\}}\R^n}}=\emptyset$, 
    and the fact that 
        $f(1,2, \ldots, d)\in(f(\R^d)\cap\R^{\delta})$ 
    hence 
ensure that 
\begin{equation}
    f(\R^d)\cap\R^{\delta}=f(\R^d)
    \qquad\text{and}\qquad
    f\pr{\R^d}\cap\prb{\textstyle\bigcup_{n \in \N\setminus\{\delta\}}\R^n}=\emptyset
    .
\end{equation}
The proof of \cref{Lemma:codomain} is thus complete. 
\end{proof}

\begin{definition}
[Image dimensions]
\label{def:smart_restriction}\cfadd{Lemma:codomain}\cfconsiderloaded{def:smart_restriction}
Let $f \in C\prb{\bigcup_{n \in \N}\R^n, \bigcup_{n \in \N}\R^n}$, $d \in \N$. Then we denote by $\imdim{f}{d} \in \N$ the unique natural number which satisfies $f(\R^d) \subseteq \R^{\imdim{f}{d}}$ \cfload.
\end{definition}
\cfclear

\begin{definition}
[Ceiling of real numbers]
\label{def:ceiling}
We denote by $\ceil{\cdot} \colon \R \to \Z$ the function which satisfies for all $x \in \R$ that $\ceil{x} = \min(\Z \cap [x, \infty))$.
\end{definition}

\begingroup
\newcommand{\const}{C}
\cfclear
\begin{athm}{lemma}{Lemma:extension_on_R_poly_a}
Let 
    $a\in C(\R,\R)$,
	$m,n\in\N$,
	$f \in C( \R^m,\R^n)$, 
	$c,\const ,r \in [0,\infty)$,
	$\eps\in(0,1]$
	satisfy for all 
		$\radius \in \N \cap [r, \infty)$
	that 
		$\const =\max\{1,c(\max\{2,2\ceil{r}\})^{\max\{1,c\}}\}$ 
	and 
    \begin{equation}
    \label{extension_on_R_poly_a:1}
        \CostLipA{a}{[-\radius,\radius]^m }
		{f|_{[-\radius,\radius]^m}}
		{cm^c \! \radius^c}
		{\varepsilon} 
		\leq 
		cm^c \! \radius^c \varepsilon^{-c}
    \end{equation}
\cfload.
Then it holds for all 
	$\radius \in [1, \infty)$
that
\begin{equation}
    \CostLipA{a}{[-\radius,\radius]^m }
	{f|_{[-\radius,\radius]^m}}
	{\const m^\const  \! \radius^\const }
	{\varepsilon} 
	\leq 
	\const m^\const  \! \radius^\const  \varepsilon^{-\const }.
\end{equation}
\end{athm}
\begin{aproof}
Throughout this proof assume w.l.o.g.\ that $r\geq 1$.
Note that 
    the fact that $C\geq c2\ceil r^{\max\{1,c\}}$
ensures that for all 
	$\radius \in [1, \infty)$ 
it holds that 
	$\radius\leq\ceil{r\radius}\in \N \cap [r, \infty)$ and
\begin{equation}
	cm^c\ceil{r\radius}^c
	\leq 
	c\ceil{r}^{\max\{1,c\}}m^c\ceil{\radius}^{\max\{1,c\}}
	\leq 
	c\pr{2\ceil{r}}^{\max\{1,c\}}m^c\!\radius^{\max\{1,c\}}
	\leq
	\const m^{\const }\!\radius^{\const }
	.
\end{equation}
	This, 
	\eqref{extension_on_R_poly_a:1}, and 
	\cref{Lemma:Monotonicity_of_Cost} 
show that for all 
	$\radius \in [1, \infty)$
it holds that
\begin{equation}
\begin{split}
    \CostLipA{a}{[-\radius,\radius]^m }
	{f|_{[-\radius,\radius]^m}}
	{\const m^\const  \! \radius^\const }
	{\varepsilon}
    &\leq
    \CostLipA{a}{[-\ceil{r\radius},\ceil{r\radius}]^m }
	{f|_{[-\ceil{r\radius},\ceil{r\radius}]^m}}
	{cm^c\ceil{r\radius}^c}
	{\varepsilon}
    \\&\leq 
	cm^c\ceil{r\radius}^c \varepsilon^{-c}
	\leq 
	\const m^\const \!\radius^{\const }\varepsilon^{-\const }.
\end{split}
\end{equation}
\finishproofthis
\end{aproof}
\endgroup

\cfclear
\begin{athm}{lemma}{Lemma:extension_on_R_poly}
Let $a\in C(\R,\R)$, $f \in C\prb{\bigcup_{d \in \N}\R^d, \bigcup_{d \in \N}\R^d}$.
Then the following three statements are equivalent:
\begin{enumerate}[(i)]
    \item \label{extension_on_R_poly:item_1}
    There exists $c \in [0,\infty)$ such that for all $d,\radius \in \N$, $\varepsilon \in (0,1]$ it holds that
    \begin{equation}
        \CostLipA{a}{[-\radius,\radius]^d }{\pr{[-\radius,\radius]^d \ni x \mapsto f(x) \in \R^{\imdim{f}{d}}}}{cd^c \! \radius^c}{\varepsilon} \leq c d^c \! \radius^c \varepsilon^{-c}
    \end{equation}
    \cfload.
    
    \item \label{extension_on_R_poly:item_2}
    There exist $c,r \in [0,\infty)$ such that for all $d \in \N$, $\radius \in \N \cap [r, \infty)$, $\varepsilon \in (0,1]$ it holds that
    \begin{equation}
        \CostLipA{a}{[-\radius,\radius]^d }{\pr{[-\radius,\radius]^d \ni x \mapsto f(x) \in \R^{\imdim{f}{d}}}}{cd^c \! \radius^c}{\varepsilon} \leq c d^c \! \radius^c \varepsilon^{-c}.
    \end{equation}
    
    \item \label{extension_on_R_poly:item_3}
    There exists $c \in [0,\infty)$ such that for all $d \in \N$, $\radius \in [1, \infty)$, $\varepsilon \in (0,1]$ it holds that
    \begin{equation}
        \CostLipA{a}{[-\radius,\radius]^d }{\pr{[-\radius,\radius]^d \ni x \mapsto f(x) \in \R^{\imdim{f}{d}}}}{cd^c \! \radius^c}{\varepsilon} \leq c d^c \! \radius^c \varepsilon^{-c}.
    \end{equation}
\end{enumerate}
\end{athm}

\begin{aproof}
First observe that 
	\cref{extension_on_R_poly:item_3} implies \cref{extension_on_R_poly:item_1} 
and 
	\cref{extension_on_R_poly:item_1} implies \cref{extension_on_R_poly:item_2}. 
Moreover, note that 
	\cref{Lemma:extension_on_R_poly_a}
	(applied with
		$m\is d$,
		$n\is \imdim fd$,
		$f\is (\R^d\ni x\mapsto f(x)\in\R^{\imdim fd})$,
		$c\is c$,
		$r\is r$,
		$\eps\is\eps$,
		for $d\in\N$, $R\in[1,\infty)$, $\eps\in(0,1]$
	in the notation of \cref{Lemma:extension_on_R_poly_a})
shows that 
	\cref{extension_on_R_poly:item_2} implies \cref{extension_on_R_poly:item_3}.
The proof of \cref{Lemma:extension_on_R_poly} is thus complete.
\end{aproof}

\cfclear
\begin{definition}
[ANN approximation space for multi-dimensional functions]
\label{def:polyD_mult}
\cfconsiderloaded{def:polyD_mult}
We denote by $\Dapprox$ the set given by
\begin{equation}
\begin{split}
\label{polyD_mult:1}
    &\Dapprox=
		\pRbbb{f \in C\prb{ \textstyle{\bigcup_{d \in \N}\R^d, \bigcup_{d \in \N}\R^d}} \colon
		\\
		&\PR*{
		\begin{array}{c}
		\exists \, c \in [0, \infty)\colon \forall \, d, \radius \in \N, \varepsilon \in (0,1]  \colon \PRb{
				\prb{f([-\radius,\radius]^d)\subseteq[-cd^c \! \radius^c,cd^c \! \radius^c]^{\imdim{f}{d}}}\land{}
				\\
				(\imdim{f}{d}\leq cd^c)\land
				\prb{\CostLipA{\ReLU}{[-\radius,\radius]^d }{\pr{[-\radius,\radius]^d \ni x \mapsto f(x) \in \R^{\imdim{f}{d}}}}{cd^c \! \radius^c}{\varepsilon}{}
				\leq{} c d^c \! \radius^c \varepsilon^{-c}}}
		\end{array}
    }{}}
		\end{split}
\end{equation}
\cfload.
\end{definition}

\cfclear

\subsection{ANN approximations for compositions of functions}
\label{Subsection:3.4}

In this section, we use the upper bounds on the number of parameters of the composition of ANNs established in \cref{subsec:compositions} above to derive upper bounds on the cost of approximating the composition of two functions by ANNs in terms of the costs of approximating each of the two functions by ANNs.

\begin{athm}{prop}{Prop:composition_cost_ReLU}
Let $d_1,d_2,d_3 \in \N$, 
$\varepsilon, L_1, L_2, R_1, R_2 \in [0, \infty)$ and
let
$f_1 \colon [-R_1,R_1]^{d_1} \to \R^{d_2}$ and
$f_2 \colon [-R_2,R_2]^{d_2} \to \R^{d_3}$ 
satisfy $f_1([-R_1,R_1]^{d_1})\subseteq[-R_2,R_2]^{d_2}$.
Then\footnote{Note that for all sets $A, B, C, D$ and all
functions $f\colon A\to B$ and $g\colon C\to D$
such that $f(A)\subseteq C$ we denote by $g\circ f\colon A\to D$
the function which satisfies for all $a\in A$ that
$(g\circ f)(a)=g(f(a))$.}
\begin{equation}
\begin{split}
\label{Prop:composition_cost_ReLU:eq}
    &\CostLip{\ReLU}{A}{f_2 \circ f_1}{L_2L_1}{\varepsilon}
\leq
\CostLipB{\ReLU}{B}{f_1}{L_1}{\tfrac{\varepsilon}{2L_2}}
+
\CostLipB{\ReLU}{B}{f_2}{L_2}{\tfrac{\varepsilon}{2}}
+
4d_2(d_2+1)
\end{split}
\end{equation}
\cfout.

\end{athm}

\begin{proof}
[Proof of \cref{Prop:composition_cost_ReLU}]
Throughout this proof let $a_1,a_2 \in \R$ satisfy $a_1 = 2L_2$ and $a_2 = 2$ and 
assume w.l.o.g.\ that 
$\CostLipB{\ReLU}{A}{f_2}{L_2}{\tfrac{\varepsilon}{a_2}} + \CostLipB{\ReLU}{B}{f_1}{L_1}{\tfrac{\varepsilon}{a_1}} < \infty$ \cfload. 
Observe that \cref{cor_cost_of_Lip_approx_set_equivalence}, \cref{lem_cost_of_Lip_approx_set_equivalence_incodomain}, 
and the assumption that $\CostLipB{\ReLU}{A}{f_2}{L_2}{\tfrac{\varepsilon}{a_2}} + \CostLipB{\ReLU}{B}{f_1}{L_1}{\tfrac{\varepsilon}{a_1}} < \infty$ 
imply that there exist
$\mathscr{f}_1, \mathscr{f}_2 \in \ANNs$ which satisfy that
\begin{enumerate}[(I)]
\item \label{composition_cost_ReLU:item1} it holds for all $k \in \{1,2\}$ that $ \realisation_\ReLU(\mathscr{f}_k) \in C(\R^{d_k},\R^{d_{k+1}})$,

\item \label{composition_cost_ReLU:item1a} it holds for all $x \in [-R_1,R_1]^{d_1}$ that $ (\realisation_\ReLU(\mathscr{f}_1))(x) \in [-R_2,R_2]^{d_2}$,

\item \label{composition_cost_ReLU:item2} it holds for all $k \in \{1,2\}$, $x,y \in [-R_k,R_k]^{d_k}$ that $\mednorm{\functionnbReLUANN{\mathscr{f}_k} (x) - \functionnbReLUANN{\mathscr{f}_k} (y)} \leq L_k \mednorm{x-y}$,

\item \label{composition_cost_ReLU:item3} it holds for all $k \in \{1,2\}$, $x \in [-R_k,R_k]^{d_k}$ that $\mednorm{f_k(x) - \functionnbReLUANN{\mathscr{f}_k} (x)} \leq \frac{\varepsilon}{a_k}$, and

\item \label{composition_cost_ReLU:item4} it holds for all $k \in \{1,2\}$ that $\paramANN(\mathscr{f}_k) \leq \CostLipB{\ReLU}{A}{f_k}{L_k}{\tfrac{\varepsilon}{a_k}}+(2-k)2d_2(d_2+1)$
\end{enumerate}
\cfload. Note that \cref{composition_cost_ReLU:item1}, \cref{Prop:identity_representation}, and \cref{Lemma:PropertiesOfCompositions_n2} ensure that for all $x \in \R^{d_1}$ it holds that $\realisation_{\ReLU}(\mathscr{f}_2\bullet\mathbb{I}_{d_2}\bullet \mathscr{f}_1)\in C(\R^{d_1},\R^{d_{3}})$ and 
\begin{equation}
\label{composition_cost_ReLU:a}
    \begin{split}
        (\realisation_{\ReLU}(\compANN{\mathscr{f}_2 }{\ReLUidANN{d_2}}\bullet \mathscr{f}_1))(x) = \prb{[\realisation_{\ReLU}(\mathscr{f}_2)] \circ [\realisation_{\ReLU}(\mathscr{f}_1)]}(x)
    \end{split}
\end{equation}
\cfload. 
    This,
    the triangle inequality,
    \cref{composition_cost_ReLU:item1a}, 
    \cref{composition_cost_ReLU:item2}, and
    \cref{composition_cost_ReLU:item3}
show that
\begin{equation}
\label{composition_cost_ReLU:1}
\begin{split}
    &      \sup_{x \in [-R_1,R_1]^{d_1}} \mednorm{(f_2 \circ f_1)(x) - (\realisation_{\ReLU}(\mathscr{f}_2 \bullet\mathbb{I}_{d_2}\bullet \mathscr{f}_1))(x) }
    \\&= \sup_{x \in [-R_1,R_1]^{d_1}} \bignorm{\prb{[\realisation_{\ReLU}(\mathscr{f}_2)] \circ [\realisation_{\ReLU}(\mathscr{f}_1)]}(x)-(f_2 \circ f_1)(x) }
		\\&\leq \PRbbbb{\sup_{x \in [-R_1,R_1]^{d_1}} \bignorm{\prb{[\realisation_{\ReLU}(\mathscr{f}_2)] \circ [\realisation_{\ReLU}(\mathscr{f}_1)]}(x)-([\realisation_{\ReLU}(\mathscr{f}_2)] \circ f_1)(x) }}
		\\&\quad +
		\PRbbbb{\sup_{x \in [-R_1,R_1]^{d_1}} \bignorm{\prb{[\realisation_{\ReLU}(\mathscr{f}_2)] \circ f_1}(x)-(f_2 \circ f_1)(x) }	}	
    \\&\leq \PRbbbb{\sup_{x \in [-R_1,R_1]^{d_1}} L_2\bignorm{\pr{\realisation_{\ReLU}(\mathscr{f}_1)}(x)-f_1(x) }}
		+
		\PRbbbb{\sup_{y \in [-R_2,R_2]^{d_2}} \bignorm{\pr{\realisation_{\ReLU}(\mathscr{f}_2)} (y)-f_2 (y) }	}
    \\&\leq \tfrac{L_2\varepsilon}{a_1} + \tfrac{\varepsilon}{a_2}
    \leq \tfrac{\varepsilon}{2} + \tfrac{\varepsilon}{2}
    = \varepsilon
    .
\end{split}
\end{equation}
Furthermore, note that 
    \cref{composition_cost_ReLU:item2} and 
    \eqref{composition_cost_ReLU:a} 
imply that for all 
    $x,y \in [-R_1,R_1]^{d_1}$ 
it holds that
\begin{equation}
\label{composition_cost_ReLU:2}
\begin{split}
    &\mednorm{(\realisation_{\ReLU}(\mathscr{f}_2 \bullet\mathbb{I}_{d_2}\bullet \mathscr{f}_1))(x) - (\realisation_{\ReLU}(\mathscr{f}_2 \bullet\mathbb{I}_{d_2}\bullet \mathscr{f}_1))(y)}
    \\&= 
    \bignorm{\prb{[\realisation_{\ReLU}(\mathscr{f}_2)] \circ [\realisation_{\ReLU}(\mathscr{f}_1)]}(x)
    -\prb{[\realisation_{\ReLU}(\mathscr{f}_2)] \circ [\realisation_{\ReLU}(\mathscr{f}_1)]}(y)}
    \\&\leq
    L_2\mednorm{(\realisation_{\ReLU}(\mathscr{f}_1))(x)-(\realisation_{\ReLU}(\mathscr{f}_1))(y)} \leq L_2L_1 \mednorm{x-y}.
\end{split}
\end{equation}
Combining this with 
    \eqref{Cost_Lipschitz_set}
    and \eqref{composition_cost_ReLU:1} 
ensures that
\begin{equation}
\begin{split}
\label{cost_bounded_by_example}
    \CostLip{\ReLU}{A}{f_2 \circ f_1}{L_2L_1}{\varepsilon} 
     &\leq \paramANN(\mathscr{f}_2 \bullet\mathbb{I}_{d_2}\bullet \mathscr{f}_1) 
		.
\end{split}
\end{equation}
Moreover, observe that \cref{PropertiesOfCompositions_id_n:Params2} of \cref{Lemma:PropertiesOfCompositions_n3} (applied with $n \curvearrowleft 2$, $\mathscr{f}_1 \curvearrowleft \mathscr{f}_2$, $\mathscr{f}_2 \curvearrowleft \mathscr{f}_1$ in the notation of \cref{Lemma:PropertiesOfCompositions_n3}) and \cref{composition_cost_ReLU:item4} show that
\begin{equation}
\begin{split}
    \paramANN(\mathscr{f}_2 \bullet\mathbb{I}_{d_2}\bullet \mathscr{f}_1) 
    &\leq 
    2\,\paramANN(\mathscr{f}_1)+2\,\paramANN(\mathscr{f}_2)
    \\&\leq 
    2 \, \CostLipB{\ReLU}{B}{f_1}{L_1}{\tfrac{\varepsilon}{a_1}} 
    + 2\prb{\CostLipB{\ReLU}{A}{f_2}{L_2}{\tfrac{\varepsilon}{a_2}}+2d_2(d_2+1)}
    \\&= 
    2\, \CostLipB{\ReLU}{B}{f_1}{L_1}{\tfrac{\varepsilon}{a_1}}
		+2 \, \CostLipB{\ReLU}{A}{f_2}{L_2}{\tfrac{\varepsilon}{a_2}}
		+4d_2(d_2+1).
\end{split}
\end{equation}
This and \eqref{cost_bounded_by_example} establish \eqref{Prop:composition_cost_ReLU:eq}. The proof of \cref{Prop:composition_cost_ReLU} is thus complete.
\end{proof}

\subsection{Compositions of ANN approximable functions}
\label{Subsection:3.5}

In this section, we employ the upper bounds on the cost of a composition of two functions from the previous subsection to establish that the ANN approximation space for multi-dimensional functions (see \cref{def:polyD_mult}) is closed under composition. This central result in \cref{lemma:D_closed} is proved with the help of two auxiliary lemmas, \cref{Lemma:Comp_theory_1} and \cref{lemma:catching_image}.

\cfclear
\begingroup
\newcommand{\df}{{m}}
\newcommand{\dfp}{{m}}
\newcommand{\dg}{{n}}
\newcommand{\dgp}{{n}}
\newcommand{\ig}{{o}}
\newcommand{\f}{f}
\newcommand{\g}{g}
\renewcommand{\c}{c}
\newcommand{\C}{C}
\newcommand{\cf}{{\mathfrak f}}
\newcommand{\cg}{{\mathfrak g}}
\renewcommand{\r}{R}
\begin{athm}{lemma}{Lemma:Comp_theory_1}
    Let
        $\df,\dg,\ig,\c,\C\in\N$,
        $\cf,\cg\in(0,\infty)$,
        $\eps\in(0,1]$,
    let
        $\f\colon \R^{\df}\to\R^{\dg}$
        and $\g\colon \R^{\dg}\to\R^{\ig}$
        be functions,
    and assume for all $\r \in \N$, $\delta\in(0,\eps]$ that 
        $\c=\ceil{\max\{\cf,\cg,2\}}$, 
        $\C =2^{\c+2}\c^{(2\c^2+\c+1)}$,
        $\dg\leq \cf\dfp^\cf$,
        $\f([-\r,\r]^\df)\subseteq [-\cf \dfp^{\cf }\!\r^{\cf },\cf \dfp^{\cf }\!\r^{\cf }]^{\dg}$,
        $\CostLipB{\ReLU}{[-\r,\r]^\df }{\f|_{[-\r,\r]^\df}}{\cf \dfp^{\cf } \! \r^{\cf }}{\delta}\leq \cf \dfp^{\cf } \! \r^{\cf }\delta^{-\cf }$, 
        and $\CostLipB{\ReLU}{[-\r,\r]^d }{\g|_{[-\r,\r]^\dg}}{\cg \dgp^{\cg } \! \r^{\cg }}{\delta}\leq \cg \dgp^{\cg } \! \r^{\cg }\delta^{-\cg }$
    \cfload.
    Then it holds for all $\r \in \N$ that 
    \begin{equation}
    \label{Lemma:Comp_theory_1:statement}
        \CostLipA{\ReLU}{[-\r,\r]^d }
        {(\g \circ \f)|_{[-\r,\r]^\df}}
        {\C \dfp^{\C} \! \r^{\C}}{\varepsilon}\leq \C \dfp^\C \! \r^\C\varepsilon^{-\C}.
    \end{equation}
\end{athm}
\begin{aproof}
First, \nobs that 
    \cref{Lemma:Monotonicity_of_Cost} 
ensures that for all 
    $\r \in \N$,
	$\delta\in(0,\eps]$
it holds that
\begin{equation}
\label{Comp_theory_1:a}
\begin{split}
    \CostLipA{\ReLU}{[-\r,\r]^\df}{\f|_{[-\r,\r]^\df}}{\c\dfp^\c \! \r^\c}{\delta} 
	\leq 
	\CostLipA{\ReLU}{[-\r,\r]^\df}{\f|_{[-\r,\r]^\df}}{\cf\dfp^\cf \! \r^\cf}{\delta} 
	\leq 
	\cf \dfp^\cf \! \r^\cf \delta^{-\cf}
	\leq 
	\c \dfp^\c \! \r^\c \delta^{-\c}\ifnocf.
\end{split}
\end{equation}
\cfload[.]%
In addition, \nobs that 
	\cref{Lemma:Monotonicity_of_Cost} 
implies that for all 
    $\r \in \N$,
	$\delta\in(0,\eps]$
it holds that
\begin{equation}
\label{Comp_theory_1:b}
\begin{split}
    \CostLipA{\ReLU}{[-\r,\r]^\dg }{\g|_{[-\r,\r]^\dg}}{\c\dgp^\c \! \r^\c}{\delta} 
	\leq
    \CostLipA{\ReLU}{[-\r,\r]^\dg }{\g|_{[-\r,\r]^\dg}}{\cg\dgp^\cg \! \r^\cg}{\delta} 
	\leq 
	\cg \dgp^\cg \! \r^\cg \delta^{-\cg}
	\leq
	\c \dgp^\c \! \r^\c \delta^{-\c}.
\end{split}
\end{equation}
Next \nobs that 
    the assumption that $\C = 2^{\c+2}\c^{(2\c^2+\c+1)}$ 
	and the fact that $c\geq 2$
show that
\begin{equation}
	\max\{\c^{2\c+2},2\c^2+\c,\c^2+\c\}\leq 4\c^{4\c}\leq\C.
\end{equation}
	This 
	and the fact that $\dg\leq\cf\dfp^\cf\leq\c\dfp^\c$
demonstrate that for all 
	$\r \in \N$ 
it holds that
\begin{equation}
\label{Comp_theory_1:1}
\begin{split}
    \c\dgp^\c(\c\dfp^\c\!\r^\c)^\c
    ( \c\dfp^\c \! \r^\c)
    \leq
    \c(\c\dfp^\c)^\c (\c\dfp^\c\!\r^\c)^\c
    ( \c\dfp^\c \! \r^\c)
    =
	\c^{2\c+2}\dfp^{(2\c^2+\c)}\!\r^{(\c^2+\c)}
	\leq 
    \C \dfp^{\C} \! \r^{\C}
    .
\end{split}
\end{equation}
    \Cref{Lemma:Monotonicity_of_Cost} 
    hence
ensures that for all 
    $\r \in \N$
it holds that
\begin{equation}
\label{Comp_theory_1:3}
\begin{split}
    \CostLipA{\ReLU}{[-\r,\r]^\df }{(\g \circ \f)|_{[-\r,\r]^{\df}}}
    {\C \dfp^{\C} \! \r^{\C}}{\varepsilon}
    \leq
    \CostLipA{\ReLU}{[-\r,\r]^\df }{(\g \circ \f)|_{[-\r,\r]^\df}}
    {\c\dgp^\c (\c\dfp^\c\!\r^\c)^\c  (\c\dfp^\c \! \r^\c)}{\varepsilon}.
\end{split}
\end{equation}
	This,
	the fact that $f([-\r,\r]^{\df})\subseteq[-\cf\dfp^\cf\!\r^\cf,\cf\dfp^\cf\!\r^\cf]^\dg\subseteq[-\c\dfp^\c\!\r^\c,\c\dfp^\c\!\r^\c]^\dg$,
    and \cref{Prop:composition_cost_ReLU} 
        (applied with
            $d_1 \curvearrowleft \df$,
            $d_2 \curvearrowleft \dg$,
            $d_3 \curvearrowleft \ig$,
            $\varepsilon \curvearrowleft \eps$,
            $L_1 \curvearrowleft \c \dfp^\c \! \r^\c$,
            $L_2 \curvearrowleft \c\dgp^\c (\c\dfp^\c\!\r^\c)^\c$,
            $R_1 \curvearrowleft \r$,
            $R_2 \curvearrowleft \c\dfp^\c\!\r^\c$,
            $f_1 \curvearrowleft \f|_{[-\r,\r]^{\df}}$,
            $f_2 \curvearrowleft \g|_{[-\c\dfp^\c\!\r^\c,\c\dfp^\c\!\r^\c]^{\dg}}$
            for $\r \in \N$ 
        in the notation of \cref{Prop:composition_cost_ReLU})
prove that for all 
    $\r \in \N$
it holds that
\begin{equation}
\label{Comp_theory_1:4}
\begin{split}
    &
    \CostLipA{\ReLU}{[-\r,\r]^\df }
    {(\g \circ \f)|_{[-\r,\r]^{\df}}}
    {\C \dfp^{\C} \! \r^{\C}}{\eps}
    \\&\leq
	2\,\CostLipB{\ReLU}{}
	{f|_{[-\r,\r]^\df}}
	{\c\dfp^\c\!\r^\c}
	{\tfrac{\eps}{2\c\dgp^\c (\c\dfp^\c\!\r^\c)^\c}}
	\\&\quad+2\,\CostLipA{\ReLU}{}
	{g|_{[-\c\dfp^\c\!\r^\c ,\c\dfp^\c\!\r^\c ]^{\dg}}}
	{\c\dgp^\c (\c\dfp^\c\!\r^\c)^\c}
	{\tfrac{\eps}{2}}
	+4\dg(\dg+1)
    .
\end{split}
\end{equation}
Combining 
    this
with
    \eqref{Comp_theory_1:a}
    and \eqref{Comp_theory_1:b} 
shows that for all 
    $\r \in \N$ 
it holds that
\begin{equation}
\begin{split}
    &
    \CostLipB{\ReLU}{[-\r,\r]^d }
    {(\g \circ \f)|_{[-\r,\r]^{\df}}}
    {\C \dfp^{\C} \! \r^{\C}}{\eps}
    \\&\leq 
    2\c\dfp^\c\!\r^\c\PR*{\tfrac{2\c\dgp^\c (\c\dfp^\c\!\r^\c)^\c}{\eps}}^\c%
    +2\c\dgp^\c  (\c\dfp^\c\!\r^\c)^\c \PR*{\tfrac{2}{\eps}}^c%
    +4\dg(\dg+1)
    \\&=
    2\c\dfp^\c\!\r^\c2^\c\c^\c\dgp^{(\c^2)}(\c\dfp^\c\!\r^\c)^{(\c^2)}\eps^{-\c}%
    +2\c\dgp^\c  (\c\dfp^\c\!\r^\c)^\c 2^\c\eps^{-\c}%
    +4\dg(\dg+1)
    \\&=
    2^{\c+1}\c^{(\c^2+\c+1)}\dfp^{(\c^3+\c)}\!\r^{(\c^3+\c)}\dgp^{(\c^2)}\eps^{-\c}%
    +2^{\c+1}\c^{\c+1}\dgp^\c\dfp^{(\c^2)}\!\r^{(\c^2)} \eps^{-\c}%
    +4\dg(\dg+1)
    \\&\leq
    2^{\c+1}\c^{(\c^2+\c+1)}\dfp^{(\c^3+\c)}\!\r^{(\c^3+\c)}\dgp^{(\c^2)}\eps^{-\c}%
    +2^{\c+1}\c^{\c+1}\dgp^\c\dfp^{(\c^2)}\!\r^{(\c^2)} \eps^{-\c}%
    +8\dgp^2
    \\&\leq
    2^{\c+2}\c^{(\c^2+\c+1)}\dfp^{(\c^3+\c)}\!\r^{(\c^3+\c)}\dgp^{(\c^2)}\eps^{-\c}
	.
\end{split}
\end{equation}
	The fact that $\dg\leq \cf\dfp^\cf\leq \c\dfp^\c$
	therefore
implies that for all
    $\r\in\N$
it holds that
\begin{equation}
\begin{split}
	\CostLipB{\ReLU}{[-\r,\r]^d }
    {(\g \circ \f)|_{[-\r,\r]^{\df}}}
    {\C \dfp^{\C} \! \r^{\C}}{\eps}
	&\leq
    2^{\c+2}\c^{(\c^2+\c+1)}\dfp^{(\c^3+\c)}\!\r^{(\c^3+\c)}(\c\dfp^\c)^{(\c^2)}\eps^{-\c}
	\\&=
    2^{\c+2}\c^{(2\c^2+\c+1)}\dfp^{(2\c^3+\c)}\!\r^{(\c^3+\c)}\eps^{-\c}
	.
\end{split}
\end{equation}
Combining
    this
with
    the fact that 
    \begin{equation}
        \max\{2^{(\c+2)}\c^{(2\c^2+\c+1)},2\c^3+\c,\c^3+\c,\c\}
		= 
		2^{\c+2}\c^{(2\c^2+\c+1)}
		=
		\C
    \end{equation}
establishes that for all 
    $\r \in \N$
it holds that
\begin{equation}
	\CostLipB{\ReLU}{[-\r,\r]^d }
    {(\g \circ \f)|_{[-\r,\r]^{\df}}}
    {\C \dfp^{\C} \! \r^{\C}}{\eps}
	\leq
    \C\dfp^\C\!\r^\C\eps^{-\C}
	.
\end{equation}
\finishproofthis
\end{aproof}
\endgroup

\begingroup
\newcommand{\const}{c}
\cfclear
\begin{athm}{lemma}{lemma:catching_image}
Let $f,g \in C\prb{\bigcup_{d \in \N}\R^d, \bigcup_{d \in \N}\R^d}$, $\const , \fconstant ,\gconstant \in(0,\infty)$ satisfy for all $d,\radius\in\N$ that 
$f([-R,R]^d)\subseteq [-\fconstant d^{\fconstant }\!\radius^{\fconstant },\fconstant d^{\fconstant }\!\radius^{\fconstant }]^{\imdim{f}{d}}$, 
$g([-R,R]^d)\subseteq [-\gconstant d^{\gconstant }\!\radius^{\gconstant },\gconstant d^{\gconstant }\!\radius^{\gconstant }]^{\imdim{g}{d}}$,
$\imdim{f}{d}\leq \fconstant d^{\fconstant }$,
$\imdim{g}{d}\leq \gconstant d^{\gconstant }$, 
and $\const =\max\{2\gconstant\ceil{\fconstant },\allowbreak \gconstant \ceil{\fconstant }^{2\gconstant } \}$ \cfload.
Then it holds for all $d,\radius\in\N$ that 
$\imdim{g\circ f}{d}\leq \const  d^{\const }$ and
$(g\circ f)([-R,R]^d)\subseteq[-\const d^\const \!\radius^\const ,\const d^\const \!\radius^\const ]^{\imdim{g\circ f}{d}}$.
\end{athm}

\begin{proof}
[Proof of \cref{lemma:catching_image}]
Note that the fact that for all $x\in[0,\infty)$ it holds that $x\leq\ceil{x}$ and the assumption that for all $d,\radius\in\N$ it holds that $f([-R,R]^d)\subseteq [-\fconstant d^{\fconstant }\!\radius^{\fconstant },\fconstant d^{\fconstant }\!\radius^{\fconstant }]^{\imdim{f}{d}}$
show that for all $d,\radius\in\N$ it holds that
\begin{equation}
\begin{split}
\label{Lemma:catching_image:eq2}
(g\circ f)([-R,R]^d)&\subseteq g\prb{\PRb{-\ceil{\fconstant d^{\fconstant }\!\radius^{\fconstant }},\ceil{\fconstant d^{\fconstant }\!\radius^{\fconstant }}}^{\imdim{f}{d}}}
.
\end{split}
\end{equation}
This, the assumption that for all $d,\radius\in\N$ it holds that $g([-R,R]^d)\subseteq [-\gconstant d^{\gconstant }\!\radius^{\gconstant },\gconstant d^{\gconstant }\!\radius^{\gconstant }]^{\imdim{g}{d}}$, and the fact that for all $d\in\N$ it holds that $\imdim{g}{\imdim{f}{d}}=\imdim{g\circ f}{d}$ demonstrate that for all $d,\radius\in\N$ it holds that
\begin{equation}
\begin{split}
\label{Lemma:catching_image:eq3}
(g\circ f)([-R,R]^d)&\subseteq  \PRb{-\gconstant (\imdim{f}{d})^{\gconstant }\ceil{\fconstant d^{\fconstant }\!\radius^{\fconstant }}^{\gconstant },\gconstant (\imdim{f}{d})^{\gconstant }\ceil{\fconstant d^{\fconstant }\!\radius^{\fconstant }}^{\gconstant }}^{\imdim{g\circ f}{d}}
.
\end{split}
\end{equation}
In addition, observe that 
the assumption that for all $d\in\N$ it holds that $\imdim{f}{d}\leq\fconstant d^{\fconstant }$, 
the fact that for all $x,y\in(0,\infty)$ it holds that $\ceil{xy}\leq\ceil{\ceil{x}\ceil{y}}=\ceil{x}\ceil{y}$,
 and the fact that for all $n\in\N$, $x\in(0,\infty)$ it holds that $\ceil{n^x}=\ceil{n^{\ceil{x}}n^{x-\ceil{x}}}\leq \ceil{n^{\ceil{x}}}=n^{\ceil{x}}$ ensure that for all $d,\radius\in\N$ it holds that
\begin{equation}
\begin{split}
\gconstant (\imdim{f}{d})^{\gconstant }\ceil{\fconstant d^{\fconstant }\!\radius^{\fconstant }}^{\gconstant }
\leq &~\gconstant (\fconstant d^{\fconstant })^{\gconstant }\ceil{\fconstant }^{\gconstant }\ceil{d^{\fconstant }}^{\gconstant }\ceil{\radius^{\fconstant }}^{\gconstant }\\
\leq &~\gconstant (\ceil{\fconstant} d^{\fconstant })^{\gconstant }\ceil{\fconstant }^{\gconstant }d^{\ceil{\fconstant }\gconstant }\radius^{\ceil{\fconstant }\gconstant }\\
\leq &~\gconstant \ceil{\fconstant }^{2\gconstant }d^{2\ceil{\fconstant }\gconstant}\radius^{\ceil{\fconstant }\gconstant }\\
\leq &~\const d^{\const }\!\radius^{\const }.
\end{split}
\end{equation}
This and \eqref{Lemma:catching_image:eq3} imply that for all $d,\radius\in\N$ it holds that
\begin{equation}
\label{lemma:catching_image:eq1}
\begin{split}
(g\circ f)([-R,R]^d)&\subseteq  [-\const d^{\const }\!\radius^{\const },\const d^{\const }\!\radius^{\const }]^{\imdim{g\circ f}{d}}.
\end{split}
\end{equation}
Furthermore, note that the assumption that for all $d\in\N$ it holds that $\imdim{f}{d}\leq\fconstant d^{\fconstant }$ and $\imdim{g}{d}\leq\gconstant d^{\gconstant }$ demonstrates that for all $d\in\N$ it holds that
\begin{equation}
\begin{split}
\imdim{g\circ f}{d}=\imdim{g}{\imdim{f}{d}}
\leq \gconstant (\imdim{f}{d})^{\gconstant } 
\leq \gconstant (\fconstant d^{\fconstant })^{\gconstant }
= \gconstant \fconstant^{\gconstant }d^{\fconstant\gconstant }
\leq\const  d^{\const }.
\end{split}
\end{equation}
The proof of \cref{lemma:catching_image} is thus complete.
\end{proof}
\endgroup

\begingroup
\newcommand{\const}{C}
\newcommand{\cONST}{c}
\newcommand{\constf}{\mathfrak f}
\newcommand{\constg}{\mathfrak g}
\cfclear
\begin{athm}{lemma}{lemma:D_closed}
Let $f,g\in\Dapprox$ \cfload. Then $(g\circ f)\in\Dapprox$.
\end{athm}
\begin{aproof}
    \Nobs that 
        the assumption that $f,g\in\Dapprox$
    shows that there exist
        $\constf ,\constg \in[0,\infty)$
    which satisfy that
    \begin{enumerate}[(I)]
        \item\llabel{it1}
        for all 
            $d,R\in\N$
        it holds that
            $f([-R,R]^d)\subseteq[-\constf d^{\constf }\!R^{\constf },\constf d^{\constf }\!R^{\constf }]^{\imdim fd}$
            and $g([-R,R]^d)\subseteq[-\constg d^{\constg }\!R^{\constg },\constg d^{\constg }\!R^{\constg }]^{\imdim gd}$,
        \item\llabel{it2}
        for all
            $d\in\N$
        it holds that
            $\imdim fd\leq \constf d^{\constf }$
            and $\imdim gd\leq \constg d^{\constg }$,
        and
        \item\llabel{it3}
        for all
            $d,R\in\N$,
            $\eps\in(0,1]$
        it holds that
        \begin{equation}
        \begin{split}
            \CostLipA{\ReLU}{}{\pr{[-R,R]^d\ni x\mapsto f(x)\in\R^{\imdim fd}}}{\constf d^{\constf }\!R^{\constf }}{\eps}
            &\leq 
            \constf d^{\constf }\!R^{\constf }\eps^{-\constf }
            \\\text{and}\qquad 
            \CostLipA{\ReLU}{}{\pr{[-R,R]^d\ni x\mapsto g(x)\in\R^{\imdim gd}}}{\constg d^{\constg }\!R^{\constg }}{\eps}
            &\leq 
            \constg d^{\constg }\!R^{\constg }\eps^{-\constg }
        \end{split}
        \end{equation}
    \end{enumerate}
    \cfload.
    \Nobs that
        \lref{it1},
        \lref{it2},
        \lref{it3}, and 
        \cref{Lemma:Comp_theory_1}
        (applied with
            $m\is d$,
            $n\is \imdim fd$,
            $o\is \imdim{g\circ f}d$,
            $\mathfrak f\is \constf $,
            $\mathfrak g\is \constg $,
            $\eps\is\eps$,
            $f\is(\R^d\ni x\mapsto f(x)\in\R^{\imdim fd})$,
            $g\is(\R^{\imdim fd}\ni x\mapsto g(x)\in\R^{\imdim {g\circ f}d})$
            for $d\in\N$, $\eps\in(0,1]$
        in the notation of \cref{Lemma:Comp_theory_1})
    prove that there exists
        $\cONST\in(0,\infty)$
    such that for all
        $d,R\in\N$,
        $\eps\in(0,1]$
    it holds that
    \begin{equation}
        \eqlabel{1}
        \CostLipA{\ReLU}{}
        {\pr{[-R,R]^d\ni x\mapsto (g\circ f)(x)\in \R^{\imdim{g\circ f}d}}}
        {\cONST d^\cONST\!R^\cONST}
        \eps
        \leq
        \cONST d^\cONST\!R^\cONST\eps^{-\cONST}
        .
    \end{equation}
    Moreover, \nobs that 
        \cref{lemma:catching_image}
    ensures that there exists 
        $\const \in(0,\infty)$ 
    such that for all 
        $d,\radius\in\N$ 
    it holds that 
        $\imdim{g\circ f}{d}\leq\const  d^{\const }$ 
        and $(g\circ f)([-R,R]^d)\subseteq[-\const d^\const \!\radius^\const ,\const d^\const \!\radius^\const ]^{\imdim{g\circ f}{d}}$.
    Combining 
        this, 
        \eqqref{1}, and 
        \cref{Lemma:Monotonicity_of_Cost} 
    with 
        \eqref{polyD_mult:1} 
    establishes that 
        $(g\circ f)\in\Dapprox$. 
    This completes the proof of \cref{lemma:D_closed}.
\end{aproof}
\endgroup

\subsection{Parallelizations of ANN approximable functions}

\label{Subsection:3.6}

The main result of this section, \cref{Lemma:Comp_theory_3} establishes, roughly speaking, that certain sequences of multi-dimensional functions which arise as parallelizations of univariate functions from the ANN approximation spaces introduced in \cref{def:polyC} are themselves in the ANN approximation space for multi-dimensional functions introduced in \cref{def:polyD_mult}
and can in that sense be approximated by ANNs without the curse of dimensionality.

\cfclear
\begin{athm}{lemma}{lemma:Lipschitz_error}
Let $L \in \R$, $d \in \N$, $m_1,m_2,\ldots,m_d\in\N$, let $g_i \in C(\R^{m_i},\R)$, $i \in \{1, 2, \ldots, d\}$, satisfy for all $i \in \{1, 2, \ldots, d\}$, $x,y\in\R^{m_i}$ that $\vass{g_i(x) - g_i(y)} \leq L\norm{x - y}$, and let $f \in C\prb{\R^{\PR{\sum_{i=1}^dm_i}},\R^d}$ satisfy for all $x=(x_1,x_2,\ldots,x_d)\in\prb{\bigtimes_{i=1}^d\R^{m_i}}$ that $f(x_1,x_2,\ldots,x_d)=(g_1(x_1), g_2(x_2), \ldots, g_d(x_d))$.
Then it holds for all $x,y\in\R^{\PR{\sum_{i=1}^dm_i}}$ that
\begin{equation}
    \mednorm{f(x) - f(y)} \leq L \mednorm{x-y} 
\end{equation}
\cfload.
\end{athm}

\begin{aproof}
Observe that for all $x=(x_1,x_2,\ldots,x_d)$, $y=(y_1,y_2,\ldots,y_d)\in\prb{\bigtimes_{i=1}^d\R^{m_i}}$ it holds that
\begin{equation}
\begin{split}
    \mednorm{f(x) - f(y)}  
    = \PR*{\ssum_{i=1}^d \vass*{g_i(x_i) - g_i(y_i)}^2}^{\frac{1}{2}}
    \leq 
    L \PR*{\ssum_{i=1}^d \norm{x_i-y_i}^2}^{\frac{1}{2}}
    \leq L \mednorm{x - y}.
\end{split}
\end{equation}
This completes the proof of \cref{lemma:Lipschitz_error}.
\end{aproof}

\begingroup
\newcommand{\const}{\kappa}
\newcommand{\cONST}{C}
\cfclear
\begin{athm}{lemma}{Lemma:Comp_theory_2.1}
    Let 
        $d,\radius\in\N$, 
        $c,\const ,\cONST ,r \in [0, \infty)$, 
        $\varepsilon \in (0,1]$, 
        $f_1,f_2,\ldots,f_d \in \Capprox{c}{r}$, 
        $F \in C(\R^d, \R^d)$ 
    satisfy for all 
        $x = (x_1, x_2, \ldots, x_d) \in \R^d$ 
        that $F(x) = (f_1(x_1), f_2(x_2), \ldots, f_d(x_d))$, 
        $\const =\max\{c,cr^c\}$, and 
        $\cONST  = \max\{61\const ^2, 2\const +2\}$ 
\cfload.
Then
\begin{equation}
\begin{split}
    \CostLipA{\ReLU}{\indicator{[-\radius,\radius]^d}}{F|_{[-\radius,\radius]^d}}{\cONST  d^{\cONST } \! \radius^{\cONST }}{\varepsilon} 
    \leq
    \cONST  d^{\cONST } \! \radius^{\cONST } \varepsilon^{-\cONST }
\end{split}
\end{equation}
\cfout.
\end{athm}

\begin{aproof}
\Nobs that \cref{cor:polyC} ensures that 
$\Capprox{c}{r}\subseteq \Capprox{\const }{1}$.
    The assumption that $f_1,f_2,\ldots,f_d \in \Capprox{c}{r}$ 
    therefore 
implies that for all 
    $m\in\{1,2,\ldots,d\}$ 
it holds that
\begin{equation}
    \CostLipB{\ReLU}{\indicator{[-\radius,\radius]^d}}{(f_m)|_{[-\radius,\radius]}}{\const   \radius^{\const }}{d^{-\frac{1}{2}}\varepsilon}
    \leq 
    \const   \radius^{\const }\prb{d^{-\frac{1}{2}}\varepsilon}^{-\const }
    =
    \const   \radius^{\const }d^{\frac{\const }{2}}\varepsilon^{-\const }
    \ifnocf.
\end{equation}
\cfload[.]%
    This 
    and \cref{cor_cost_of_Lip_approx_set_equivalence} 
assure that there exist 
    $\mathscr{f}_1,\mathscr{f}_2,\ldots,\mathscr{f}_d\in\ANNs$
which satisfy that
\begin{enumerate}[(I)] 
\item \label{Comp_theory_3:item_1} 
    it holds for all 
        $m\in\{1,2,\ldots,d\}$ 
    that 
        $\realisation_{\ReLU}( \mathscr{f}_m) \in C(\R,\R)$,
\item \label{Comp_theory_3:item_2} 
    it holds for all 
        $m\in\{1,2,\ldots,d\}$, 
        $x \in [-\radius,\radius]$ 
    that 
        $\vass{f_m (x) - (\realisation_{\ReLU}(\mathscr{f}_m))(x)}
        \leq d^{-\frac{1}{2}}\varepsilon$, 
\item \label{Comp_theory_3:item_3} 
    it holds for all 
        $m\in\{1,2,\ldots,d\}$, 
        $x,y \in [-\radius,\radius]$ 
    that
    \begin{equation}
        \vass{(\realisation_{\ReLU}( \mathscr{f}_m ))(x) - (\realisation_{\ReLU}( \mathscr{f}_m ))(y)} \leq \const \radius^\const \vass{x-y},
    \end{equation}
    and
\item \label{Comp_theory_3:item_4} 
    it holds for all 
        $m\in\{1,2,\ldots,d\}$ 
    that
        $\paramANN(\mathscr{f}_m) \leq   \const \radius^\const  d^{\frac{\const }{2}}\varepsilon^{-\const }$
\end{enumerate}
\cfload. 
Next let $\mathscr{g} \in \ANNs$ satisfy
\begin{equation}
\label{Comp_theory_3:1}
    \mathscr{g} 
    = 
    \parallelization_{d, (\ReLUidANN{1},\ReLUidANN{1},\ldots,\ReLUidANN{1})}\prb{\mathscr{f}_{1}, \mathscr{f}_{2}, \ldots, \mathscr{f}_{d} }
\end{equation}
\cfload. 
\Nobs that
    \eqref{Comp_theory_3:1},
    \cref{Comp_theory_3:item_1},
    \cref{Comp_theory_3:item_2},
    \cref{Prop:identity_representation}, and 
    \cref{Lemma:PropertiesOfParallelizationRealization}
ensure that for all 
    $x=(x_1,x_2,\ldots,x_d)\in[-\radius,\radius]^d$ 
it holds that 
    $\realisation_{\ReLU}(\mathscr{g}) \in C(\R^d,\R^d)$ and 
\begin{equation}
\label{Lemma:Comp_theory_2.1:prop:approx}
\begin{split}
    \norm{F(x) - \pr{\realisation_{\ReLU}(\mathscr{g})}(x)}
    &=
    \PR*{\ssum_{m=1}^{d}\vass*{f_m(x_m) - \pr{\realisation_{\ReLU}(\mathscr{f}_m)}(x_m)}^2}^\frac{1}{2}
    \leq
    \PRb{d\prb{d^{-1}\varepsilon^2}}^\frac{1}{2}
    =
    \varepsilon
\end{split}
\end{equation}
\cfload. 
In addition, \nobs that 
    \eqref{Comp_theory_3:1}, 
    \cref{Comp_theory_3:item_3}, 
    \cref{Lemma:PropertiesOfParallelizationRealization}, and 
    \cref{lemma:Lipschitz_error} 
demonstrate that for all 
    $x,y\in[-\radius,\radius]^d$ 
it holds that
\begin{equation}
\label{Lemma:Comp_theory_2.1:prop:lip}
\norm{\pr{\realisation_{\ReLU}(\mathscr{g})}(x) - \pr{\realisation_{\ReLU}(\mathscr{g})}(y)}
\leq \const  \radius^\const \norm{x-y}.
\end{equation}
In the next step, \nobs that
    \cref{Prop:identity_representation}
implies that
    $\dims(\ReLUidANN{1})=(1,2,1)$.
Combining
    this,
    \eqref{Comp_theory_3:1}, and
    \cref{Comp_theory_3:item_4}
with
    \cref{Lemma:PropertiesOfParallelization} 
        (applied with
        $n \curvearrowleft d$,
		$(\mathscr{g}_1,\mathscr{g}_2,\allowbreak\dots, \mathscr{g}_n) \curvearrowleft (\ReLUidANN{1},\ReLUidANN{1},\ldots,\ReLUidANN{1})$,
		$\mathscr{f} \curvearrowleft (\mathscr{f}_1,\mathscr{f}_2,\allowbreak\dots, \mathscr{f}_d)$ 
        in the notation of \cref{Lemma:PropertiesOfParallelization}), and 
    the fact that $\max_{m\in\{1,2,\dots,d\}}\lengthANN(\mathscr{f}_m)\leq \max_{m\in\{1,2,\dots,d\}}\param(\mathscr{f}_m)$ 
shows that
\begin{equation}
\begin{split}
\label{Comp_theory_3:2}
    \paramANN\pr{\mathscr{g}}
 	&\leq 
    \frac{1}{2} \PRbbb{\ssum_{j=1}^d
	\prb{2\,\paramANN(\mathscr{f}_j)
    +
	\prb{2\PRb{\max\nolimits_{m\in\{1,2,\dots,d\}} \lengthANN(\mathscr f_m)}+1}(2+1)}}^2
	\\&\leq 
    \tfrac{1}{2}
	\prb{2d \const    \radius^\const d^{\frac{\const }{2}} \varepsilon^{-\const }  
	+9d  \const    \radius^\const d^{\frac{\const }{2}} \varepsilon^{-\const }
	}^{2}
	\\&\leq
	61\const ^2 d^{2+\const } \! \radius^{2\const } \varepsilon^{-2\const }\leq
            \cONST  d^{\cONST } \! \radius^{\cONST } \varepsilon^{-\cONST }.
\end{split}
\end{equation}
This, \eqref{Lemma:Comp_theory_2.1:prop:approx}, \eqref{Lemma:Comp_theory_2.1:prop:lip}, \cref{Lemma:Monotonicity_of_Cost}, and the fact that $\cONST \geq\const  $ establish that
\begin{equation}
    \begin{split}
        &
            \CostLipA{\ReLU}{\indicator{[-\radius,\radius]^d}}{F|_{[-\radius,\radius]^d}}{\cONST  d^{\cONST } \! \radius^{\cONST }}{\varepsilon} 
        \leq
            \CostLipB{\ReLU}{\indicator{[-\radius,\radius]^d}}{F|_{[-\radius,\radius]^d}}{\const   \radius^\const }{\varepsilon} 
        \leq
            \cONST  d^{\cONST } \! \radius^{\cONST } \varepsilon^{-\cONST }.
    \end{split}
\end{equation}
This completes the proof of \cref{Lemma:Comp_theory_2.1}\cfload.
\end{aproof}
\endgroup

\begingroup
\newcommand{\const}{C}
\cfclear
\begin{athm}{lemma}{Lemma:Comp_theory_3}
    Let 
        $c,r \in [0, \infty)$, 
        $(f_n)_{n\in \N} \subseteq \Capprox{c}{r}$, 
        $F \in C\prb{\bigcup_{d \in \N}\R^d, \bigcup_{d \in \N}\R^d}$ 
    satisfy for all 
        $d,\radius \in \N$, 
        $x = (x_1, x_2, \ldots, x_d) \in \R^d$, 
        $y\in[-\radius,\radius]$ 
    that 
        $F(x) = (f_1(x_1), f_2(x_2), \ldots, f_d(x_d))$ 
        and $\vass{f_{d}(y)} \leq cd^c\!R^c$
    \cfload.
    Then $F \in \Dapprox$
    \cfout.
\end{athm}
\begin{proof}[Proof of \cref{Lemma:Comp_theory_3}]
    \setnote
    \Nobs that 
        \cref{Lemma:Comp_theory_2.1} 
    ensures that there exists 
        $\const \in[\max\{1,c\},\infty)$ 
    which satisfies that for all 
        $d, \radius \in \N$, 
        $\varepsilon \in (0,1]$ 
    it holds that
    \begin{equation}
    \label{Lemma:Comp_theory_3:eq1}
    \begin{split}
        \CostLipA{\ReLU}{\indicator{[-\radius,\radius]^d}}{\pr{[-\radius,\radius]^d \ni x \mapsto F(x) \in \R^{d}}}{\const  d^{\const } \! \radius^{\const }}{\varepsilon} 
        \leq
        \const  d^{\const } \! \radius^{\const } \varepsilon^{-\const }
    \end{split}
    \end{equation}
    \cfload. 
    Moreover, \nobs that 
        the assumption that 
            for all 
                $d,\radius \in \N$, 
                $y\in[-R,R]$ 
            it holds that 
                $\vass{f_{d}(y)} \leq cd^c\!R^c$ 
        and the fact that 
            $c\leq \const $ 
    imply that for all 
        $d,\radius \in \N$ 
    it holds that 
    \begin{equation}
    \label{Lemma:Comp_theory_3:eq2}
        F([-\radius,\radius]^d)
        \subseteq
        \prb{\sssbigtimes_{k=1}^{d}[-ck^c\! R^c,ck^c\! R^c]}
        \subseteq
        [-cd^c\! R^c,cd^c\! R^c]^d
        \subseteq
        [-\const  d^{\const } \! \radius^{\const },\const  d^{\const } \! \radius^{\const }]^d
        . 
    \end{equation}
    Furthermore, \nobs that 
        the fact that 
            $\const \geq 1$ 
    implies that for all 
        $d\in\N$ 
    it holds that 
        $\imdim{F}{d}=d\leq \const d^\const $ 
    \cfload. 
        This, 
        \eqref{Lemma:Comp_theory_3:eq1}, and 
        \eqref{Lemma:Comp_theory_3:eq2} 
    ensure that 
        $F \in \Dapprox$ \cfload.
    This completes the proof of \cref{Lemma:Comp_theory_3}.
\end{proof}
\endgroup

\section{ANN approximations for Lipschitz continuous functions}
\label{Section:4}

In this section we present results on the approximation of
locally Lipschitz continuous functions from $\R$ to $\R$
by ANNs with the ReLU activation function.
The results in this section are essentially well-known
and only for completeness we also include detailed proofs here.

\subsection{Linear interpolations with ANNs}
\label{Subsection:4.1}

In this section, we introduce ANNs for piecewise linear interpolation, which will be used in the following section to obtain results on the approximability of locally Lipschitz continuous functions from $\R$ to $\R$.

\cfclear
\begin{definition}[Interpolating ANNs]
\label{def:interpolatingDNN}
\cfconsiderloaded{def:interpolatingANN}
Let 
	$N \in \N$, 
	$h_0,h_1,\ldots, h_{N},\gp_0,\gp_1,\ldots,\gp_N \in \R$ 
satisfy $\gp_0<\gp_1<\ldots<\gp_N$.
Then we denote by 
$\interpolatingNN{h_0,h_1,\ldots, h_{N}}{\gp_0,\gp_1,\ldots,\gp_N } \in \ANNs$ the \cfadd{def:neuralnetwork}ANN given by
\begin{equation}
\label{def:interpolatingDNN:eq1}
\begin{split}
	&\interpolatingNN{h_0,h_1,\ldots, h_{N}}{\gp_0,\gp_1,\ldots,\gp_N } = 
	\mleft(\!
		\pr*{\!
			\begin{pmatrix} 1\\1\\\vdots\\1 \end{pmatrix},\!
			\begin{pmatrix} -\gp_0\\-\gp_1\\\vdots\\-\gp_N \end{pmatrix} 
		\!}\mright.
		,\\&\qquad\mleft.\vphantom{\begin{pmatrix} 1\\1\\\vdots\\1 \end{pmatrix}}
		\pr*{
			 \begin{pmatrix} 
			 	\tfrac{h_1-h_0}{\gp_1 - \gp_0} &
				\tfrac{h_2-h_1}{\gp_2-\gp_1} - \tfrac{h_1-h_0}{\gp_1-\gp_0} &
				\ldots &
				\tfrac{h_N-h_{N-1}}{\gp_N-\gp_{N-1}} - \tfrac{h_{N-1}-h_{N-2}}{\gp_{N-1}-\gp_{N-2}} &
				-\tfrac{h_N-h_{N-1}}{\gp_N-\gp_{N-1}}
			 \end{pmatrix}
			 ,
			  h_0 
		}
	\!\mright) \\
& \in (\R^{(N+1) \times 1} \times \R^{N+1}) \times (\R^{1 \times (N+1)} \times \R)
\end{split}
\end{equation}
\cfload.
\end{definition}

\cfclear
\begin{athm}{lemma}{DNN_interpolation}
Let 
	$N \in \N$, 
	$h_0,h_1,\ldots, h_{N},\gp_0,\gp_1,\ldots,\gp_N \in \R$ 
satisfy $\gp_0<\gp_1<\ldots\allowbreak<\gp_N$.
Then 
\begin{enumerate}[(i)]
\item \label{DNN_interpolation:item1}
it holds that
	$\dims\prb{\interpolatingNN{h_0,h_1,\ldots, h_{N}}{\gp_0,\gp_1,\ldots,\gp_N } } 
= 
	(1, N+1, 1)$,
\item \label{DNN_interpolation:realization}
it holds that
	$\functionReLUANN\prb{\interpolatingNN{h_0,h_1,\ldots, h_{N}}{\gp_0,\gp_1,\ldots,\gp_N } }
    \in C(\R,\R)$,
\item \label{DNN_interpolation:item2}
it holds for all 
	$n\in \{0, 1, \ldots, N\}$ 
that $\prb{ \functionReLUANN \prb{ \interpolatingNN{h_0,h_1,\ldots, h_{N}}{\gp_0,\gp_1,\ldots,\gp_N }}} (\gp_n)=h_n$,

\item \label{DNN_interpolation:item3}
it holds for all 
$n\in \{0, 1, \ldots, N-1 \}$, $x \in (-\infty,\xi_0]$, $y\in\PR{\xi_n,\xi_{n+1}}$, $z\in[\xi_N,\infty)$
that $\prb{ \functionReLUANN \prb{\interpolatingNN{h_0,h_1,\ldots, h_{N}}{\gp_0,\gp_1,\ldots,\gp_N }}}(x)=h_0$, $\prb{ \functionReLUANN \prb{\interpolatingNN{h_0,h_1,\ldots, h_{N}}{\gp_0,\gp_1,\ldots,\gp_N }}}(z)=h_N$, and
\begin{equation}
\prb{ \functionReLUANN \prb{\interpolatingNN{h_0,h_1,\ldots, h_{N}}{\gp_0,\gp_1,\ldots,\gp_N }}}(y)=h_n + \prb{\tfrac{h_{n+1} - h_n}{\gp_{n+1} - \gp_n} }(y - \gp_n),\quad
\end{equation}
\item \label{DNN_interpolation:item4}
it holds for all 
	$n\in \{ 1,2, \ldots, N \}$, $x \in [\gp_{n-1}, \gp_n]$
that 
\begin{equation}
    \prb{  \functionReLUANN \prb{ \interpolatingNN{h_0,h_1,\ldots, h_{N}}{\gp_0,\gp_1,\ldots,\gp_N }}}(x) \in \PRb{ \min\{h_{n-1}, h_n\}, \max\{h_{n-1}, h_n\} },
\end{equation}
and
\item \label{DNN_interpolation:item5}
it holds for all 
	$x \in \R$
that
\begin{equation}
\prb{ \functionReLUANN \prb{ \interpolatingNN{h_0,h_1,\ldots, h_{N}}{\gp_0,\gp_1,\ldots,\gp_N } }}(x)
\in \PRb{ \min\nolimits_{n \in \{0, 1, \ldots, N \}} h_n, \max\nolimits_{n \in \{0, 1, \ldots, N \}} h_n }
\end{equation}
\end{enumerate}
\cfout.
\end{athm}

\begin{proof}[Proof of \cref{DNN_interpolation}]
Throughout this proof
let $c_n \in \R$, $n\in \{0, 1, \ldots, N\}$, satisfy 
for all
$
 n\in \N\cap(0,N)$ that $c_0=\tfrac{h_1-h_0}{\gp_1 - \gp_0}$, $c_N=-\tfrac{h_N-h_{N-1}}{\gp_N-\gp_{N-1}}$, and
\begin{equation}
	c_n=		\tfrac{h_{n+1}-h_n}{\gp_{n+1}-\gp_n} - \tfrac{h_n-h_{n-1}}{\gp_n-\gp_{n-1}}
\ifnocf.
\end{equation}
\cfload[.]%
Observe that 
    \cref{def:interpolatingDNN:eq1} 
implies
    \cref{DNN_interpolation:item1}.
Note that
    \cref{ANNrealization:ass2} and
    \eqref{def:interpolatingDNN:eq1} 
prove that for all
	$x \in \R$ 
it holds that
    $\functionReLUANN\prb{\interpolatingNN{h_0,h_1,\ldots, h_{N}}{\gp_0,\gp_1,\ldots,\gp_N } }
    \in C(\R,\R)$
and
\begin{equation}
\begin{split}
\label{DNN_interpolation_explicit:eq}
	\prb{\functionReLUANN \prb{\interpolatingNN{h_0,h_1,\ldots, h_{N}}{\gp_0,\gp_1,\ldots,\gp_N }}}(x)
&= 
	\begin{pmatrix}
		c_0&c_1&\ldots&c_N
	\end{pmatrix}
	\begin{pmatrix}
		\max\{x-\gp_0,0\}\\
		\max\{x-\gp_1,0\}\\
		\vdots\\
		\max\{x-\gp_N,0\}
	\end{pmatrix}
	+
	h_0
\\&=
	h_0+\sum_{k=0}^{N}c_k\max\{x-\gp_k,0\}\ifnocf.
\end{split}
\end{equation}
\cfload[.]%
    Hence, 
we obtain 
    \cref{DNN_interpolation:realization}.
    This and 
    the assumption that 
        $\xi_0<\xi_1<\ldots<\xi_N$ 
ensure that for all 
	$x \in (-\infty,\gp_0]$
it holds that
\begin{equation}
\label{DNN_interpolation:eq1}
\begin{split}
	\prb{\functionReLUANN \prb{\interpolatingNN{h_0,h_1,\ldots, h_{N}}{\gp_0,\gp_1,\ldots,\gp_N }}}(x)
=
	h_0+0
=
	h_0.
\end{split}
\end{equation}
In addition, observe that the assumption that $\xi_0<\xi_1<\ldots<\xi_N$ and the fact that
for all
 $n\in \{0, 1, \ldots, N-1\}$ it holds that $\sum_{k=0}^{n}c_k = \tfrac{h_{n+1}-h_n}{\gp_{n+1}-\gp_n}$
show that for all
	$n\in \{0, 1, \ldots, N-1 \}$, 
	$x\in [\gp_n,\gp_{n+1}]$ 
it holds that
\begin{equation}
\label{DNN_interpolation:eq2}
\begin{split}
    &
    \prb{\functionReLUANN \prb{\interpolatingNN{h_0,h_1,\ldots,h_{N}}{\gp_0,\gp_1,\ldots,\gp_N }}}(x)
	-
	\prb{\functionReLUANN \prb{\interpolatingNN{h_0,h_1,\ldots, h_{N}}{\gp_0,\gp_1,\ldots,\gp_N }}}(\gp_n)
    \\&= 
	\sum_{k=0}^{N}c_k\PRb{\max\{x-\gp_k,0\}-\max\{\gp_{n}-\gp_k,0\}}
    \\&=
    \sum_{k=0}^{n}c_k \PR{(x-\gp_k)-(\gp_{n}-\gp_k)} 
    \\&= 
    \sum_{k=0}^{n}c_k (x-\gp_{n})
    =
	\prb{\tfrac{h_{n+1}-h_n}{\gp_{n+1}-\gp_n}}(x-\gp_n)
    .
\end{split}
\end{equation}
This and \eqref{DNN_interpolation:eq1} demonstrate that for all
	$x\in [\gp_0,\gp_{1}]$
it holds that
\begin{equation}
\label{DNN_interpolation:eqn0}
\begin{split}
&	\prb{\functionReLUANN \prb{\interpolatingNN{h_0,h_1,\ldots, h_{N}}{\gp_0,\gp_1,\ldots,\gp_N }}}(x)\\
&=\prb{\functionReLUANN \prb{\interpolatingNN{h_0,h_1,\ldots, h_{N}}{\gp_0,\gp_1,\ldots,\gp_N }}}(\gp_0)
	+
	\prb{\functionReLUANN \prb{\interpolatingNN{h_0,h_1,\ldots, h_{N}}{\gp_0,\gp_1,\ldots,\gp_N }}}(x) - \prb{\functionReLUANN \prb{\interpolatingNN{h_0,h_1,\ldots, h_{N}}{\gp_0,\gp_1,\ldots,\gp_N }}}(\gp_0) \\
&=h_0
	+
	\prb{\tfrac{h_{1}-h_0}{\gp_{1}-\gp_0}}(x-\gp_0)\ifnocf.
\end{split}
\end{equation}
Moreover, note that \eqref{DNN_interpolation:eq2} implies that
for all 
	$n \in \N\cap(0,N)$, 
	$x\in [\gp_n,\gp_{n+1}]$ 
with
$\forall \, y\in [\gp_{n-1},\gp_{n}] \colon\allowbreak
	\prb{ \functionReLUANN \allowbreak \prb{\interpolatingNN{h_0,h_1,\ldots, h_{N}}{\gp_0,\gp_1,\ldots,\gp_N }}}(y)
\allowbreak =
	h_{n-1} + \prb{\tfrac{h_{n} - h_{n-1}}{\gp_{n} - \gp_{n-1}} }(y- \gp_{n-1})
$
it holds that
\begin{equation}
\begin{split}
\label{DNN_interpolation:eqinductionstep}
	&\prb{\functionReLUANN \prb{\interpolatingNN{h_0,h_1,\ldots, h_{N}}{\gp_0,\gp_1,\ldots,\gp_N }}}(x)
\\&=
	\prb{\functionReLUANN \prb{\interpolatingNN{h_0,h_1,\ldots, h_{N}}{\gp_0,\gp_1,\ldots,\gp_N }}}(\gp_n)
	+
	\prb{\functionReLUANN \prb{\interpolatingNN{h_0,h_1,\ldots, h_{N}}{\gp_0,\gp_1,\ldots,\gp_N }}}(x) - \prb{\functionReLUANN \prb{\interpolatingNN{h_0,h_1,\ldots, h_{N}}{\gp_0,\gp_1,\ldots,\gp_N }}}(\gp_n) \\
&=
	h_{n-1} + \prb{\tfrac{h_{n} - h_{n-1}}{\gp_{n} - \gp_{n-1}} }(\gp_{n}- \gp_{n-1})
	+
	\prb{\tfrac{h_{n+1}-h_n}{\gp_{n+1}-\gp_n}}(x-\gp_n)\\
&=
	h_{n} 
	+
	\prb{\tfrac{h_{n+1}-h_n}{\gp_{n+1}-\gp_n}}(x-\gp_n)\ifnocf.
\end{split}
\end{equation}
\cfload[.]%
Combining this and \eqref{DNN_interpolation:eqn0} with induction proves that for all
	$n\in \{0, 1, \ldots, N-1 \}$,
	$x\in [\gp_n,\gp_{n+1}]$
it holds that
\begin{equation}
\label{DNN_interpolation:eq3}
	\prb{\functionReLUANN \prb{\interpolatingNN{h_0,h_1,\ldots, h_{N}}{\gp_0,\gp_1,\ldots,\gp_N }}}(x)
=
	h_n + \prb{\tfrac{h_{n+1} - h_n}{\gp_{n+1} - \gp_n} }(x - \gp_n).
\end{equation}
The fact that for all
$n\in \{0, 1, \ldots, N\}$ it holds that $\gp_n\leq \gp_{N}$, the fact that
$\sum_{k=0}^{N}c_k=0$, and \eqref{DNN_interpolation_explicit:eq} therefore
imply that for all 
	$x\in [\gp_N,\infty)$ 
it holds that
\begin{align}\begin{split}
&	\prb{\functionReLUANN \prb{\interpolatingNN{h_0,h_1,\ldots, h_{N}}{\gp_0,\gp_1,\ldots,\gp_N }}}(x)
	-
	\prb{\functionReLUANN \prb{\interpolatingNN{h_0,h_1,\ldots, h_{N}}{\gp_0,\gp_1,\ldots,\gp_N }}}(\gp_N)\\
&= \sum_{k=0}^{N}c_k\PR{\max\{x-\gp_k,0\}-\max\{\gp_{N}-\gp_k,0\}}\\
&=\sum_{k=0}^{N}c_k [(x-\gp_k)-(\gp_{N}-\gp_k)]= \sum_{k=0}^{N}c_k (x-\gp_{N})=0.\end{split}
\end{align}
This and \eqref{DNN_interpolation:eq3} show that for all 
	$x\in [\gp_N,\infty)$
it holds that 
\begin{equation}
\begin{split}
\label{DNN_interpolation:eqlast}
	\prb{\functionReLUANN \prb{\interpolatingNN{h_0,h_1,\ldots, h_{N}}{\gp_0,\gp_1,\ldots,\gp_N }}}(x)
&=
	\prb{\functionReLUANN \prb{\interpolatingNN{h_0,h_1,\ldots, h_{N}}{\gp_0,\gp_1,\ldots,\gp_N }}}(\gp_N)\\
&=
	h_{N-1} + \prb{\tfrac{h_{N} - h_{N-1}}{\gp_{N} - \gp_{N-1}} }(\gp_N - \gp_{N-1})
=
	h_N.    
\end{split}
\end{equation}
Combining this, \eqref{DNN_interpolation:eq1}, and \eqref{DNN_interpolation:eq3} establishes \cref{DNN_interpolation:item3}.
Note that \cref{DNN_interpolation:item3} implies \cref{DNN_interpolation:item2,DNN_interpolation:item4}. Observe that \eqref{DNN_interpolation:eq1}, \eqref{DNN_interpolation:eqlast}, and \cref{DNN_interpolation:item4} ensure \cref{DNN_interpolation:item5}.
The proof of \cref{DNN_interpolation} is thus complete.
\end{proof}

\subsection{ANN approximations for locally Lipschitz continuous functions}
\label{Subsection:4.3}

The main result of this section, \cref{Lemma:loc_Comp_theory_15_sum}, shows, roughly speaking, that locally Lipschitz continuous functions from $\R$ to $\R$ belong to the ANN approximation spaces introduced in \cref{def:polyC} above and in that sense can be approximated by ANNs with the ReLU activation function without the curse of dimensionality.
This is achieved by means of linear interpolations and the corresponding ANNs introduced in the preceding subsection.

\begingroup
\newcommand{\newconstant}{a}
\cfclear
\begin{athm}{lemma}{Prop:interpolation_3}
    Let 
        $\radius, c,\varepsilon \in (0, \infty)$, 
        $\newconstant \in [0, \infty)$, 
        $N \in \N$, 
        $(\gp_{n})_{n \in \Z} \subseteq \R$ 
    satisfy for all 
        $n \in \Z$ 
    that 
        $\gp_{n} = \clip{-R}{R}\pr{n(2c)^{-1}(1+2\radius)^{-\newconstant} \varepsilon}$ and 
        $N = \min\{k \in \N \colon \gp_{k} = \radius\}$,
    let
        $f \in C([-\radius,\radius],\R)$ 
    satisfy for all 
        $x,y\in[-\radius,\radius]$ 
    that
    \begin{align}\label{interpolation_302}
        \abs{f(x)-f(y)}\leq c(1+\abs{x}+\abs{y})^{\newconstant}\vass{x-y},
    \end{align}
    and let 
        $g \colon \R \to \R$ 
    satisfy for all 
        $n\in\Z\cap[-N,N)$, 
        $x\in (\gp_{n},\gp_{n+1}]$, 
        $y\in[R,\infty)$ 
    that 
        $g(y)=f(R)$, 
        $g(-y)=f(-R)$, 
        and
    \begin{equation}\label{interpolation_303}
        g(x)
        = 
        f(\gp_{n}) + \pr*{\frac{f(\gp_{n+1}) - f(\gp_{n})}{\gp_{n+1}-\gp_{n}}} (x-\gp_{n})
    \end{equation}
    \cfload. 
    Then 
    \begin{enumerate}[(i)]
        \item \label{interpolation_309} 
            it holds for all 
                $n\in\Z$ 
            that 
            $
                \gp_{n+1}-\gp_{n} 
                \leq 
                \frac{\varepsilon}{2c (1 + 2\radius)^{\newconstant}}
            $,
        \item \label{interpolation_308b} 
            it holds for all 
                $n\in\Z$ 
            that 
                $g (\gp_{n})=f(\gp_{n})$,
        \item \label{interpolation_307} 
            it holds for all 
                $x,y\in\R$ 
            that 
            $
                \abs{g(x)-g(y)}
                \leq 
                c(1+2\radius)^{\newconstant}\abs{x-y}$, 
            and %
        \item \label{interpolation_308} 
            it holds that 
                $\sup_{x\in[-\radius,\radius]}\abs{f(x)-g(x)}\leq \varepsilon$.
    \end{enumerate}
\end{athm}

\begin{proof}[Proof of \cref{Prop:interpolation_3}]
    \Nobs that
        the fact that 
            for all
                $x,y\in\R$
            it holds that
                $\abs{\clip{-R}R(x)-\clip{-R}R(y)}\leq\abs{x-y}$
    shows that for all
        $n\in\Z$
    it holds that
    \begin{equation}
    \begin{split}
        \gp_{n+1}-\gp_n
        &=
        \abs{\gp_{n+1}-\gp_n}
        \leq
        \abs{(n+1)(2c)^{-1}(1+2\radius)^{-\newconstant} \varepsilon-n(2c)^{-1}(1+2\radius)^{-\newconstant} \varepsilon}
        \\&\leq
        \abs{(2c)^{-1}(1+2\radius)^{-\newconstant} \varepsilon}
        =
        (2c)^{-1}(1+2\radius)^{-\newconstant} \varepsilon
        .
    \end{split}
    \end{equation}
        This
    establishes
        \cref{interpolation_309}.
    \Nobs that 
        the fact that 
            for all 
                $x\in\R$ 
            it holds that 
                $\clip{-R}{R}(-x)=-\clip{-R}{R}(x)$ 
    implies that for all 
        $n\in\Z$ 
    it holds that 
        $\gp_{-n}=-\gp_n$ 
    \cfload. 
        This 
    ensures that for all 
        $M\in\Z\cap[N,\infty)$ 
    it holds that 
    \begin{equation}
        \gp_{M}
        =
        R
        \qandq
        \gp_{-M}
        =
        -R
        .
    \end{equation}
        Hence, 
    we obtain that for all 
        $M\in\Z\cap[N,\infty)$ 
    it holds that
    \begin{equation}
        g(\gp_{M})
        =
        g(R)
        =
        f(R)
        =
        f(\gp_{M})
        \qandq
        g(\gp_{-M})
        =
        g(-R)
        =
        f(-R)
        =
        f(\gp_{-M})
        .
    \end{equation}
    Combining 
        this 
    with 
        \eqref{interpolation_303} 
    proves that for all 
        $n\in\Z\cap[-N,N)$, 
        $x\in[\gp_{n},\gp_{n+1}]$ 
    it holds that
    \begin{equation}
    \label{interpolation_303_boundary}
        g(x)
        =
        f(\gp_{n}) + \pr*{\frac{f(\gp_{n+1}) - f(\gp_{n})}{\gp_{n+1}-\gp_{n}}} (x-\gp_{n})
        .
    \end{equation}
        This 
    establishes 
        \cref{interpolation_308b}.
    \Nobs that 
        \eqref{interpolation_302} and 
        \eqref{interpolation_303_boundary} 
    imply that for all
        $n\in\Z\cap[-N,N)$, 
        $x,y \in [\gp_{n}, \gp_{n+1}]$
    it holds that
    \begin{equation}
    \label{interpolation_306}
    \begin{split}
        \abs{g(x)-g(y)}
        &=
        \vass*{\pr*{\frac{f(\gp_{n+1}) - f(\gp_{n})}{\gp_{n+1}-\gp_{n}}} (x-y)}
        \\&\leq 
        c(1 + \vass*{\gp_{n+1}} + \vass*{\gp_{n}})^{\newconstant} \abs{x-y}
        \\&\leq
        c(1+2\radius)^{\newconstant} \vass{x-y}
        .
    \end{split}
    \end{equation}
        This, 
        \cref{interpolation_308b}, and 
        \eqref{interpolation_302} 
    demonstrate that for all
        $n,m \in \Z\cap[-N,N)$, 
        $x \in [\gp_{n}, \gp_{n+1}]$, 
        $y \in [\gp_m, \gp_{m+1}]$
        with $n < m$
    it holds that
    \begin{equation}\begin{split}
    \label{interpolation_310}
        \vass{g(x)-g(y)}
        &\leq 
        \vass{g(x)-g(\gp_{n+1})}+\vass{g(\gp_{n+1})-g(\gp_m)}+\vass{g(\gp_m)-g(y)}
        \\&=
        \vass{g(x)-g(\gp_{n+1})}+\vass{f(\gp_{n+1})-f(\gp_m)}+\vass{g(\gp_m)-g(y)}
        \\&\leq
        c(1 + \vass*{\gp_{n+1}} + \vass*{\gp_{n}})^{\newconstant} (\gp_{n+1}-x)
            \\&\quad+ c(1 + \vass*{\gp_{n+1}} + \vass*{\gp_{m}})^{\newconstant}(\gp_m-\gp_{n+1})
            \\&\quad  + c(1 + \vass*{\gp_{m}} + \vass*{\gp_{m+1}})^{\newconstant} (y-\gp_m)
        \\&\leq 
        c(1+2\radius)^{\newconstant} (y-x)
        = 
        c(1+2\radius)^{\newconstant} \vass{y-x}
        .
    \end{split}
    \end{equation}
    Combining 
        this 
    with 
        \eqref{interpolation_306} 
    demonstrates that for all 
        $x,y\in[-R,R]$ 
    it holds that
    \begin{equation}
    \label{interpolation_310.1}
        \vass{g(x)-g(y)}
        \leq 
        c (1+2\radius)^{\newconstant}\vass{x-y}
        .
    \end{equation}
    Moreover, \nobs that 
        the assumption that 
            for all 
                $y\in[R,\infty)$ 
            it holds that 
                $g(y)=f(R)$ and 
                $g(-y)=f(-R)$ 
    ensures that for all 
        $x\in\R$ 
    it holds that 
        $g(x)=g\pr*{\clip{-R}{R}(x)}$. 
        This and 
        \eqref{interpolation_310.1} 
    show that for all 
        $x,y\in\R$ 
    it holds that
    \begin{equation}
    \label{interpolation_310.2}
    \begin{split}
        \vass{g(x)-g(y)}
        &= 
        \vass*{g\pr*{\clip{-R}{R}(x)}-g\pr*{\clip{-R}{R}(y)}}
        \\&\leq 
        c (1+2\radius)^{\newconstant}\vass{\clip{-R}{R}(x)-\clip{-R}{R}(y)}
        \\&\leq 
        c (1+2\radius)^{\newconstant}\vass{x-y}
        .
    \end{split}
    \end{equation}
        This 
    establishes 
        \cref{interpolation_307}.
    \Nobs that 
        \eqref{interpolation_302}, 
        \cref{interpolation_308b}, and 
        \cref{interpolation_307} 
    imply that for all
        $n\in \Z\cap[-N,N)$, 
        $x\in [\gp_{n}, \gp_{n+1}]$
    it holds that
    \begin{equation}
    \label{Prop:interpolation_3:eq1}
    \begin{split}
        \abs{f(x) - g(x)}
    &= 
        \abs{f(x)- f(\gp_{n})+f(\gp_{n})-g(x)}
    \\&= 
        \abs{f(x)-f(\gp_{n})+g(\gp_{n})-g(x)}
    \\&\leq 
        \abs{f(x)-f(\gp_{n})}-\abs{g(\gp_{n})-g(x)}
    \\&\leq
        c((1 + \vass*{\gp_{n}} + \vass{x})^{\newconstant}+(1 + 2\radius)^{\newconstant}) \abs{x-\gp_{n}} 
    \\&\leq
        2c(1 + 2\radius)^{\newconstant} \abs{\gp_{n+1}-\gp_{n}}
    \\&\leq
        \frac{2c(1 + 2\radius)^{\newconstant}\varepsilon}{2c (1+2\radius)^{\newconstant}}
    = 
        \varepsilon
    .
    \end{split}
    \end{equation}
        This 
    establishes 
        \cref{interpolation_308}.
    The proof of \cref{Prop:interpolation_3} is thus complete.
\end{proof}
\endgroup

\begingroup
\cfclear
\newcommand{\newconstant}{a}
\begin{athm}{prop}{Prop:locally_Lipschitz_approx_2}
    Let 
        $\radius \in [1, \infty)$, 
        $c\in(0,\infty)$, 
        $\newconstant \in [0, \infty)$, 
        $\varepsilon \in (0, 1]$, 
        $N \in \N$, 
        $(\gp_{n})_{n \in \Z} \subseteq \R$ 
    satisfy for all 
        $n \in \Z$ 
    that 
        $\gp_{n} = \clip{-R}{R}\pr{n(2c)^{-1}(1+2\radius)^{-\newconstant} \varepsilon}$ and 
        $N = \min\{k \in \N \colon \gp_{k} = \radius\}$,
    let
        $f \in C([-\radius,\radius],\R)$ 
    satisfy for all 
        $x,y\in[-\radius,\radius]$ 
    that
    \begin{align}\label{locally_Lipschitz_approx_2_02}
        \abs{f(x)-f(y)}
        \leq 
        c(1+\abs{x}+\abs{y})^{\newconstant}\abs{x-y}
        ,
    \end{align}
    and let 
        $\mathscr{f}\in\ANNs$ 
    satisfy 
        $\mathscr{f} = \interpolatingNN{f(\gp_{-N}),f(\gp_{-N+1}),\ldots, f(\gp_{N})}{\gp_{-N}, \gp_{-N+1},\ldots,\gp_{N} }$
    \cfload.
    Then
    \begin{enumerate}[(i)]
    \item \label{locally_Lipschitz_approx_2:item0}
        it holds that
        $
            \dims(\mathscr{f}) 
        = 
            (1, 2N+1, 1)
        $,
    \item \label{locally_Lipschitz_approx_2:realization}
        it holds that
        $
            \functionReLUANN(\mathscr{f}) 
        \in
            C(\R,\R)
        $,
    \item \label{locally_Lipschitz_approx_2:item6} 
        it holds for all 
            $n\in \Z\cap[-N,N]$ 
        that 
            $(\functionReLUANN (\mathscr{f}) )(\gp_n)=f(\gp_n)$,
    \item \label{locally_Lipschitz_approx_2:item1}
        it holds for all 
            $x,y\in\R$ 
        that 
        $
            \vass{(\functionReLUANN (\mathscr{f}) )(x)-(\functionReLUANN (\mathscr{f}) )(y)}
        \leq 
            c(1+2\radius)^{\newconstant} \vass{x-y}
        $, 
    \item \label{locally_Lipschitz_approx_2:item2}
        it holds that
        $ 
            \sup_{x \in [-\radius,\radius]}\vass{f(x)-(\functionReLUANN (\mathscr{f}) )(x)}
        \leq 
            \varepsilon
        $, and
    \item \label{locally_Lipschitz_approx_2:item3}
        it holds that
        $
            \paramANN(\mathscr{f}) 
        \leq
            12c  (1+2\radius)^{\newconstant} \radius \varepsilon^{-1} + 10
        \leq
            (12c + 10) 3^{\newconstant} \!\radius^{\newconstant+1} \varepsilon^{-1}
        $
    \end{enumerate}
    \cfout.
\end{athm}
\begin{proof}[Proof of \cref{Prop:locally_Lipschitz_approx_2}]
    \Nobs that 
        \cref{DNN_interpolation} and 
        \cref{Prop:interpolation_3} 
    prove 
        \cref{locally_Lipschitz_approx_2:item0,%
        locally_Lipschitz_approx_2:realization,%
        locally_Lipschitz_approx_2:item6,%
        locally_Lipschitz_approx_2:item2,%
        locally_Lipschitz_approx_2:item1}.
    \Nobs that 
        the assumption that 
            $N=\min\{k \in \N \colon \gp_{k} = \radius\}$ 
    ensures that 
        $\gp_{N-1}<\radius$. 
        This 
    implies that 
        $N-1 \leq \frac{2c(1+2\radius)^{\newconstant}\! \radius}{\varepsilon}$. 
        \Cref{locally_Lipschitz_approx_2:item0} 
        hence 
    assures that
    \begin{equation}
    \begin{split}
        \paramANN(\mathscr{f})
    &=
        3(2N) + 4
    \leq 
        6\prb{\tfrac{2c (1+2\radius)^{\newconstant}\! \radius}{\varepsilon} + 1}+4
    \\&=
        12c (1+2\radius)^{\newconstant}\! \radius\varepsilon^{-1} + 10
    \\&\leq
        12c  3^{\newconstant} \!\radius^{\newconstant+1} \varepsilon^{-1} + 10
    \leq
        (12c + 10) 3^{\newconstant} \!\radius^{\newconstant+1} \varepsilon^{-1}
    .
    \end{split}
    \end{equation}
        This 
    establishes 
        \cref{locally_Lipschitz_approx_2:item3}.
    The proof of \cref{Prop:locally_Lipschitz_approx_2} is thus complete.
\end{proof}
\endgroup

\begingroup
\newcommand{\newconstant}{a}
\cfclear
\begin{athm}{cor}{Coro:loc_Lipschitz_cost_constant_targetfct}
Let ${c}\in\R$, $f \in C(\R, \R)$ satisfy for all $x\in \R$ that $f(x)={c}$\cfload.
Then
\begin{enumerate}[(i)]
    \item\label{Coro:loc_Lipschitz_cost_constant_targetfct_item_0} 
        it holds that 
            $\realisation_\ReLU(\affineANN_{0,{c}}) \in C(\R,\R)$,
    \item\label{Coro:loc_Lipschitz_cost_constant_targetfct_item_1} 
        it holds for all 
            $x \in \R$ 
        that 
            $ (\realisation_\ReLU(\affineANN_{0,{c}}))(x)=f(x)$,
    \item\label{Coro:loc_Lipschitz_cost_constant_targetfct_item_2} 
        it holds for all 
            $x,y \in \R$ 
        that 
            $\vass{ (\realisation_\ReLU(\affineANN_{0,{c}}))(x)-(\realisation_\ReLU(\affineANN_{0,{c}}))(y)} =0$, 
        and
    \item\label{Coro:loc_Lipschitz_cost_constant_targetfct_item_3} 
        it holds that 
            $\paramANN(\affineANN_{0,{c}}) = 2$ 
\end{enumerate}
\cfout.
\end{athm}

\begin{proof}[Proof of \cref{Coro:loc_Lipschitz_cost_constant_targetfct}]
    \Nobs[observe] that 
        \eqref{ANNrealization:ass2} and 
        the fact that 
            $\dims(\affineANN_{0,{c}})=(1,1)$ 
    establish 
        \cref{Coro:loc_Lipschitz_cost_constant_targetfct_item_0,%
        Coro:loc_Lipschitz_cost_constant_targetfct_item_1,%
        Coro:loc_Lipschitz_cost_constant_targetfct_item_2,%
        Coro:loc_Lipschitz_cost_constant_targetfct_item_3}.
    The proof of \cref{Coro:loc_Lipschitz_cost_constant_targetfct} is thus complete.
\end{proof}
\endgroup

\begingroup
\newcommand{\newconstant}{a}
\cfclear
\begin{athm}{cor}{Coro:loc_Lipschitz_cost_bis}
    Let 
        $c, \newconstant \in [0, \infty)$, 
        $f \in C(\R,\R)$ 
    satisfy for all 
        $x,y\in \R$ 
    that 
        $\vass{f(x)-f(y)} \leq c(1+\abs{x}+\abs{y})^{\newconstant}  \vass{x-y}$.
    Then it holds for all 
        $\radius \in [1, \infty)$, 
        $\varepsilon \in (0,1]$ 
    that 
    \begin{equation}
        \CostLip{\ReLU}{[-\radius, \radius]}{f|_{[-\radius, \radius]} }{c(1+2\radius)^{\newconstant} }{\varepsilon} 
    \leq 
        12c \radius (1+2\radius)^{\newconstant} \varepsilon^{-1} + 10
    \end{equation}
\cfout.
\end{athm}

\begin{proof}[Proof of \cref{Coro:loc_Lipschitz_cost_bis}]
    \Nobs that 
        \cref{Prop:locally_Lipschitz_approx_2} and 
        \cref{Coro:loc_Lipschitz_cost_constant_targetfct} 
    establish that for all 
        $\radius \in [1, \infty)$, 
        $\varepsilon \in (0,1]$ 
    there exists 
        $\mathscr{f} \in \ANNs$ 
    such that 
    \begin{enumerate}[(I)]
        \item
            it holds that
                $\functionReLUANN (\mathscr{f})\in C(\R,\R)$,
        \item 
            it holds for all 
                $x,y\in\R$ 
            that 
            \begin{equation}
                \vass{(\functionReLUANN (\mathscr{f}) )(x)-(\functionReLUANN (\mathscr{f}) )(y)}
            \leq 
                c(1+2\radius)^{\newconstant} \vass{x-y}
            ,
            \end{equation}
        \item 
            it holds that
            $ 
                \sup_{x \in [-\radius,\radius]}\vass{f(x)-(\functionReLUANN (\mathscr{f}) )(x)}
            \leq 
                \varepsilon
            $, 
            and
        \item
        it holds that
        $
            \paramANN(\mathscr{f})
        \leq
            12c (1+2\radius)^{\newconstant}\!\radius \varepsilon^{-1} + 10
        $
    \end{enumerate}
    \cfload. 
        This and 
        \cref{lem_cost_of_Lip_approx_set_equivalence} 
    establish that for all 
        $\radius \in [1, \infty)$, 
        $\varepsilon \in (0,1]$ 
    it holds that 
        $
            \CostLip{\ReLU}{[-\radius, \radius]}{f}{c(1+2\radius)^{\newconstant} }{\varepsilon} 
        \leq 
            12c (1+2\radius)^{\newconstant}\!\radius \varepsilon^{-1} + 10
        $ \cfload.
    This completes the proof of \cref{Coro:loc_Lipschitz_cost_bis}.
\end{proof}
\endgroup

\begingroup
\newcommand{\newconstant}{a}
\cfclear
\begin{athm}{cor}{Lemma:loc_Comp_theory_15_sum}
    Let 
        $c, \newconstant \in [0, \infty)$, 
        $f \in C(\R,\R)$ 
    satisfy for all 
        $x,y\in \R$ 
    that 
        $
            \vass{f(x)-f(y)} 
        \leq 
            c(1+\abs{x}+\abs{y})^{\newconstant} \vass{x-y}
        $.
    Then 
        $f \in \Capprox{(12c+10) 3^{\newconstant}}{1}$ \cfout.
\end{athm}

\begin{aproof}
    Note that 
        \cref{Lemma:Monotonicity_of_Cost} and 
        \cref{Coro:loc_Lipschitz_cost_bis} 
    establish that for all 
        $\radius \in [1, \infty)$, $\varepsilon \in (0,1]$ 
    it holds that 
    \begin{equation} 
    \begin{split}
        \CostLip{\ReLU}{[-\radius, \radius]}{f_{[-\radius, \radius]}}{c 3^{\newconstant} \!\radius^{\newconstant} }{\varepsilon}
        &\leq
        \CostLip{\ReLU}{[-\radius, \radius]}{f_{[-\radius, \radius]}}{c (1+2\radius)^{\newconstant} }{\varepsilon}
        \\&\leq 
        12c (1+2\radius)^{\newconstant}\!\radius \varepsilon^{-1} + 10
        \\&\leq 
        (12c+10)3^{\newconstant}  \!\radius^{\newconstant+1} \varepsilon^{-1}
    \end{split}
    \end{equation}
    \cfload. 
        This, 
        the fact that 
            $\max\{c 3^{\newconstant}, \newconstant, (12c + 10) 3^{\newconstant}, \newconstant+1, 1\} = \max\{(12c + 10) 3^{\newconstant}, \newconstant+1\} = (12c + 10) 3^{\newconstant}$, 
        and \cref{Lemma:Monotonicity_of_Cost} 
    ensure that for all 
        $\radius \in [1, \infty)$, 
        $\varepsilon \in (0,1]$ 
    it holds that
    \begin{equation} 
    \begin{split}
        \CostLipB{\ReLU}{[-\radius, \radius]}{f|_{[-\radius, \radius]}}{(12c + 10) 3^{\newconstant} \!\radius^{(12c + 10) 3^{\newconstant}} }{\varepsilon}
    &\leq 
        \CostLip{\ReLU}{[-\radius, \radius]}{f|_{[-\radius, \radius]}}{c 3^{\newconstant} \!\radius^{\newconstant} }{\varepsilon}
    \\&\leq 
        (12c + 10) 3^{\newconstant}  \!\radius^{(12c + 10) 3^{\newconstant}} \varepsilon^{-(12c + 10) 3^{\newconstant}}
    .
    \end{split}
    \end{equation}
    The proof of \cref{Lemma:loc_Comp_theory_15_sum} is thus complete.
\end{aproof}
\endgroup

\section{ANN representations for maximum functions}
\label{Section:5}

In this section we present for every $d\in\N$ an explicit exact representation of the $d$-dimensional maximum function as the realization of a ReLU ANN with a bound on the number of parameters that is quadratic in $d$.
As a consequence, we prove in \cref{Coro:max_d_cost} in \cref{Subsection:5.2}
below that certain vector-valued multi-dimensional maximum functions
are in the approximation space for multi-dimensional functions defined in
\cref{Section:3} above.

The results in this section are essentially well-known.
In particular, we refer, e.g., to Beck et al.~\cite[Subsection~3.1.2]{BeckJentzenKuckuck2019},
to Cheridito et al.~\cite[Section IV]{cheridito2021efficient},
and to Jentzen \& Riekert \cite[Subsection~3.2]{jentzen2020strong2}
for closely related results.

\subsection{Explicit ANN representations for maximum functions}
\label{Subsection:5.1}

In this subsection we establish that for every $d\in\N$ the $d$-dimensional maximum function
which maps every $(x_1, x_2, \ldots, x_d) \in \R^d$ to $\max\{x_1, x_2, \ldots, x_d\}\in\R$
can be represented exactly as a realization of an ANN with the ReLU activation function.
We also obtain an upper bound  on the number of parameters (quadratic in the input dimension) of the concrete ANN representation
we employ, which will be instrumental in establishing, in the following subsection, that the maximum functions are in the 
ANN approximation spaces defined in \cref{Section:3} above.

\cfclear
\begin{definition}
[ANN representations of maximum functions]
\label{def:max_d}\cfconsiderloaded{def:max_d}
We denote by $\maxANN_d\in\ANNs$, $d \in \{2,3,\ldots\}$, the \cfadd{def:neuralnetwork}ANNs which satisfy that 
\begin{enumerate}[(i)]
\item \label{max_d_def:1} it holds for all $d \in \{2,3,\ldots\}$ that
$\mathcal{I}(\maxANN_d) = d$,

\item \label{max_d_def:2} it holds for all $d \in \{2,3,\ldots\}$ that
$\mathcal{O}(\maxANN_d) = 1$,  and

\item \label{max_d_def:4} it holds for all $d \in \{2,3,\ldots\}$ that
\begin{equation}
\label{max_d_def:4eq}
    \maxANN_d = 
    \begin{dcases}
       \!\! \mleft( \!\! \mleft( \!\!
	\begin{pmatrix}
		1 & -1 \\
		0 & 1 \\
		0 & -1
	\end{pmatrix},
	\begin{pmatrix}
		0 \\
		0 \\
		0
	\end{pmatrix} \! \! \mright), \mleft(
	\begin{pmatrix}
		1 & 1 & -1
	\end{pmatrix}, 0 \mright) \! \! \mright)
	& \colon d=2 \\
	\maxANN_{\frac{d+1}{2}} \bullet \prb{\parallelizationSpecial_{\frac{d+1}{2}}( \maxANN_2, \maxANN_2 , \ldots, \maxANN_2, \ReLUidANN{1}) }
	& \colon d \in \{3, 5, 7, \ldots\} \\
	\compANN{\maxANN_{\frac{d}{2}} }{ \prb{\parallelizationSpecial_{\frac{d}{2}}( \maxANN_2, \maxANN_2, \ldots, \maxANN_2) }}
	& \colon d \in \{4, 6, 8, \ldots\}
    \end{dcases}
\end{equation}
\end{enumerate}
\cfadd{Lemma:max_d_welldefined}
(cf.\ \cref{def:ANN,def:neuralnetwork,def:simpleParallelization,def:ANNcomposition,def:ReLU_identity} and, e.g., \cite[Lemma~3.7]{jentzen2020strong2}).
\end{definition}

\cfclear
\begin{athm}{lemma}{lem:max_d_lipschitz}
    Let 
        $d \in \N$, 
        $x=(x_1,x_2,\ldots,x_d),\,\allowbreak y=(y_1,y_2,\ldots,y_d)\in\R^d$. 
    Then 
    \begin{equation}
    \begin{split}
        &\vass{\max\{x_1,x_2,\ldots,x_d\}-\max\{y_1,y_2,\ldots,y_d\}}
    \\&\leq
        \max\{\vass{x_1-y_1},\vass{x_2-y_2},\ldots,\vass{x_d-y_d}\}
    \leq
        \norm{x-y}
    \end{split}
    \end{equation}
    \cfout.
\end{athm}

\begin{proof}[Proof of \cref{lem:max_d_lipschitz}]
    Throughout this proof 
        assume w.l.o.g.\ that 
            $\max\{x_1,x_2,\ldots,x_d\}\geq\max\{y_1,y_2,\ldots,y_d\}$ and 
        let 
            $m\in\{1,2,\ldots,d\}$ 
        satisfy 
            $x_m=\max\{x_1,x_2,\ldots,x_d\}$. 
    \Nobs that
    \begin{equation}
    \begin{split}
        &\vass{\max\{x_1,x_2,\ldots,x_d\}-\max\{y_1,y_2,\ldots,y_d\}}
        \\&=
        x_m-\max\{y_1,y_2,\ldots,y_d\}
        \leq 
        x_m-y_m
        =
        \vass{x_m-y_m}
        \\&\leq
        \max\{\vass{x_1-y_1},\vass{x_2-y_2},\ldots,\vass{x_d-y_d}\}
        \\&=
        \PRb{\max\pRb{\vass{x_1-y_1}^2,\vass{x_2-y_2}^2,\ldots,\vass{x_d-y_d}^2}}^{\frac12}
        \\&\leq
        \PRbbb{\ssum_{i=1}^d\abs{x_i-y_i}^2}^{\frac12}
        =
        \norm{x-y}
    \end{split}
    \end{equation}
    \cfload.
    This completes the proof of \cref{lem:max_d_lipschitz}.
\end{proof}

\cfclear
\begin{athm}{prop}{Prop:max_d}
    Let 
        $d \in \{2,3,\ldots\}$.
    Then 
    \begin{enumerate}[(i)]
    \item \label{max_d:item_1} 
        it holds that 
            $\realisation_{\ReLU}\pr*{{\maxANN_d}} \in C(\R^d,\R)$,
    \item \label{max_d:item_2} 
        it holds for all 
            $x = (x_1,x_2, \ldots, x_d)$, 
            $y=(y_1,y_2, \ldots, y_d) \in \R^d$ 
        that 
        $
            \abs{(\realisation_{\ReLU}({\maxANN_d}))(x) 
            - (\realisation_{\ReLU}({\maxANN_d}))(y)}
        \leq 
            \max\{\abs{x_1-y_1}, \abs{x_2-y_2}, \ldots, \abs{x_d-y_d}\} 
        \leq 
            \mednorm{x-y}
        $,
    \item \label{max_d:item_3} 
        it holds for all 
            $x = (x_1,x_2, \ldots, x_d) \in \R^d$ 
        that 
            $(\realisation_{\ReLU}({\maxANN_d}))(x) = \max\{x_1, x_2, \ldots, x_d\}$, 
    \item \label{max_d:item_4} 
        it holds that
            $\lengthANN(\maxANN_d) = \ceil{\log_2(d)}+1$,
    \item \label{max_d:item_5} 
        it holds for all 
            $i \in \N_0$ 
        that
            $\mathbb D_i(\maxANN_d) \leq 3\ceil[\big]{\tfrac{d}{2^i}}$,
        and
    \item \label{max_d:item_6} 
        it holds that
        $\paramANN(\maxANN_d) \leq 3d^2 + 18d + 12\ceil{\log_2(d)} - \tfrac{13}{2}$
    \end{enumerate}
\cfout.
\end{athm}

\begin{aproof}
\Nobs that 
    the fact that $\inDimANN(\maxANN_d)=d$ and $\outDimANN(\maxANN_d)=1$
implies
    \cref{max_d:item_1} \cfload.

\Nobs that
\cite[Proposition~3.10]{jentzen2020strong2}
(see also \cite[Proposition~4.2.7]{jentzen2023mathematical})
establishes
\cref{max_d:item_3,max_d:item_4,max_d:item_5}.
\Nobs that 
    \cref{max_d:item_3} and
    \cref{lem:max_d_lipschitz}
imply
    \cref{max_d:item_2}.

\Nobs that 
    the fact that 
        $\sum_{i=1}^{\infty}\tfrac{1}{4^i}=\tfrac{1}{3}$ 
    and \cref{max_d:item_5} 
show that for all 
    $s \in \{3,4,\ldots\}$ 
it holds that
\begin{equation}
\begin{split}
	\paramANN(\maxANN_s) 
&=
    \sum_{i=0}^{\lengthANN\pr*{\maxANN_s}-1}\!\!\pr*{\singledims_{i}\pr*{\maxANN_s}+1}\singledims_{i+1}\pr*{\maxANN_s}
\\&=
    \PR*{\sum_{i=1}^{\lengthANN\pr*{\maxANN_s}-2}\!\!\pr*{\singledims_{i}\pr*{\maxANN_s}+1}\singledims_{i+1}\pr*{\maxANN_s}}+\pr*{\singledims_{0}\pr*{\maxANN_s}+1}\singledims_{1}\pr*{\maxANN_s}
	\\&\quad+\pr*{\singledims_{\lengthANN\pr*{\maxANN_s}-1}\pr*{\maxANN_s}+1}\singledims_{\lengthANN\pr*{\maxANN_s}}\pr*{\maxANN_s}
\\&\leq 
    \PR*{\sum_{i=1}^{\ceil{\log_2(s)}-1} \!\!\!\! \pr*{ 3 \ceil[\big]{\tfrac{s}{2^i}}+1} 3 \ceil[\big]{\tfrac{s}{2^{i+1}}}}
	+ (s+1) 3 \ceil[\big]{\tfrac{s}{2}} + 4
\\&\leq 
    \PR*{ \sum_{i=1}^{\ceil{\log_2(s)}-1} \!\!\!\! \pr*{ 3 \prb{\tfrac{s}{2^i}+1}+1} 3 \pr*{\tfrac{s}{2^{i+1}}+1}} + (s+1)\prb{\tfrac{3s}{2}+\tfrac{3}{2}} + 4
\\&= 
    \PR*{\sum_{i=1}^{\ceil{\log_2(s)}-1} \!\!\!\! 
    \prb{\tfrac{9s^2}{2^{2i+1}}+\tfrac{15s}{2^i}+12}}
	+ \tfrac{3s^2}{2}+3s + \tfrac{11}{2}
\\&\leq 
    \tfrac{9s^2}{2}\PR*{\sum_{i=1}^{\ceil{\log_2(s)}-1} \!\!\!\! \tfrac{1}{4^{i}}}
    + 15s\PR*{\sum_{i=1}^{\ceil{\log_2(s)}-1} \!\!\!\! \tfrac{1}{2^{i}}}
    + \PR*{\sum_{i=1}^{\ceil{\log_2(s)}-1} \!\!\! 12}
	+ \tfrac{3s^2}{2}+3s + \tfrac{11}{2}
\\&\leq 
    \tfrac{9s^2}{6}
    + 15s
    + 12\ceil{\log_2(s)}- 12
	+ \tfrac{3s^2}{2}+3s + \tfrac{11}{2}
\\&=
    3s^2+ 18s+ 12\ceil{\log_2(s)}- \tfrac{13}{2}
    \ifnocf.
\end{split}
\end{equation}
\cfload[.]%
Combining 
    this 
with 
    the fact that 
        $\paramANN(\maxANN_2)=13\leq 12+ 36+ 12- \tfrac{13}{2}$ 
establishes 
    \cref{max_d:item_6}.
\finishproofthus
\end{aproof}

\subsection{Maximum functions in the spaces of ANN approximable functions}
\label{Subsection:5.2}

The main result of this section, \cref{Coro:max_d_1bis}, establishes, roughly speaking, that a certain sequence of vector-valued
multi-dimensional maximum functions is in the ANN approximation space defined in \cref{def:polyD_mult}.
This  will be 
a major ingredient for the main results of this article.

\cfclear
\begin{athm}{cor}{Coro:max_d_cost}
Let $d \in \N$ and let $F \colon \R^d \to \R^d$ satisfy for all $x = (x_1, x_2, \ldots, x_d) \in \R^d$ that $F(x) = \prb{x_1, \max\{x_1,x_2\}, \ldots,  \max\{x_1,x_2,\dots,x_d\}}$.
Then $\CostLipB{\ReLU}{\R^d}{F}{\sqrt{d}}{0} \leq 3d^4+30d^3$
\cfout.
\end{athm}

\begin{proof}[Proof of \cref{Coro:max_d_cost}]
    Throughout this proof 
        assume w.l.o.g.\ that 
            $d\geq 2$, 
        let 
            $A \in \R^{d^2 \times d}$ 
        satisfy 
            for all 
                $x=(x_1,x_2, \ldots, x_d) \in \R^d$ that 
            \begin{equation}
            \begin{split}
                Ax 
                = 
                (\underbrace{x_1, x_1, \ldots, x_1}_d, 
                    \underbrace{x_1,x_2, x_2, \ldots, x_2}_d, 
                    \ldots, 
                    \underbrace{x_1, x_2, \ldots, x_d}_d)
                ,
            \end{split}
            \end{equation}
        and let 
            $\mathscr{g} \in \ANNs$ 
        satisfy 
            $\mathscr{g} = \compANN{\PR{\parallelizationSpecial_{d}(\maxANN_d, \maxANN_d, \ldots, \maxANN_d)}}{ \affineANN_{A,0}}$ 
            \cfadd{Prop:identity_representation}\cfload.
    \Nobs that 
        \cref{Lemma:PropertiesOfParallelizationEqualLength},
        \cref{Lemma:PropertiesOfCompositions_n2}, and
        \cref{Prop:max_d} 
    ensure that for all 
        $x \in \R^d$ 
    it holds that 
    \begin{equation}\eqlabel{eq:F}
        \realisation_{\ReLU}(\mathscr{g})\in C(\R^d,\R^d)
        \qquad\text{and}\qquad
        (\realisation_{\ReLU}(\mathscr{g}))(x) = F(x)
        \ifnocf.
    \end{equation}
    \cfload[.]%
        This and 
        \cref{lem:max_d_lipschitz} 
    imply that for all 
        $x=(x_1,x_2,\ldots,x_d)$, $y=(y_1,y_2,\ldots,y_d) \in \R^d$ 
    it holds that
    \begin{equation}
    \label{cor:max_d_lipschitz}
    \begin{split}
        \norm{(\realisation_{\ReLU}(\mathscr{g}))(x) - (\realisation_{\ReLU}(\mathscr{g}))(y)}
        &=
        \PR*{\ssum_{i=1}^{d} \vass{\max\{x_1,x_2,\dots,x_i\}-\max\{y_1,y_2,\dots,y_i\}}^2}^{\frac{1}{2}}
        \\&\leq
        \PR*{\ssum_{i=1}^d \norm{(x_1,x_2,\dots,x_i)-(y_1,y_2,\dots,y_i)}^2}^{\frac12}
        \\&\leq
        \PR*{\ssum_{i=1}^{d} \norm{x-y}^2}^{\frac{1}{2}}
        =
        \sqrt{d} \norm{x-y}
    \end{split}
    \end{equation}
    \cfload.
    Observe that 
        \cref{lem:dimcomp} 
    demonstrates that for all 
        $i\in\N$ 
    it holds that 
        $
            \singledims_0(\mathscr{g})
            =
            d\leq d^2
            =
            \singledims_0\pr*{\parallelizationSpecial_{d}(\maxANN_d, \maxANN_d, \ldots, \maxANN_d)}
        $ and 
    \begin{equation}
        \singledims_i(\mathscr{g})
        =
        \singledims_i\pr*{\parallelizationSpecial_{d}(\maxANN_d, \maxANN_d, \ldots, \maxANN_d)}
        .
    \end{equation}
    Combining 
        this,
        \cref{Lemma:ParallelizationImprovedBoundsOne}, and 
        \cref{Prop:max_d} 
    with 
        the fact that 
            for all 
                $n\in\N$ 
            it holds that 
                $\log_2(n)\leq n$ 
    ensures that
    \begin{equation}
    \begin{split}
        \paramANN(\mathscr{g})
        &\leq 
        \paramANN\pr*{\parallelizationSpecial_{d}(\maxANN_d, \maxANN_d, \ldots, \maxANN_d)} 
        \\&\leq 
        d^2\,\paramANN(\maxANN_d)
        \\&\leq 
        d^2\prb{3d^2 + 18d + 12\ceil{\log_2(d)} - \tfrac{13}{2}}
        \\&\leq
        d^2\prb{3d^2 + 18d + 12d}
        \\&=
        3d^4+30d^3
    \end{split}
    \end{equation}
    \cfload. 
        This,
        \eqref{cor:max_d_lipschitz}, and 
        \eqqref{eq:F}
    prove that 
        $\CostLipB{\ReLU}{{\R^d}}{F}{\sqrt{d}}{0} \leq 3d^4+30d^3$ \cfload.
    This completes the proof of \cref{Coro:max_d_cost}.
\end{proof}

\cfclear
\begin{athm}{cor}{Coro:max_d_classes}
    Let 
        $F \colon \prb{\bigcup_{d \in \N}\R^d} \to \prb{\bigcup_{d \in \N}\R^d}$ 
    satisfy for all 
        $d \in \N$, 
        $x = (x_1, x_2, \allowbreak\ldots, x_d) \in \R^d$ 
    that 
        $F(x) = \prb{x_1, \max\{x_1,x_2\}, \ldots,  \max\{x_1,x_2,\dots,x_d\}}$ 
    \cfload.
Then 
    $F \in \Dapprox$
\cfout.
\end{athm}
\begin{proof}[Proof of \cref{Coro:max_d_classes}]
    \Nobs that 
        \cref{Lemma:Monotonicity_of_Cost} and 
        \cref{Coro:max_d_cost} 
    ensure that for all 
        $d,\radius \in \N$, 
        $\varepsilon \in (0,1]$ 
    it holds that 
    \begin{equation}
    \begin{split}
        &\CostLipB{\ReLU}{{[-\radius,\radius]^d}}{\pr{[-\radius,\radius]^d \ni x \mapsto F(x) \in \R^d}}{33d^{33}\! \radius^{33}}{\varepsilon}
        \\&\leq
        \CostLipB{\ReLU}{{[-\radius,\radius]^d}}{\pr{\R^d \ni x \mapsto F(x) \in \R^d}}{\sqrt{d}}{0}
        \\&\leq 
        3d^4+30d^3
        \leq 
        33d^{4}\!\radius \varepsilon^{-1}
        \leq 
        33d^{33}\!\radius^{33} \varepsilon^{-33}
        \ifnocf.
    \end{split}
    \end{equation}
    \cfload[.]%
    Combining 
        this and 
        the fact that 
            for all 
                $d,\radius\in\N$ 
            it holds that 
                $F([-\radius,\radius]^d)\subseteq[-\radius,\radius]^d$ 
    with 
        \eqref{polyD_mult:1} 
    assures that 
        $F \in \Dapprox$ \cfload.
    This completes the proof of \cref{Coro:max_d_classes}.
\end{proof}

\cfclear
\begin{athm}{cor}{Coro:max_d_1bis}
    Let 
        $d \in \N$.
    Then 
        $\CostLip{\ReLU}{\R^d}{\Clip{-1}1d}{1}{0} \leq 3d^2+3d$
    \cfout.
\end{athm}

\begin{aproof}
    \Nobs that
        \cref{lem:clipping_function}
    implies that
        $\CostLip{\ReLU}{\R^d}{\Clip{-1}1d}{1}{0} \leq 3d^2+3d$ \cfload.
    \finishproofthus
\end{aproof}

\section{ANN approximations for product functions}
\label{Section:6}

In this section we establish
in \cref{Prop:new_product_d}
in \cref{Subsection:6.3} below and
in \cref{Coro:prod_d_cost}
in \cref{Subsection:6.4} 
new bounds for the cost
of deep ReLU ANN approximations of the multi-dimensional product functions.
Results related to these findings can be found, e.g., in
Yarotsky~\cite[Subsection~3.1]{yarotsky2017error},
Schwab \& Zech~\cite[Proposition~3.3]{SchwabZech2019}, and
Cheridito et al.~\cite[Proposition~33]{cheridito2021efficient}.

In \cref{Subsection:6.1} we introduce ReLU ANN approximations of the square
function and in \cref{Subsection:6.2} we present ReLU ANN approximations
for the two-dimensional product function. These results are essentially
well-known. 
In particular, \cref{Lemma:square_01} is
proved in Grohs et al.~\cite[Section~5]{GrohsIbrgimovJentzen2021}
(see also Yarotsky \cite[Proposition~2]{yarotsky2017error},
Grohs et al.~\cite[Section~3]{grohs2019space},
and Cheridito et al.~\cite[Proposition~16]{cheridito2021efficient}).
Finally, in \cref{Subsection:6.4} we use these results 
to show that certain sequences of functions involving products
are in the approximation space for multi-dimensional functions introduced in 
\cref{Section:3} above.

\subsection{ANN approximations for the square function on compact intervals}
\label{Subsection:6.1}

In this subsection we establish, in \cref{Lemma:square_r} below, a result on approximations of the square function by ReLU ANNs with
concrete bounds on the number of parameters needed for the approximating ANNs, which will be used in the following subsection to obtain approximation results for the two-dimensional product function.

\begin{athm}{lemma}{Lemma:square_01}
	Let $\varepsilon \in (0,\infty)$.
Then there exists $\mathscr{f} \in \ANNs$ such that
\begin{enumerate}[(i)]
\item\label{square_01:item_1} it holds that $\functionReLUANN{(\mathscr{f})} \in C(\R,\R)$,
\item\label{square_01:item_4} it holds for all $x,y \in \R$ that $\vass{\functionnbReLUANN{\mathscr{f}}(x) - \functionnbReLUANN{\mathscr{f}}(y)} \leq 2\vass{x-y}$,
\item\label{square_01:item_2} it holds that $\sup_{x \in [0,1]} \vass{x^2 - \functionnbReLUANN{\mathscr{f}}(x)} \leq \varepsilon$,
\item\label{square_01:item_7} it holds that $\lengthANN(\mathscr{f}) = \max\pRb{2, \ceil[\big]{\tfrac{1}{2} \log_2(\varepsilon^{-1})} }$,
\item\label{square_01:item_6} it holds that $\dims(\mathscr{f}) = (1, 4, 4, \ldots, 4, 1)$,
and
\item\label{square_01:item_5} it holds that $\paramANN(\mathscr{f}) \leq \max\{13,10\log_2(\varepsilon^{-1}) - 7\}$
\end{enumerate}
(cf.\ \cref{def:ANN,def:ANNrealization,def:ReLU}).
\end{athm}

\begin{aproof}
\Nobs that 
    Lemmas 5.1 and 5.2 in Grohs et al.\ \cite{GrohsIbrgimovJentzen2021} 
prove \cref{square_01:item_1,square_01:item_2,%
square_01:item_5,square_01:item_6,square_01:item_7,square_01:item_4}.
\finishproofthus
\end{aproof}

\cfclear
\begin{athm}{lemma}{Lemma:square_r}
	Let 
        $\radius, \varepsilon \in (0,\infty)$.
    Then there exists 
        $\mathscr{f} \in \ANNs$ 
    such that
    \begin{enumerate}[(i)]
        \item\label{Square_r:item_1} 
            it holds that 
                $\functionReLUANN{(\mathscr{f})} \in C(\R,\R)$,
        \item\label{Square_r:item_4} 
            it holds for all 
                $x,y \in \R$ 
            that 
            $
                \vass{\functionnbReLUANN{\mathscr{f}}(x) - \functionnbReLUANN{\mathscr{f}}(y)} 
                \leq 
                2\radius \vass{x-y}
            $,
        \item\label{Square_r:item_2} 
            it holds that 
            $
                \sup_{x \in [-\radius, \radius]} \vass{x^2 - \functionnbReLUANN{\mathscr{f}}(x)} 
                \leq 
                \varepsilon
            $,
        \item\label{Square_r:item_7} 
            it holds that 
            $
                \lengthANN(\mathscr{f}) 
                =
                \max\pRb{3, \ceil[\big]{ \log_2 (\radius)-\tfrac{1}{2} \log_2(\varepsilon)  + 1}}
            $,
        \item\label{Square_r:item_6} 
            it holds that 
                $\dims(\mathscr{f}) = (1, 2, 4, 4, \ldots, 4, 1)$,
            and
        \item\label{Square_r:item_5} 
            it holds that 
            $
                \paramANN(\mathscr{f}) 
                \leq 
                \max\{21, 20\log_2 (\radius) - 10\log_2(\varepsilon) +1\}
            $
    \end{enumerate}
    \cfout.
\end{athm}

\begin{aproof}
    Throughout this proof 
        let 
            $\varphi_1\in\pr*{\pr*{\R^{2\times 1}\times\R^2}\times\pr*{\R^{1\times 2}\times\R}}\subseteq\ANNs$, 
            $\varphi_2 \in\pr*{\R^{1\times 1}\times\R}\subseteq \ANNs$ 
        satisfy
        \begin{equation}
        \label{Lemma:square_r:networks}
            \varphi_1
            =
            \pr*{\!\pr*{\!
            \begin{pmatrix}
                    \radius^{-1}\\
                    -\radius^{-1}
            \end{pmatrix}
            ,
            \begin{pmatrix}
                    0\\
                    0
            \end{pmatrix}
            \!}
            ,
            \pr*{
            \begin{pmatrix}
                    1 & 1
            \end{pmatrix}
            ,
            0
            }\!}\qandq
            \varphi_2
            =
            \affineANN_{\radius^2,0}
        \end{equation}
    \cfload.
    \Nobs that 
        \eqref{Lemma:square_r:networks} 
    implies that 
    \begin{equation}
    \llabel{eq:dims}
        \dims(\varphi_1) 
        = 
        (1,2,1)
        \qandq
        \dims(\varphi_2) 
        =
        (1,1)
        .
    \end{equation}
    \Moreover
        \eqref{Lemma:square_r:networks} 
    ensures that for all 
        $x \in \R$ 
    it holds that 
    \begin{equation}
    \label{Lemma:square_r:eq:varphi}
        \functionReLUANN(\varphi_1), \functionReLUANN{(\varphi_2)} \in C(\R,\R),
        \quad
        \functionnbReLUANN{\varphi_1}(x) 
        = 
        \frac{\vass{x}}{\radius},
        \quad\text{and}\quad
        \functionnbReLUANN{\varphi_2}(x) 
        = 
        \radius^2x
    \end{equation}
    \cfload. 
    \Moreover
        \cref{Lemma:square_01}  (applied with
		$\varepsilon \curvearrowleft  \frac{\varepsilon}{\radius^2}$
in the notation of \cref{Lemma:square_01})
    shows that there exists 
        $\mathscr{g} \in \ANNs$ 
    which satisfies that
    \begin{enumerate}[(I)]
    \item 
        it holds that 
            $\functionReLUANN{(\mathscr{g})} \in C(\R,\R)$,
		\item\label{Square_r:Lip} 
        it holds for all 
            $x,y \in \R$ 
        that 
        $
            \abs{\functionnbReLUANN{\mathscr{g}}(x) - \functionnbReLUANN{\mathscr{g}}(y)} 
            \leq 
            2\abs{x-y}
        $,
    \item\label{Square_r:approx} 
        it holds that 
        $
            \sup_{x \in [0,1]} \vass{x^2 - \functionnbReLUANN{\mathscr{g}}(x)} 
            \leq 
            \frac{\varepsilon}{\radius^2}
        $,
    \item\label{square_t:7} 
        it holds that 
        $
            \lengthANN(\mathscr{f}) 
            = 
            \max\pRb{2, \ceil[\big]{ \log_2 (\radius)-\tfrac{1}{2} \log_2(\varepsilon)}}
        $,
    \item\label{Square_r:6} 
        it holds that 
            $\dims(\mathscr{g}) = (1, 4, 4, \ldots, 4, 1)$,
        and
    \item\label{Square_r:param} 
        it holds that 
        $
            \paramANN(\mathscr{g}) 
            \leq 
            \max\{13,20\log_2(\radius)-10\log_2(\varepsilon)-7\}
        $\ifnocf.
    \end{enumerate}
    \cfload[.]%
    Next let $\mathscr{f} \in \ANNs$ satisfy
    \begin{equation}
    \label{Square_r:phi}
        \mathscr{f} 
        = 
        \compANN{\compANN{\varphi_2}{\mathscr{g}}}{\varphi_1}
    \end{equation}
    \cfload.
    \Nobs that 
        \lref{eq:dims}, 
        \eqref{Square_r:phi}, 
        \cref{square_t:7}, 
        \cref{Square_r:6}, and 
        \cref{lem:dimcomp} 
    assure that 
    \begin{equation} 
    \label{Square_r:eq_pr_-1}
        \lengthANN(\mathscr{f}) 
        = 
        \max\pRb{3, \ceil[\big]{ \log_2 (\radius)-\tfrac{1}{2} \log_2(\varepsilon)  + 1}}
        \qandq
        \dims(\mathscr{f}) 
        = 
        (1, 2, 4, 4, \ldots, 4, 1)
        .
    \end{equation}
    \Moreover 
        \eqref{Lemma:square_r:eq:varphi}, 
        \eqref{Square_r:phi}, 
        and 
        \cref{Lemma:PropertiesOfCompositions_n2} 
    prove that for all 
        $x \in \R$ 
    it holds that
    \begin{equation}
    \label{Square_r:phi_reduced}
        \functionReLUANN{(\mathscr{f})} \in C(\R,\R)
        \qandq
        \functionnbReLUANN{\mathscr{f}}(x) 
        = 
        \radius^2 \PRb{\functionnbReLUANN{\mathscr{g}}\prb{\tfrac{\vass{x}}{\radius}}}
        \ifnocf.
    \end{equation}
    \cfload[.]%
    Combining 
        this 
    with 
        \cref{Square_r:approx} 
    demonstrates that
    \begin{equation}
    \begin{split}
    \label{Square_r:eq_pr_0}
        &
        \sup\nolimits_{x \in [-\radius, \radius]} \vass{x^2 - \functionnbReLUANN{\mathscr{f}}(x)} = \sup\nolimits_{x \in [-\radius, \radius]} \vass[\big]{\radius^2\PRb{\tfrac{\vass{x}}{\radius}}^2 - \radius^2 \PRb{\functionnbReLUANN{\mathscr{g}}\prb{\tfrac{\vass{x}}{\radius}}}} 
        \\&= 
        \sup\nolimits_{y \in [0,1]} \vass[\big]{\radius^2 y^2 - \radius^2 \PRb{\functionnbReLUANN{\mathscr{g}}(y)}} = \radius^2 \PRb{ \sup\nolimits_{x \in [0,1]} \vass{x^2 - \functionnbReLUANN{\mathscr{g}}(x)} } 
        \leq 
        \varepsilon
        .
    \end{split}
    \end{equation}
    \Moreover 
        \eqref{Square_r:phi_reduced} and 
        \cref{Square_r:Lip} 
    imply that for all 
        $x,y\in\R$ 
    it holds that
    \begin{equation}
    \begin{split}
    \label{Square_r:eq_pr_2}
        &
        \vass*{\functionnbReLUANN{\mathscr{f}}(x)-\functionnbReLUANN{\mathscr{f}}(y) }
        \\&=
        \radius^2\vass[\big]{\functionnbReLUANN{\mathscr{g}}\prb{\tfrac{\vass{x}}{\radius}}-\functionnbReLUANN{\mathscr{g}}\prb{\tfrac{\vass{y}}{\radius}}}
        \leq 
        2\radius^2\vass[\big]{\tfrac{\vass{x}}{\radius}-\tfrac{\vass{y}}{\radius}}
        \leq 
        2\radius\vass{x-y}.
    \end{split}
    \end{equation}
    \Moreover 
        \eqref{Lemma:square_r:eq:varphi}, 
        \cref{Square_r:param}, and 
        \cref{Lemma:PropertiesOfCompositions_n1} 
    show that
    \begin{equation}
    \begin{split}
        \paramANN(\mathscr{f}) 
        &\leq 
        \paramANN(\mathscr{g}) + 2(1+1)+1(2+1)+4(2+1) - 1(2+1) - 4(1+1)
        \\&\leq 
        \max\{13+8,20\log_2 (\radius) - 10\log_2(\varepsilon) -7 + 8\}
        \\&\leq 
        \max\{21,20\log_2 (\radius) - 10\log_2(\varepsilon) +1\}
        \ifnocf.
    \end{split}
    \end{equation}
    \cfload[.]%
    Combining 
        this 
    with 
        \eqref{Square_r:eq_pr_-1}, 
        \eqref{Square_r:phi_reduced}, 
        \eqref{Square_r:eq_pr_0}, and
        \eqref{Square_r:eq_pr_2} 
    establishes 
        \cref{Square_r:item_1,Square_r:item_2,%
        Square_r:item_4,Square_r:item_5,Square_r:item_6,Square_r:item_7}.
    \finishproofthus
\end{aproof}

\subsection{ANN approximations for two-dimensional product functions}
\label{Subsection:6.2}

In the following lemma, we establish a result on approximations by ReLU ANNs of the two-dimensio\-nal product function with concrete upper bounds on the number of parameters needed for the approximating ANNs.

\cfclear
\begin{athm}{lemma}{Lemma:prod_r}
	Let $\radius, \varepsilon \in (0,\infty)$.
Then there exists $\mathscr{f} \in \ANNs$ such that
\begin{enumerate}[(i)]
\item\label{prod_r:item_well_def} 
it holds that 
    $\functionReLUANN{(\mathscr{f})} \in C(\R^2,\R)$,
\item\label{prod_r:item_Lip} 
it holds for all 
    $x,y \in \R^2$ 
that 
    $\vass{\functionnbReLUANN{\mathscr{f}}(x) - \functionnbReLUANN{\mathscr{f}}(y)} \leq \sqrt{32}\radius \norm{x-y}$,
\item\label{prod_r:item_approx_int} 
it holds that 
    $\sup_{x,y \in [-\radius, \radius]} \vass{xy - \functionnbReLUANN{\mathscr{f}}(x,y)} \leq \varepsilon$,
\item\label{prod_r:item_7} 
it holds that 
    $\lengthANN(\mathscr{f}) = \max\pRb{3, \ceil[\big]{\log_2 (\radius)-\tfrac{1}{2} \log_2(\varepsilon) + \tfrac{3}{2} + \tfrac{1}{2}\log_2(3)}}$,
\item\label{prod_r:item_6} 
it holds that 
    $\dims(\mathscr{f}) = (2, 6, 12, 12, \ldots, 12, 1)$,
and
\item\label{prod_r:item_cost} 
it holds that 
    $\paramANN(\mathscr{f}) \leq \max\{157 , 211 + 180\log_2 (\radius) - 90\log_2(\varepsilon)\}$
\end{enumerate}
\cfout.
\end{athm}

\begin{proof}
[Proof of \cref{Lemma:prod_r}.]
Throughout this proof let $\varphi_1 \in\pr*{\R^{3\times 2}\times\R^3}\subseteq \ANNs$, $\varphi_2 \in\pr*{\R^{1\times 3}\times\R}\subseteq \ANNs$ satisfy
\begin{equation}
\label{Lemma:prod_r:eq:networks}
    \varphi_1 = \pr*{ \!\begin{pmatrix} 1 & 1 \\ 1 & 0 \\ 0 & 1   \end{pmatrix}, \begin{pmatrix} 0 \\ 0 \\ 0    \end{pmatrix}\!}\qandq
		\varphi_2 = \pr*{ \tfrac{1}{2}\begin{pmatrix} 1 & -1 & -1    \end{pmatrix}, 0}
\end{equation}
\cfload. \Nobs that \eqref{Lemma:prod_r:eq:networks} ensures that
\begin{equation}
\label{Lemma:prod_r:eq:varphi_new}
\dims(\varphi_1) = (2, 3),\;\;
\dims(\varphi_2) = (3,1),\;\;
\functionReLUANN{(\varphi_1)} \in C(\R^2,\R^3),
\;\;\text{and}\;\;
\functionReLUANN{(\varphi_2)} \in C(\R^3,\R)
\ifnocf. 
\end{equation}
\cfload[.]%
\Moreover \cref{Lemma:prod_r:eq:networks} implies that for all $x,y,z \in \R$ it holds that
\begin{equation}
\label{Lemma:prod_r:eq:varphi}
\functionnbReLUANN{\varphi_1}(x,y) = (x+y,x,y),\qandq\functionnbReLUANN{\varphi_2}(x,y,z) = \frac{x-y-z}{2}
\ifnocf. 
\end{equation}
\cfload[.]\Nobs that \cref{Lemma:square_r} (applied with
		$\radius \curvearrowleft 2\radius$,
		$\varepsilon \curvearrowleft \tfrac{2\varepsilon}{3}$
in the notation of \cref{Lemma:square_r}) shows that there exists $\mathscr{g} \in \ANNs$ which satisfies that
\begin{enumerate}[(I)]
\item \label{prod_r:1} it holds that $\functionReLUANN{(\mathscr{g})} \in C(\R,\R)$,
\item \label{prod_r:3} it holds for all $x,y \in \R$ that $ \vass{ \functionnbReLUANN{\mathscr{g}}(x) - \functionnbReLUANN{\mathscr{g}}(y) } \leq 4\radius  \vass{x-y}$,
\item \label{prod_r:2} it holds that $\sup_{x \in [-2\radius,2\radius]} \vass{x^2 - \functionnbReLUANN{\mathscr{g}}(x)} \leq \frac{2\varepsilon}{3}$,
\item \label{prod_r:4} it holds that $\lengthANN(\mathscr{g}) = \max\pRb{3, \ceil[\big]{ \log_2 (\radius)-\tfrac{1}{2} \log_2(\varepsilon) + \tfrac{3}{2} + \tfrac{1}{2}\log_2(3)}}$,
\item \label{prod_r:6} it holds that $\dims(\mathscr{g}) = (1,2,4, 4, \ldots, 4, 1)$,
and
\item\label{prod_r:param} it holds that $\paramANN(\mathscr{g}) \leq \max\{21,20\log_2 (\radius) - 10\log_2(\varepsilon) + 10\log_2(3) +11\}$.
\end{enumerate}
Next let $\mathscr{f} \in \ANNs$ satisfy
\begin{equation}
\label{prod_r:phi}
	\mathscr{f} = \compANN{\compANN{\varphi_2}{(\parallelizationSpecial_{3}(\mathscr{g},\mathscr{g},\mathscr{g}))}}{\varphi_1}%
\end{equation}
\cfload. 
\Nobs that 
    \eqref{Lemma:prod_r:eq:varphi_new}, 
    \eqref{prod_r:phi},
    \cref{prod_r:6}, 
    \cref{Lemma:PropertiesOfParallelizationEqualLengthDims}, 
    \cref{Lemma:PropertiesOfCompositions_n2}, and
    \cref{lem:dimcomp}
ensure that 
\begin{equation}
\label{prod_r:phi:eq1}
    \lengthANN(\mathscr{f}) = \lengthANN(\mathscr{g})
    \qandq
    \dims(\mathscr{f}) = (2, 6, 12, 12, \ldots, 12, 1) %
    .
\end{equation}
Next \nobs that \cref{Lemma:PropertiesOfCompositions_n2} and \eqref{Lemma:prod_r:eq:varphi} prove that for all $x,y\in\R$ it holds that $\functionReLUANN{(\mathscr{f})} \in C(\R^2,\R)$ and
\begin{equation}
\label{prod_r:phi:prep}
    \functionnbReLUANN{\mathscr{f}}(x,y)
    =
    \tfrac{1}{2}\PRb{\functionnbReLUANN{\mathscr{g}}(x+y)
    -\functionnbReLUANN{\mathscr{g}}(x)
    -\functionnbReLUANN{\mathscr{g}}(y)}.
\end{equation}
This, \cref{prod_r:2}, and \cref{Lemma:PropertiesOfCompositions_n2} demonstrate that for all $x,y \in [-\radius, \radius]$ it holds that
\begin{equation}
\begin{split}
\label{prod_r:phi:eq2}
	&\vass{xy - \functionnbReLUANN{\mathscr{f}}(x,y)}
	\\&= \tfrac{1}{2}\abs[\big]{(x+y)^2 - x^2 - y^2 - \functionnbReLUANN{\mathscr{g}}(x+y) + \functionnbReLUANN{\mathscr{g}}(x) + \functionnbReLUANN{\mathscr{g}}(y)} \\
	&\leq \tfrac{1}{2} \vass{(x+y)^2 - \functionnbReLUANN{\mathscr{g}}(x+y)} + \tfrac{1}{2}\vass{x^2 - \functionnbReLUANN{\mathscr{g}}(x)} 
	+ \tfrac{1}{2}\vass{y^2 - \functionnbReLUANN{\mathscr{g}}(y)} \leq \varepsilon.
\end{split}
\end{equation}
\Moreover 
    \eqref{prod_r:phi:prep} and 
    \cref{prod_r:3} 
show that for all 
    $x_1,x_2,y_1,y_2\in \R$ 
it holds that
\begin{equation}
\begin{split}
\label{prod_r:phi:eq3}
	&\vass{\functionnbReLUANN{\mathscr{f}}(x_1,x_2) - \functionnbReLUANN{\mathscr{f}}(y_1,y_2)}
	\\&\leq
    \tfrac{1}{2}\prb{ \vass{\functionnbReLUANN{\mathscr{g}}(x_1+x_2) - \functionnbReLUANN{\mathscr{g}}(y_1+y_2)} \\
	&\quad+ \vass{\functionnbReLUANN{\mathscr{g}}(x_1) - \functionnbReLUANN{\mathscr{g}}(y_1)} + \vass{\functionnbReLUANN{\mathscr{g}}(x_2) - \functionnbReLUANN{\mathscr{g}}(y_2)} } \\
	&\leq
    2\radius \prb{ \vass{(x_1+x_2)-(y_1+y_2)} + \vass{x_1-y_1} + \vass{x_2-y_2} } \\
	&\leq
    4\radius \pr*{ \vass{x_1-y_1} + \vass{x_2-y_2} } 
    \leq
    \sqrt{32}\radius \norm{(x_1-y_1,x_2-y_2)}
\end{split}
\end{equation}
\cfload. 
\Moreover 
    the fact that 
        $3\leq 256^{\nicefrac{1}{5}}=2^{\nicefrac{8}{5}}$ 
implies that
    $\log_2(3)\leq\tfrac{8}{5}$. 
Combining 
    \eqref{prod_r:phi},
    \cref{prod_r:param},
    and \cref{Lemma:ParallelizationImprovedBoundsOne} 
with 
    \cref{Lemma:PropertiesOfCompositions_n1} 
    therefore 
ensures that
\begin{equation}
\begin{split}
    \paramANN({\mathscr{f}}) 
    &\leq 
    \paramANN(\varphi_2)+ 9 \paramANN(\mathscr{g}) +\paramANN(\varphi_1)+1(12+1)+6(2+1) - 1(3+1)-6(3+1)
        \\&\quad  -3(12+1)-3(2+1) 
        \\&= 9 \paramANN(\mathscr{g}) +44  -76
    \\&\leq 
    \max\{189-32, 9(11 + 10\log_2(3)) -32 + 180\log_2 (\radius) - 90\log_2(\varepsilon)\}
    \\&\leq 
    \max\{157 , 211 + 180\log_2 (\radius) - 90\log_2(\varepsilon)\}
    .
\end{split}
\end{equation}
    This, 
    \eqref{prod_r:phi:eq1}, 
    \eqref{prod_r:phi:prep}, 
    \eqref{prod_r:phi:eq2},  
    \eqref{prod_r:phi:eq3}, and
    \cref{prod_r:4}
establish 
    \cref{prod_r:item_Lip,prod_r:item_7,prod_r:item_cost,prod_r:item_6,prod_r:item_well_def,prod_r:item_approx_int}.
The proof of \cref{Lemma:prod_r} is thus complete\cfload.
\end{proof}

\subsection{ANN approximations for multi-dimensional product functions}
\label{Subsection:6.3}

In this subsection, we employ the results from \cref{Subsection:6.2} on ANN approximations of the two-dimensional product function to obtain results on ANN approximations of multi-dimensional product functions. After a number of preparatory lemmas, we establish the main result of this subsection, \cref{Prop:new_product_d}, which asserts, for every $d\in\N$ and $R\in[1,\infty)$, the existence of ReLU ANNs that approximate the product function $[-R,R]^d\ni(x_1,x_2,\dots,x_d)\mapsto \prod_{i=1}^d x_i\in\R$ up to an arbitrarily small error $\varepsilon\in(0,\infty)$ with concrete upper bounds (polynomial in the input dimension $d$, in the reciprocal of the prescribed approximation error $\varepsilon$, and in the diameter $R$ of the compact set on which the product function is to be approximated) on the number of parameters needed for the approximating ANNs.

\cfclear
\begingroup
\begin{athm}{lemma}{Lem:product_d}
    Let 
        $d \in \N $, 
        $\radius, \varepsilon \in (0, \infty)$.
    Then there exists 
        $\mathscr{f} \in \ANNs$ 
    such that
    \begin{enumerate}[(i)]
        \item \llabel{it:realization} 
            it holds that 
                $\realisation_\ReLU(\mathscr{f}) \in C(\R^{2d},\R^d)$,
        \item \llabel{it:lipschitz} 
            it holds for all 
                $x,y \in \R^{2d}$ 
            that 
            $
                \norm{\functionnbReLUANN{\mathscr{f}} (x) - \functionnbReLUANN{\mathscr{f}} (y)}
                \leq 
                \sqrt{32}\radius\mednorm{x-y}
            $,
        \item \llabel{it:approx} 
            it holds for all 
                $x = (x_1,x_2, \ldots, x_{2d}) \in [-\radius, \radius]^{2d}$ 
            that
            \begin{equation}
                \norm*{\pr*{x_1x_2,x_3x_4,\ldots,x_{2d-1}x_{2d}} - \functionnbReLUANN{\mathscr{f}} (x)} 
                \leq 
                \varepsilon,
            \end{equation}
        \item \llabel{it:length} 
            it holds that 
            $
                \lengthANN(\mathscr{f}) 
                = 
                \max\pRb{3, \ceil[\big]{\log_2(\radius)-\tfrac12 \log_2(\varepsilon) + \tfrac32 + \tfrac12\log_2(3)}}
            $,        
        \item \llabel{it:dims} 
            it holds that
                $\dims(\mathscr{f}) = (2d, 6d, 12d, 12d, \ldots, 12d, d)$, and
        \item \llabel{it:params} 
            it holds that
            $
                \paramANN(\mathscr{f}) 
                \leq 
                d^2\max\{157 , 211 + 45\pr{4\log_2(R) - 2\log_2(\varepsilon)+\log_2(d)}\}
            $
    \end{enumerate}
    \cfout.
\end{athm}

\newcommand{\g}{\mathscr{g}}
\begin{aproof}
    \Nobs that 
        \cref{Lemma:prod_r}
				(applied with
				$\radius \curvearrowleft \radius$,
		$\varepsilon \curvearrowleft d^{-\frac{1}{2}}\varepsilon$	
in the notation of \cref{Lemma:prod_r})
    proves that there exists
        $\g\in \ANNs$ 
    which satisfies that
    \begin{enumerate}[(I)]
    \item \llabel{it:single_realization} 
        it holds that 
            $\realisation_{\ReLU}({\g}) \in C(\R^2,\R)$,
    \item \llabel{it:single_lipschitz} 
        it holds for all 
            $x,y\in\R^{2}$ 
        that 
        $
            \vass{\functionnbReLUANN{\g}(x) - \functionnbReLUANN{\g}(y)}
            \leq 
            \sqrt{32} R \mednorm{x-y}
        $,
    \item \llabel{it:single_approx} 
        it holds that 
        $
            \sup_{x,y \in [-R, R] } \vass{xy - \functionnbReLUANN{\g}(x,y)} 
            \leq 
            d^{-\frac{1}{2}}\varepsilon
        $,
    \item \llabel{it:single_length} 
        it holds that 
        $
            \lengthANN(\g) 
            = 
            \max\pRb{3, \ceil[\big]{\log_2(\radius)-\tfrac12 \log_2(\varepsilon) + \tfrac32 + \tfrac12\log_2(3)}}
        $,   
    \item \llabel{it:single_dims} 
        it holds that 
            $\dims(\g) = (2, 6, 12, 12, \ldots, 12, 1)$,
        and
    \item \llabel{it:single_params} 
        it holds that
        $
            \paramANN(\g) 
            \leq 
            \max\{157 , 211 + 180\log_2(R) - 90\log_2(\varepsilon)+45\log_2(d)\}
        $
    \end{enumerate}
    \cfload. 
    Next let 
        $\mathscr{f} \in \ANNs$ 
    satisfy 
    \begin{equation}
    \begin{split}
    \llabel{eq:fdef}
        \mathscr{f}
        &=
        \parallelizationSpecial_{d} \pr*{\g,\g,\ldots,\g}
        \ifnocf.
    \end{split}
    \end{equation}
		\cfload[.]%
    \Nobs that 
        \lref{eq:fdef}, 
        \lref{it:single_realization}, 
        \lref{it:single_dims}, and 
        \cref{Lemma:PropertiesOfParallelizationEqualLengthDims} 
    ensure that
    \begin{equation}
    \begin{split}
    \llabel{eq:realdims}
        \lengthANN(\mathscr f)=\lengthANN(\g)
        \qandq
        \dims(\mathscr{f}) 
        = 
        d\dims(\g)=(2d, 6d, 12d, 12d, \ldots, 12d, d)
        .
    \end{split}
    \end{equation}
    \Moreover 
        \lref{eq:fdef},
		\lref{it:single_params}, and 
        \cref{Lemma:ParallelizationImprovedBoundsOne} 
    imply that
    \begin{equation}
    \begin{split}
    \llabel{eq:params}
        \paramANN(\mathscr{f})
        \leq 
        d^2\paramANN(\g) 
        \leq 
        d^2\max\{157 , 211 + 180\log_2(R) - 90\log_2(\varepsilon)+45\log_2(d)\}
        .
    \end{split}
    \end{equation}
    \Moreover 
        \lref{eq:fdef},
        \lref{it:single_approx}, and
        \cref{Lemma:PropertiesOfParallelizationEqualLength}
    show that for all 
        $x=(x_1,x_2,\ldots,\allowbreak x_{2d})\in[-R,R]^{2d}$ 
    it holds that
        $\realisation_\ReLU(\mathscr{f}) \in C(\R^{2d},\R^d)$
        and
    \begin{equation}
    \begin{split}
    \llabel{eq:approx}
        &
        \norm*{\pr*{x_1x_2,x_3x_4,\ldots,\allowbreak x_{2d-1}x_{2d}} - \functionnbReLUANN{\mathscr{f}} (x)}
        \\&=
        \PR*{\ssum_{i=1}^d\vass*{x_{2i-1}x_{2i} - \functionnbReLUANN{\g} (x_{2i-1},x_{2i})}^2}^{\frac{1}{2}}
        \leq
        \PR*{\ssum_{i=1}^d d^{-1}\varepsilon^2}^{\frac{1}{2}}
        =
        \varepsilon
        .
    \end{split}
    \end{equation}
    Next we combine 
        \lref{eq:fdef}, 
        \lref{it:single_lipschitz}, and 
        \cref{Lemma:PropertiesOfParallelizationEqualLength}
    with 
        \cref{lemma:Lipschitz_error}
            (applied with
                $L\is \sqrt{32}R$,
                $d\is d$,
                $(g_1,g_2,\dots,g_d)\is(\realisation_{\ReLU}({\g}),\realisation_{\ReLU}({\g}),\dots,\realisation_{\ReLU}({\g}))$,
                $f\is \realisation_{\ReLU}(\mathscr f)$
            in the notation of \cref{lemma:Lipschitz_error})
    to obtain that for all 
        $x,y\in\R^{2d}$ 
    it holds that
    \begin{equation}
        \norm{\functionnbReLUANN{\mathscr f}(x) - \functionnbReLUANN{\mathscr f}(y)}
        \leq 
        \sqrt{32} R \mednorm{x-y}
        .
    \end{equation}
        This,
        \lref{eq:realdims},
        \lref{eq:params},
        \lref{eq:approx}, and
        \lref{it:single_length}        
    establish 
        \cref{\loc.it:realization,\loc.it:lipschitz,\loc.it:approx,\loc.it:length,\loc.it:dims,\loc.it:params}. 
    \finishproofthus
\end{aproof}
\endgroup

\cfclear
\begin{athm}{lemma}{lem:comp_lipschitz_approx}
Let 
    $n\in\N$, 
    $d_0,d_1,\ldots,d_n\in\N$, 
    $L_1,L_2,\ldots,L_n,\varepsilon_1,\varepsilon_2,\ldots,\varepsilon_n\in[0,\infty)$,
for every $i\in\{1,2,\ldots,n\}$
    let $D_i\subseteq\R^{d_{i-1}}$,
for every $i\in\{1,2,\ldots,n\}$ let
    $f_i\colon D_i\to\R^{d_i}$ and 
    $g_i\colon \R^{d_{i-1}}\to\R^{d_i}$ 
satisfy for all $x\in D_i$ that 
\begin{equation}
\label{lem:comp_lipschitz_approx:eq1}
\norm{f_i(x)-g_i(x)}\leq\varepsilon_i,
\end{equation}
and assume for all $j\in\N\cap(0,n)$, $x,y\in\R^{d_j}$ that 
\begin{equation}
\label{lem:comp_lipschitz_approx:eq1.1}
f_j(D_j)\subseteq D_{j+1}
\qquad\text{and}\qquad
\norm{g_{j+1}(x)-g_{j+1}(y)}\leq L_{j+1}\norm{x-y}
\ifnocf.
\end{equation}
\cfload[.]Then it holds for all $x\in D_1$ that
\begin{equation}
\begin{split}
\label{lem:comp_lipschitz_approx:eq4}
&\norm{\pr*{f_n\circ f_{n-1}\circ\ldots\circ f_1}(x)-\pr*{g_n\circ g_{n-1}\circ\ldots\circ g_1}(x)}\leq\ssum_{i=1}^n\PRbbb{\prbbb{\sprod_{j=i+1}^n L_j}\varepsilon_i}.
\end{split}
\end{equation}
\end{athm}

\begin{proof}[Proof of \cref{lem:comp_lipschitz_approx}]
Throughout this proof 
    assume w.l.o.g.\ that 
        $n\geq 2$, 
    let 
        $F_i\colon D_1\to\R^{d_{i}}$, $i\in\{1,2,\ldots,n\}$, 
    satisfy 
        for all 
            $i\in\{1,2,\ldots,n\}$, $x\in D_1$ 
        that
        \begin{equation}
        \label{lem:comp_lipschitz_approx:eq2}
            F_i(x)=\pr*{f_i\circ f_{i-1}\circ\ldots\circ f_1}(x),
        \end{equation}
    and let 
        $G_i\colon \R^{d_{i-1}}\to\R^{d_n}$, $i\in\{1,2,\ldots,n\}$, 
    satisfy 
        for all 
            $i\in\{1,2,\ldots,n\}$, 
            $x\in \R^{d_{i-1}}$ 
        that
        \begin{equation}
        \label{lem:comp_lipschitz_approx:eq3}
            G_i(x)=\pr*{g_n\circ g_{n-1}\circ\ldots\circ g_i}(x).
        \end{equation}
\Nobs that 
    \cref{lem:comp_lipschitz_approx:eq1.1} and
    \eqref{lem:comp_lipschitz_approx:eq3}
imply that for all 
    $i\in\{2,3,\ldots,n\}$, 
    $x,y\in\R^{d_{i-1}}$ 
it holds that
\begin{equation}
\label{lem:comp_lipschitz_approx:eq5}
    \norm{G_i(x)-G_i(y)}
    \leq 
    \prbbb{\sprod_{j=i}^n L_j}\norm{x-y}
    .
\end{equation}
Combining this with \cref{lem:comp_lipschitz_approx:eq1}, \cref{lem:comp_lipschitz_approx:eq2}, and \cref{lem:comp_lipschitz_approx:eq3} ensures that for all $x\in D_1$ it holds that
\begin{align}
&\nonumber\norm{F_n(x)-G_1(x)}
\\&\nonumber=\norm*{F_n(x)-G_n(F_{n-1}(x))+\PR*{\ssum_{i=2}^{n-1}\pr*{G_{i+1}\pr*{F_{i}(x)}-G_{i}\pr*{F_{i-1}(x)}}}+G_2(F_{1}(x))-G_{1}(x)}
\\&\nonumber\leq\norm*{F_n(x)-G_n(F_{n-1}(x))}+\PR*{\ssum_{i=2}^{n-1}\norm*{G_{i+1}\pr*{F_{i}(x)}-G_{i}\pr*{F_{i-1}(x)}}}+\norm*{G_2(F_{1}(x))-G_{1}(x)}
\\&\nonumber=\norm*{f_n(F_{n-1}(x))-g_n(F_{n-1}(x))}+\PR*{\ssum_{i=2}^{n-1}\norm*{G_{i+1}\pr*{f_i\pr*{F_{i-1}(x)}}-G_{i+1}\pr*{g_i\pr*{F_{i-1}(x)}}}}
\\&\quad+\norm*{G_2(f_{1}(x))-G_{2}(g_{1}(x))}
\\&\nonumber\leq\varepsilon_n+\PRbbb{\ssum_{i=2}^{n-1}\prbbb{\sprod_{j=i+1}^n L_j}\norm*{f_i\pr*{F_{i-1}(x)}-g_i\pr*{F_{i-1}(x)}}}+\prbbb{\sprod_{j=2}^n L_j}\norm*{f_{1}(x)-g_{1}(x)}
\\&\nonumber\leq\varepsilon_n+\PRbbb{\ssum_{i=2}^{n-1}\prbbb{\sprod_{j=i+1}^n L_i}\varepsilon_i}+\prbbb{\sprod_{j=2}^n L_j}\varepsilon_1
.
\end{align}
This completes the proof of \cref{lem:comp_lipschitz_approx}.
\end{proof}

\cfclear
\begingroup
\begin{athm}{lemma}{lem:log2_est}
Let 
    $d\in\N$. 
Then 
    $\log_2\pr{d}\leq \frac{5d}{9}$.
\end{athm}

\begin{aproof}
    Throughout this proof let 
        $f\in C((0,\infty),\R)$ 
    satisfy for all 
        $x\in(0,\infty)$ 
    that 
        $f(x)=\frac{x}{\log_2(x)}=\frac{\ln(2)x}{\ln(x)}$. 
    \Nobs that 
        the fact that 
            for all 
                $x\in[3,\infty)$ 
            it holds that 
                $\ln(x)\geq 1$ 
    ensures that for all 
        $x\in[3,\infty)$ 
    it holds that
    \begin{equation}
    \llabel{log2_estimate:0}
    \begin{split}
        f'(x)
        =
        \ln(2)\PR*{\frac{1}{\ln(x)}-\frac{x}{x[\ln(x)]^2}}
        =
        \frac{\ln(2)}{\ln(x)}\PR*{1-\frac{1}{\ln(x)}}
        \geq 
        0
        .
    \end{split}
    \end{equation}
    Combining 
        this 
    with 
        the fact that 
            $\log_2(3)\leq\log_2(32^{\nicefrac{1}{3}})=\frac{5}{3}$ 
        and the fundamental theorem of calculus 
    shows that for all 
        $x\in[3,\infty)$ 
    it holds that
    \begin{equation}
    \llabel{log2_estimate:1}
    \begin{split}
        f(x)
        =
        f(3)+\int_3^x f'(y)\,\mathrm{d}y
				\geq f(3)
				= 
        \frac{3}{\log_2(3)}
        &\geq 
        \frac{9}{5}
        .
    \end{split}
    \end{equation}
        This,
        the fact that 
            $\log_2(1)=0$,
        and the fact that
            $\log_2(2)=1\leq\frac{10}{9}$
    establish that
    $
        \log_2\pr{d}
        \leq 
        \frac{5d}{9}
    $.
\finishproofthus
\end{aproof}
\endgroup

\cfclear
\begingroup
\newcommand{\f}{\mathscr f}
\begin{athm}{lemma}{lem:prod_pow2}
    Let
        $d\in\N$,
        $R,\eps\in(0,\infty)$
    \cfload.
    Then there exists $\f\in\ANNs$ such that
    \begin{enumerate}[(i)]
        \item 
				\label{functiontoshow}
            it holds that
                $\realisation_\ReLU(\f)\in C\prb{\R^{(2^d)},\R}$,
        \item
				\label{lipschitztoshow}
            it holds for all
                $x,y\in\R^{(2^d)}$
            that
            $
                \abs{(\realisation_\ReLU(\f))(x)-(\realisation_\ReLU(\f))(y)}
                \leq
                2^{\frac{5d}{2}}\! R^{(2^{d}-1)}
                \norm{x-y}
            $,
        \item
				\label{approxtoshow}
            it holds for all
                $x=(x_1,x_2,\dots,x_{2^d})\in[-R,R]^{(2^d)}$
            that
            \begin{equation}
                \abs[\bigg]{\PRbbb{\sprod_{i=1}^{{2^d}} x_i}-(\realisation_\ReLU(\f))(x)}
                \leq
                \eps,
            \end{equation}
            and
        \item
				\label{paramtoshow}
            it holds that
                $\paramANN(\f)\leq  426d4^{d}+90\ReLU\pr*{\log_2(R)}8^{d}+90\ReLU\pr*{\log_2\pr*{\eps^{-1}}}4^{d}$
    \end{enumerate}
    \cfout.
\end{athm}
\newcommand{\h}[1]{\mathscr{h}_{#1}}
\begin{aproof}
    Throughout this proof 
        assume w.l.o.g.\ that 
            $d\geq 2$, 
						for every $i\in\{1,2,\ldots,d\}$
        let $D_i\subseteq\R^{\pr{2^{d-i+1}}}$ 
        satisfy 
            $D_i=\PRb{-R^{\pr{2^{i-1}}},R^{\pr{2^{i-1}}}}^{\pr{2^{d-i+1}}}$, and 
        for every $i\in\{1,2,\ldots,d\}$ 
				let $p_i\colon D_i\to\R^{\pr{2^{d-i}}}$
        satisfy for all 
            $x=(x_1,x_2,\dots,x_{2^{d-i+1}})\in D_i$ 
        that 
        \begin{equation}
            p_i(x)
            =
            (x_1x_2,x_3x_4,\ldots,x_{2^{d-i+1}-1}x_{2^{d-i+1}}).
        \end{equation}
    \Nobs that
        \cref{Lem:product_d}
			(applied with
                $d \curvearrowleft 2^{d-i}$,
                $R \curvearrowleft R^{(2^{i-1})}$,
                $\varepsilon \curvearrowleft 2^{\frac{5i-5d}{2}}\! R^{(2^{i}-2^{d})}d^{-1}\eps$
                for $i\in\{1,2,\ldots,d\}$
            in the notation of \cref{Lem:product_d}) 				
    shows that for every $i\in\{1,2,\dots,d\}$ there exists
        $\h i\in\ANNs$ 
    which satisfies that
    \begin{enumerate}[(I)]
        \item 
			\llabel{function}
            it holds 
            that
                $\realisation_\ReLU(\h i)\in C\prb{\R^{(2^{d-i+1})},\R^{(2^{d-i})}}$,
        \item
			\llabel{lipschitz}
            it holds for all
                      $x,y\in\R^{(2^{d-i+1})}$
            that
            $
                \norm{(\realisation_\ReLU(\h i))(x)-(\realisation_\ReLU(\h i))(y)}
                \leq
                \sqrt{32}R^{(2^{i-1})}\norm{x-y}
            $,
        \item
            \llabel{approx}
            it holds for all
                $x\in D_i$
            that
            \begin{equation}
                \norm{p_i(x)-(\realisation_\ReLU(\h i))(x)}
                \leq
                2^{\frac{5i-5d}{2}}\! R^{(2^{i}-2^{d})}d^{-1}\eps
                ,
            \end{equation}
			and
        \item
			\llabel{param}
			it holds
            that
		    \begin{align}
                &\paramANN(\h i)
                \\&\leq\nonumber
                4^{d-i}\max\pRb{157 , 211 + 45\prb{4\log_2\prb{R^{(2^{i-1})}} - 2\log_2\prb{2^{\frac{5i-5d}{2}}\! R^{(2^{i}-2^{d})}d^{-1}\eps}\allowbreak+d-i}}
                \ifnocf.
            \end{align}
    \end{enumerate}
    \cfload[.]%
    Next let 
        $\f\in\ANNs$ 
    satisfy
    \begin{equation}
    \llabel{fdef}
        \f
        =
        \compANN{\mathscr{h}_d}{\ReLUidANN{2}}\bullet\compANN{\mathscr{h}_{d-1}}{\ReLUidANN{2^{2}}}\bullet\ldots\bullet\mathscr{h}_2\bullet\ReLUidANN{2^{d-1}}\bullet\mathscr{h}_1
        \ifnocf.
    \end{equation}
    \cfload[.]%
    \Nobs that 
        \lref{fdef}, 
        \lref{function},
        \cref{Prop:identity_representation}, and
        \cref{Lemma:PropertiesOfCompositions_n2}
    ensure that 
    \begin{equation}
    \llabel{functioncheck}
        \realisation_\ReLU(\f)
        =
        \PR*{\realisation_\ReLU(\mathscr{h}_d)}\circ\PR*{\realisation_\ReLU(\mathscr{h}_{d-1})}\circ\ldots\circ\PR*{\realisation_\ReLU(\mathscr{h}_1)}
        \in C\prb{\R^{(2^d)},\R}
        \ifnocf.
    \end{equation}
    \cfload[.]%
        \Lref{lipschitz}, and 
        induction 
        therefore
    imply that for all 
        $x,y\in\R^{(2^d)}$ 
    it holds that
    \begin{equation}
    \llabel{eq:lipschitz}
    \begin{split}
        \vass[\big]{
            (\realisation_\ReLU(\f))(x)
            -
            (\realisation_\ReLU(\f))(y)
        }
        &\leq 
        \prb{\sqrt{32}}^{d} R^{(2^{d-1}+2^{d-2}+\ldots+2^0)}\norm{x-y}
        \\&= 
        2^{\frac{5d}{2}}\! R^{(2^{d}-1)}\norm{x-y}
        .
    \end{split}
    \end{equation}
    Next \nobs that
        the fact that
            for all
                $i\in\{1,2,\dots,d\}$,
                $x,y\in[-R^{(2^{i-1})},R^{(2^{i-1})}]$
            it holds that
                $xy\in[-R^{(2^{i})},R^{(2^{i})}]$
    demonstrates that for all
        $i\in\N\cap(0,d)$
    it holds that
        $p_i(D_i)\subseteq D_{i+1}$.
    Combining
        this,
        \lref{functioncheck},
        \lref{lipschitz}, and
        \lref{approx}
    with
        \cref{lem:comp_lipschitz_approx} 
            (applied with
                $n \is d$,
                $(d_0,\allowbreak d_1,\dots,d_n)\is(2^d,\allowbreak 2^{d-1},\dots,2^0)$,
                $(L_i)_{i\in\{1,2,\ldots,n\}} \is \prb{\sqrt{32} R^{(2^{i-1})}}{}_{i\in\{1,2,\ldots,d\}}$,
                $(\varepsilon_i)_{i\in\{1,2,\ldots,n\}} \is \prb{2^{\frac{5i-5d}{2}}\! R^{(2^{i}-2^{d})}d^{-1}\eps}{}_{i\in\{1,2,\ldots,d\}}$,
                $(D_1,\allowbreak D_2,\dots,D_n)\is(D_1,\allowbreak D_2,\dots,D_d)$,
                $(f_1,\allowbreak f_2,\dots,f_n)\is(p_1,\allowbreak p_2,\dots,p_d)$,
                $(g_1,\allowbreak g_2,\dots,g_n)\is \pr{\realisation_\ReLU(\mathscr{h}_1),\allowbreak \realisation_\ReLU(\mathscr{h}_2),\dots,\realisation_\ReLU(\mathscr{h}_d)}$
            in the notation of \cref{lem:comp_lipschitz_approx})
    ensures that for all
        $x=(x_1,x_2,\dots,x_{2^d})\in[-R,R]^{(2^d)}$
    it holds that
    \begin{equation}
    \begin{split}
    \llabel{approxcheck}
        &\abs[\bigg]{\PRbbb{\sprod_{i=1}^{{2^d}} x_i} - (\realisation_\ReLU(\f))(x)}
        \\&=
        \abs[\big]{
            \pr*{p_d\circ p_{d-1}\circ \ldots\circ p_1}(x)
            -
            \prb{\PR*{\realisation_\ReLU(\mathscr{h}_d)}\circ\PR*{\realisation_\ReLU(\mathscr{h}_{d-1})}\circ\ldots\circ\PR*{\realisation_\ReLU(\mathscr{h}_1)}}(x)
        }
        \\&\leq
        \ssum_{i=1}^d\PRbbb{\prbbb{\sprod_{j=i+1}^d \sqrt{32} R^{(2^{j-1})}}2^{\frac{5i-5d}{2}}\! R^{(2^{i}-2^{d})}d^{-1}\eps}
        \\&=
        d^{-1}\varepsilon\PRbbb{\ssum_{i=1}^{d}\pr*{2^{\frac{5d-5i}{2}}\! R^{(2^{d}-2^{i})}2^{\frac{5i-5d}{2}}R^{(2^{i}-2^{d})}}}
        =
        \eps
        .
    \end{split}
    \end{equation}
    \Moreover 
        \lref{param} and 
        the fact that
            for all
                $i\in\{1,2,\dots,d\}$ 
            it holds that
    \begin{equation}
    \begin{split}
        &4\log_2\prb{R^{(2^{i-1})}} - 2\log_2\prb{2^{\frac{5i-5d}{2}}\! R^{(2^{i}-2^{d})}d^{-1}\eps}+d-i
        \\&=
        2^{i+1}\log_2(R)-2\prb{\tfrac{5i-5d}2+(2^i-2^d)\log_2(R)-\log_2(d)+\log_2(\eps)}+d-i
        \\&=
        2\pr*{2^{d}\log_2(R) - \log_2\pr{\eps}+\log_2\pr{d}+3d-3i}
        \\&\leq 
        2\pr*{2^{d}\ReLU\pr{\log_2(R)} + \ReLU\pr{\log_2\pr{\eps^{-1}}}+\log_2\pr{d}+3d-3i}
    \end{split}
    \end{equation}
    imply that for all
        $i\in\{1,2,\dots,d\}$ 
    it holds that
    \begin{equation}
    \begin{split}
    \llabel{param_h_i}
        &\paramANN(\h i)
        \\&\leq
        4^{d-i}\max\pRb{157 , 211 + 90\pr*{2^{d}\ReLU\pr{\log_2(R)} + \ReLU\pr{\log_2\pr{\eps^{-1}}}+\log_2\pr{d}+3d-3i}}
        \\&\leq 
        4^{d-i}\PRb{ 90\prb{2^{d}\ReLU\pr{\log_2(R)} + \ReLU\pr{\log_2\pr{\eps^{-1}}}+\log_2\pr{d}}+211+270d}
        .
    \end{split}
    \end{equation}
    Combining 
        this,
        \cref{Lemma:PropertiesOfCompositions_n3}, and 
        \cref{lem:log2_est} 
    with 
        the fact that 
            $\sum_{i=1}^{d}4^{d-i}=\sum_{i=0}^{d-1}4^{i}= \frac{4^d-1}{3}$ 
    shows that
    \begin{equation}
    \begin{split}
    \llabel{paramcheck}
        \paramANN(\f)
        &\leq
        3\PR*{\ssum_{i=1}^d\paramANN(\h i)} %
        \\&\leq
        (4^{d}-1)\PRb{ 90\pr*{2^{d}\ReLU\pr*{\log_2(R)} + \ReLU\pr{\log_2\pr{\eps^{-1}}}+\log_2\pr{d}}+211+270d}
        \\&\leq
        4^d\PRb{ 90\ReLU\pr*{\log_2(R)}2^d + 90\ReLU\pr{\log_2\pr{\eps^{-1}}}+50d+106d+270d}
        \\&\leq
        426d4^{d}+90\ReLU\pr*{\log_2(R)}8^{d}+90\ReLU\pr{\log_2\pr{\eps^{-1}}}4^{d}
        .
    \end{split}
    \end{equation}
    Combining 
        this 
    with 
        \lref{functioncheck}, 
        \lref{eq:lipschitz}, and
        \lref{approxcheck}
    establishes 
        \cref{functiontoshow,lipschitztoshow,approxtoshow,paramtoshow}\cfload. 
    \finishproofthus
    \end{aproof}
    \endgroup

    \cfclear
    \begingroup
    \begin{athm}{prop}{Prop:new_product_d}
    Let $d \in \N$, $\radius\in[1,\infty)$, $\varepsilon \in (0, \infty)$.
    Then there exists $\mathscr{f} \in \ANNs$ such that

    \begin{enumerate}[(i)]
    \item \label{new_continuous_d} it holds 
    that $\realisation_\ReLU(\mathscr{f}) \in C(\R^d,\R)$,

    \item \label{new_Lipschitz_d} it holds for all
    $x,y \in \R^d$ that 
        $\vass{\functionnbReLUANN{\mathscr{f}} (x) 
    - \functionnbReLUANN{\mathscr{f}} (y)} 
    \leq \sqrt{32}d^{\frac{5}{2}}\! R^{2d-1} \mednorm{x-y}$,

    \item\label{new_approx_d} it holds for all
    $x = (x_1,x_2, \ldots, x_d) \in [-\radius, \radius]^d$ that 
    $\vass[\big]{\PRb{\sprod_{i = 1}^d x_i} - \functionnbReLUANN{\mathscr{f}} (x)} \leq \varepsilon$, and

    \item\label{new_cost_d} it holds 
    that
    $\paramANN(\mathscr{f}) \leq 1896  d^3+720 \log_2(R)d^3 + 360\ReLU\pr{\log_2\pr{\eps^{-1}}}d^2 
    $
    \end{enumerate}
    \cfout. %
\end{athm}

\newcommand{\D}{D}
\begin{proof}
[Proof of \cref{Prop:new_product_d}]
Throughout this proof 
assume w.l.o.g.\ that\cfadd{Prop:identity_representation}
    $d\geq 2$,
let 
    $\D\in\N$ 
satisfy 
    $\D=2^{\ceil{\log_2(d)}}$ 
and let 
    $A\in\R^{\D\times d}$, $B\in\R^{\D}$ 
satisfy for all 
    $x=(x_1,x_2,\ldots,x_d)\in\R^d$ 
that
\begin{equation}
\llabel{affinedef}
    Ax+B=(x_1,x_2,\ldots,x_d,1,1,\ldots,1)
\end{equation}
\cfload.
\Nobs that \cref{lem:prod_pow2} (applied with 
$d \curvearrowleft \ceil{\log_2(d)}$,
$R \curvearrowleft R$,
$\eps \curvearrowleft \eps$
 in the notation of \cref{lem:prod_pow2}) ensures that there exists $\mathscr{g}\in\ANNs$ which satisfies that
\begin{enumerate}[(I)]
        \item 
				\llabel{functionprop}
            it holds that
                $\realisation_\ReLU(\mathscr{g})\in C(\R^{\D},\R)$,
        \item
				\llabel{lipschitzprop}
            it holds for all
                $x,y\in\R^{\D}$
            that
            $
                \abs{(\realisation_\ReLU(\mathscr{g}))(x)-(\realisation_\ReLU(\mathscr{g}))(y)}
                \leq
                D^{\frac{5}{2}}\! R^{\D-1}
                \norm{x-y}
            $,
        \item
				\llabel{approxprop}
            it holds for all
                $x=(x_1,x_2,\dots,x_{\D})\in[-R,R]^{\D}$
            that
            $
                \abs[\big]{\PRb{\sprod_{i=1}^{\D} x_i} - (\realisation_\ReLU(\mathscr{g}))(x)}
                \leq
                \eps
            $,
            and
        \item
		\llabel{paramprop}
            it holds that
            \begin{equation}
                \paramANN(\mathscr{g})
                \leq 426\log_2(\D)\D^2+ 90\ReLU\pr*{\log_2(R)}\D^3+ 90\ReLU\pr{\log_2\pr{\eps^{-1}}}\D^2
            \end{equation}
    \end{enumerate}
\cfload[.]%
\Nobs that 
    \cref{Lemma:PropertiesOfCompositions_n2}, 
    \lref{functionprop}, 
    and the fact that 
        $\realisation_\ReLU(\affineANN_{A,B})\in C(\R^d,\R^\D)$ 
imply that 
\begin{equation}
\llabel{functioncheck}
\realisation_\ReLU(\compANN{\mathscr{g}}{\affineANN_{A,B}})=\PR{\realisation_\ReLU(\mathscr{g})}\circ\PR{\realisation_\ReLU(\affineANN_{A,B})}\in C(\R^{d},\R)\ifnocf.
\end{equation}
\cfload[.]%
\Moreover 
    \lref{lipschitzprop},
    \lref{affinedef},
    the fact that 
        $D\leq 2d$,
		and the assumption that $R\geq 1$
show that for all 
    $x,y\in\R^d$ 
it holds that
\begin{equation}
\llabel{lipschitzcheck}
\begin{split}
    &\vass[\big]{(\realisation_\ReLU(\compANN{\mathscr{g}}{\affineANN_{A,B}}))(x)-(\realisation_\ReLU(\compANN{\mathscr{g}}{\affineANN_{A,B}}))(y)}
    \\&=
    \vass*{\prb{\PR{\realisation_\ReLU(\mathscr{g})}\circ\PR{\realisation_\ReLU(\affineANN_{A,B})}}(x)-\prb{\PR{\realisation_\ReLU(\mathscr{g})}\circ\PR{\realisation_\ReLU(\affineANN_{A,B})}}(y)}
    \\&\leq 
    D^{\frac{5}{2}}\! R^{\D-1}\norm{\pr{\realisation_\ReLU(\affineANN_{A,B})}(x)-\pr{\realisation_\ReLU(\affineANN_{A,B})}(y)}
    \\&=
    D^{\frac{5}{2}}\! R^{\D-1}\norm{Ax+B-(Ay+B)}
    \\&=
    D^{\frac{5}{2}}\! R^{\D-1}\norm{x-y}
    \leq 
    \sqrt{32}d^{\frac{5}{2}}\! R^{2d-1}\norm{x-y}
    .
\end{split}
\end{equation}
\Moreover
    \lref{affinedef} and 
    the assumption that 
        $R\geq 1$
ensure that for all
    $x\in[-R,R]^d$ 
it holds that 
    $Ax+B\in[-R,R]^D$.
    \Lref{approxprop}
    therefore
demonstrates that for all 
    $x=(x_1,x_2,\ldots,x_d)\in[-R,R]^d$ 
it holds that
\begin{equation}
\llabel{approxcheck}
\begin{split}
    &\vass*{\PRbbb{\sprod_{i=1}^{d} x_i} - (\realisation_\ReLU(\compANN{\mathscr{g}}{\affineANN_{A,B}}))(x)}
    \\&=
    \vass*{\PRbbb{\sprod_{i=1}^{d} x_i}- \prb{\PR{\realisation_\ReLU(\mathscr{g})}\circ\PR{\realisation_\ReLU(\affineANN_{A,B})}}(x)}
    \\&=
    \vass*{\PRbbb{\sprod_{i=1}^{d} x_i} - \pr{\realisation_\ReLU(\mathscr{g})}(x_1,x_2,\ldots,x_d,1,1,\ldots,1)}\leq\varepsilon
    .
\end{split}
\end{equation}
In the next step \nobs that
    \cref{Lemma:PropertiesOfCompositions_n2} and
    \cref{lem:dimcomp}
imply that  for all 
    $i\in\{1,2,\ldots,\lengthANN(\mathscr{g})\}$ 
it holds that
    $\lengthANN(\compANN{\mathscr{g}}{\affineANN_{A,B}})=\lengthANN(\mathscr g)$,
    $\singledims_0(\compANN{\mathscr{g}}{\affineANN_{A,B}})=\singledims_0(\affineANN_{A,B})=d\leq D=\singledims_0(\mathscr{g})$, and
    $\singledims_i(\compANN{\mathscr{g}}{\affineANN_{A,B}})=\singledims_i(\mathscr{g})$.
Hence, we obtain that
$
\paramANN(\compANN{\mathscr{g}}{\affineANN_{A,B}})
\leq 
\paramANN(\mathscr{g})
$.
Combining 
    this,
    \lref{paramprop},
    and \cref{lem:log2_est} 
with 
    the fact that 
        $D\leq 2d$
and
		the assumption that $R\geq 1$
shows that
\begin{equation}
\begin{split}
\llabel{paramcheck}
\paramANN(\compANN{\mathscr{g}}{\affineANN_{A,B}})
&\leq 
426\log_2(\D)\D^2+ 90\log_2(R)\D^3+ 90\ReLU\pr{\log_2\pr{\eps^{-1}}}\D^2
\\&\leq 
237 \D^3+90 \log_2(R)\D^3 + 90\ReLU\pr{\log_2\pr{\eps^{-1}}}\D^2
\\&\leq 
1896 d^3+720 \log_2(R)d^3 + 360\ReLU\pr{\log_2\pr{\eps^{-1}}}d^2
.
\end{split}
\end{equation}
This,  \lref{functioncheck}, \lref{lipschitzcheck}, and \lref{approxcheck} establish \cref{new_cost_d,new_approx_d,new_Lipschitz_d,new_continuous_d}. The proof of \cref{Prop:new_product_d} is thus complete\cfload.
\end{proof}
\endgroup

\subsection{Product functions in the spaces of ANN approximable functions}
\label{Subsection:6.4}

In this subsection, we use the upper bound on the number of parameters needed for ReLU ANN approximations of the product functions obtained in \cref{Prop:new_product_d} above, to obtain upper bounds on the cost for ReLU ANN approximations of the function $\R^d\ni(x_1,x_2,\ldots,x_d)\mapsto\pr{x_1, x_1x_2, \ldots, x_1x_2 \cdots  x_d}\in\R^d$ for every $d\in \N$ and then establish, in the main result of this subsection, \cref{Coro:prod_d_classes_1}, that this sequence of functions (when suitably clipped) is in the ANN approximation space for multi-dimensional functions defined in \cref{def:polyD_mult} and can in that sense be approximated by ReLU ANNs without the curse of dimensionality. This result will play a major role in obtaining the main results of this article.

\cfclear
\begin{athm}{cor}{Coro:prod_d_cost}
Let 
    $d \in \N$,
    $R\in[1,\infty)$,
    $\eps\in(0,\infty)$
and let 
    $F \colon \R^d \to \R^d$ 
satisfy for all 
    $x = (x_1, x_2, \ldots, x_d) \in \R^d$ 
that 
    $F(x) = \prb{x_1, x_1x_2, \ldots, x_1x_2 \cdots  x_d}$.
Then
\begin{equation}
\label{Coro:prod_d_cost_ineq}
\begin{split}
    &\CostLipB{\ReLU}{[-\radius, \radius]^d}{F|_{[-\radius, \radius]^d}}{\sqrt{32}d^{3}\! R^{2d-1}}{\varepsilon}
    \leq
    2296 d^5+720 \log_2(R)d^5 + 360\ReLU\pr{\log_2\pr{\eps^{-1}}}d^4
\end{split}
\end{equation}
\cfout.
\end{athm}

\begin{aproof}
Throughout this proof 
    assume w.l.o.g.\ that\cfadd{Prop:identity_representation}
        $d\geq 2$
    and let 
        $A \in \R^{d^2 \times d}$, 
        $B \in \R^{d^2}$ 
    satisfy 
        for all 
            $x=(x_1,x_2, \ldots, x_d) \in \R^d$ 
        that
        \begin{equation}
        \label{Coro:prod_d_cost:eq:linear}
            Ax+B 
            =
            (\underbrace{x_1, 1, 1, \ldots, 1}_d, \underbrace{x_1,x_2, 1, \ldots, 1}_d, \ldots, \underbrace{x_1, x_2, \ldots, x_d}_d)
            \ifnocf.
        \end{equation}
\cfload[.]%
\Nobs that 
    \eqref{Coro:prod_d_cost:eq:linear} 
implies that for all 
    $x=(x_1,x_2, \ldots, x_d)$, $y=(y_1,y_2, \ldots, y_d)\in\R^d$ 
it holds that
\begin{equation}
\label{Comp_theory_14:1.1}
\begin{split}
    &\norm*{\pr*{\realisation_{\ReLU}\pr*{\affineANN_{A,B}}}(x)-\pr*{\realisation_{\ReLU}\pr*{\affineANN_{A,B}}}(y)}
    =
    \norm*{Ax+B-(Ay+B)}
    \\&=
    \PR*{\sum_{i=1}^d\sum_{j=1}^i\vass{x_j-y_j}^2}^{\frac{1}{2}}
    \leq
    \PR*{\sum_{i=1}^d\sum_{j=1}^d\vass{x_j-y_j}^2}^{\frac{1}{2}}
    =
    \sqrt{d}\norm{x-y}
    \ifnocf.
\end{split}
\end{equation}
\cfload[.]%
\Moreover 
    \cref{Prop:new_product_d} (applied with
		$d \curvearrowleft d$,
		$\radius \curvearrowleft \radius$,
		$\varepsilon \curvearrowleft d^{-\frac12}\varepsilon$ 
		in the notation of \cref{Prop:new_product_d}) 
proves that there exists
    $\mathscr f\in \ANNs$ 
which satisfies that
\begin{enumerate}[(I)]
\item \label{Coro:prod_d_cost:i1}
    it holds that
        $\realisation_\ReLU{(\mathscr{f})} \in C(\R^d,\R)$,
\item \label{Coro:prod_d_cost:i2}
    it holds for all $x,y\in\R^d$ that
    $
        \vass*{
            \functionnbReLUANN{\mathscr{f}}(x) 
            - 
            \functionnbReLUANN{\mathscr{f}}(y)
        } 
        \leq 
        \sqrt{32}d^{\frac{5}{2}}\! R^{2d-1}\mednorm{x-y}
    $,
\item \label{Coro:prod_d_cost:i3}
    it holds for all $x=(x_1,x_2,\ldots,x_d)\in[-R,R]^d$ that 
    $
    \abs[\big]{
        \PRb{\sprod_{i = 1}^d x_i}
        - 
        \functionnbReLUANN{\mathscr{f}}(x) 
    } 
    \leq 
    d^{-\frac12}\varepsilon
    $,
    and
\item \label{Coro:prod_d_cost:i4}
    it holds that
    $
        \paramANN(\mathscr{f}) 
        \leq
        1896  d^3+720 \log_2(R)d^3 + 360\ReLU\pr{\log_2\pr{d^2\eps^{-1}}}d^2 
    $
\end{enumerate}
\cfload.
Next let 
    $\mathscr g\in \ANNs$ 
satisfy 
\begin{equation}
\label{Comp_theory_14:1}
    \mathscr{g} = \compANN{\parallelizationSpecial_{d}\pr*{\mathscr{f}, \mathscr{f}, \ldots, \mathscr{f} } }{ \affineANN_{A,B}}
\end{equation}
\cfload.
\Nobs that
    \eqref{Comp_theory_14:1}, 
		\cref{Coro:prod_d_cost:i1},
    \cref{Lemma:PropertiesOfParallelizationEqualLength}, and
    \cref{Lemma:PropertiesOfCompositions_n2} ensure that
\begin{equation}
\begin{split}
\label{Comp_theory_14:item_bis_1}  
    \realisation_{\ReLU}( \mathscr{g})=\PR*{\realisation_\ReLU\pr*{\parallelizationSpecial_{d}\pr*{\mathscr{f}, \mathscr{f}, \ldots, \mathscr{f} }}}\circ\PR{\realisation_\ReLU(\affineANN_{A,B})} \in C(\R^d,\R^d)
    \ifnocf.
\end{split}
\end{equation}
\cfload[.]%
    This,
    \eqref{Coro:prod_d_cost:eq:linear},
    \cref{Coro:prod_d_cost:i3},
    \cref{Lemma:PropertiesOfParallelizationEqualLength}, and
    the assumption that $R\geq 1$
show that for all 
    $x=(x_1,x_2, \ldots, x_d) \in [-\radius, \radius]^d$ 
it holds that
\begin{equation}
\label{Comp_theory_14:item_bis_2}
\begin{split}
    \mednorm{F(x) - (\realisation_{\ReLU}( \mathscr{g}))(x)}&= \PR*{\sum_{j=1}^d\vass*{\PR[\big]{\sprod\nolimits_{i = 1}^j x_i}
    - 
    \functionnbReLUANN{\mathscr{f}}(x_1,x_2, \ldots, x_j,1,1,\ldots,1) }^2}^{\frac{1}{2}}
    \\&\leq 
    \PRb{d\prb{d^{-\frac12}\varepsilon}^2}^{\frac{1}{2}}=\varepsilon.
\end{split}
\end{equation}
\Moreover
		\eqref{Comp_theory_14:1.1},
		\eqref{Comp_theory_14:item_bis_1},
		\cref{Coro:prod_d_cost:i2},
		\cref{Lemma:PropertiesOfParallelizationEqualLength}, and
    \cref{lemma:Lipschitz_error} demonstrate that for all 
        $x ,y \in \R^d$ 
    it holds that
\begin{equation}
\label{Comp_theory_14:item_bis_3} 
        \mednorm{(\realisation_{\ReLU}( \mathscr{g}))(x) - (\realisation_{\ReLU}( \mathscr{g}))(y)}
        \leq 
        \sqrt{32}d^{\frac{5}{2}}\! \radius^{2d -1}\prb{\sqrt{d}\mednorm{x-y}}
				=
				\sqrt{32}d^{3}\! R^{2d-1}\mednorm{x-y}
        .
    \end{equation}
In the next step \nobs that
    \cref{Lemma:PropertiesOfCompositions_n2} and
    \cref{lem:dimcomp}
imply that for all 
    $i\in\{1,2,\ldots,\lengthANN(\mathscr{g})\}$ 
it holds that
    $\lengthANN(\mathscr g)=\lengthANN(\parallelizationSpecial_{d}\pr*{\mathscr{f}, \mathscr{f}, \ldots, \mathscr{f} })$,
    $\singledims_0(\mathscr{g})=\singledims_0(\affineANN_{A,B})=d\leq d^2=\singledims_0(\parallelizationSpecial_{d}\pr*{\mathscr{f}, \mathscr{f}, \ldots, \mathscr{f} })$, and
    $\singledims_i(\mathscr{g})=\singledims_i(\parallelizationSpecial_{d}\pr*{\mathscr{f}, \mathscr{f}, \ldots, \mathscr{f} })$.
Hence, we obtain that
$
\paramANN(\mathscr{g})
\leq 
\paramANN(\parallelizationSpecial_{d}\pr*{\mathscr{f}, \mathscr{f}, \ldots, \mathscr{f} })
$.
This,  
    \cref{Coro:prod_d_cost:i4}, 
		and 
    \cref{Lemma:ParallelizationImprovedBoundsOne} 
demonstrate that
\begin{equation}
\label{Comp_theory_14:0}
\begin{split}
    \paramANN(\mathscr{g})
    &\leq 
    d^2\prb{1896  d^3+720 \log_2(R)d^3 + 360\ReLU\pr{\log_2\pr{d^2\eps^{-1}}}d^2}
    \\&=
    d^2\prb{1896  d^3+720 \log_2(R)d^3 + 360\ReLU\pr{2\log_2(d)+\log_2\pr{\eps^{-1}}}d^2}
    \\&\leq
    d^2\prb{1896  d^3+720 \log_2(R)d^3 + 360\ReLU\pr{\log_2\pr{\eps^{-1}}}d^2+720\log_2(d)d^2}
    .
\end{split}
\end{equation}
    \Cref{lem:log2_est}
    hence
shows that
\begin{equation}
\begin{split}
    \paramANN(\mathscr{g})
    &\leq
    d^2\prb{1896  d^3+720 \log_2(R)d^3 + 360\ReLU\pr{\log_2\pr{\eps^{-1}}}d^2+400d^3}
    \\&=
    2296  d^5+720 \log_2(R)d^5 + 360\ReLU\pr{\log_2\pr{\eps^{-1}}}d^4
    .
\end{split}
\end{equation}
Combining 
    this, 
    \cref{Comp_theory_14:item_bis_1}, 
    \cref{Comp_theory_14:item_bis_2}, and 
    \cref{Comp_theory_14:item_bis_3} 
ensures that 
\begin{equation}
\begin{split}
    &\CostLipB{\ReLU}{[-\radius, \radius]^d}{F|_{[-\radius, \radius]^d}}{\sqrt{32}d^3\! \radius^{2d-1}}{\varepsilon} 
    \leq  
    2296  d^5+720 \log_2(R)d^5 + 360\ReLU\pr{\log_2\pr{\eps^{-1}}}d^4
\end{split}
\end{equation}
\cfload[.]%
\finishproofthis
\end{aproof}

\cfclear
\begin{athm}{cor}{Coro:prod_d_classes_1}
Let 
    $F \in C\prb{\bigcup_{d \in \N}\R^d, \bigcup_{d \in \N}\R^d}$ 
satisfy
    for all 
        $d \in \N$, 
        $x = (x_1, x_2, \ldots, x_d) \allowbreak \in \R^d$ 
    that 
        $F(x) = \prb{\clip{-1}{1}(x_1), \allowbreak \prod_{i = 1}^2 \clip{-1}{1}(x_i), \allowbreak \ldots, \allowbreak \prod_{i = 1}^d \clip{-1}{1}(x_i)}$ 
    \cfload.
Then 
    $F \in \Dapprox$
\cfout.
\end{athm}

\begin{aproof}
Observe that 
    \cref{Lemma:Monotonicity_of_Cost} and 
    \cref{Coro:max_d_1bis} 
ensure that for all 
    $d \in \N$, 
    $\radius \in [1,\infty)$ 
it holds that
\begin{equation}
\label{Coro:prod_d_classes_1:eq1}
\begin{split}
    &\CostLipA{\ReLU}{\R^d}{
    (\Clip{-1}1d)|_{[-\radius,\radius]^d}}{1}{0} 
    \leq 
    3d^2+3d
\end{split}
\end{equation}
\cfload. 
    This, 
    \cref{Lemma:Monotonicity_of_Cost}, and 
    \cref{Prop:composition_cost_ReLU} 
        (applied with
            $d_1 \is d$,
            $d_2 \is d$,
            $d_3 \is d$,
            $\varepsilon \is \varepsilon$,
            $L_1 \is 1$,
            $L_2 \is \sqrt{32}d^3$,
            $R_1 \is R$,
            $R_2 \is 1$,
            $f_1 \is (\Clip{-1}1d)|_{[-R,R]^d}$,
            $f_2 \is \prb{[-1,1]^d\ni x\mapsto F(x)\in\R^d}$
        for 
            $d \in \{2,3,\ldots\}$, 
            $\radius \in[1, \infty)$, 
            $\varepsilon \in (0,1]$ 
        in the notation of \cref{Prop:composition_cost_ReLU}) 
show that for all 
    $d \in \{2,3,\ldots\}$, 
    $\radius \in[1, \infty)$, 
    $\varepsilon \in (0,1]$ 
it holds that
\begin{equation}
\begin{split}
    &\CostLipB{\ReLU}{[-\radius, \radius]^d}{\pr{[-\radius, \radius]^d \ni x \mapsto F(x) \in \R^d}}{\sqrt{32}d^3}{\varepsilon}
    \\&=
    \CostLipB{\ReLU}{[-\radius, \radius]^d}{\pr{[-\radius, \radius]^d \ni x \mapsto (F\circ \Clip{-1}1d)(x) \in \R^d}}{\sqrt{32}d^3}{\varepsilon}
    \\&
    \leq 
    4d(d+1)+2\,\CostLipB{\ReLU}{[-\radius, \radius]^d}{(\Clip{-1}1d)|_{[-\radius,\radius]^d}}{1}{0}
    \\&
    \quad
    +2\,\CostLipB{\ReLU}{[-1, 1]^d}{\pr{[-1, 1]^d \ni x \mapsto F(x) \in \R^d}}{\sqrt{32}d^3}{\tfrac{\varepsilon}{2}}
    \\&
    \leq 
    4d(d+1)+2(3d^2+3d)+2\,\CostLipB{\ReLU}{[-1, 1]^d}{\pr{[-1, 1]^d \ni x \mapsto F(x) \in \R^d}}{\sqrt{32}d^3}{\tfrac{\varepsilon}{2}}
		\\&
		=10d^2+10d+2\,\CostLipB{\ReLU}{[-1, 1]^d}{\pr{[-1, 1]^d \ni x \mapsto F(x) \in \R^d}}{\sqrt{32}d^3}{\tfrac{\varepsilon}{2}}
    .
\end{split}
\end{equation}
    This, 
    \cref{Lemma:Monotonicity_of_Cost},
    \cref{Coro:prod_d_cost}, 
    the fact that 
        for all
            $d\in\N$,
            $x=(x_1,x_2,\ldots,x_d)\in[-1,1]^d$ 
        it holds that 
            $F(x)=\pr{x_1, x_1x_2, \ldots, x_1x_2 \cdots  x_d}$,		
    and the fact that 
        for all 
            $\varepsilon \in (0,1]$ 
        it holds that 
            $\ReLU(\log_2(2\varepsilon^{-1}))=\max\{1+\log_2(\varepsilon^{-1}),0\}\leq 1+\varepsilon^{-1}$ 
imply that for all 
    $d \in \{2,3,\ldots\}$, 
    $\radius \in[1, \infty)$, 
    $\varepsilon \in (0,1]$ 
it holds that
\begin{equation}
\label{Coro:prod_d_classes_1:eq5}
\begin{split}
		&\CostLipB{\ReLU}{[-\radius, \radius]^d}{\pr{[-\radius, \radius]^d \ni x \mapsto F(x) \in \R^d}}{14 d^{14}\!\radius^{14}}{\varepsilon}
    \\&\leq\CostLipB{\ReLU}{[-\radius, \radius]^d}{\pr{[-\radius, \radius]^d \ni x \mapsto F(x) \in \R^d}}{\sqrt{32}d^3}{\varepsilon}
    \\&\leq 
    10d^2+10d + 2\prb{2296 d^5 + 360\ReLU\pr{\log_2\pr{2\eps^{-1}}}d^4}
    \\&\leq 
    10d^2+10d + 4592 d^5 + 720d^4+720d^4\varepsilon^{-1}
    \\&\leq 
    4955 d^5 + 360d^5\varepsilon^{-1}
    \\&\leq 
    5315 d^5\varepsilon^{-1}
    \leq 
    14 d^{14}\!\radius^{14}\varepsilon^{-14}
    .
\end{split}
\end{equation}
\Moreover 
    \eqref{Coro:prod_d_classes_1:eq1}, 
    \cref{Lemma:Monotonicity_of_Cost}, and 
    the fact that 
        for all 
            $x\in\R$ 
    it holds that 
        $F(x)=\Clip{-1}11(x)$
ensure that for all     
    $\radius \in[1, \infty)$, 
    $\varepsilon \in (0,1]$ 
it holds that
\begin{equation}
\begin{split}
    \CostLipA{\ReLU}{\R^1}{
    [-\radius,\radius]\ni x\mapsto F(x)\in \R}{14\radius^{14}}{\varepsilon}
		&\leq 
		\CostLipA{\ReLU}{\R^1}{
    (\Clip{-1}11)|_{[-\radius,\radius]}}{1}{0}
		\\&\leq
		6
		\leq 14 \radius^{14}\varepsilon^{-14}.
\end{split}
\end{equation}
Combining
    this and
    \eqref{Coro:prod_d_classes_1:eq5}
with
    the fact that 
        for all 
            $d,\radius\in\N$
    it holds that 
        $\imdim{F}{d}=d$ and
        $F([-\radius,\radius]^d)\subseteq[-1,1]^d$ 
establishes that 
    $F \in \Dapprox$ \cfload.
This completes the proof of \cref{Coro:prod_d_classes_1}\cfload.
\end{aproof}

\section{ANN approximations for high-dimensional functions}
\label{Section:7}
In this section we combine the results of
\cref{Section:4,Section:5,Section:6} above
on the approximation capabilities of ANNs regarding certain concrete classes
of functions with the properties of approximation spaces
proved in \cref{Section:3} above to establish
the main result of this article,
\cref{Theo:example_multiple_composition_loclip} in \cref{Subsection:7.2} below, as well
as several corollaries, which directly imply \cref{Theo:introduction} in the introduction
and the examples given in the introduction.

\subsection{Compositions involving maxima, products, and approx\-i\-ma\-ble functions}
\label{Subsection:7.1}

In this subsection, we establish, roughly speaking, that certain sequences of functions that are obtained as compositions of multi-dimensional product functions, maximum functions, and parallelizations of univariate functions from the approximation spaces for one-dimensional functions defined in \cref{def:polyC} are themselves in the approximation space for multi-dimensional functions defined in \cref{def:polyD_mult} and can in that sense be approximated by ReLU ANNs without the curse of dimensionality. These results are essentially direct consequences of the main results from \cref{Section:5,Section:6,Subsection:3.6} above.

\cfclear
\begin{athm}{prop}{Prop:example_multiple_composition}
	Let $c,r \in [0, \infty)$, $n \in \N$, $a_1, a_2, \ldots, a_n \in \N_0 \cup \{-1\}$, $(f_{k,d})_{(k,d) \in \N^2} \subseteq \Capprox{c}{r}$,
	and let
	$F_{k} \colon \prb{\bigcup_{d \in \N} \R^d} \to \prb{\bigcup_{d \in \N} \R^d}$, $k \in \N_0 \cup \{-1\}$,
	satisfy for all $k,d,\radius \in \N$, $x\in[-R,R]$, $v = (v_1, v_2, \ldots, v_d) \in \R^d$ that
	$\vass{f_{k,d}(x)} \allowbreak\leq cd^c\! R^c$,
	$F_{-1}(v) = \prb{v_1, \allowbreak \max\{v_1,v_2\}, \allowbreak \ldots, \allowbreak \max\{v_1, v_2, \ldots, v_d\} }$, 
	$F_{0} (v) = \prb{\clip{-1}{1}(v_1), \allowbreak \prod_{i = 1}^2 \clip{-1}{1}(v_i), \allowbreak \ldots, \allowbreak \prod_{i = 1}^d \clip{-1}{1}(v_i)}$, and
	$F_{k}(v) = \prb{f_{k,1}(v_1), f_{k,2}(v_2), \allowbreak \ldots,\allowbreak f_{k,d}(v_d)}$
	\cfload.
	Then $
	(F_{a_n} \circ \ldots \circ F_{a_2} \circ F_{a_1}) \in \Dapprox
	$
	\cfout.
\end{athm}

\begin{proof}
	[Proof of \cref{Prop:example_multiple_composition}]
	Note that \cref{Lemma:Comp_theory_3} ensures that
	$\{ F_1, F_2, \ldots \}\subseteq \Dapprox$ \cfload.
	Moreover, observe that \cref{Coro:max_d_classes} shows that $F_{-1} \in \Dapprox$.
	Furthermore, note that \cref{Coro:prod_d_classes_1} proves that $F_{0} \in \Dapprox$.
	Combining this, the fact that $F_{-1} \in \Dapprox$, and the fact that $\{ F_1, F_2, \ldots\} \subseteq \Dapprox$ with \cref{lemma:D_closed} and induction establishes that
	$
	(F_{a_n} \circ \ldots \circ F_{a_2} \circ F_{a_1}) \in \Dapprox
	$.
	This completes the proof of \cref{Prop:example_multiple_composition}.
\end{proof}

\cfclear
\begin{athm}{cor}{Coro:example_multiple_composition_bis}
	Let 
        $c,r \in [0, \infty)$, 
        $n \in \N$, 
    let 
        $a \colon \{1,2, \ldots, n\} \to (\N_0 \cup \{-1\})$ 
    satisfy 
        $\inf\prb{[a^{-1}(\N)]\cup \{\infty\}} \leq \inf\prb{[a^{-1}(\{0\})]\cup \{\infty\}}$,
	let 
        $F_{k} \colon \prb{\bigcup_{d \in \N} \R^d} \to \prb{\bigcup_{d \in \N} \R^d}$, $k \in \N_0 \cup \{-1\}$,
	and let 
        $(f_{k,d})_{(k,d) \in \N^2} \subseteq \Capprox{c}{r}$,
	satisfy 
        for all 
            $k, d \in \N$, 
            $x \in \R$, 
            $v = (v_1, v_2, \ldots, v_d) \in \R^d$ 
        that 
            $\vass{f_{k,d}(x)} \leq 1$,
	        $F_{-1}(v) = \prb{v_1, \allowbreak \max\{v_1,v_2\}, \allowbreak \ldots, \allowbreak \max\{v_1, v_2, \ldots, v_d\} }$, 
	        $F_{0} (v) = \prb{v_1, \allowbreak v_1 v_2, \allowbreak \ldots, \allowbreak v_1 v_2\cdots v_d }$,
	        and $F_{k}(v) = \prb{f_{k,1}(v_1), f_{k,2}(v_2), \allowbreak \ldots, f_{k,d}(v_d)}$
	\cfload.
	Then $
	(F_{a(n)} \circ \ldots \allowbreak \circ F_{a(2)} \circ F_{a(1)}) \in \Dapprox
	$
	\cfout.
\end{athm}
\begin{aproof}
    Throughout this proof 
        let 
            $G_k\colon\prb{\bigcup_{d \in \N} \R^d} \to \prb{\bigcup_{d \in \N} \R^d}$, $k \in \N_0 \cup \{-1\}$,
        satisfy 
            for all
                $k\in\N\cup\{-1\}$,
                $d\in\N$,
                $x=(x_1,x_2,\dots,x_d)\in\R^d$
            that
                $G_k=F_k$
                and
                \begin{equation}
                    \textstyle
                    G_0(x)=\prb{\clip{-1}{1}(x_1), \allowbreak \prod_{i = 1}^2 \clip{-1}{1}(x_i), \allowbreak \ldots, \allowbreak \prod_{i = 1}^d \clip{-1}{1}(x_i)}
                \end{equation}
    \cfload.
    \Nobs that 
        the assumption that
            for all 
                $k,d\in\N$,
                $x\in\R$
            it holds that
                $\abs{f_{k,d}(x)} \leq 1$
    ensures that for all
        $k,d\in\N$,
        $v\in\R^d$
    it holds that
    \begin{equation}
        \llabel{eq1}
        F_k(v)\in[-1,1]^d
        .
    \end{equation}
    Furthermore, \nobs that for all
        $d\in\N$,
        $v\in[-1,1]^d$
    it holds that
    \begin{equation}
        G_{-1}(v) = F_{-1}(v)\in[-1,1]^d
        .
    \end{equation}
        This,
        \lref{eq1}
        and the assumption that
            $\min\prb{[a^{-1}(\N)]\cup \{\infty\}} \leq \min\prb{[a^{-1}(\{0\})]\cup \{\infty\}}$, 
    show that 
        $(F_{a(n)} \circ \ldots \circ F_{a(2)} \circ F_{a(1)}) = (G_{a(n)} \circ \ldots \circ G_{a(2)} \circ G_{a(1)})$.
    Combining this with
        \cref{Prop:example_multiple_composition} 
    establishes that
        $(F_{a(n)} \circ \ldots \circ F_{a(2)} \circ F_{a(1)}) \in \Dapprox$ \cfload.
	\finishproofthis
\end{aproof}

\subsection{Compositions involving maxima, products, and regular functions}
\label{Subsection:7.2}

Combining the results obtained in the previous subsection with the results on ANN approximations of Lipschitz continuous functions from \cref{Section:4}, in this subsection we establish the central result of this article, \cref{Theo:example_multiple_composition_loclip} below, which shows, roughly speaking, that certain sequences of functions that arise as compositions of maximum functions, product functions, and parallelizations of one-dimensional locally Lipschitz functions are in the approximation space for multi-dimensional functions defined in \cref{def:polyD_mult} and can in this sense be approximated by ReLU ANNs without the curse of dimensionality. This conclusion is recast in more basic terms in \cref{Coro:example_multiple_composition_loclip_2}, while \cref{Coro:example_multiple_composition_loclip_2bis} provides a simpler statement with a weaker conclusion, where the upper bound on the number of parameters of the approximating ANNs depends in an arbitrary way -- instead of polynomially -- on the diameter of the region where the approximation property holds. \cref{Theo:introduction} in the introduction is a direct consequence of \cref{Coro:example_multiple_composition_loclip_2bis}.
Finally, \cref{Coro:example_multiple_composition_loclip_proj}, which follows straightforwardly from \cref{Coro:example_multiple_composition_loclip_2bis}, concerns projections of the functions considered in \cref{Theo:example_multiple_composition_loclip} onto their last component, which serves, e.g., to establish the examples \cref{eq:intro_ex1}--\cref{eq:intro_ex12} from the introduction.

\cfclear
\begin{athm}{theorem}{Theo:example_multiple_composition_loclip}
	Let $r \in [0, \infty)$, $n \in \N$, let $a \colon \{1,2, \ldots, n\} \to (\N_0 \cup \{-1\})$ satisfy $\min\prb{[a^{-1}(\N)]\cup \{\infty\}} \leq \min\prb{[a^{-1}(\{0\})]\cup \{\infty\}}$, let
	$f_{k,d} \colon \R \to \R$, $k,d \in \N$, 
	satisfy for all $k, d \in \N$, $x,y \in \R$ that $\vass{f_{k,d}(x)} \leq 1$ and 
	$\vass{f_{k,d}(x) - f_{k,d}(y)} \leq r(1 + \abs{x} + \abs{y})^{r} \abs{x-y}$,
	and let 
	$F_{k} \colon \prb{\bigcup_{d \in \N} \R^d} \to \prb{\bigcup_{d \in \N} \R^d}$, $k \in \N_0 \cup \{-1\}$,
	satisfy for all $d,k \in \N$, $v = (v_1, v_2, \ldots, v_d) \in \R^d$ that
	$F_{-1}(v) = \prb{v_1, \allowbreak \max\{v_1,v_2\}, \allowbreak \ldots, \allowbreak \max\{v_1, v_2, \ldots, v_d\} }$, 
	$F_{0} (v) = \prb{v_1, \allowbreak v_1 v_2, \allowbreak \ldots, \allowbreak v_1 v_2 \cdots  v_d }$,
	and
	$F_{k}(v) = \prb{f_{k,1}(v_1), f_{k,2}(v_2), \allowbreak \ldots, f_{k,d}(v_d)}$\cfload.
	Then $
	(F_{a(n)} \circ \ldots \allowbreak\circ F_{a(2)} \circ F_{a(1)}) \in \Dapprox
	$
	\cfout.
\end{athm}

\begin{proof}
[Proof of \cref{Theo:example_multiple_composition_loclip}]
Note that \cref{Lemma:loc_Comp_theory_15_sum} ensures that for all $k,d \in \N$ it holds that $f_{k,d} \in \Capprox{(12r+10) 3^{r}}{1}$ \cfload.
Combining this with \cref{Coro:example_multiple_composition_bis} establishes that
$ (F_{a(n)} \circ \ldots \circ F_{a(2)} \circ F_{a(1)}) \in \Dapprox$ \cfload.
This completes the proof of \cref{Theo:example_multiple_composition_loclip}.
\end{proof}

\cfclear
\begin{athm}{cor}{Coro:example_multiple_composition_loclip_2}
	Let 
        $r \in [0, \infty)$, 
        $n \in \N$, 
    let 
        $a \colon \{1,2, \ldots, n\} \to (\N_0 \cup \{-1\})$ 
    satisfy 
        $\min\prb{[a^{-1}(\N)]\cup \{\infty\}} \leq \min\prb{[a^{-1}(\{0\})]\cup \{\infty\}}$, 
    let
	    $f_{k,d} \colon \R \to \R$, $k,d \in \N$, 
	satisfy 
        for all 
            $k, d \in \N$, 
            $x,y \in \R$ 
        that 
            $\vass{f_{k,d}(x)} \leq 1$ and 
	        $\vass{f_{k,d}(x) - f_{k,d}(y)} \leq r(1 + \abs{x} + \abs{y})^{r} \abs{x-y}$,
	and let 
	    $F_{k} \colon \prb{\bigcup_{d \in \N} \R^d} \to \prb{\bigcup_{d \in \N} \R^d}$, $k \in \N_0 \cup \{-1\}$,
	satisfy for all 
        $d,k \in \N$, 
        $v = (v_1, v_2, \ldots, v_d) \in \R^d$ 
    that
	    $F_{-1}(v) = \prb{v_1, \allowbreak \max\{v_1,v_2\}, \allowbreak \ldots, \allowbreak \max\{v_1, v_2, \ldots, v_d\} }$, 
	    $F_{0} (v) = \prb{v_1, \allowbreak v_1 v_2, \allowbreak \ldots, \allowbreak v_1 v_2 \cdots  v_d }$,
    	and	$F_{k}(v) = \prb{f_{k,1}(v_1), f_{k,2}(v_2), \allowbreak \ldots, f_{k,d}(v_d)}$
    \cfload.
	Then there exist 
        $(\mathscr{F}_{d,\radius,\varepsilon})_{(d,\radius, \varepsilon) \in \N \times [1, \infty) \times (0,1]} \subseteq \ANNs$ 
        and $c \in \R$ 
    such that for all 
        $d \in \N$, 
        $\radius \in [1, \infty)$, 
        $\varepsilon \in (0,1]$ 
    it holds that
	    $\paramANN(\mathscr{F}_{d,\radius,\varepsilon}) \leq cd^c\!\radius^c\varepsilon^{-c}$, 
        $\realisation_{\ReLU}(\mathscr{F}_{d,\radius,\varepsilon}) \in C( \R^d, \R^d)$, and
	\begin{equation}
	    \sup_{x \in [-\radius,\radius]^d}\norm{ (F_{a(n)} \circ \ldots \circ F_{a(2)} \circ F_{a(1)})(x) - (\realisation_{\ReLU}(\mathscr{F}_{d,\radius,\varepsilon}))(x) } \leq \varepsilon
	\end{equation}
	\cfout.
\end{athm}

\begin{proof}
	[Proof of \cref{Coro:example_multiple_composition_loclip_2}]
	Note that \cref{Theo:example_multiple_composition_loclip} assures that
    $
	(F_{a(n)} \circ \ldots \circ F_{a(2)} \circ F_{a(1)}) \in \Dapprox\ifnocf.
    $
    \cfload[.]This and \cref{Lemma:extension_on_R_poly} imply that 
    there exists $c \in [0, \infty)$ which satisfies that for all $d \in \N$, $\radius \in [1, \infty)$, $\varepsilon \in (0, 1]$
    it holds that 
    $\CostLipB{\ReLU}{[-\radius,\radius]^d}{\pr{[-\radius,\radius]^d \ni x \mapsto (F_{a(n)} \circ \ldots \circ F_{a(2)} \circ F_{a(1)})(x) \in \R^d}}{cd^c \! \radius^c \allowbreak}{\varepsilon} \leq c d^c \! \radius^c \varepsilon^{-c}$ \cfload.
	Combining this with \cref{cor_cost_of_Lip_approx_set_equivalence} establishes that there exist $(\mathscr{F}_{d,\radius,\varepsilon})_{(d,\radius, \varepsilon) \in \N \times [1, \infty) \times (0,1]} \subseteq \ANNs$ such that
	\begin{enumerate}[(I)]
	
    \item\label{example_multiple_composition_loclip_2:item_1} it holds for all $d \in \N$, $\radius \in [1, \infty)$, $\varepsilon \in (0,1]$ that $\realisation_{\ReLU}(\mathscr{F}_{d,\radius,\varepsilon}) \in C(\R^d,\R^d)$,
    
    \item\label{example_multiple_composition_loclip_2:item_2} it holds for all $d \in \N$, $\radius \in [1, \infty)$, $\varepsilon \in (0,1]$ that 
    \begin{equation}
        \sup_{x \in [-R,R]^d}\norm{ (F_{a(n)} \circ \ldots \circ F_{a(2)} \circ F_{a(1)})(x) - (\realisation_{\ReLU}(\mathscr{F}_{d,\radius,\varepsilon}))(x)} \leq \varepsilon,
    \end{equation}

    \item\label{example_multiple_composition_loclip_2:item_3} it holds for all $d \in \N$, $\radius \in [1, \infty)$, $\varepsilon \in (0,1]$, $x,y \in \R^d$ with $x \ne y$ that $\mednorm{(\realisation_{\ReLU}(\mathscr{F}_{d,\radius,\varepsilon}))(x) - (\realisation_{\ReLU}(\mathscr{F}_{d,\radius,\varepsilon}))(y)} \leq cd^c \!\radius^c \mednorm{x-y}$, and
    \item\label{example_multiple_composition_loclip_2:item_4} it holds for all $d \in \N$, $\radius \in [1, \infty)$, $\varepsilon \in (0,1]$ that $\paramANN(\mathscr{F}_{d,\radius,\varepsilon}) \leq cd^c \!\radius^c \varepsilon^{-c}$\ifnocf.
    \end{enumerate}
\cfload[.]This completes the proof of \cref{Coro:example_multiple_composition_loclip_2}.
\end{proof}

\cfclear
\begin{athm}{cor}{Coro:example_multiple_composition_loclip_2bis}
	Let $r,R \in [0, \infty)$, $n \in \N$, $a_1,a_2,\ldots,a_n \in \N_0 \cup \{-1\}$, let
	$f_{k,d} \colon \R \to \R$, $k,d \in \N$, satisfy for all $k, d \in \N$, $x,y \in \R$ that $\vass{f_{k,d}(x)} \leq 1 \leq a_1$ and 
	$\vass{f_{k,d}(x) - f_{k,d}(y)} \leq r(1 + \abs{x} + \abs{y})^{r} \abs{x-y}$,
	and let 
	$F_{k} \colon \prb{\bigcup_{d \in \N} \R^d} \to \prb{\bigcup_{d \in \N} \R^d}$, $k \in \N_0 \cup \{-1\}$,
	satisfy for all $k,d \in \N$, $v = (v_1, v_2, \ldots, v_d) \in \R^d$ that
	$F_{-1}(v) = \prb{v_1, \allowbreak \max\{v_1,v_2\}, \allowbreak \ldots, \allowbreak \max\{v_1, v_2, \ldots, v_d\} }$, 
	$F_{0} (v) = \prb{v_1, \allowbreak v_1 v_2, \allowbreak \ldots, \allowbreak v_1 v_2{}\cdots{} v_d }$,
	and
	$F_{k}(v) = \prb{f_{k,1}(v_1), f_{k,2}(v_2), \allowbreak \ldots, f_{k,d}(v_d)}$\cfload.
Then there exist $(\mathscr{F}_{d,\varepsilon})_{(d, \varepsilon) \in \N \times (0,1]} \subseteq \ANNs$ and $c \in \R$ such that for all $d \in \N$, $\varepsilon \in (0,1]$ it holds that
	$\paramANN(\mathscr{F}_{d,\varepsilon}) \leq cd^c\varepsilon^{-c}$, $\realisation_{\ReLU}(\mathscr{F}_{d,\varepsilon}) \in C( \R^d, \R^d)$, and
	\begin{equation}
	    \sup_{x \in [-R,R]^d}\norm{ (F_{a_n} \circ \ldots \circ F_{a_2} \circ F_{a_1})(x) - (\realisation_{\ReLU}(\mathscr{F}_{d,\varepsilon}))(x) } \leq \varepsilon
	\end{equation}
	\cfout.
\end{athm}

\begin{proof}
	[Proof of \cref{Coro:example_multiple_composition_loclip_2bis}]
	Note that the assumption that $a_1\geq 1$ implies that 
	$
	\inf\pr{\{i\in\{1,2,\ldots,n\}\colon a_i\in\N\}\cup\{\infty\}}=1\leq \inf\pr{\{i\in\{1,2,\ldots,n\}\colon a_i=0\}\cup\{\infty\}}.
	$
	\cref{Coro:example_multiple_composition_loclip_2} hence ensures that there exist $(\mathscr{F}_{d,\varepsilon})_{(d, \varepsilon) \in \N \times (0,1]} \subseteq \ANNs$ and $c \in \R$ such that for all $d \in \N$, $\varepsilon \in (0,1]$ it holds that
	$\paramANN(\mathscr{F}_{d,\varepsilon}) \leq cd^c\varepsilon^{-c}$, $\realisation_{\ReLU}(\mathscr{F}_{d,\varepsilon}) \in C( \R^d, \R^d)$, and
	\begin{equation}
	    \sup_{x \in [-R,R]^d}\norm{ (F_{a_n} \circ \ldots \circ F_{a_2} \circ F_{a_1})(x) - (\realisation_{\ReLU}(\mathscr{F}_{d,\varepsilon}))(x) } \leq \varepsilon\ifnocf.
	\end{equation}
\cfload[.]This completes the proof of \cref{Coro:example_multiple_composition_loclip_2bis}.
\end{proof}

\cfclear
\begin{athm}{cor}{Coro:example_multiple_composition_loclip_proj}
Let 
    $r,R \in [0, \infty)$, 
    $n \in \N$, 
    $a_1,a_2,\ldots,a_n \in \N_0 \cup \{-1\}$,
let
    $p_d\colon\R^d\to\R$, $d\in\N$,
satisfy 
    for all
        $d\in\N$,
        $x=(x_1,x_2,\dots,x_d)\in\R^d$
    that
        $p_d(x)=x_d$,
let
	$f_{k,d} \colon \R \to \R$, $k,d \in \N$, 
satisfy
    for all 
        $k, d \in \N$, 
        $x,y \in \R$ 
    that 
        $\vass{f_{k,d}(x)} \leq 1 \leq a_1$ and 
	    $\vass{f_{k,d}(x) - f_{k,d}(y)} \leq r(1 + \abs{x} + \abs{y})^{r} \abs{x-y}$,
and let 
	$F_{k} \colon \prb{\bigcup_{d \in \N} \R^d} \to \prb{\bigcup_{d \in \N} \R^d}$, $k \in \N_0 \cup \{-1\}$,
satisfy 
    for all 
        $k,d \in \N$, 
        $v = (v_1, v_2, \ldots, v_d) \in \R^d$ 
    that
	    $F_{-1}(v) = \prb{v_1, \allowbreak \max\{v_1,v_2\}, \allowbreak \ldots, \allowbreak \max\{v_1, v_2, \ldots, v_d\} }$, 
	    $F_{0} (v) = \prb{v_1, \allowbreak v_1 v_2, \allowbreak \ldots, \allowbreak v_1 v_2{}\cdots{} v_d }$,
	    and
	    $F_{k}(v) = \prb{f_{k,1}(v_1), f_{k,2}(v_2), \allowbreak \ldots, f_{k,d}(v_d)}$\cfload.
Then there exist 
    $(\mathscr{F}_{d,\varepsilon})_{(d, \varepsilon) \in \N \times (0,1]} \subseteq \ANNs$ 
    and $c \in \R$ 
such that 
    for all 
        $d \in \N$, 
        $\varepsilon \in (0,1]$ 
    it holds that
	    $\paramANN(\mathscr{F}_{d,\varepsilon}) \leq cd^c\varepsilon^{-c}$, 
        $\realisation_{\ReLU}(\mathscr{F}_{d,\varepsilon}) \in C( \R^d, \R)$, 
        and
        \begin{equation}
            \llabel{claim}
            \sup_{x \in [-R,R]^d}\abs{ (p_d\circ F_{a_n} \circ \ldots \circ F_{a_2} \circ F_{a_1})(x) - (\realisation_{\ReLU}(\mathscr{F}_{d,\varepsilon}))(x) } \leq \varepsilon
        \end{equation}
\cfout.
\end{athm}
\begin{aproof}
    Throughout this proof
        let 
            $P_d\in\R^{1\times d}$, $d\in\N$, 
        satisfy for all 
            $d\in\{2,3,\dots\}$ 
        that
            $P_1 = (1)\in\R^{1\times 1}$ and
            $P_d=\begin{pmatrix}0&0&\cdots&0&1\end{pmatrix}\in\R^{1\times d}$.
	\Nobs that
        \cref{Coro:example_multiple_composition_loclip_2bis}
    establishes that there exist 
        $(\mathscr{G}_{d,\varepsilon})_{(d, \varepsilon) \in \N \times (0,1]} \subseteq \ANNs$ 
        and 
        $c \in \R$ 
    which satisfy that for all 
        $d \in \N$, 
        $\varepsilon \in (0,1]$ 
    it holds that
	    $\paramANN(\mathscr{G}_{d,\varepsilon}) \leq cd^c\varepsilon^{-c}$, 
        $\realisation_{\ReLU}(\mathscr{G}_{d,\varepsilon}) \in C( \R^d, \R^d)$, and
	\begin{equation}
        \llabel{eq:approx}
	    \sup_{x \in [-R,R]^d}\norm{ (F_{a_n} \circ \ldots \circ F_{a_2} \circ F_{a_1})(x) - (\realisation_{\ReLU}(\mathscr{G}_{d,\varepsilon}))(x) } 
        \leq 
        \varepsilon
	\end{equation}
    \cfload[.]%
    \Moreover
        \cref{lem:dimcomp}
    implies that for all
        $d\in\N$,
        $\eps\in(0,1]$,
        $n\in\{0,1,\dots,\lengthANN(\mathscr{G}_{d,\varepsilon})-1\}$
    it holds that
        $\lengthANN(\compANN{\affineANN_{P_d,0}}{\mathscr{G}_{d,\varepsilon}})=\lengthANN(\mathscr{G}_{d,\varepsilon})$,
        $\singledims_n(\compANN{\affineANN_{P_d,0}}{\mathscr{G}_{d,\varepsilon}})=\singledims_n(\mathscr{G}_{d,\varepsilon})$,
        and 
        $\singledims_{\lengthANN(\mathscr{G}_{d,\varepsilon})}(\compANN{\affineANN_{P_d,0}}{\mathscr{G}_{d,\varepsilon}})=1\leq \singledims_{\lengthANN(\mathscr{G}_{d,\varepsilon})}(\mathscr{G}_{d,\varepsilon})$
    \cfload.
        The fact that
            for all
                $d\in\N$,
                $\eps\in(0,1]$
        it holds that
            $\paramANN(\mathscr{G}_{d,\varepsilon}) \leq cd^c\varepsilon^{-c}$
        \hence 
    demonstrates that for all
        $d\in\N$,
        $\eps\in(0,1]$
    it holds that
    \begin{equation}
        \llabel{eq:1}
        \paramANN(\compANN{\affineANN_{P_d,0}}{\mathscr{G}_{d,\varepsilon}})
        \leq
        \paramANN(\mathscr{G}_{d,\varepsilon})
        \leq 
        cd^c\varepsilon^{-c}
        .
    \end{equation}
    Next \nobs that
        \cref{Lemma:PropertiesOfCompositions_n2} and
        the fact that
            for all
                $d\in\N$,
                $x\in\R^d$
            it holds that
                $p_d(x)=P_dx=(\realisation_\ReLU(\affineANN_{P_d,0}))(x)$
    shows that for all
        $d\in\N$,
        $\varepsilon \in (0,1]$ 
    it holds that
    \begin{equation}
        \llabel{eq:2}
        \realisation_{\ReLU}(\compANN{\affineANN_{P_d,0}}{\mathscr{G}_{d,\varepsilon}})
        =
        [\realisation_\ReLU(\affineANN_{P_d,0})]\circ [\realisation_{\ReLU}(\mathscr{G}_{d,\varepsilon})]
        =
        p_d\circ [\realisation_{\ReLU}(\mathscr{G}_{d,\varepsilon})]
        \in C(\R^d,\R)
        .
    \end{equation}
        This,
        \lref{eq:approx},
        and the fact that for all
            $d\in\N$,
            $x=(x_1,x_2,\dots,x_d)\in\R^d$
        it holds that
            $\abs{p_d(x)}=\abs{x_d}\leq \norm{x}$
    prove that for all
        $d \in \N$, 
        $\varepsilon \in (0,1]$ 
    it holds that
    \begin{equation}
    \begin{split}
	    &\sup_{x \in [-R,R]^d}\abs[\big]{ (p_d\circ F_{a_n} \circ \ldots \circ F_{a_2} \circ F_{a_1})(x) - (\realisation_{\ReLU}(\compANN{\affineANN_{P_d,0}}{\mathscr{G}_{d,\varepsilon}}))(x) } 
        \\&=
	    \sup_{x \in [-R,R]^d}\abs[\big]{ p_d\prb{(F_{a_n} \circ \ldots \circ F_{a_2} \circ F_{a_1})(x)} - p_d\prb{(\realisation_{\ReLU}(\mathscr{G}_{d,\varepsilon}))(x)} } 
        \\&=
	    \sup_{x \in [-R,R]^d}\abs[\big]{ p_d\prb{(F_{a_n} \circ \ldots \circ F_{a_2} \circ F_{a_1})(x) - (\realisation_{\ReLU}(\mathscr{G}_{d,\varepsilon}))(x)} } 
        \\&\leq
	    \sup_{x \in [-R,R]^d}\norm{ (F_{a_n} \circ \ldots \circ F_{a_2} \circ F_{a_1})(x) - (\realisation_{\ReLU}(\mathscr{G}_{d,\varepsilon}))(x) } 
        \leq 
        \varepsilon
        .
    \end{split}
	\end{equation}
    Combining
        this
    with
        \lref{eq:1} and 
        \lref{eq:2}
    establishes
        \lref{claim}.
    \finishproofthis
\end{aproof}

\subsection*{Acknowledgements}
This project is based on the master thesis of PB written from
January 2020 to June 2020 at ETH Zurich under the supervision 
of AJ and PC.
This work has been supported by the Ministry of Culture and Science NRW as part of the Lamarr Fellow Network.
Moreover, this work has been funded by the Deutsche Forschungsgemeinschaft (DFG, German Research
Foundation) under Germany's Excellence Strategy EXC 2044-390685587, Mathematics Münster:
Dynamics--Geometry--Structure.

\bibliographystyle{acm}
 \bibliography{bibfile}

\begin{thebibliography}{100}

\bibitem{almira2021negative}
{\sc Almira, J., de~Teruel, P.~L., Romero-López, D., and Voigtlaender, F.}
\newblock Negative results for approximation using single layer and multilayer feedforward neural networks.
\newblock {\em J. Math. Anal. Appl. 494}, 1 (2021), 124584.

\bibitem{bach2017breaking}
{\sc Bach, F.}
\newblock Breaking the curse of dimensionality with convex neural networks.
\newblock {\em J. Mach. Learn. Res. 18}, 19 (2017), 53 pages.

\bibitem{barron1992neural}
{\sc Barron, A.}
\newblock Neural net approximation.
\newblock In {\em Proceedings of the 7th Yale Workshop on Adaptive and Learning Systems\/} (1992), pp.~69--72.

\bibitem{barron1993universal}
{\sc Barron, A.~R.}
\newblock Universal approximation bounds for superpositions of a sigmoidal function.
\newblock {\em IEEE Trans. Inform. Theory 39}, 3 (1993), 930--945.

\bibitem{barron1994approximation}
{\sc Barron, A.~R.}
\newblock Approximation and estimation bounds for artificial neural networks.
\newblock {\em Mach. Learn. 14}, 1 (1994), 115--133.

\bibitem{BeckBeckerCheriditoJentzenNeufeld2019}
{\sc Beck, C., Becker, S., Cheridito, P., Jentzen, A., and Neufeld, A.}
\newblock Deep {S}plitting {M}ethod for {P}arabolic {PDE}s.
\newblock {\em SIAM J. Sci. Comput. 43}, 5 (2021), A3135--A3154.

\bibitem{Kuckuck2020overview}
{\sc Beck, C., Hutzenthaler, M., Jentzen, A., and Kuckuck, B.}
\newblock An overview on deep learning-based approximation methods for partial differential equations.
\newblock {\em Discrete Contin. Dyn. Syst. Ser. B 28}, 6 (2023), 3697--3746.

\bibitem{BeckJentzenKuckuck2019}
{\sc Beck, C., Jentzen, A., and Kuckuck, B.}
\newblock Full error analysis for the training of deep neural networks.
\newblock {\em Infin. Dimens. Anal. Quantum Probab. Relat. Top. 25}, 2 (2022), Paper no.~2150020, 77~pp.

\bibitem{Bellman1957}
{\sc Bellman, R.}
\newblock {\em Dynamic Programming}.
\newblock Princeton Landmarks in Mathematics. Princeton University Press, Princeton, NJ, 2010.
\newblock Reprint of the 1957 edition.

\bibitem{beneventano2020highdimensional}
{\sc Beneventano, P., Cheridito, P., Jentzen, A., and von Wurstemberger, P.}
\newblock {High-dimensional approximation spaces of artificial neural networks and applications to partial differential equations}.
\newblock {\em arXiv:2012.04326\/} (2020), 32 pages.

\bibitem{BernerGrohsJentzen2018}
{\sc Berner, J., Grohs, P., and Jentzen, A.}
\newblock Analysis of the generalization error: empirical risk minimization over deep artificial neural networks overcomes the curse of dimensionality in the numerical approximation of {B}lack-{S}choles partial differential equations.
\newblock {\em SIAM J. Math. Data Sci. 2}, 3 (2020), 631--657.

\bibitem{blum1991approximation}
{\sc Blum, E.~K., and Li, L.~K.}
\newblock Approximation theory and feedforward networks.
\newblock {\em Neural Netw. 4}, 4 (1991), 511--515.

\bibitem{BolcskeiGrohsKutyniokPetersen2019OptimalApproximation}
{\sc B\"{o}lcskei, H., Grohs, P., Kutyniok, G., and Petersen, P.}
\newblock Optimal approximation with sparsely connected deep neural networks.
\newblock {\em SIAM J. Math. Data Sci. 1}, 1 (2019), 8--45.

\bibitem{breiman1993hinging}
{\sc Breiman, L.}
\newblock Hinging hyperplanes for regression, classification, and function approximation.
\newblock {\em IEEE Trans. Inform. Theory 39}, 3 (1993), 999--1013.

\bibitem{BurgerNeubauer2001}
{\sc Burger, M., and Neubauer, A.}
\newblock Error bounds for approximation with neural networks.
\newblock {\em J. Approx. Theory 112}, 2 (2001), 235--250.

\bibitem{candes1998ridgelets}
{\sc Candes, E.~J.}
\newblock {\em Ridgelets: {T}heory and applications}.
\newblock ProQuest LLC, Ann Arbor, MI, 1998.
\newblock Ph.D. Thesis, Stanford University.

\bibitem{caragea2020neural}
{\sc Caragea, A., Petersen, P., and Voigtlaender, F.}
\newblock {Neural network approximation and estimation of classifiers with classification boundary in a Barron class}.
\newblock {\em Ann. Appl. Probab. 33}, 4 (2023), 3039 -- 3079.

\bibitem{carroll1989construction}
{\sc Carroll, and Dickinson}.
\newblock Construction of neural nets using the radon transform.
\newblock In {\em International 1989 Joint Conference on Neural Networks\/} (1989), vol.~1, pp.~607--611.

\bibitem{chen1995approximation}
{\sc Chen, T., and Chen, H.}
\newblock Approximation capability to functions of several variables, nonlinear functionals, and operators by radial basis function neural networks.
\newblock {\em IEEE Trans. on Neural Networks 6}, 4 (1995), 904--910.

\bibitem{cheridito2021efficient}
{\sc Cheridito, P., Jentzen, A., and Rossmannek, F.}
\newblock Efficient approximation of high-dimensional functions with neural networks.
\newblock {\em IEEE Trans. Neural Netw. Learn. Syst. 33}, 7 (2022), 3079--3093.

\bibitem{ChuiMhaskar1994}
{\sc Chui, C.~K., Li, X., and Mhaskar, H.~N.}
\newblock Neural networks for localized approximation.
\newblock {\em Math. Comp. 63}, 208 (1994), 607--623.

\bibitem{chui2019deep}
{\sc Chui, C.~K., Lin, S.-B., and Zhou, D.-X.}
\newblock Deep neural networks for rotation-invariance approximation and learning.
\newblock {\em Anal. Appl. (Singap.) 17}, 5 (2019), 737--772.

\bibitem{cohen2016expressive}
{\sc Cohen, N., Sharir, O., and Shashua, A.}
\newblock On the expressive power of deep learning: A tensor analysis.
\newblock In {\em 29th Annual Conference on Learning Theory\/} (23--26 Jun 2016), V.~Feldman, A.~Rakhlin, and O.~Shamir, Eds., vol.~49 of {\em Proceedings of Machine Learning Research}, PMLR, pp.~698--728.

\bibitem{Cybenko1989}
{\sc Cybenko, G.}
\newblock Approximation by superpositions of a sigmoidal function.
\newblock {\em Math. Control Signals Systems 2}, 4 (1989), 303--314.

\bibitem{daniely2017depth}
{\sc Daniely, A.}
\newblock Depth separation for neural networks.
\newblock In {\em Proceedings of the 2017 Conference on Learning Theory\/} (07--10 Jul 2017), S.~Kale and O.~Shamir, Eds., vol.~65 of {\em Proceedings of Machine Learning Research}, PMLR, pp.~690--696.

\bibitem{devore1997approximation}
{\sc DeVore, R.~A., Oskolkov, K.~I., and Petrushev, P.~P.}
\newblock Approximation by feed-forward neural networks.
\newblock {\em Ann. Numer. Math. 4}, 1--4 (1997), 261--287.
\newblock The heritage of P. L. Chebyshev: a Festschrift in honor of the 70th birthday of T. J. Rivlin.

\bibitem{donahue1997rates}
{\sc Donahue, M.~J., Gurvits, L., Darken, C., and Sontag, E.}
\newblock Rates of convex approximation in non-{H}ilbert spaces.
\newblock {\em Constr. Approx. 13}, 2 (1997), 187--220.

\bibitem{EHanJentzen2017}
{\sc E, W., Han, J., and Jentzen, A.}
\newblock Deep learning-based numerical methods for high-dimensional parabolic partial differential equations and backward stochastic differential equations.
\newblock {\em Commun. Math. Stat. 5}, 4 (2017), 349--380.

\bibitem{Jiequn2020AlgorithmsPDEs}
{\sc E, W., Han, J., and Jentzen, A.}
\newblock Algorithms for solving high dimensional {PDEs}: from nonlinear monte carlo to machine learning.
\newblock {\em Nonlinearity 35}, 1 (2021), 278--310.

\bibitem{wang2018exponential}
{\sc E, W., and Wang, Q.}
\newblock Exponential convergence of the deep neural network approximation for analytic functions.
\newblock {\em Science China Mathematics 61}, 10 (2018), 1733--1740.

\bibitem{ElbraechterSchwab2018}
{\sc Elbr\"{a}chter, D., Grohs, P., Jentzen, A., and Schwab, C.}
\newblock {DNN} expression rate analysis of high-dimensional {PDE}s: Application to option pricing.
\newblock {\em Constr. Approx. 55\/} (2022), 3--71.

\bibitem{elbraechter2021deep}
{\sc Elbr\"{a}chter, D., Perekrestenko, D., Grohs, P., and B\"{o}lcskei, H.}
\newblock Deep neural network approximation theory.
\newblock {\em IEEE Trans. Inform. Theory 67}, 5 (2021), 2581--2623.

\bibitem{eldan2016power}
{\sc Eldan, R., and Shamir, O.}
\newblock The power of depth for feedforward neural networks.
\newblock In {\em 29th Annual Conference on Learning Theory\/} (23--26 Jun 2016), V.~Feldman, A.~Rakhlin, and O.~Shamir, Eds., vol.~49 of {\em Proceedings of Machine Learning Research}, PMLR, pp.~907--940.

\bibitem{Ellacott1994}
{\sc Ellacott, S.~W.}
\newblock Aspects of the numerical analysis of neural networks.
\newblock {\em Acta Numer. 3\/} (1994), 145--202.

\bibitem{funahashi1989approximate}
{\sc Funahashi, K.-I.}
\newblock On the approximate realization of continuous mappings by neural networks.
\newblock {\em Neural Netw. 2}, 3 (1989), 183--192.

\bibitem{gallant1988there}
{\sc Gallant, and White}.
\newblock There exists a neural network that does not make avoidable mistakes.
\newblock In {\em IEEE 1988 International Conference on Neural Networks\/} (1988), vol.~1, pp.~657--664.

\bibitem{GononGrohsEtAl2019}
{\sc Gonon, L., Grohs, P., Jentzen, A., Kofler, D., and Šiška, D.}
\newblock Uniform error estimates for artificial neural network approximations for heat equations.
\newblock {\em IMA J. Numer. Anal. 42}, 3 (2022), 1991--2054.

\bibitem{gonon2020deep}
{\sc Gonon, L., and Schwab, C.}
\newblock {Deep ReLU network expression rates for option prices in high-dimensional, exponential Lévy models}.
\newblock Tech. Rep. 2020-52, Seminar for Applied Mathematics, ETH Z{\"u}rich, Switzerland, 2020.

\bibitem{gribonval2019approximation}
{\sc Gribonval, R., Kutyniok, G., Nielsen, M., and Voigtlaender, F.}
\newblock Approximation spaces of deep neural networks.
\newblock {\em Constr. Approx. 55\/} (2022), 259--367.

\bibitem{GrohsHerrmann2020}
{\sc Grohs, P., and Herrmann, L.}
\newblock Deep neural network approximation for high-dimensional elliptic {PDE}s with boundary conditions.
\newblock {\em IMA J. Numer. Anal. 42}, 3 (05 2021), 2055--2082.

\bibitem{GrohsWurstemberger2018}
{\sc Grohs, P., Hornung, F., Jentzen, A., and von Wurstemberger, P.}
\newblock {A proof that artificial neural networks overcome the curse of dimensionality in the numerical approximation of Black--Scholes partial differential equations}.
\newblock {\em Mem. Amer. Math. Soc. 284}, 1410 (2023), 106~pp.

\bibitem{grohs2019space}
{\sc Grohs, P., Hornung, F., Jentzen, A., and Zimmermann, P.}
\newblock Space-time error estimates for deep neural network approximations for differential equations.
\newblock {\em Adv. Comput. Math. 49}, 1 (2023), Paper no.~4, 78~pp.

\bibitem{GrohsIbrgimovJentzen2021}
{\sc Grohs, P., Ibragimov, S., Jentzen, A., and Koppensteiner, S.}
\newblock Lower bounds for artificial neural network approximations: A proof that shallow neural networks fail to overcome the curse of dimensionality.
\newblock {\em J. Complexity 77\/} (2023), Paper no.~101746, 53~pp.

\bibitem{GrohsJentzenSalimova2019}
{\sc Grohs, P., Jentzen, A., and Salimova, D.}
\newblock {Deep neural network approximations for solutions of PDEs based on Monte Carlo algorithms}.
\newblock {\em Partial Differ. Equ. Appl. 3}, 4 (2022), Paper no.~45, 41~pp.

\bibitem{grohs2021proof}
{\sc Grohs, P., and Voigtlaender, F.}
\newblock Proof of the theory-to-practice gap in deep learning via sampling complexity bounds for neural network approximation spaces.
\newblock {\em Found. Comput. Math. 24\/} (2024), 1085--1143.

\bibitem{guhring2019error}
{\sc G\"{u}hring, I., Kutyniok, G., and Petersen, P.}
\newblock Error bounds for approximations with deep {ReLU} neural networks in {$W^{s,p}$} norms.
\newblock {\em Anal. Appl. (Singap.) 18}, 5 (2020), 803--859.

\bibitem{guehring2020expressivity}
{\sc G\"uhring, I., Raslan, M., and Kutyniok, G.}
\newblock Expressivity of deep neural networks.
\newblock In {\em Mathematical Aspects of Deep Learning}, P.~Grohs and G.~Kutyniok, Eds. Cambridge University Press, Cambridge, 2023, pp.~149--199.

\bibitem{GuliIsm2018a}
{\sc Guliyev, N.~J., and Ismailov, V.~E.}
\newblock Approximation capability of two hidden layer feedforward neural networks with fixed weights.
\newblock {\em Neurocomputing 316\/} (2018), 262--269.

\bibitem{GULIYEV2018296}
{\sc Guliyev, N.~J., and Ismailov, V.~E.}
\newblock {On the approximation by single hidden layer feedforward neural networks with fixed weights}.
\newblock {\em Neural Netw. 98\/} (2018), 296--304.

\bibitem{HanJentzenE2018}
{\sc Han, J., Jentzen, A., and E, W.}
\newblock Solving high-dimensional partial differential equations using deep learning.
\newblock {\em Proc. Natl. Acad. Sci. USA 115}, 34 (2018), 8505--8510.

\bibitem{Hanin2017}
{\sc Hanin, B.}
\newblock Universal function approximation by deep neural nets with bounded width and {ReLU} activations.
\newblock {\em Mathematics 7}, 10 (2019).

\bibitem{heinrich2006randomized}
{\sc Heinrich, S.}
\newblock The randomized information complexity of elliptic {PDE}.
\newblock {\em J. Complexity 22}, 2 (2006), 220--249.

\bibitem{heinrich1999monte}
{\sc Heinrich, S., and Sindambiwe, E.}
\newblock {M}onte {C}arlo complexity of parametric integration.
\newblock {\em J. Complexity 15}, 3 (1999), 317--341.

\bibitem{hornik1991approximation}
{\sc Hornik, K.}
\newblock Approximation capabilities of multilayer feedforward networks.
\newblock {\em Neural Netw. 4}, 2 (1991), 251--257.

\bibitem{hornik1993some}
{\sc Hornik, K.}
\newblock Some new results on neural network approximation.
\newblock {\em Neural Netw. 6}, 8 (1993), 1069--1072.

\bibitem{hornik1989multilayer}
{\sc Hornik, K., Stinchcombe, M., and White, H.}
\newblock Multilayer feedforward networks are universal approximators.
\newblock {\em Neural Netw. 2}, 5 (1989), 359--366.

\bibitem{hornik1990universal}
{\sc Hornik, K., Stinchcombe, M., and White, H.}
\newblock Universal approximation of an unknown mapping and its derivatives using multilayer feedforward networks.
\newblock {\em Neural Netw. 3}, 5 (1990), 551--560.

\bibitem{HornungJentzenSalimova2020}
{\sc Hornung, F., Jentzen, A., and Salimova, D.}
\newblock Space-time deep neural network approximations for high-dimensional partial differential equations.
\newblock {\em arXiv:2006.02199\/} (2020), 52 pages.

\bibitem{HutzenthalerJentzenKruseNguyen2019}
{\sc Hutzenthaler, M., Jentzen, A., Kruse, T., and Nguyen, T.~A.}
\newblock A proof that rectified deep neural networks overcome the curse of dimensionality in the numerical approximation of semilinear heat equations.
\newblock {\em Partial Differ. Equ. Appl. 1}, 2 (2020), Paper No. 10, 34 pages.

\bibitem{irie1988capabilities}
{\sc Irie, and Miyake}.
\newblock Capabilities of three-layered perceptrons.
\newblock In {\em IEEE 1988 International Conference on Neural Networks\/} (1988), vol.~1, pp.~641--648.

\bibitem{jentzen2023mathematical}
{\sc Jentzen, A., Kuckuck, B., and von Wurstemberger, P.}
\newblock {\em {Mathematical Introduction to Deep Learning: Methods, Implementations, and Theory}}.
\newblock arXiv:2310.20360v2, 2023.
\newblock Preprint, \url{https://arxiv.org/abs/2310.20360v2}.

\bibitem{jentzen2020strong2}
{\sc Jentzen, A., and Riekert, A.}
\newblock Strong overall error analysis for the training of artificial neural networks via random initializations.
\newblock {\em Commun. Math. Stat. 12\/} (2024), 385--434.

\bibitem{JentzenSalimovaWelti2018}
{\sc Jentzen, A., Salimova, D., and Welti, T.}
\newblock A proof that deep artificial neural networks overcome the curse of dimensionality in the numerical approximation of {K}olmogorov partial differential equations with constant diffusion and nonlinear drift coefficients.
\newblock {\em Commun. Math. Sci. 19}, 5 (2021), 1167--1205.

\bibitem{jones1992simple}
{\sc Jones, L.~K.}
\newblock A simple lemma on greedy approximation in {H}ilbert space and convergence rates for projection pursuit regression and neural network training.
\newblock {\em Ann. Statist. 20}, 1 (1992), 608--613.

\bibitem{kainen2009complexity}
{\sc Kainen, P.~C., K\r{u}rkov\'{a}, V., and Sanguineti, M.}
\newblock Complexity of {G}aussian-radial-basis networks approximating smooth functions.
\newblock {\em J. Complexity 25}, 1 (2009), 63--74.

\bibitem{KaiKurSang2012}
{\sc Kainen, P.~C., K\r{u}rkov\'{a}, V., and Sanguineti, M.}
\newblock Dependence of computational models on input dimension: tractability of approximation and optimization tasks.
\newblock {\em IEEE Trans. Inform. Theory 58}, 2 (2012), 1203--1214.

\bibitem{kidger2020universal}
{\sc Kidger, P., and Lyons, T.}
\newblock {Universal Approximation with Deep Narrow Networks}.
\newblock In {\em Proceedings of Thirty Third Conference on Learning Theory\/} (09--12 Jul 2020), J.~Abernethy and S.~Agarwal, Eds., vol.~125 of {\em Proceedings of Machine Learning Research}, PMLR, pp.~2306--2327.

\bibitem{klusowski2018approximation}
{\sc Klusowski, J.~M., and Barron, A.~R.}
\newblock Approximation by combinations of {R}e{LU} and squared {R}e{LU} ridge functions with {$\ell^1$} and {$\ell^0$} controls.
\newblock {\em IEEE Trans. Inform. Theory 64}, 12 (2018), 7649--7656.

\bibitem{kurkova2002comparison}
{\sc K\r{u}rkov\'{a}, V., and Sanguineti, M.}
\newblock Comparison of worst case errors in linear and neural network approximation.
\newblock {\em IEEE Trans. Inform. Theory 48}, 1 (2002), 264--275.

\bibitem{kurkova2008geometric}
{\sc K\r{u}rkov\'{a}, V., and Sanguineti, M.}
\newblock Geometric upper bounds on rates of variable-basis approximation.
\newblock {\em IEEE Trans. Inform. Theory 54}, 12 (2008), 5681--5688.

\bibitem{KurKaiKre1997}
{\sc K{\r{u}}rkov{\'a}, V., Kainen, P.~C., and Kreinovich, V.}
\newblock {Estimates of the Number of Hidden Units and Variation with Respect to Half-Spaces}.
\newblock {\em Neural Netw. 10}, 6 (1997), 1061--1068.

\bibitem{lavretsky2002geometric}
{\sc Lavretsky, E.}
\newblock On the geometric convergence of neural approximations.
\newblock {\em IEEE Trans. on Neural Networks 13}, 2 (2002), 274--282.

\bibitem{lee2017ability}
{\sc Lee, H., Ge, R., Ma, T., Risteski, A., and Arora, S.}
\newblock On the ability of neural nets to express distributions.
\newblock In {\em Proceedings of the 2017 Conference on Learning Theory\/} (07--10 Jul 2017), S.~Kale and O.~Shamir, Eds., vol.~65 of {\em Proceedings of Machine Learning Research}, PMLR, pp.~1271--1296.

\bibitem{leshno1993multilayer}
{\sc Leshno, M., Lin, V.~Y., Pinkus, A., and Schocken, S.}
\newblock Multilayer feedforward networks with a nonpolynomial activation function can approximate any function.
\newblock {\em Neural Netw. 6}, 6 (1993), 861--867.

\bibitem{li2020better}
{\sc Li, B., Tang, S., and Yu, H.}
\newblock Better approximations of high dimensional smooth functions by deep neural networks with rectified power units.
\newblock {\em Commun. Comput. Phys. 27}, 2 (2020), 379--411.

\bibitem{ChenWu2019}
{\sc Liang, C., and Wu, C.}
\newblock A note on the expressive power of deep rectified linear unit networks in high-dimensional spaces.
\newblock {\em Math. Methods Appl. Sci. 42}, 9 (2019), 3400--3404.

\bibitem{LuShenYangZhang2020}
{\sc Lu, J., Shen, Z., Yang, H., and Zhang, S.}
\newblock Deep network approximation for smooth functions.
\newblock {\em SIAM J. Math. Anal. 53}, 5 (2021), 5465--5506.

\bibitem{MaioPinkus1999}
{\sc Maiorov, V., and Pinkus, A.}
\newblock Lower bounds for approximation by {MLP} neural networks.
\newblock {\em Neurocomputing 25}, 1 (1999), 81--91.

\bibitem{maiorov2000near}
{\sc Maiorov, V.~E., and Meir, R.}
\newblock On the near optimality of the stochastic approximation of smooth functions by neural networks.
\newblock {\em Adv. Comput. Math. 13}, 1 (2000), 79--103.

\bibitem{makovoz1996random}
{\sc Makovoz, Y.}
\newblock Random approximants and neural networks.
\newblock {\em J. Approx. Theory 85}, 1 (1996), 98--109.

\bibitem{makovoz1998uniform}
{\sc Makovoz, Y.}
\newblock Uniform approximation by neural networks.
\newblock {\em J. Approx. Theory 95}, 2 (1998), 215--228.

\bibitem{mhaskar1996neural}
{\sc Mhaskar, H.~N.}
\newblock Neural networks for optimal approximation of smooth and analytic functions.
\newblock {\em Neural Comput. 8}, 1 (1996), 164--177.

\bibitem{mhaskar1993}
{\sc Mhaskar, H.~N.}
\newblock Approximation properties of a multilayered feedforward artificial neural network.
\newblock {\em Adv. Comput. Math. 1\/} (2017), 61--80.

\bibitem{MhaskarMicchelli1994}
{\sc Mhaskar, H.~N., and Micchelli, C.~A.}
\newblock Dimension-independent bounds on the degree of approximation by neural networks.
\newblock {\em IBM J. Res. Dev. 38}, 3 (1994), pp. 227--284.

\bibitem{MhaskarMicchelli1995}
{\sc Mhaskar, H.~N., and Micchelli, C.~A.}
\newblock Degree of approximation by neural and translation networks with a single hidden layer.
\newblock {\em Adv. in Appl. Math. 16}, 2 (1995), 151--183.

\bibitem{MhaskarPoggio2016}
{\sc Mhaskar, H.~N., and Poggio, T.}
\newblock Deep vs. shallow networks: an approximation theory perspective.
\newblock {\em Anal. Appl. (Singap.) 14}, 6 (2016), 829--848.

\bibitem{nguyen1999approximation}
{\sc Nguyen-Thien, T., and Tran-Cong, T.}
\newblock Approximation of functions and their derivatives: A neural network implementation with applications.
\newblock {\em Applied Math. Model. 23}, 9 (1999), 687--704.

\bibitem{novak1997curse}
{\sc Novak, E., and Ritter, K.}
\newblock The curse of dimension and a universal method for numerical integration.
\newblock In {\em Multivariate Approximation and Splines\/} (Basel, 1997), G.~N{\"u}rnberger, J.~W. Schmidt, and G.~Walz, Eds., Birkh{\"a}user Basel, pp.~177--187.

\bibitem{NovakWozniakowski2008}
{\sc Novak, E., and Wo\'{z}niakowski, H.}
\newblock {\em Tractability of multivariate problems. {V}olume I: {L}inear information}, vol.~6 of {\em EMS Tracts in Mathematics}.
\newblock European Mathematical Society (EMS), Z\"{u}rich, 2008.

\bibitem{NovakWozniakowski2010}
{\sc Novak, E., and Wo\'{z}niakowski, H.}
\newblock {\em Tractability of multivariate problems. {V}olume {II}: {S}tandard information for functionals}, vol.~12 of {\em EMS Tracts in Mathematics}.
\newblock European Mathematical Society (EMS), Z\"{u}rich, 2010.

\bibitem{park1991universal}
{\sc Park, J., and Sandberg, I.~W.}
\newblock Universal approximation using radial-basis-function networks.
\newblock {\em Neural Comput. 3}, 2 (1991), 246--257.

\bibitem{perekrestenko2018universal}
{\sc Perekrestenko, D., Grohs, P., Elbr{\"a}chter, D., and B{\"o}lcskei, H.}
\newblock The universal approximation power of finite-width deep {R}e{LU} networks.
\newblock {\em arXiv:1806.01528\/} (2018), 16 pages.

\bibitem{petersen2018optimal}
{\sc Petersen, P., and Voigtlaender, F.}
\newblock Optimal approximation of piecewise smooth functions using deep {ReLU} neural networks.
\newblock {\em Neural Netw. 108\/} (2018), 296--330.

\bibitem{pinkus1999approximation}
{\sc Pinkus, A.}
\newblock Approximation theory of the {MLP} model in neural networks.
\newblock {\em Acta Numer. 8\/} (1999), 143--195.

\bibitem{poggio2017why}
{\sc Poggio, T., Mhaskar, H., Rosasco, L., Miranda, B., and Liao, Q.}
\newblock Why and when can deep--but not shallow--networks avoid the curse of dimensionality: {A} review.
\newblock {\em Int. J. Autom. Comput. 14}, 5 (2017), 503--519.

\bibitem{ReisingerZhang2019}
{\sc Reisinger, C., and Zhang, Y.}
\newblock Rectified deep neural networks overcome the curse of dimensionality for nonsmooth value functions in zero-sum games of nonlinear stiff systems.
\newblock {\em Anal. Appl. (Singap.) 18}, 06 (2020), 951--999.

\bibitem{safran17a}
{\sc Safran, I., and Shamir, O.}
\newblock Depth-width tradeoffs in approximating natural functions with neural networks.
\newblock In {\em Proceedings of the 34th International Conference on Machine Learning\/} (06--11 Aug 2017), D.~Precup and Y.~W. Teh, Eds., vol.~70 of {\em Proceedings of Machine Learning Research}, PMLR, pp.~2979--2987.

\bibitem{schmitt2000lower}
{\sc Schmitt, M.}
\newblock Lower bounds on the complexity of approximating continuous functions by sigmoidal neural networks.
\newblock In {\em Advances in Neural Information Processing Systems\/} (2000), S.~Solla, T.~Leen, and K.~M\"{u}ller, Eds., vol.~12, MIT Press, pp.~328--334.

\bibitem{SchwabZech2019}
{\sc Schwab, C., and Zech, J.}
\newblock Deep learning in high dimension: {N}eural network expression rates for generalized polynomial chaos expansions in {UQ}.
\newblock {\em Anal. Appl. (Singap.) 17}, 1 (2019), 19--55.

\bibitem{ShahamCloningerCoifman2018}
{\sc Shaham, U., Cloninger, A., and Coifman, R.~R.}
\newblock Provable approximation properties for deep neural networks.
\newblock {\em Appl. Comput. Harmon. Anal. 44}, 3 (2018), 537--557.

\bibitem{shen2019nonlinear}
{\sc Shen, Z., Yang, H., and Zhang, S.}
\newblock Nonlinear approximation via compositions.
\newblock {\em Neural Netw. 119\/} (2019), 74--84.

\bibitem{shen2020deep}
{\sc Shen, Z., Yang, H., and Zhang, S.}
\newblock Deep network approximation characterized by number of neurons.
\newblock {\em Commun. Comput. Phys. 28}, 5 (2020), 1768--1811.

\bibitem{voigtlaender2019approximation}
{\sc Voigtlaender, F., and Petersen, P.}
\newblock Approximation in {$L^p(\mu)$} with deep {ReLU} neural networks.
\newblock In {\em 2019 13th International conference on Sampling Theory and Applications (SampTA)\/} (2019), IEEE, pp.~1--4.

\bibitem{yarotsky2017error}
{\sc Yarotsky, D.}
\newblock Error bounds for approximations with deep {ReLU} networks.
\newblock {\em Neural Netw. 94\/} (2017), 103--114.

\bibitem{JMLR:v25:23-0912}
{\sc Zhang, S., Lu, J., and Zhao, H.}
\newblock Deep network approximation: Beyond relu to diverse activation functions.
\newblock {\em Journal of Machine Learning Research 25}, 35 (2024), 1--39.

\end{thebibliography}
 
\end{document}